%% file: Lobus_-_April_2023.tex
\documentclass[12pt,amsfonts]{article}

\usepackage{latexsym}
\usepackage{amsfonts,amssymb}
\usepackage{mathrsfs} 
 \usepackage{hyperref}

\newtheorem{theorem}{Theorem}[section]
\newtheorem{lemma}[theorem]{Lemma}
\newtheorem{remark}
{Remark}
\newtheorem{corollary}[theorem]{Corollary}
\newtheorem{proposition}[theorem]{Proposition}

\newtheorem{definition}[theorem]{Definition} 
\newcommand{\br}[1]{\begin{remark}\label{#1}\rm}
\newcommand{\er}{\end{remark}}
\newcommand{\bd}[1]{\begin{definition}\label{#1}\rm}
\newcommand{\ed}{\end{definition}}
\newcommand{\bt}[1]{\begin{theorem}\label{#1}}
\newcommand{\et}{\end{theorem}}
\newcommand{\bprop}[1]{\begin{proposition}\label{#1}}
\newcommand{\eprop}{\end{proposition}}
\newcommand{\bcor}[1]{\begin{corollary}\label{#1}}
\newcommand{\ecor}{\end{corollary}}

\newcommand{\D}{\displaystyle}
\newcommand{\T}{\textstyle}
\newcommand{\lra}{\longrightarrow}
\newcommand{\Ra}{\Longrightarrow}

\newcommand{\stack}[2]{\raisebox{-2pt} 
{\renewcommand{\arraystretch}{.01} 
\begin{tabular}{c} 
$#2$\\$\scriptscriptstyle #1$ 
\end{tabular} 
}}

\renewcommand{\l}{{\it leb}}

\newcommand{\vp}{\varphi}
\newcommand{\ve}{\varepsilon}
\newcommand{\nid}{\noindent}
\newcommand{\qed}{\hfill$\Box$} 

\def\1{\, {\rm I}\mskip-10mu 1}

\renewcommand{\t}[1]{\tilde{#1}} 
\newcommand{\bmu}{\mbox{\boldmath${\mu}$}}
\newcommand{\bnu}{\mbox{\boldmath${\nu}$}} 
\newcommand{\sbmu}{\mbox{\tiny\boldmath${\mu}$}}

\begin{document}
\title{Quasi-Invariance under Flows Generated by Non-Linear PDEs} 
\par
\author{J\"org-Uwe L\"obus
\\  Department of Mathematics, University of Link\"oping  \\ 
SE-581 83 Link\"oping \\ 
Sweden \\ {\tt jorg-uwe.lobus@liu.se} \\ 
{\tt https://orcid.org/0000-0001-5646-1277}
}
\date{}
\maketitle
{\footnotesize
\noindent
\begin{quote}
{\bf Abstract} 

The paper is concerned with the change of probability measures $\bmu$
along non-random probability measure valued trajectories $\nu_t$, $t\in 
[-1,1]$. Typically solutions to non-linear PDEs, modeling spatial 
development as time progresses, generate such trajectories. Depending on 
in which direction the map $\nu\equiv\nu_0\mapsto\nu_t$ does not exit 
the state space, for $t\in [-1,0]$ or for $t\in [0,1]$, 
the Radon-Nikodym derivative $d\bmu\circ\nu_t/d\bmu$ 
is determined. It is also investigated how  Fr\'echet differentiability 
of the solution map of the PDE can contribute to the existence of this 
Radon-Nikodym derivative. The first application is a certain Boltzmann 
type equation. Here the Fr\'echet derivative of the solution map is 
calculated explicitly and quasi-invariance is established. The second 
application is a PDE related to the asymptotic behavior of a Fleming-Viot 
type particle system. Here it is 
demonstrated how quasi-invariance can be used in order to derive a corresponding 
integration by parts formula. 
\noindent 



\noindent
{\bf Running head} Quasi-Invariance

\end{quote}
}

\input paper-April-2023
\end{document}

%% file: paper-April-2023.tex
\section{Introduction}\label{sec:1}
\setcounter{equation}{0}

{\bf Equivalence of families of differentiable measures} and related 
Radon-Nikodym derivatives have been studied from different points 
of view in \cite{Be85}, \cite{Be06}, \cite{Bo97}, \cite{BM-W99}, 
\cite{Cr83}, \cite{DS88}, \cite{SvW93}, \cite{UZ00}. A very 
comprehensive presentation of the actual state of the subject has been 
given in \cite{Bo10} where the particular importance of the paper 
\cite{SvW93} is emphasized. Moreover, a mathematically elegant approach 
to a particular situation has been carried out in \cite{Be06}, Appendix. 

In the present paper we consider a probability measure valued solution 
$t\mapsto\nu_t$ to a non-linear PDE of type $\frac{d}{dt}\nu_t=A^f 
\nu_t$ where $A^f$ is an operator independent of time $t$. We assume 
that $\nu_t=h_t\cdot\lambda$ where $\lambda$ is a Borel-measure on 
some locally compact metric space $D$ and $h_t\in L^v(D)$ for some 
$1\le v\le\infty$. 

Let $\bmu$ be a certain probability measure over the space of 
probability measures $\nu$ on $D$. Moreover, let $\bmu\circ\nu_{t}$ 
be the image measures relative to the maps $\nu\equiv\nu_0\mapsto 
\nu_t$. We consider two cases, $t\in [-1,0]$ (Case 1) and 
$t \in [0,1]$ (Case 2). These cases correspond to the direction of 
time for which the maps $\nu\equiv\nu_0\mapsto\nu_t$ do not exit 
the state space. Here the state space is the set of initial measures 
$\nu_0\equiv h_0\cdot \lambda$ where the set of all $h_0$ is an 
appropriately chosen subset of $L^v(D)$. We are interested in the 
equivalence of $\bmu\circ \nu_{t}$ with respect to $\bmu$. In the 
paper we develop a certain concept of flows on bounded intervals 
of time and we call the measure $\bmu$ {\it quasi-invariant} with 
respect to the flow generated by the maps $\nu\mapsto\nu_t$, $t\in 
[-1,0]$ in Case 1 or $t\in [0,1]$ in Case 2, if for all such $t$ 
the measure $\bmu \circ\nu_{t}$ is equivalent to $\bmu$. 

At first glance one may be attempted to adjust the situation to the 
known theory in \cite{SvW93}, \cite{BM-W99}, and \cite{Bo10}, Section 
11.2, or \cite{Be06}, Appendix. However, as it holds for a number of 
PDEs, the solutions $\nu_t$ may no longer belong to some meaningful 
state space of measures for large $t>1$ in Case 1, or small $t<-1$ in 
Case 2. In this context, it may not be possible to extend the map 
$[-1,0]\ni t\mapsto\nu_t$ to a map $[-1,\infty)\ni t\mapsto\nu_t$ in 
Case 1. Similarly, it may not be possible to extend the map $[0,1] 
\ni t\mapsto\nu_t$ to a map $(-\infty,1]\ni t\mapsto\nu_t$ in Case 2. 
This observation requires a completely new framework relative to 
the quasi-invariance problem. As a consequence, the divergence term in 
the Radon-Nikodym derivative ${d\bmu\circ\nu_{t}}/{d\bmu}$ for $t\in 
[-1,0]$ in Case 1 or $t\in [0,1]$ in Case 2 includes an additional 
new component reflecting the observation. 

In addition, the 
verification of one of the required forms of differentiability of $t 
\mapsto\bmu\circ\nu_{t}$ in \cite{SvW93} and \cite{Bo10}, Section 11.2, 
is a major part of our direct proof of equivalence of the family $\bmu 
\circ\nu_{t}$, $t\in [-1,0]$ or $t\in [0,1]$. For a more detailed 
discussion of the work in \cite{SvW93} as well as \cite{Bo10}, Section 
11.2, and in the Appendix of \cite{Be06}, we would like to refer to 
Subsection \ref{sec:1:2}. 

We provide a direct treatment. This allows to focus on conditions which 
come up naturally in the proof of equivalence of probability  measures 
$\bmu$ under the flow generated by the maps $\nu\mapsto\nu_t$ relative to  
a non-linear PDE of type $\frac{d}{dt}\nu_t=A^f\nu_t$, $t\in [-1,0]$ or 
$t\in [0,1]$. In particular, these are the general conditions (j)-(jjj) 
of Subsection \ref{sec:2:1} on the vector field $A^f\cdot$ and the 
differentiability of the measure $\bmu$. Moreover hypotheses (jv) and 
(v) of Subsection \ref{sec:2:12} are conditions on the differentiability 
of the flow. The discussion in Subsection \ref{sec:2:12} refers to the 
question in which sense these sufficient conditions for the change of 
measure result, Theorem \ref{Theorem2.3}, are even necessary and hence 
optimal, cf. Remark (7) in Subsection \ref{sec:2:12} and Proposition 
\ref{Proposition2.5}. Except for the preliminaries for the 
applications in Sections \ref{sec:3} and \ref{sec:4} and some 
definitions in Subsection \ref{sec:2:11} which are included for 
better readability and comparability with the existing theory, all 
work in the paper is original. 
\bigskip 

For an effective {\bf presentation of the main results}, we shall now 
focus on descriptions that avoid technical details as much as possible. 
We will highlight the main results, among other things by the labels 
(\ref{1.5})-(\ref{R.4}) and refer to the corresponding theorems etc.~in 
the paper.
\bigskip

{\bf The Radon-Nikodym derivative ${d\bmu\circ \nu_{t}}/{d\bmu}$} for 
the two cases $t\in [-1,0]$ or $t\in [0,1]$ is the main objective of 
the paper. The quasi-invariance result of Theorem \ref{Theorem2.3} 
provides a Radon-Nikodym derivative which is measurable with respect 
to the $\sigma$-algebra generated by the open sets relative to the 
Prokhorov distance $\pi$ on the space $E$ of all probability measure 
on $D$. On an intuitive level, Theorem \ref{Theorem2.3} reads for 
$t\in [0,1]$ (Case 2) as
\begin{eqnarray}\label{1.5}
\frac{d\bmu\circ\nu_{t}}{d\bmu}=\exp\left\{\int_{s=0}^t\delta (A^f) 
(\nu_{s})\, ds\right\}\quad\bmu\mbox{\rm -a.e.} 
\end{eqnarray}
where $\delta (A^f)$ can be interpreted as the {\it divergence} of the 
measure $\bmu$ relative to the vector field $A^f\cdot$ corresponding to 
the PDE $\frac{d}{dt}\nu_t=A^f\nu_t$. To verify Fr\'echet 
differentiability of the map $\nu\to\nu_t$, embedded in suitable Banach 
spaces, can be beneficial in applications. For this consider conditions 
(jv) and (v) and Proposition \ref{Proposition2.4} of Subsection 
\ref{sec:2:12} and observe how Proposition \ref{Proposition2.4} is used 
in the application of Section \ref{sec:3}. Since the Fr\'echet 
differentiability of a continuous linear map $\nu\to\nu_t$ is obvious, 
the present paper is mainly directed to non-linear PDEs. 
\medskip 

{\bf The existence of the divergence $\delta (A^f)$ of the measure 
$\bmu$} along the solution paths $\nu\equiv\nu_0\mapsto\nu_t$ of the 
measure valued PDE is the form of differentiability of $\bmu$ we 
require for the quasi-invariance result of Theorem \ref{Theorem2.3}. 
More precisely, we suppose for $A^f\nu\equiv a^f\nu\cdot\lambda$ 
and a suitable class of test functions $f$ that 
\begin{eqnarray*}
\int Df(a^f)\, d\bmu=-\int f\, \delta (A^f)\, d\bmu+b\int f\, d 
\bmu_\gamma\, ,  
\end{eqnarray*}
where $\bmu_\gamma$ is a certain measure concentrated on the part of 
the state space which is exited by the trajectories immediately after 
time zero, and $b=1$ in Case 1 and $b=-1$ in Case 2. This is a 
mathematically appropriate condition. In concrete situations however 
one might be interested in separating hypotheses on $\bmu$ from 
hypotheses on the PDE and its solutions. 

In Subsection \ref{sec:2:11}, Paragraph \ref{sec:2:2:2}, we construct 
a certain class of measures $\bmu$. We show that there are divergences 
$\t \delta (b_n)$ of such measures $\bmu$ along a certain collection of 
basis vectors $b_n$, $n\ge 2$, such that $\int (Df,b_n)\, d\bmu=-\int f 
\, \t \delta (b_n)\, d\bmu$ for a suitable class of test functions $f$, 
cf. Proposition \ref{Proposition2.3}. Note that the basis vector $b_1$ 
is a constant function which plays an extra role. Let $a_n$ be the 
coefficients of the vector field $a^f\cdot$ relative to the basis vectors 
$b_n$, $n\ge 2$. In Proposition \ref{Proposition2.2} we show that in the 
specified setting of Subsection \ref{sec:2:11} we have 
\begin{eqnarray}\label{R.2}
\delta(A^f)=\sum_{n=2}^\infty\left((Da_n,b_n)+a_n\t \delta(b_n)\right)\, . 
\end{eqnarray} 
Interpreting the left-hand side as a Skorokhod type stochastic integral 
we provide arguments to interpret the right-hand side as an Ogawa type 
stochastic integral minus the trace of the gradient of the integrand. 
These and several other relations between differentiability of measures 
and non-Gaussian stochastic calculus are investigated in Subsection 
\ref{sec:2:11}. We even discover several parallels to the stochastic 
analysis the on Wiener space. 
\medskip 

{\bf Quasi-invariance of measures along the flow of particular 
non-linear PDEs} is a relevant topic in contemporary analysis. For 
example for several types of non-linear Schr\"odinger equations, 
quasi-invariance of certain Gaussian measures has been demonstrated in 
the recent works \cite{DT21}, \cite{FT19}, \cite{OST18}, \cite{OT17}, 
and \cite{PTV20}. In addition, the paper \cite{Tz15} addresses to the 
problem of quasi-invariance along flows generated by certain non-linear 
Hamiltonian PDEs such as the Benjamin-Bona-Mahony equation. In 
\cite{GOTW22}, \cite{OT20}, \cite{STX20}, and the references therein, 
the problem is treated relative to different forms of the defocusing 
wave equation. Here, the solutions to the corresponding PDEs are paths 
taking values in certain function spaces.  

We add the two applications of the theory of the present paper to 
the list of flows along which quasi-invariance of certain not 
necessarily Gaussian measures has been established. These are the 
flows generated by a class of Boltzmann type equations and by a 
PDE arising from a certain Fleming-Viot type particle system. 
\medskip 

{\bf The Boltzmann type equation} in its integrated (mild) form,  
\begin{eqnarray*} 
p_t=S(t)\, p_0+\lambda\int_0^t S(t-s)\, Q(p,p)(s)\, ds\, ,\quad 
t\in [0,T], 
\end{eqnarray*} 
is one of the two applications investigated in the paper. Here $p_0$ 
is a probability density at time zero, $S(t)$, $t\ge 0$, denotes the 
semigroup describing the dynamics with no collisions inside the 
physical space $\Omega$, and $Q$ is the collision kernel. The solution 
$p\equiv p(r,v,t)$ is interpreted as a particle density at the physical 
state $r\in\Omega$, the velocity $v$ belonging to the velocity space 
$V$, and the time $t\in [0,T]$. The parameter $\lambda>0$ is meaningful 
for the construction of local solutions for given probability densities 
$p_0$ belonging to certain open subsets of $L^1(\Omega\times V)$. In 
order to distinguish it from the measure $\lambda$ of Sections 
\ref{sec:1} and \ref{sec:2} we may call it parameter $\lambda$ in 
the Boltzmann type equation. 

After introducing the particular Boltzmann type equation, we collect 
the results of \cite{Lo20} which we use in the present paper, cf. 
Subsection \ref{sec:3:1} and Proposition \ref{Proposition3.2} (a) of 
Subsection \ref{sec:3:2}. Here, we refer to local and global existence 
and uniqueness of solutions to the Boltzmann type equation which, among 
other things, lead to injectivity of the map $p_0\cdot\l\mapsto p_t 
\cdot\l$ and verification of Case 2, see Lemma \ref{Lemma3.4} 
(a),(c) and its proof. In Subsection \ref{sec:3:2} we examine the 
map 
\begin{eqnarray*}
\Psi(p_0,q)(t):=S(t)\, p_0+\lambda\int_0^t S(t-s)\, Q(q,q)\, (s) 
\, ds\, ,\quad t\in [0,T],     
\end{eqnarray*} 
in order to derive some regularity properties of solutions to 
Boltzmann type equations. In particular, we verify Fr\'echet 
differentiability of the solution map $p_0\mapsto p\equiv (p_t 
)_{t\in [0,T]}$ in Proposition \ref{Proposition3.3} (c) and 
Fr\'echet differentiability of the map $p_0\mapsto p_t$ in the 
proof of Corollary \ref{Corollary3.5}. We provide an explicit 
representation of the Fr\'echet derivative of the solution map 
$p_0\mapsto p$, 
\begin{eqnarray}\label{R.3}
\nabla p(\cdot )=\left(\sum_{k=0}^\infty\left(\nabla_2\Psi\left(p_0,p  
\vphantom{l^2}\right)\right)^k\circ\nabla_1\Psi\left(p_0,p\vphantom 
{l^2}\right)\right)(\cdot ) 
\end{eqnarray} 
where the exponent $k$ refers to $k$-fold composition and the right-hand side 
converges absolutely in the operator norm. Here, $\nabla_1$ is the Fr\'echet 
derivative of the map $p_0\mapsto\Psi(p_0,q)$ and $\nabla_2$ is the Fr\'echet 
derivative of the map $q\mapsto\Psi(p_0,q)$. To the author's best knowledge, 
Fr\'echet derivatives of solution maps of similar Boltzmann type equations 
have not been determined before. 

In Subsection \ref{sec:3:3} we are concerned with compatibility with the 
general framework of Sections \ref{sec:1} and \ref{sec:2}. In particular 
Fr\'echet differentiability in the sense of Proposition \ref{Proposition2.4} 
is derived in order to verify conditions (jv) and (v) of Subsection 
\ref{sec:2:12}. In other words, {\it Quasi-invariance in the sense of 
the results of Subsection \ref{sec:2:12} is obtained}. This is summarized 
in Remark (2) of Section \ref{sec:3}. 

We mention that the verification of conditions (jv) and (v) of Subsection 
\ref{sec:2:12} is a result concerning the sensitivity of the trajectory 
$t\mapsto p_t\cdot\l$ of a Boltzmann type equation with respect to the 
initial datum. Such a result is of independent interest. Sensitivity of 
solution paths relative to initial values has been investigated for a 
variety of PDEs. For a comprehensive overview we refer to \cite{Ko19}. 
\medskip 

{\bf The asymptotic particle density of a Fleming-Viot type system} 
generates a flow under which we aim to state quasi-invariance in 
the sense of Section \ref{sec:2}. It should be mentioned that the 
extensive and comprehensive examination of the underlying particle 
system and similar ones is an actual topic in probability and 
analysis for about two decades. For this let us refer to \cite{BHM00}, 
\cite{GK04}, \cite{Lo051}, \cite{FM07}, \cite{Vi11}, \cite{Vi14}, 
\cite{LPR18}, \cite{CDGR20}, \cite{JM22} and the sources cited in these papers. 
This system has its roots in the sciences, cf. for example \cite{CBH05}. 

The particle density we investigate the present paper is the 
$L^1$-normalized solution to the heat equation with Dirichlet boundary 
conditions. For a bounded domain $D\subset {\Bbb R}^d$ with smooth boundary 
$\partial D$ and $h\in L^2(D)$ with $h\ge 0$ and $\|h\|_{L^1(D)}=1$ let 
$g_t$ denote the unique solution to $\frac12\Delta g=\frac{\partial} 
{\partial t}g$ with boundary condition $g_t=0$ on $\partial D$, $t\ge 0$, 
and initial datum $g_0=h$. Introduce $z(t):=\int_{y\in D}g_t(y)\, dy$. 
Then $h_t:=g_t/z(t)$ satisfies uniquely 
\begin{eqnarray*} 
\left\{
\begin{array}{l}
\D\frac{\partial}{\partial t}h_t(x)=\frac12\Delta h_t(x)-\frac{z'(t)}{z(t)} 
h_t(x)  \\
\D h_t(x)|_{x\in \partial D}=0,\ h_0=h 
\vphantom{\int}
\end{array} 
\right.
\end{eqnarray*}
where we keep in mind that $z$ depends on the initial datum $h$. We 
consider the map $\nu\equiv h\cdot\l\mapsto h_t\cdot\l\equiv\nu_t$ where 
$\l$ denotes the Lebesgue measure on $(D,{\cal B}(D))$. As above, let 
$\bmu$ be a probability measure on the metric space $E$ of all probability 
measures over $(D,{\cal B}(D))$ endowed with  the Prokhorov distance. 
Subsection \ref{sec:4:1} is devoted to compatibility with Sections 
\ref{sec:1} and \ref{sec:2} which includes {\it quasi-invariance in the sense 
of the results of Subsection \ref{sec:2:12}}. The latter is is summarized 
in Remark (4) of Section \ref{sec:4}. However, we are also interested in a 
more detailed analysis of the map $\nu\mapsto\nu_t$ for which we apply the 
quasi-invariance result. 

For $f\in L^2(E,\bmu)$ we introduce $T_tf(\nu):=f(\nu_t)$. In Subsection 
\ref{sec:4:2} we first examine properties of the infinitesimal generator 
$(A,D(A))$ of the semigroup $T_t$, $t\ge 0$. Representations of $Af$ are 
provided, for example for cylindrical functions of the form $f(\nu)=\vp 
((h_1,\nu),\ldots,(h_r,\nu))$, $\nu\in E$, $\vp\in C_b^1({\Bbb R}^r)$, $r 
\in {\Bbb N}$. Here $h_1,h_2,\ldots\, $ are eigenfunctions of the Dirichlet 
Laplacian on $D$ corresponding to the eigenvalues $0>2\lambda_1>2\lambda_2 
\ge\ldots\ $. Indicating the dependence of $z$ on the initial value by 
using the notation $z\equiv z(\nu,t)$, for such cylindrical functions it 
holds that
\begin{eqnarray*}
Af=\sum_{i=1}^r\frac{\partial\vp}{\partial x_i}\cdot (h_i,\cdot)\left( 
\lambda_i-z'(\cdot ,0)\right)\, .  
\end{eqnarray*}
As application of the quasi-invariance results of Subsection 
\ref{sec:2:12}, in Subsection \ref{sec:4:2} we are mainly interested 
in establishing integration by parts relative to the infinitesimal 
generator $(A,D(A))$. It is shown in Theorem \ref{Theorem4.6} that 
\begin{eqnarray}\label{R.4} 
\langle -Af,g\rangle_{L^2}+\langle -Ag,f\rangle_{L^2}=\left\langle 
\delta (A^f)\cdot f,g\right\rangle_{L^2}+\int fg\, d\bmu_\gamma\, , 
\quad f,g\in D(A). 
\end{eqnarray} 
\medskip 

The notation used in the paper follows the common standard. However since 
we just consider flows with time independent generators we ease the notation 
of flows by just giving the time difference as subindex or argument. The 
superscript $+$ indicates the general direction in time. 
\medskip 

Although the labels of equations, theorems, propositions, and lemmas 
refer to the sections where they are introduced, the labels of frequently 
used conditions do not. Labels of conditions are placed on the left-hand 
side. The reader may use the following list to locate them. 
\begin{itemize}
\item{} (a1), (a2) in Paragraph \ref{sec:1:1:5}, (a3) in Paragraph \ref{sec:1:1:6} 
\item{} (i), (ii), (iii), (a4) in the introductory text to Section \ref{sec:2} 
\item{} (j'), (jj'), (jjj') in Paragraph \ref{sec:2:2:5} 
\item{} (jv), (jv'), (v), (v') in the beginning of Subsection \ref{sec:2:12} 
\item{} (vj), (vj') in Theorem \ref{Theorem2.7} of Subsection \ref{sec:2:12} 
\item{} (jjj''), (vjj) in Subsection \ref{sec:4:2} 
\item{} (k), (kk), (kkk), (kw), (w), (wk), (wkk) in the introductory text to Section \ref{sec:3} 
\end{itemize}
%

\subsection{Related Work}\label{sec:1:2}

The change of measure formula (\ref{1.5}) is known from \cite{SvW93}, 
\cite{Be06}, and \cite{Bo10}. Furthermore, it appears already in the 
note \cite{DS88} and is stated in \cite{BM-W99}. In a Gaussian setting, 
formula (\ref{1.5}) has been established in \cite{Bo97}, \cite{Cr83}, 
and \cite{UZ00}. In this subsection we shall discuss why the hypotheses 
under which formula (\ref{1.5}) is derived in these references are not 
or not directly verifiable under the assumptions of the present paper.  
\medskip 

{\bf References \cite{Cr83} and \cite{UZ00}. } 
These works use arguments from Wiener space calculus and stochastic 
analysis. The constructions lead to to the formula (\ref{1.5}) relative 
to the change of the Wiener measure along differentiable flows with 
increments in the Cameron-Martin space. 
\medskip

{\bf Reference \cite{Bo97}. } 
In the survey paper \cite{Bo97}, the theorems Theorem 9.2.2 -- Theorem 
9.2.4 focus on Gaussian measures $\bmu$. The conditions under which the 
theorems are derived refer to this concept. However the proof of Theorem 
9.2.4 provides one of the main ideas how to prove our Theorem 
\ref{Theorem2.3}, namely to establish a first order ODE of the form 
(\ref{2.19}) whose solution is (\ref{1.5}). 
\medskip

{\bf References \cite{Bo10} and \cite{SvW93}. } 
A major part of Subsection 11.2 in \cite{Bo10} is based on the paper 
\cite{SvW93}. A class of topologies is introduced under which the 
family of measures (in our terminology) $\bmu\circ\nu_t$ is assumed to 
be differentiable with respect to $t$. In addition, Theorem 11.2.7 in 
\cite{Bo10} and Theorem 3.3 in \cite{SvW93} require a certain integrability 
of the derivative $\frac{d}{dt}\bmu\circ\nu_t$ in advance, while Theorem 
11.2.13 in \cite{Bo10} assumes a certain invariance of a class of test 
functions under the map $\nu\to\nu_t$.  

Here in the present paper, the only form of differentiability of the 
measure $\bmu$ under which we will follow the flow is (jj) in Subsection 
\ref{sec:2:1}. It will take an effort to establish the derivative 
$\frac{d}{dt}\bmu\circ\nu_t$, cf. Step 3 in the proof of our Theorem 
\ref{Theorem2.3}. Part (b) of our Theorem \ref{Theorem2.3} together with  
condition (jj) of Section \ref{sec:2}, namely $\delta (A^f)\in L^1(E, 
\bmu)$, says that $t\mapsto\bmu\circ\nu_t$ is, in the terminology of 
\cite{SvW93} and \cite{Bo10}, $\tau_{{\cal B}_b}$-differentiable. In 
other words, to some extend we get compatibility with \cite{SvW93} and 
\cite{Bo10}, Subsection 11.2, as one result of our work. 
\medskip

{\bf Reference \cite{DS88}. } 
Proposition 1 and Theorem 1 in the paper \cite{DS88} suppose that $t 
\mapsto\bmu\circ\nu_t$ is continuously $\tau_{{\cal B}_b}$-differentiable. 
This is what we obtain but not assume. The change of measure formula in 
Theorem 2 of \cite{DS88} is not compatible with our formula (\ref{1.5}). 
This formula investigates change of measure under a certain class of 
maps in the state space independent of time. Theorem 2 of \cite{DS88} 
is remarkable since the crucial condition is Fr\'echet differentiability 
of a map similar to our $\nu\mapsto\nu_t$. 
\medskip

{\bf Reference \cite{BM-W99}. } 
The first goal of the paper \cite{BM-W99} is to give conditions that 
provide the existence of flows along a certain class of vector fields 
over a Hausdorff locally convex space $X$. The second goal is the 
quasi-invariance of a class of differentiable measures under such 
flows with state space $X$. For these purposes a sophisticated set of 
conditions is given referring particularly to a certain finite 
dimensional projection technique which is applied. The construction 
of the flows and the proof of quasi-invariance are carried out 
simultaneously. In particular it is shown that the Radon-Nikodym 
process has a measurable modification with continuous paths.

The flows exhibit a certain property, (4.9) in \cite{BM-W99}, which 
for $t=0$ is comparable with our Theorem \ref{Theorem2.3} (b). Based 
on (4.9) in \cite{BM-W99} and the already proved quasi-invariance, 
including the just mentioned property of the Radon-Nikodym process, 
the formula (\ref{1.5}) is derived. 

In contrast, we assume the existence of a flow which is, for example, 
the flow relative to the solution map of some non-linear PDE modeling 
a particle density as in Sections \ref{sec:3} and \ref{sec:4}. In the 
present paper, the state space of the flow is some subset of $E$, the 
space of all probability measures on the separable locally compact 
metric space $D$. Since $E$ is not a locally convex space, our 
framework is not directly compatible with that of \cite{BM-W99}. 

Moreover, once we have demonstrated quasi-invariance including 
measurability and continuity of the Radon-Nikodym process by 
the construction  (\ref{2.16}) and (\ref{2.17}) below, the proof of 
the formula (\ref{1.5}) is quite similar to the proof of Theorem 
4.4 (i) in \cite{BM-W99}. In \cite{BM-W99}, the underlying idea is 
referred to as {\it Bell's method}. 

\medskip
{\bf Reference \cite{Be85}. } 
In a simpler linear setting formula (\ref{1.5}) has already been 
derived in \cite{Be85}. D. Bell's paper \cite{Be85} may have been 
the first contribution toward the formula (\ref{1.5}). 

\medskip
{\bf Reference \cite{Be06}. } 
In the new appendix to this reprint from 1987, D. Bell derives 
(\ref{1.5}) if the maps (in our notation) $\nu\mapsto\nu_t$, $t\in 
{\Bbb R}$, either act on smooth finite-dimensional compact manifolds 
without boundary or in  Banach spaces. He uses assumptions as for 
example continuous (Fr\'echet) differentiability of the maps $\nu 
\mapsto\nu_t$ and differentiability of the vector field (in our 
notation) $A^f\cdot$ in order to establish a first order ODE for 
${\Bbb R}\ni t\mapsto d\bmu\circ\nu_{-t}/d\bmu$. Strengthening the 
setup by more conditions on the vector field $A^f\cdot$ and its 
divergence $\delta(A^f)$, the solution to this ODE turns out to be 
(\ref{1.5}). 

Determined by the applications, our approach has to use assumptions 
on the maps $\nu\mapsto\nu_t$ and the vector field $A^f\cdot$ that are 
weak in comparison with the ones of D. Bell. For this see Subsections 
\ref{sec:2:1}, \ref{sec:2:12} below and especially Propositions 
\ref{Proposition2.4} as well as \ref{Proposition2.5}. In particular  
differentiability of the vector field $A^f\cdot$ as used in \cite{Be06} 
excludes differential expressions as in our applications of Sections 
\ref{sec:3} and \ref{sec:4} below.

\subsection{Notation and Basic Hypotheses}\label{sec:1:1}

Let $(D,\rho)$ be a separable locally compact metric space. Let 
${\cal B}$ denote the Borel $\sigma$-algebra on $D$. Let $C_0(D)$ 
denote the space of all continuous functions on $D$ such that for 
each $\ve >0$ the set $\{x\in D:|f(x)|\ge\ve\}$ is compact. One  
reason for the particular choice of $(D,\rho)$ is that, in this  
setting, $C_0(D)$ is separable. Another reason is that $(D,\rho)$ 
is $\sigma$-compact. See also Remarks (1) and (2) of Section 
\ref{sec:2}. 

Let $\lambda$ be a finite measure on $(D,{\cal B})$. If not 
otherwise stated, a {\it measure} is always non-negative. Denote 
$L^v(D)\equiv L^v(D,{\cal B},\lambda)$, $1\le v\le\infty$. Let $F$ 
denote the set of all finite signed measures on $(D,{\cal B})$ endowed 
with the topology of {\it narrow convergence}. That is, a sequence 
$\nu_n\in F$, $n\in {\Bbb N}$, is said to (narrowly) converge to $\nu 
\in F$ if 
\begin{eqnarray*} 
\lim_{n\to\infty}\int f\, d\nu_n=\int f\, d\nu\quad\mbox{\rm for all } 
f\in C_b(D) 
\end{eqnarray*} 
where $C_b(D)$ is the set of all bounded and continuous real functions 
$f$ on $(D,\rho)$. To avoid confusion we will solely use the notation 
narrow convergence in the paper, although for probability measures 
$\nu_n$, $n\in {\Bbb N}$, and $\nu$, the above definition is often 
referred to as {\it weak convergence}, see \cite{Re13}. 

Let $E\subset F$ denote the subset of all probability measures, 
and ${\cal B}(E)$ the Borel $\sigma$-algebra on $E$ relative to 
the Prokhorov metric $\pi$ on $E$. Since $(D,\rho)$ is a separable 
metric space, convergence in $(E,\pi)$ is equivalent to narrow 
convergence of probability measures, cf. the Portmanteau theorem 
in the form of \cite{EK86}, Theorem 3.3.1. 
\medskip

For $\nu\in F$ with $\nu=h\cdot\lambda$ and $h\in L^v(D)$ for some 
$1\le v\le\infty$ and $g\in L^w(D)$ where $1/v+1/w=1$ we may write 
$(\nu,g)$, $(g,\nu)$, $(h,g)$, or $(g,h)$ instead of $\int g\, d\nu$ 
or $\int gh\, d\lambda$. Here we set $w:=\infty$ if $v=1$. The following 
is a useful detail for the setup of the paper. It will be proved in 
Subsection \ref{sec:1:4}. 
\begin{lemma}\label{Lemma1.0} 
The sets $\, {\cal S}^v:=E\cap\{h\cdot\lambda:h\in L^v(D)\}$, $1\le v< 
\infty$, belong to ${\cal B}(E)$. 
\end{lemma}

\subsubsection{Formal definition of flows and trajectories generalizing 
solutions to PDEs - a motivation}\label{sec:1:1:1} Below, we introduce 
two flows, a function valued one $V_t^+$, $t\in [0,S+1]$, and a measure 
valued one $W^+(t,\cdot)$, $t\in [0,S+1]$, where $S\ge 1$. 

Consider a PDE of type $\frac{d}{dt}h_t=Ah_t$ where $A$ is an operator 
independent of the time $t\ge 0$. Assume that it is known that this PDE 
has a unique solution on some time interval $t\in [0,S+1]$ for every 
initial value $h_0=h$ belonging to a certain set ${\cal V}$. Let such a 
solution be denoted by $h_t\equiv\t V^+_t(h)$, $t\in [0,S+1]$. 

At some fixed time $t\in [0,S+1]$ there may exist more solution paths 
of the PDE than just those given by $\{h_s\equiv\t V^+_s(h):h\in {\cal 
V},\, s\in [0,t]\}$. Indeed, if $\{h_u\equiv\t V^+_u(h):u\le S+1-t,\ h 
\in {\cal V}\}$ is not a subset of ${\cal V}$ then at time $t$ there 
may also exist solution paths with initial values belonging to $\{h_u 
\equiv\t V^+_u(h):u\le S+1-t,\ h\in {\cal V}\} \setminus {\cal V}$. 
Those solution paths may not exist on the whole time interval $[0,S+1]$. 
But at least they exist on the time interval $[0,t]$.

To indicate that we now also consider solution paths $h_u$ that may 
just exist on time intervals $[0,t]$ with $t\in (0,S+1]$ we drop the 
$~\t{}~$ from the notation, i. e. we use the symbol $V^+_u(h)$.  

In order to state a flow property of type $V_{t}^+\circ V_{s}^+ 
=V_{t+s}^+$, $s,t\in [0,S+1]$ with $s+t\in [0,S+1]$, we need to 
carefully define the domains of the operators $V_t^+$, $t\in 
[0,S+1]$. The just outlined observation is the motivation for the 
definition of $D(V_t^+)$, $t\in [0,S+1]$, below. See also Remark 
(1) below. Adequately we define a measure valued flow by $W^+(t, 
\nu):=V^+_t(h)\cdot\lambda\, ,\quad t\in [0,S+1]$, $\nu=h\cdot 
\lambda$ and select the respective domains of definition. 
\medskip

Furthermore, we recall that, at time $t=1$, we observe or wish 
to control a probability measure $\bmu$ over probability measures 
$V^+_1(h)\cdot \lambda$. Recall in addition that the 
quasi-invariance formula of Theorem \ref{Theorem2.3} looks backward 
in time. In contrast, to prove this formula, one must also look 
forward in time, cf. conditions (jv) and (v) or the calculations in 
Step 1 of the proof of Theorem \ref{Theorem2.3}. This gives rise  
to introduce the trajectory $\nu_t\equiv U^+(t,\nu)$ on reasonably 
selected domains $D(U^+(t,\cdot))$, $t\in [-1,1]$, following $W^+( 
1,\cdot)$ back and forth $t$ units of time. With this understanding, 
we are examining probability measures $\bmu$ over the set of 
probability measures $\nu\equiv\nu_0\equiv U^+(0,\nu)$ at time $t=0$. 

\subsubsection{The flow $V_t^+$, $t\in [0,S+1]$}\label{sec:1:1:2}  In 
this paragraph we provide a certain formalization / generalization 
of the flow generated by a non-linear PDE of the form $\frac{d}{dt}h_t 
=Ah_t$ where $A$ is a time independent operator.  

Let $1\le v<\infty$, $1\le S<\infty$, and $r>1$. Denote $B^v(r):= 
\{f\in L^v(D):\|f\|_{L^v(D)}<r\}$. Suppose we are given an open 
set ${\cal V}\subseteq L^v(D)$ with $B^v(r)\subseteq {\cal V}$ and, 
for all $t\in [0,S+1]$, an injective map $\t V_t^+:{\cal V}\mapsto 
L^v(D)$. Suppose also that, for all $t\in [0,S+1]$, there is a 
unique injective extension $V_t^+:D(V_t^+):=\bigcup_{s\in [0,S+1-t 
]}\{\t V_s^+(h):h\in {\cal V}\}\mapsto L^v(D)$ of $(\t V_t^+,{\cal 
V})$ such that 
\begin{itemize}
\item[$\bullet$] $V^+_0$ is the identity in $D(V_0^+)$ and 
$\vphantom{\D\int}$
\item[$\bullet$] for all $s,t\in [0,S+1]$ with $s+t\in [0,S+1]$ and 
$h\in {\cal V}$, which yields $\t V_s^+(h)\in D(V_t^+)$, 
we have 
\begin{eqnarray}\label{1.1} 
V_t^+\circ\t V_s^+(h)=\t V_{s+t}^+(h)\, .
\end{eqnarray} 
\end{itemize} 
We note that 
\begin{eqnarray}\label{1.2} 
D(V^+_s)\supseteq D(V^+_t)\supseteq D(V^+_{S+1})={\cal V}\supseteq 
B^v(r)\, ,\quad 0\le s\le t\le S+1.
\end{eqnarray} 
%
By (\ref{1.1}) we have for $s\in [0,S+1]$ with $s+t\in [0,S+1]$ 
\begin{eqnarray}\label{1.2.1} 
\{h\in D(V_s^+):V_s^+(h)\in D(V_t^+)\}\supseteq D(V^+_{s+t})\, . 
\end{eqnarray} 
Any $\t h\in D(V_s^+)$ has the representation $\t h=\t V_u^+(h)$ for 
some $h\in {\cal V}$ and some $u\in [0,S+1]$ with $u+s\in [0,S+1]$. 
From (\ref{1.1}) we obtain for all $u,s,t\in [0,S+1]$ with $u+s+t\in 
[0,S+1]$ and $\t h=\t V_u^+(h)\in D(V_{s+t}^+)$ for some $h\in {\cal 
V}$ the identity 
\begin{eqnarray}\label{1.2.2}
V_t^+\circ V_s^+(\t h)=\t V_{u+s+t}^+(h)=V_{s+t}^+(\t h)\, . 
\end{eqnarray} 
Conversely, relations (\ref{1.2.1}) and (\ref{1.2.2}) {\it define the 
type of function valued flow}  $\left(V^+_t,D(V^+_t)\right)$, $t\in 
[0,S+1]$, we are dealing with in the paper. 
\medskip 

\nid 
{\bf Remark} (1) For a more figurative description of the set 
$D(V_t^+)$, $t\in [0,S+1]$, we look at all solution paths $\t V_s^+h$, 
$s\in [0,S+1]$, with values in $L^v(D)$ where $h\in {\cal V}$. 
Then $D(V_t^+)$ is the set of all functions $\t V_s^+h\in L^v(D)$ 
such that the corresponding function valued solution path is, after 
time $s$, at least $t$ units of time ``alive" before terminating 
at time $S+1$. 

\subsubsection{The flow $W^+(t,\cdot)$, $t\in [0,S+1]$, and the trajectory 
$U^+(t,\cdot)$, $t\in [-1,1]$}\label{sec:1:1:3}  The goal is to follow 
non-random measure valued trajectories $\nu_t\equiv U^+(t,\nu)$, $t\in 
[-1,1]$, describing the dynamics of certain phenomena over time. 
The time range has been chosen to be $t\in [-1,1]$ since a change 
of measure formula will be established which naturally requires 
to look backward in time. 
\medskip

Set $D(W^+(t,\cdot)):=\{h\cdot\lambda:h\in D (V_t^+)\}$, $t\in [0,S+1]$, 
and introduce $W^+(t,\cdot)$ by 
\begin{eqnarray*}
W^+(t,\nu):=V^+_t(h)\cdot\lambda\, ,\quad \nu\equiv h\cdot\lambda\in 
D(W^+(t,\cdot)),\  t\in [0,S+1]. 
\end{eqnarray*} 
Since $V_t^+:D(V_t^+)\mapsto L^v(D)$ is injective, $W^+(t,\cdot): 
D(W^+(t,\cdot))\mapsto \{h\cdot\lambda:h\in L^v(D)\}$ is injective 
as well for all $t\in [0,S+1]$. In addition, if $s,t\in [0,S+1]$ with 
$s+t\in [0,S+1]$ we have the relation 
\begin{eqnarray}\label{1.2.3} 
&&\left\{\nu\in D(W^+(s,\cdot)):W^+(s,\nu)\in D(W^+(t,\cdot))\right\}
\supseteq D\left(W^+(s+t,\cdot)\right)\qquad
\end{eqnarray} 
and for all $\nu\in D(W^+(s+t,\cdot))$ it holds that 
\begin{eqnarray}\label{1.2.4} 
W^+(t,\cdot)\circ W^+(s,\nu)=W^+(s+t,\nu)\, . 
\end{eqnarray} 

Introduce $D(U^+(t,\cdot)):= \{W^+(1,\nu):\nu\in D(W^+(2,\cdot))\}$, $t\in 
[-1,1]$, and 
\begin{eqnarray}\label{1.2.5} 
U^+(t,\nu):= W^+(t,\nu)\, ,\quad \nu\in D(U^+(t,\cdot))\, ,\quad t\in [0,1].
\end{eqnarray} 
Recalling the just mentioned injectivity, for $t\in [-1,0]$ and $\nu\in 
D(U^+(t,\cdot))$, let $U^+(t,\nu)$ be the unique element in $D(W^+(-t, 
\cdot))$ satisfying 
\begin{eqnarray}\label{1.2.6} 
W^+(-t,\cdot)\circ U^+(t,\nu)=\nu\, . 
\end{eqnarray} 

We remark that relations (\ref{1.2.1})-(\ref{1.2.6}) are the {\it defining 
properties of the type of measure valued trajectory}  $\left(U^+(t,\cdot), 
D(U^+(t,\cdot))\right)$, $t\in [-1,1]$, we are working with in the paper. 
Below, we will also use the notation $\nu_t\equiv U^+(t,\nu)$, $\nu\in 
D(U^+(t,\cdot))$, $t\in [-1,1]$, and $\nu_t\equiv W^+(t,\nu)$, $\nu\in 
D(W^+(t,\cdot))$, $t\in [0,S+1]$ where we stress compatibility. 

\begin{lemma}\label{Lemma1.2} 
(a) Let $v=1$. For any sequence $m\equiv(m_1,m_2,\ldots\, )$ of 
strictly decreasing positive real numbers with $\lim_{n\to\infty}m_n 
=0$, there exists a set ${\cal S}_m^1\subseteq E\cap\{h\cdot\lambda: 
h\in L^1(D)\}$ which is compact in $(E,\pi)$ with the following 
properties. We have ${\cal S}_m^1\subseteq {\cal S}_{\t m}^1$ if 
$\, \t m_n\le m_n$ for all $n\in {\Bbb N}$. Moreover, 
\begin{eqnarray*}
{\cal S}^1&&\hspace{-.5cm}=E\cap\{h\cdot\lambda:h\in L^1(D)\} \\ 
&&\hspace{-.4cm}\equiv\left\{h\cdot\lambda:h\in L^1(D),\ \left\|h 
\right\|_{L^1(D)}=1,\ h\ge 0\right\}=\bigcup_{\mbox{\rm\footnotesize 
all sequences $m$}}{\cal S}_m^1\, .
\end{eqnarray*}
Furthermore, for any finite measure $\bnu$ on $(E,{\cal B}(E))$ with 
$\bnu({\cal S}^1)=\bnu(E)$ and any $\ve >0$ there exists a sequence 
$m(\ve)$ with the above properties satisfying 
\begin{eqnarray*}
\bnu\left({\cal S}_{m(\ve)}^1\right)\ge\bnu(E)-\ve. 
\end{eqnarray*}
(b) Let $1<v<\infty$. There exists an increasing sequence ${\cal 
S}_m^v$, $m\in {\Bbb N}$, of compact subsets of $(E,\pi)$ such that 
\begin{eqnarray*}
{\cal S}^v&&\hspace{-.5cm}=E\cap\{h\cdot\lambda:h\in L^v(D)\} \\ 
&&\hspace{-.4cm}\equiv\left\{h\cdot\lambda:h\in L^v(D),\ \left\|h 
\right\|_{L^1(D)}=1,\ h\ge 0\right\}=\bigcup_{m\in {\Bbb N}}{\cal 
S}_m^v\, .
\end{eqnarray*}
For any finite measure $\bnu$ on $(E,{\cal B}(E))$ with $\bnu({\cal 
S}^v)=\bnu(E)$ and any $\ve >0$ there exists $m(\ve)\in {\Bbb N}$ with 
the above properties such that 
\begin{eqnarray*}
\bnu\left({\cal S}_{m(\ve)}^v\right)\ge\bnu(E)-\ve. 
\end{eqnarray*}
\end{lemma} 

For the proof see Subsection \ref{sec:1:4}. As a consequence, each 
${\cal S}_m^v$ is closed in $(E,\pi)$; note that, for $v=1$, $m$ is 
a certain sequence and, for $1<v<\infty$, we have $m\in {\Bbb N}$. 
This yields 
\begin{eqnarray*}
{\cal S}^v_m\in {\cal B}(E)\quad\mbox{\rm for each }m
\end{eqnarray*}
and $1\le v<\infty$. Recall also Lemma \ref{Lemma1.0}. 

\subsubsection{State space and direction of inclusion} 
\label{sec:1:1:5} Let $1\le v<\infty$. Identify $\{V_1^+(h):h\in 
D(V_2^+)$ with $V_t^+h\ge 0$ and $\|V_t^+h\|_{L^1(D)}=1$ for all 
$t\in [0,2]\}$ with ${\cal G}:=\{\nu\in D(U^+(0,\cdot)):\nu_t\in 
E$ for all $t\in [-1,1]\}$. Choose the {\it state space} $G\subseteq 
{\cal G}$ with $G\in {\cal B}(E)$. 
\medskip

The following hypotheses will be important in the paper. We emphasize 
that $G\subseteq {\cal G}$ implies that $\nu_{-t}\in E$, $t\in [-1,1]$, 
whenever $\nu\in G$. 
\begin{itemize}
\item[(a1)] If $B\in{\cal B}(E)$, $B\subseteq G$, and $t\in [-1,1]$ 
suppose $\{\nu_{-t}:\nu\in B\}\in {\cal B}(E)$. Furthermore, assume 
for $t\in [-1,1]$ that the map $G\ni\nu\mapsto\nu_{-t}$ is ${\cal B} 
(E)$-measurable. For $\nu\equiv h\cdot\lambda\in G$ the function 
$[-1,1]\ni t\mapsto h_{-t}$ is continuous in $L^v(D)$.
\end{itemize}

We consider the following two types of maps. 
\medskip 

In \textbf{\textit{ Case 1}} we suppose in addition to (a1) that for 
all $t\in [0,1]$ the map $G\ni\nu\mapsto\nu_{-t}$ is injective and we 
have $\{\nu_{-t}:\nu\in G\}\subseteq G$. In \textbf{\textit{ Case 2}} 
we assume in addition to (a1) that for all $t\in [0,1]$ the map $G\ni 
\nu\mapsto\nu_t$ is injective and it holds that $\{\nu_t:\nu\in G\} 
\subseteq G$. We note that Case 1 and Case 2 are not necessarily 
mutually exclusive. 
\medskip 

For $t\in [-1,1]$ introduce the set $G_t:=\{\nu_t:\nu\in G\}$. 
Observe that in Case 1 the family of maps $G\ni\nu\mapsto\nu_{-s}$, 
$s\in [0,1]$, extends by iteration to $s\in [0,\infty)$. The set  
$G_{-s}:=\{\nu_{-s}:\nu\in G\}$ is therefore well-defined for $s\in 
[-1,\infty)$. We have $G_{-s}\subseteq G$ bijectively for $s\in 
[0,\infty)$. This and the injectivity of $V_\cdot^+$ supposed in 
Paragraph \ref{sec:1:1:2} implies even $G_{-s_2}\subseteq G_{-s_1}$ 
for $-1\le s_1\le s_2$. For $t\in [-1,\infty)$ introduce $G^-_{-t} 
:=\bigcup_{s\in (t,\infty)}G_{-s}$. By (a1) and, for $t>1$ by its 
iterated application, we have $G_{-t}\in {\cal B}(E)$ for all $t\in 
[-1,\infty)$. The last two properties yield $G^-_{-t}=\bigcup_{n\in 
{\Bbb N}}G_{-t-\frac1n}\in {\cal B}(E)$ for all $t\in [-1,\infty)$. 

In the same way we introduce in Case 2 the sets $G_s:=\{\nu_s:\nu\in 
G\}$ for $s\in [-1,\infty)$ and $G^+_t:=\bigcup_{s\in (t,\infty)} 
G_s$, $t\in [-1,\infty)$. Here we have $G_{s_1}\supseteq G_{s_2}$ for 
$-1\le s_1\le s_2$. Correspondingly, $G_t\in {\cal B}(E)$ and $G^+_t 
\in {\cal B}(E)$ for all $t\in [-1,\infty)$. 

We mention that $G_0=G$ and abbreviate the notation by $G^-\equiv 
G_0^-$ (Case 1) as well as $G^+\equiv G_0^+$ (Case 2). For the next 
hypothesis let us remind of the different meaning of the index $m$ 
in the notation of the set ${\cal S}^v_m$ depending on $v=1$ or 
$1<v<\infty$, cf. Lemma \ref{Lemma1.2}. 
\begin{itemize}
\item[(a2)] If $v=1$ then for any sequence of strictly decreasing 
positive real numbers $N\equiv (N_1,N_2,\ldots\, )$ with $N_1<1$ 
and $\lim_{n\to\infty}N_n=0$ there exists a closed set $G_N$ in 
$(E,\pi)$ such that $G_N\subseteq G_{\t N}$ if $\t N_n\le N_n$ for 
all $n\in {\Bbb N}$ and  
\begin{itemize}
\item[(1)] $G=\bigcup_{{\rm all}\ N} G_N$,
\item[(2)] for any index $N$, we have $G_N\subseteq {\cal S}^v_m$ 
for some index $m$, 
\item[(3)] for any index $N$ there is some $t_0\equiv t_0(N)>0$ and
there exists some index $\hat{N}\equiv \hat{N}(t_0,N)$ such that 
$\nu_t\in G_{\hat{N}}$ whenever $\nu\in G_N$ and $t\in [0,1+t_0]$ 
in Case 1, or $t\in [-(1+t_0),0]$ in Case 2, 
\item[(4)] for any index $N$ the map $G_N\ni\nu\mapsto\nu_t$ is 
continuous in $(E,\pi)$, in Case 1 for $t\in [-(1+t_0),0]$, in Case 
2 for $t\in [0,1+t_0]$.
\end{itemize}
Furthermore if $1<v<\infty$ then there exists a sequence of increasing 
closed sets $G_N$ in $(E,\pi)$, $N\in {\Bbb N}$, such that (1)-(4). 
Note the different meaning of $N$ and $m$ depending on $v=1$ or 
alternatively $1<v<\infty$. 
\end{itemize}

\subsubsection{The measure to follow along the trajectory $U^+$ } 
\label{sec:1:1:4} 

Let $\bmu$ be a probability measure on $(E,{\cal B}(E))$ such that 
the closure of $G$ in $(E,\pi)$ is supp $\bmu$. In particular, suppose 
\begin{eqnarray*}
\bmu(G)=1\, . 
\end{eqnarray*} 

By (a1), for $t\in [-1,1]$ the measure $\bmu\circ\nu_{-t}$ is well-defined 
on $(E,{\cal B}(E))$ by  
\begin{eqnarray*}
\bmu\circ\nu_{-t}(B):=\bmu(\{\nu_{-t}:\nu\in B\cap G\})\, ,\quad B\in 
{\cal B}(E),  
\end{eqnarray*}
which includes $\bmu\circ\nu_{-t}(E\setminus G)=0$. We mention also 
that $\bmu\circ\nu_{-t}(G)=1$ may not be true in Case 1 if $\, t\in 
(0,1]$, and in Case 2 if $\, t\in [-1,0)$. 
\medskip

Since in Case 1 we are interested in equivalence of the measures 
$\bmu\circ\nu_{-t}$, $t\in [0,1]$, it is necessary to assume that 
\begin{eqnarray*}
&&\hspace{-.5cm}\textstyle\bmu\circ\nu_{-t}(G\setminus G^-)=\bmu 
\left(\{\nu_{-t}:\nu\in G\setminus G^-\}\right)=\bmu\left(G_{-t} 
\setminus G_{-t}^-\right)
\end{eqnarray*}
is either zero for all $t\in [0,1]$ or positive for all $t\in [0,1]$. 
Here we recall $G\supseteq G^-$ from Paragraph \ref{sec:1:1:5}. This 
leads to the hypothesis 
\begin{eqnarray*}
\bmu\left(G_{-t}\setminus G_{-t}^-\right)=0\, ,\quad t\in [0,1],\ \mbox 
{\rm in Case 1}
\end{eqnarray*}
and similarly 
\begin{eqnarray*}
\bmu\left(G_t\setminus G_t^+\right)=0\, ,\quad t\in [0,1],\ \mbox{\rm 
in Case 2}. 
\end{eqnarray*}
{\bf Remark} (2) We use hypothesis (a2) in the proof of the second part 
of Lemma \ref{Lemma1.1} (a) and in the proof of Lemma \ref{Lemma1.1} (b). 
There we verify tightness of the set of measures $\{\bmu\circ\nu_{-t}:t\in 
[0,1]\}$ on $(G^-,\pi)$ in Case 1, and tightness of the set of measures 
$\{\bmu\circ\nu_t:t\in [0,1]\}$ on $(G^+,\pi)$ in Case 2. 

\subsubsection{The generator of the flow $W^+$ and the trajectory $U^+$} 
\label{sec:1:1:6} In this paragraph, for $\nu\in G$ we will use the 
notation $\nu_t\equiv h_t\cdot\lambda$, $t\in [-1,1]$. Let us assume 
\begin{itemize}
\item[(a3)] the existence of the limit 
\begin{eqnarray}\label{1.3} 
\hspace{-.3cm}\lim_{s\to 0}\frac1s\left(h_{s+t}-h_t\right)\ \mbox{\rm 
in }L^v(D)\quad \mbox{\rm for all $t\in [-1,1]$}\, ,\quad\mbox{\rm for 
$\bmu$-a.e. $\nu\equiv h\cdot\lambda\in G$}. 
\end{eqnarray} 
\end{itemize} 
Here we use the convention $s\ge 0$ if $t=-1$ and $s\le 0$ if $t=1$. 
For those $\nu\equiv h\cdot\lambda\in G$ we set $A^f\nu_t:=\lim_{s\to 0} 
\frac1s\left(h_{s+t}-h_t\right)\cdot\lambda$. The operator $A^f$ maps 
for $\bmu$-a.e. $\nu\in G$ and $t\in [-1,1]$, the element $\nu_t$ to 
some element $A^f\nu_t$ belonging to $F$. We remark that, according to 
the above definition of the flow $W^+$ and trajectory $U^+$, the operator 
$A^f$ is independent of $t\in [-1,1]$.

In other words, we shall assume that $\bmu$-a.e. 
\begin{eqnarray}\label{1.4} 
\frac{d}{dt}\nu_t\equiv\lim_{s\to 0}\frac1s\left(h_{s+t}-h_t\right) 
\cdot\lambda=A^f\nu_t\, ,\quad\nu\equiv h\cdot\lambda\in G,\ t\in 
[-1,1]. 
\end{eqnarray} 
The superscript $f$ in $A^f$ indicates the generator of the {\it 
f}~low $W^+$. Conversely, given a vector field $A^f\nu$, $\nu\in 
\bigcup_{t\in [-1,1]}\{\t\nu_t:\t\nu\in G\}$, and some initial 
value $\nu_{-1}\in \{\t\nu_{-1}:\t\nu\in G\}$ or $\nu\equiv\nu_0 
\in G$, relation (\ref{1.4}) typically becomes a {\it PDE} with 
solution $\nu_t$, $t\in [-1,1]$. 

The continuity in $L^v(D)$ of $[-1,1]\ni t\to h_t$ for $\nu_t\equiv 
h_t\cdot\lambda$, $\nu\in G$, postulated in hypothesis (a1) will be 
used in the proof of the following lemma, cf. Subsection \ref{sec:1:4}. 
We emphasize again that in Case 1, the family of maps $G\ni\nu 
\mapsto\nu_{-s}$, $s\in [0,1]$, can be extended by iteration to all 
$s\in [0,\infty)$. In the same way, in Case 2, the family of maps 
$G\ni\nu\mapsto\nu_s$, $s\in [0,1]$, can be extended to $s\in [0, 
\infty)$.
\begin{lemma}\label{Lemma1.1} 
Suppose (a1) and (a2). \\ 
(a) Let $\bnu$ be a finite measure on $(E,{\cal B}(E))$ with $\bnu 
(G)=\bnu(E)$. Then the compact sets ${\cal S}^v_m$ of Lemma 
\ref{Lemma1.2} can be chosen in such a way that for any $\eta>0$ 
there exists $m(\eta)$ with 
\begin{eqnarray*} 
\bnu\circ\nu_{-t}\left({\cal S}^v_{m(\eta)}\right)\equiv\bnu\left( 
\left\{\nu_{-t}:\nu\in {\cal S}^v_{m(\eta)}\cap G\right\}\right)\ge 
\bnu\circ\nu_{-t}(G)-\eta\quad\mbox{\rm for all}\quad t\in [-1,1] 
\, . \vphantom{\int} 
\end{eqnarray*}
In particular in Case 1, for any $\eta>0$ there exists an index 
$\hat{N}(\eta)$ and $t_0>0$ such that 
\begin{eqnarray*} 
\bnu\left(\left\{\nu_{-t}:\nu\in G_{\hat{N}(\eta)}\right\}\right) 
\ge\bnu\left(\left\{\nu_{-t}:\nu\in G\right\}\right)-\eta\quad 
\mbox{\rm for all}\quad t\in [0,1+t_0]. \vphantom{\int} 
\end{eqnarray*}
In Case 2, for any $\eta>0$ there exists an index $\hat{N}(\eta)$ 
and $t_0>0$ with 
\begin{eqnarray*} 
\bnu\left(\left\{\nu_t:\nu\in G_{\hat{N}(\eta)}\right\}\right)\ge 
\bnu\left(\left\{\nu_t:\nu\in G\right\}\right)-\eta\quad\mbox{\rm f
or all}\quad t\in [0,1+t_0]. \vphantom{\int}
\end{eqnarray*}
The statements of part (a) of the lemma hold in particular for $\bnu 
=\bmu$. \\ 
(b) Let $\eta>0$. In Case 1, there exists a set $\hat{\cal S}^-\equiv 
\hat{\cal S}^-(\eta,v)\subseteq G^-$ which is compact in $(E,\pi)$ 
such that
\begin{eqnarray*}
\bmu\circ\nu_{-t}(\hat{\cal S}^-)\ge\bmu\circ\nu_{-t}(G^-)-\eta\quad\mbox 
{\rm for all}\quad t\in [0,1]. \vphantom{\int} 
\end{eqnarray*}
In Case 2, there exists a set $\hat{\cal S}^+\equiv \hat{\cal S}^+(\eta,v) 
\subseteq G^+$ which is compact in $(E,\pi)$ such that
\begin{eqnarray*}
\bmu\circ\nu_t(\hat{\cal S}^+)\ge\bmu\circ\nu_t(G^+)-\eta\quad\mbox{ 
\rm for all}\quad t\in [0,1]. \vphantom{\int}
\end{eqnarray*}
In Case 1, there exists another compact set ${\cal S}^-\equiv {\cal S}^- 
(\eta,v)\subseteq G^-$ in $(E,\pi)$ such that
\begin{eqnarray*}
\left\{\nu_{-t}:\nu\in\hat{{\cal S}}^-,\ t\in [0,1]\right\}\subseteq{\cal 
S}^- 
\end{eqnarray*} 
and in Case 2, there is another compact set ${\cal S}^+\equiv {\cal S}^+ 
(\eta,v)\subseteq G^+$ in $(E,\pi)$ such that
\begin{eqnarray*}
\left\{\nu_t:\nu\in\hat{{\cal S}}^+,\ t\in [0,1]\right\}\subseteq{\cal 
S}^+\, . 
\end{eqnarray*} 
\end{lemma}
{\bf Remark} (3) We remind again of the different meanings of the 
subindex $m$ in the notation ${\cal S}^v_m$ depending on $v=1$ or 
$1<v<\infty$, cf. Lemma \ref{Lemma1.2}. Let us also keep in mind that 
in Case 1 we have $\bmu\circ\nu_{-t}(G\setminus G^-)=\bmu\left(G_{-t} 
\setminus G_{-t}^-\right)=0$, $t\in [0,1]$, and in Case 2 it holds 
that $\bmu\circ\nu_t(G\setminus G^+)=\bmu\left(G_t\setminus G_t^+ 
\right)=0$, $t\in [0,1]$, cf. Paragraph \ref{sec:1:1:4}. 

Lemma \ref{Lemma1.1} (a) says that {\it the set of measures} $\left 
\{\bmu\circ\nu_{-t}:t\in [-1,1]\right\}$ {\it is tight} in $(E,\pi)$, 
recall from Paragraph \ref{sec:1:1:4} that $\bmu\circ\nu_{-t}(E 
\setminus G)=0$. In addition Lemma \ref{Lemma1.1} (a) states that in 
Case 1 the set of measures $\left\{\bmu\circ\nu_{-t}:t\in [0,1+t_0] 
\right\}$ is tight in $(G,\pi)$ and in Case 2 the set of measures 
$\left\{\bmu\circ\nu_t:t\in [0,1+t_0]\right\}$ is tight in $(G,\pi)$. 
Note here that by condition (a2), $G_{\hat{N}(\eta)}$ is a compact 
subset of $G$ with respect to $\pi$. This part of Lemma \ref{Lemma1.1} 
(a) is just used in order to prove Lemma \ref{Lemma1.1} (b). 

Moreover, Lemma \ref{Lemma1.1} (b) asserts that in Case 1 the set 
of measures $\left\{\bmu\circ\nu_{-t}:t\in [0,1]\right\}$ is tight in 
$(G^-,\pi)$ and in Case 2 the set of measures $\left\{\bmu\circ\nu_t: 
t\in [0,1]\right\}$ is tight in $(G^+,\pi)$. 

The first part of Lemma \ref{Lemma1.1} (a) and Lemma \ref{Lemma1.1} (b) 
are used to prove the main result of the paper, Theorem 
\ref{Theorem2.3}.

\subsection{Proofs of Lemmas \ref{Lemma1.0}-\ref{Lemma1.1}} 
\label{sec:1:4} 

{\bf Proof of Lemma \ref{Lemma1.0}} {\it Step 1 } Let $v=1$. Suppose 
$\nu_0\in E$ is not absolutely continuous with respect to $\lambda$. 
Then there is a sequence of non-negative $g_i\in C_0(D)$ with $\|g_i 
\|_{C_0(D)}=1$ such that $\int g_i\, d\lambda\stack{i\to\infty}{\lra}0$ 
and $\liminf_{i\to\infty}\int g_i\, d\nu>0$. Use e. g. \cite{FL07}, 
Theorem 2.11, for this and take into consideration that, as mentioned 
in the beginning of Subsection \ref{sec:1:1}, $(D,\rho)$ is $\sigma 
$-compact. 

Since $C_0(D)$ is separable, there is a sequence of non-negative $h_i 
\in C_0(D)$ with $\|h_i\|_{C_0(D)}=1$, $i\in {\Bbb N}$, which is dense 
in the set of all non-negative $h\in C_0(D)$ with $\|h\|_{C_0(D)}=1$. 
Define 
\begin{eqnarray*}
f_i^{(1)}:=\int h_i\, d\lambda\quad\mbox{\rm as well as}\quad f_i^{(2)} 
(\nu):=\int h_i\, d\nu\, ,\quad \nu\in E\, ,\quad i\in {\Bbb N},  
\end{eqnarray*}
and ${\cal A}_{i,m,n}:=\left\{\nu\in E: f_i^{(2)}(\nu)\ge \frac1m\right 
\}$ if $f_i^{(1)}\le\frac1n$ as well as ${\cal A}_{i,m,n}:=\emptyset$ 
if $f_i^{(1)}>\frac1n$, $i,m,n\in {\Bbb N}$. We have ${\cal A}_{i,m,n} 
\in {\cal B}(E)$, $i,m,n\in {\Bbb N}$, and therefore the representation 
of the set $\{\nu\in E:\nu$ is not absolutely continuous with respect 
to $\lambda\}$ by 
\begin{eqnarray*}
\bigcup_{m\in {\Bbb N}}\bigcap_{n\in {\Bbb N}}\bigcup_{i\in {\Bbb N}} 
{\cal A}_{i,m,n}\in {\cal B}(E).
\end{eqnarray*}
{\it Step 2 } Let $v>1$. Here use Lemma \ref{Lemma1.2} (b) whose proof is 
carried out independently. 
\qed
\bigskip

\nid
{\bf Proof of Lemma \ref{Lemma1.2} } For a comprehensive presentation 
of the background material, we refer to \cite{Bo07}, Subsection 4.5 and 
Supplements to Chapter 4. 
\medskip

\nid
{\it Step 1 } Let $v=1$. For any sequence $m\equiv(m_1,m_2,\ldots\, )$ 
of strictly decreasing positive numbers with $\lim_{n\to\infty}m_n=0$ 
let us show that  
\begin{eqnarray*}
{\cal S}_m^1:=\left\{h\cdot\lambda\in E:\int_A h\, d\lambda\le\frac1n 
\ \ \mbox{\rm for all } A\in {\cal B}\ \mbox{\rm with }\lambda(A)<m_n 
\ \mbox {\rm and all }n\in {\Bbb N} \right\} 
\end{eqnarray*}
is a compact subset of $(E,\pi)$. Let $h_n\cdot\lambda\in {\cal 
S}_m^1$, $n\in {\Bbb N}$, be an arbitrary sequence. We shall 
demonstrate that it has a subsequence with a limit belonging to 
${\cal S}_m^1$.

In fact, $h_n$, $n\in {\Bbb N}$, is uniformly integrable by the 
definition of ${\cal S}_m^1$. Furthermore, $h_n$, $n\in {\Bbb N}$, 
converges weakly in $L^1(D)$ on some subsequence $(n_k)$ of indices 
by the just mentioned uniform integrability, $\|h_n\|_{L^1(D)}=1$, 
$n\in {\Bbb N}$, and the Dunford-Pettis theorem. Summing up, there 
exists $h\in L^1(D)$ such that $h_{n_k}\stack{k\to\infty}{\lra}h$ 
weakly in $L^1(D)$. 

The bound $0\le h$ $\lambda$-a.e. on the limiting element can be 
verified by checking $0\le\int\vp h_{n_k}\, d\lambda\stack {k\to 
\infty}{\lra}\int\vp h\, d\lambda$ for all $\vp\in L^\infty(D)$ 
with $\vp\ge 0$ $\lambda$-a.e. That $\|h\|_{L^1(D)}=1$ follows 
from $1=\int h_{n_k}\, d\lambda\stack {k\to\infty}{\lra}\int h\, d 
\lambda$. In addition, for any $A\in {\cal B}$, it holds that $\int 
\chi_A h_{n_k}\, d\lambda\stack {k\to \infty}{\lra}$ $\int\chi_A h 
\, d\lambda$ where $\chi_A$ denotes the indicator function of the 
set $A$. Thus $h\cdot\lambda\in {\cal S}_m^1$. 

Consequently, there is a subsequence $h_{n_k}\cdot\lambda$, $k\in 
{\Bbb N}$, converging narrowly to $h\cdot\lambda\in {\cal S}_m^1$ 
which is equivalent to convergence with respect to the metric 
$\pi$, cf. \cite{EK86}, Theorem 3.3.1. We have shown that ${\cal 
S}_m^1$ is a compact subset of $(E,\pi)$. 
\medskip

Now, let $\bnu$ be a finite measure on $(E,{\cal B}(E))$ with 
$\bnu({\cal S}^1)=\bnu(E)$ and let $\ve >0$. For any $n\in {\Bbb N}$, 
there is a sufficiently small $m'_n(\ve)>0$ such that 
\begin{eqnarray*}
S_n:=\left\{h\cdot\lambda:\int_A h\, d\lambda\le\frac1n\ \ \mbox{\rm 
for all } A\in {\cal B}\ \mbox{\rm with }\lambda(A)<m'_n(\ve)\right\} 
\end{eqnarray*}
satisfies $\bnu\left(S_n\right)\ge\bnu(E)-2^{-n}\ve$. Thus, for 
any strictly decreasing sequence $m(\ve)\equiv(m_1(\ve),$ $m_2 
(\ve),\ldots\, )$ with $m_n(\ve)\le m'_n(\ve)$, $n\in {\Bbb N}$, 
and $\lim_{n\to\infty}m_n(\ve)=0$ it holds that 
\begin{eqnarray*}
\bnu\left({\cal S}^1_{m(\ve)}\right)\ge\bnu\left(\bigcap_{n\in 
{\Bbb N}}S_n\right)\ge\bnu(E)-\ve\, . 
\end{eqnarray*}
{\it Step 2 } For $1<v<\infty$, set 
\begin{eqnarray*}
{\cal S}_m^v:=\left\{h\cdot\lambda\in E:\|h\|_{L^v(D)}\le m\right\}\, . 
\end{eqnarray*}
For any sequence $h_n\cdot\lambda\in {\cal S}_m^v$, $n\in {\Bbb N}$, 
there is a subsequence $h_{n_k}$, $k\in {\Bbb N}$, converging weakly 
in $L^v(D)$ to some $h\in L^v(D)$ by the boundedness of $\|h_n\|_{ 
L^v(D)}$. For the limiting element we have $\|h\|_{L^v(D)}\le\liminf_k 
\|h_{n_k}\|_{L^v(D)}\le m$. Furthermore, $h_{n_k}$, $k\in {\Bbb N}$, 
converges also weakly in $L^1(D)$ to $h$. For the remainder of Lemma 
\ref{Lemma1.2} (b) use the same arguments as in Step 1. 
\qed
\bigskip

\nid
{\bf Proof of Lemma \ref{Lemma1.1}} We keep in mind the definition 
of ${\cal S}_m^v$ from the proof of Lemma \ref{Lemma1.2}. Furthermore, 
we recall the notation of Subsection \ref{sec:1:1}. In particular, 
we use $\nu_t\equiv h_t\cdot\lambda$, $t\in [-1,1]$. We recall that 
the function $[-1,1]\ni t\mapsto h_t$ is for $\nu\equiv h_0\cdot 
\lambda\in G$ continuous in $L^v(D)$ by the hypothesis (a1). In Steps 
1 and 2 below, we prove part (a). In Step 1 we focus on the first 
statement for $t\in [-1,1]$. The additional statements of part (a) 
in Case 1 and, respectively, in Case 2 follow in Step 2. In Steps 3 
and 4 we show (b). Let $\eta>0$. 
\medskip 

\nid 
{\it Step 1 } Let $v=1$. By the just mentioned continuity of $[-1,1] 
\ni t\mapsto h_t$ in $L^1(D)$, for $\nu\equiv h_0\cdot\lambda\in G$, 
we get the equicontinuity of the set of real functions $[-1,1]\ni 
t\to\nu_t(A)$, \linebreak $A\in {\cal B}$. Therefore, for fixed $\nu$ 
and $\delta>0$, the function $[-1,1]\ni t\to\sup\{\nu_t(A):\lambda(A) 
<\delta,$ $A\in {\cal B}\}=:\vp_{\nu,\delta}(t)$ is continuous. 
Furthermore, for fixed $\nu$, $\|\vp_{\nu,\delta}\|\to 0$ as $\delta 
\to 0$ by the uniform integrability of each individual $h_t$ and Dini's 
theorem. We set 
\begin{eqnarray*}
m'_{\nu,n}:=\sup\left\{\delta:\|\vp_{\nu,\delta}\|\le\frac1n\right\} 
\, ,\quad n\in {\Bbb N}. 
\end{eqnarray*}

By this construction, we have $\nu_t\in {\cal S}_{m_\nu}^1$ for all 
$t\in [-1,1]$, for any strictly decreasing sequence $m_\nu=(m_{\nu,1}, 
m_{\nu,2},\ldots\, )$ with $m_{\nu,n}\le m'_{\nu,n}$, $n\in {\Bbb N}$. 
Moreover, for all $n\in {\Bbb N}$ there is a sufficiently small $m(\eta 
)'_n>0$ such that 
\begin{eqnarray*}
\bnu\left(\nu\in G:m'_{\nu,n}\ge m(\eta)'_n\right)\ge\bnu(G)-2^{-n}\eta 
\, . 
\end{eqnarray*}
It follows that there exists a strictly decreasing sequence $m(\eta)= 
(m(\eta)_1,$ $m(\eta)_2,\ldots\, )$ of positive real numbers with  
$\lim_{n\to\infty}m(\eta)_n=0$ such that 
\begin{eqnarray*}
\bnu\left(\left\{\nu\in G:m'_{\nu,1}\ge m(\eta)_1,\, m'_{\nu,2}\ge m 
(\eta)_2,\ldots\, \vphantom{l^1}\right\}\right)\ge\bnu(G)-\eta\, . 
\end{eqnarray*}
Accordingly, for all $t\in [-1,1]$, we have 
\begin{eqnarray*}
&&\hspace{-.5cm}\bnu\left(\nu\in G:\nu_t\in {\cal S}_{m(\eta)}^1\right) 
 \\ 
&&\hspace{.5cm}=\bnu\left(\left\{\nu\in G:\nu_t(A)\le\frac1n\ \mbox{ 
\rm if }\ \lambda(A)<m(\eta)_n ,\ A\in {\cal B},\ n\in {\Bbb N}\right\} 
\right) \\ 
&&\hspace{.5cm}\ge\bnu(G)-\eta\, . 
\end{eqnarray*}
This implies $\bnu\circ\nu_{-t}({\cal S}^1_{m(\eta)})=\bnu(\nu\in G:\nu_t 
\in {\cal S}^1_{m(\eta)}\cap G)\ge\bnu(\nu_{-t}:\nu\in G)-\eta$ for all 
$t\in [-1,1]$. 
\medskip 

Furthermore, for $1<v<\infty$, the function $[-1,1]\ni t\mapsto\|h_t\|_{L^v 
(D)}$ is continuous for any $\nu\equiv h\cdot\lambda\in G$. Thus, for those 
$\nu$, there exists $m_\nu\in {\Bbb N}$ such that $\nu_t\in {\cal S}_{m_\nu 
}^v$ for all $t\in [-1,1]$. It follows that there is a sufficiently large 
$m(\eta)\in {\Bbb N}$ such that $\bnu\left(\nu\in G:m_\nu\le m(\eta)\right) 
\ge\bnu(G)-\eta$. Thus, $\bnu(\nu\in G:\nu_t\in {\cal S}^v_{m(\eta)})\ge\bnu 
(G)-\eta$ for all $t\in [-1,1]$. This says $\bnu\circ\nu_{-t}({\cal S}^v_{m 
(\eta)})=\bnu(\nu\in G:\nu_t\in {\cal S}^v_{m(\eta)}\cap G)\ge\bnu(\nu_{-t}: 
\nu\in G)-\eta$, $t\in [-1,1]$. 
\medskip 

\nid 
{\it Step 2 } Recall condition (a2) and suppose $v=1$ and Case 1. For any 
sequence $N\equiv (N_1,N_2,\ldots\, )$ of strictly decreasing positive real 
numbers with $\lim_{n\to\infty}N_n=0$ there exist $t_0(N)$ and some (other) 
such sequence $\hat{N}$ satisfying 
\begin{eqnarray}\label{1.11}
G_N=\{\nu\in G_N,\ \nu_t\in G\}=\{\nu\in G_N,\ \nu_t \in G_{\hat{N}}\}
\end{eqnarray}
for all $t\in [0,1+t_0(N)]$. Again by (a2), for $\eta>0$ there exists a 
sequence $N(\eta)$ with 
\begin{eqnarray}\label{1.12} 
&&\hspace{-.5cm}\bnu\left(G_{-t}\right)-\bnu\left(\left\{\nu_{-t}\in G_{N 
(\eta)}:\nu\in G\right\}\right)\nonumber \\ 
&&\hspace{.5cm}=\bnu\left(\left\{\nu_{-t}\in G\setminus G_{N(\eta)}:\nu\in 
G\right\}\right)=\bnu\left(\left\{\nu\in G\setminus G_{N(\eta)}:\nu_t\in G 
\right\}\right)\nonumber \\
&&\hspace{.5cm}\le\bnu\left(G\setminus G_{N(\eta)}\right)\le\eta\, ,\quad 
t\in [0,1+t_0(N)].  
\end{eqnarray}
Recalling again condition (a2), we may choose $m(\eta)$ such that $G_{\hat 
{N}(\eta)}\subseteq {\cal S}_{m(\eta)}^1$ where $N(\eta)$ and $\hat{N}(\eta 
)$ are related by (\ref{1.11}). In particular, $G_{\hat{N}(\eta)}$ is compact.
According to (\ref{1.11}) and (\ref{1.12}) for $t\in [0,1+t_0(N)\wedge t_0 
(N(\eta))]\equiv [0,1+t_0]$ we have 
\begin{eqnarray}\label{1.14}
&&\hspace{-.5cm}\bnu\left(\left\{\nu\in G:\nu_t\in G_{\hat{N}(\eta)}\right\} 
\right)\ge\bnu\left(\left\{\nu\in G_{N(\eta)}:\nu_t\in G_{\hat{N}(\eta)} 
\right\}\right)\nonumber \\ 
&&\hspace{.5cm}=\bnu\left(\left\{\nu\in G_{N(\eta)}:\nu_t\in G\right\}\right) 
=\bnu\left(\left\{\nu_{-t}\in G_{N(\eta)}:\nu\in G\right\}\right)\nonumber \\ 
&&\hspace{.5cm}\ge\bnu\left(\left\{\nu_{-t}:\nu\in G\right\}\right)-\eta\, .  
\end{eqnarray}

Suppose now $1<v<\infty$ and Case 1. Let $N\in {\Bbb N}$. The function 
$[0,1+t_0(N)]\ni t\mapsto\|h_t\|_{L^v(D)}$ is continuous for any $\nu 
\equiv h\cdot\lambda\in G$. For $t\in [0,1+t_0(N)]$, there exist $\hat 
{N}\in {\Bbb N}$ such that we have (\ref{1.11}), where we stress the 
different meanings of $N$ and $\hat{N}$ depending on $v=1$ or $1<v< 
\infty$. Furthermore, by (a2) there exists $N(\eta)\in {\Bbb N}$ such 
that (\ref{1.12}). As in (\ref{1.14}), with $t_0\equiv t_0(N)\wedge t_0 
(N(\eta))$ we arrive at 
\begin{eqnarray*}
\bnu\left(\left\{\nu\in G:\nu_t\in G_{\hat{N}(\eta)}\right\}\right)\ge\bnu 
\left(\left\{\nu_{-t}:\nu\in G\right\}\right)-\eta\quad\mbox{\rm for}\quad 
t\in [0,1+t_0],  
\end{eqnarray*}
where also in this situation we have $G_{\hat{N}(\eta)}\subseteq{\cal 
S}_{m(\eta)}^v$ for some index $m\in\mathbb{N}$ by condition (a2), i. e. 
$G_{\hat{N}(\eta)}$ is compact. Together with (\ref{1.14}) if $v=1$, we 
obtain the remainder of part (a) and a preparation for the proof of part 
(b), namely 
\begin{eqnarray}\label{1.15}
&&\hspace{-.5cm}\bnu\left(\left\{\nu_{-t}:\nu\in G_{\hat{N}(\eta)}\right 
\}\right)=\bnu\left(\left\{\nu\in G:\nu_t\in G_{\hat{N}(\eta)}\right\} 
\right)\nonumber \\ 
&&\hspace{.5cm}\ge\bnu\left(\left\{\nu_{-t}:\nu\in G\right\}\right)-\eta 
\, ,\quad t\in [0,1+t_0],  
\end{eqnarray}
for all $1\le v<\infty$. The counterpart for Case 2 follows similarly. 
\medskip

\nid 
{\it Step 3 } Let us turn to part (b) of the lemma under the hypotheses 
of Case 1. Recall from Paragraph \ref{sec:1:1:5} the notation $G_{-s}= 
\{\nu_{-s}:\nu\in G\}$, $s\in [-1,\infty)$, and $G^-_{-t}=\bigcup_{s\in 
(t,\infty)}G_{-s}$, $t\in [-1,\infty)$. Note again that $G^-=G_0^-$. 

Because of $G^-_{-t}=\bigcup_{s\in (t,\infty)}G^-_{-s}$, $t\in [0,\infty 
)$, the function $[0,\infty)\ni t\mapsto\bmu(G^-_{-t})$ is right-continuous 
and decreasing. Furthermore, from $G_{-t}=\bigcap_{s\in (-1,t)}G^-_{-s}$, 
$t\in [0,\infty)$, and $\bmu\left(G_{-t}\setminus G_{-t}^-\right)=0$, cf. 
Paragraph \ref{sec:1:1:4}, we obtain  
\begin{eqnarray*}
\lim_{u\uparrow t}\bmu(G_{-u}^-)=\bmu\left(\textstyle{\bigcap_{s\in (-1,t) 
}}G^-_{-s}\right)=\bmu(G_{-t})=\bmu(G^-_{-t})\, ,\quad t\in [0,1]. 
\end{eqnarray*}
We may now deduce that $[0,1]\ni t\mapsto\bmu(G^-_{-t})$ is continuous 
and decreasing and even uniformly continuous. Taking into consideration 
that $t\mapsto\bmu(G^-_{-t})$ is right continuous in $t=1$ it follows now 
that, for given $\eta>0$, there is a $t_0>0$ such that 
\begin{eqnarray}\label{1.13} 
\bmu(G^-_{-t-t_0})\ge\bmu(G^-_{-t})-\frac{\eta}{2}\, ,\quad t\in [0,1].  
\end{eqnarray}
Without loss of generality we may assume that this $t_0$ is the one of 
Step 2 for $\eta$ replaced by $\eta/2$, i. e. $t_0\equiv t_0(N)\wedge 
t_0(N(\eta/2))$. Furthermore, we recall that $G_{\hat{N}(\eta/2)}$ 
appearing in (\ref{1.15}) is compact in $(E,\pi)$. Then $\hat{{\cal S}}^- 
:=\{\nu_{-t_0}:\nu\in G_{\hat{N}(\eta/2)}\}$, as the continuous image of 
the compact set $G_{\hat{N}(\eta/2)}$, is also a compact subset of $(E, 
\pi)$, cf. condition (a2). We even have $\hat{{\cal S}}^-\subseteq G_{ 
-t_0}\subseteq G^-$. For all $t\in [0,1]$ it follows now from (\ref{1.15}) 
that  
\begin{eqnarray*}
&&\hspace{-.5cm}\bmu\circ\nu_{-t}\left(\hat{{\cal S}}^-\right)=\bmu\left 
(\left\{\nu_{-t}:\nu\in\hat{{\cal S}}^-\right\}\right)=\bmu\left(\left\{ 
\nu_{-t-t_0}:\nu\in G_{\hat{N}(\eta/2)}\right\}\right) \\ 
&&\hspace{.5cm}\ge\bmu\left(\left\{\nu_{-t-t_0}:\nu\in G\right\}\right)- 
\frac{\eta}{2}=\bmu(G_{-t-t_0})-\frac{\eta}{2}\, . 
\end{eqnarray*}
Next we keep in mind $G_{-t-t_0}\supseteq G^-_{-t-t_0}$, $t\in [0,1]$, cf. 
Paragraph \ref{sec:1:1:5}. Applying (\ref{1.13}) we get 
\begin{eqnarray*}
&&\hspace{-.5cm}\bmu\circ\nu_{-t}\left(\hat{{\cal S}}^-\right)\ge 
\bmu(G_{-t-t_0})-\frac{\eta}{2}\ge\bmu(G^-_{-t-t_0})-\frac{\eta}{2} 
\ge\bmu(G^-_{-t})-\eta \\ 
&&\hspace{.5cm}=\bmu\circ\nu_{-t}(G^-)-\eta\, ,\quad t\in [0,1].  
\end{eqnarray*}

\nid
{\it Step 4 } Let now $\hat{N}(\eta/2)$ and $\hat{\hat{N}}(\eta/2)$ 
be related as $N$ and $\hat{N}$ in (\ref{1.11}) and let $\hat{t}_0>0$ 
be the counterpart of the number $t_0$ from the end of Step 3 if we 
had started the analysis in Step 2 with the index $\hat{N}(\eta/2)$ 
instead of $N$. Pretending that we had carried out the analysis of 
Steps 2 and 3 with sufficient small $t_0(N)$, $t_0(N(\eta))$, and $t_0 
(N(\eta/2))$, without loss of generality, we may suppose $\hat{t}_0= 
t_0$. In particular, the counterpart to (\ref{1.11}) is 
\begin{eqnarray*}
G_{\hat{N}(\eta/2)}=\left\{\nu\in G_{\hat{N}(\eta/2)},\ \nu_t\in 
G_{\hat{\hat{N}}(\eta/2)}\right\} 
\end{eqnarray*}
for all $t\in [0,1+t_0]$. This implies 
\begin{eqnarray*}
\left\{\nu_t:\nu\in G_{\hat{N}(\eta/2)},\ t\in [0,1+t_0]\right\} 
\subseteq G_{\hat{\hat{N}}(\eta/2)}\, .
\end{eqnarray*}
Thus for the compact set $\hat{{\cal S}}^-=\{\nu_{-t_0}:\nu\in G_{\hat 
{N}(\eta/2)}\}$ in $(E,\pi)$ it holds that 
\begin{eqnarray*}
\left\{\nu_{-t}:\nu\in\hat{{\cal S}}^-,\ t\in [0,1]\right\}\subseteq 
\left\{\nu_{-1-t_0}:\nu\in G_{\hat{\hat{N}}(\eta/2)}\right\}=:{\cal S}^- 
\subseteq G_{-1-t_0}\subseteq G^-\, . 
\end{eqnarray*}
According to condition (a2), (4) in Paragraph \ref{sec:1:1:5} the 
set ${\cal S}^-$ is the continuous image of the compact set $G_{\hat 
{\hat{N}}(\eta/2)}$ in $(E,\pi)$ and therefore compact, as well. 
\qed 

\section{The Change of Measure Formula}\label{sec:2}
\setcounter{equation}{0} 

In a first reading of this section, Subsection \ref{sec:2:11} can be 
skipped and the following definition of the spaces ${\Bbb C}^{q,1}_b 
(E,F)\equiv {\Bbb C}^{q,1}_b(E,F;X,Y)$ and ${\Bbb C}_b^{q,1}(E)\equiv 
{\Bbb C}_b^{q,1}(E;X,Y)$ can be specified to $Y={\Bbb R}$. In this 
case the set $C_b^{q,1}(E)\equiv {\Bbb C}_b^{q,1}(E;X,{\Bbb R})$ is 
precisely the set of test functions we need for the proof of Theorem 
\ref{Theorem2.3}, see Step 1 of that proof. 

In Subsection \ref{sec:2:11} the full definition of the 
spaces ${\Bbb C}^{q,1}_b(E,F)$ and ${\Bbb C}_b^{q,1}(E)$ is needed 
for a meaningful definition of the notion of the {\it trace of the 
gradient}, cf. (\ref{A.65}). Furthermore, the full definition is also 
needed for the proper formulation of conditions (j') and (jjj') of 
Subsection \ref{sec:2:11}, see Proposition \ref{Proposition2.2}. 
\medskip

Let $v$ as in Subsection \ref{sec:1:1} and let $X$ be a Banach space 
densely and continuously embedded in $L^v(D)$ by the inclusion map 
$i:X\mapsto L^v(D)$. Let $Y$ be a further real Banach space. Denote by 
$B(X,Y)$ the space of all bounded linear operators $X\mapsto Y$ endowed 
with the operator norm. Recall from Subsection \ref{sec:1:1} that $\bmu 
(G)=1$ where $G\in {\cal B}(E)$ has the property that $G\subseteq\{g 
\cdot\lambda:g\in L^v(D)\}$. For $\nu=g\cdot\lambda\in G$, $h\in X$, 
and $t\in {\Bbb R}$ let $\nu+t\cdot i\circ h\cdot\lambda:=(g+t\cdot i 
\circ h)\cdot\lambda$. 
\medskip 

For $q\in [1,\infty]$, let ${\Bbb C}^{q,1}_b(E,F)\equiv {\Bbb C}^{q,1}_b 
(E,F;X,Y)$ denote the set of all bounded and everywhere defined functions 
$f:F\to Y$ for which the following holds. 
\begin{itemize}
\item[(i)] For $\bmu$-a.e. $\nu\in E$ and all $h\in X$ the directional 
derivative 
\begin{eqnarray*}
\frac{\partial f}{\partial h}(\nu):=\lim_{t\to 0}\frac1t\left(f(\nu+ 
t\cdot i\circ h\cdot\lambda)-f(\nu)\right)
\end{eqnarray*}
exists in $Y$.  
\item[(ii)] There exists a measurable map $Df:E\to B(X,Y)$ such that   
\begin{eqnarray*}
Df(\nu)(h)=\frac{\partial f}{\partial h}(\nu)\quad\mbox{\rm for }\bmu 
\mbox{\rm -a.e. }\nu\in E\ \mbox{\rm and all }h\in X,  
\end{eqnarray*}
\item[(iii)] and $Df\in L^q(E,\bmu;B(X,Y))$. 
\end{itemize} 
If $Y={\Bbb R}$ then we may identify $B(X,Y)$ with a function space. 
\medskip

Define ${\Bbb C}_b^{q,1}(E)\equiv {\Bbb C}_b^{q,1}(E;X,Y):=\{f|_E:f\in 
{\Bbb C}_b^{q,1}(F,E)\}$, $q\in [1,\infty]$. If $X$ is unspecified and 
$Y={\Bbb R}$ then we use the notation $C_b^{q,1}(E)\equiv {\Bbb C}_b^{ 
q,1}(E;X,{\Bbb R})$ and $C_b^{q,1}(E,F)\equiv {\Bbb C}_b^{q,1}(E,F;X, 
{\Bbb R})$. 
\medskip

The next lemma is concerned with the well-definiteness of the operator 
$D$ introduced in (ii) and (iii) on ${\Bbb C}_b^{q,1}(E)$. For $f\in 
{\Bbb C}_b^{q,1}(E)$ let $\Phi(f):=\{\vp\in {\Bbb C}_b^{q,1}(F,E):\vp|_E 
=f\}$. Let $\1$ denote the constant function on the set $D$ taking the 
value 1. 
\begin{lemma}\label{LemmaA.1} 
Suppose that there is a measurable map $E\ni\nu\mapsto\mu\equiv\mu(\nu) 
\in E$ such that for $\bmu$-a.e. $\nu=k\cdot\lambda\in E$ with $k\in 
L^v(D)$ we have $\mu=i\circ g\cdot\lambda\in E$ with $g\equiv g(\nu)\in X$ 
and, for some $\ve\equiv\ve(\nu)>0$, 
\begin{itemize}
\item[{\rm (a4)}] $\qquad\ve\, i\circ g\le k\, .$ 
\end{itemize}
(a) Let $\psi\in {\Bbb C}_b^{q,1}(F,E)$ with $\psi=0$ on $E$. Then for 
$\bmu$-a.e. $\nu\in E$ we have with $g\equiv g(\nu)\in X$ introduced above  
\begin{eqnarray}\label{A.1}
D\psi(\nu)(\, \cdot\, )=D\psi(\nu)(g)\cdot(i\circ\, \cdot\, ,\1)\, . 
\end{eqnarray}
(b) Let $A\in L^p(E,\bmu;X)$ be a vector field with $(i\circ A,\1)=0$ 
$\bmu$-a.e. Then, for fixed $f\in {\Bbb C}_b^{q,1}(E)$, the expression 
$D\vp(A)\equiv D\vp(\cdot)(A(\cdot))\in L^1(E,\bmu;Y)$ is invariant 
relative to $\vp\in\Phi(f)$. 
\end{lemma} 
Proof. (a) For all $\t g\in X$ with $i\circ\t g\cdot\lambda\in E$ and 
all $h\in X\setminus\{0\}$ with $h\ge 0$ introduce 
\begin{eqnarray*}
\overline{h}(\t g):=h-(i\circ h,\1)\cdot \t g\in X\, . 
\end{eqnarray*}
For $\bmu$-a.e. $\nu=k\cdot\lambda\in E$ with $k\in L^v(D)$ and $g\equiv 
g(\nu)\in X$ we have $\ve\, i\circ g\le k$, cf. (a4) Thus for all 
$h\in X\setminus\{0\}$ with $h\ge 0$ we obtain 
\begin{eqnarray*}
\nu+t\cdot i\circ \overline{h}(g)\cdot\lambda:=\left(k+t\cdot i\circ 
h - t\cdot(i\circ h,\1)\cdot i\circ g\vphantom{l^1}\right)\cdot\lambda 
\in E\, ,\quad t\in\left[0,\ve\, (i\circ h,\1)^{-1}\right), 
\end{eqnarray*}
which implies because of the hypothesis $\psi=0$ on $E$ 
\begin{eqnarray*} 
0=\lim_{t\downarrow 0}\frac1t\left(\psi\left(\nu+t\cdot i\circ\overline 
{h}(g)\cdot\lambda\right)-\psi(\nu)\vphantom{l^1}\right)=D\psi(\nu)(h) 
-D\psi(\nu)((i\circ h,\1)\cdot g)\, . 
\end{eqnarray*}
In other words, we have (\ref{A.1}). 
\medskip 

\nid
(b) This is an immediate consequence of part (a). 
\qed 
\medskip 

Suppose (a4). Then in the sense of Lemma \ref{LemmaA.1} (b), for 
$f\in {\Bbb C}_b^{q,1}(E)$ we define the {\it gradient} $Df$ of $f$ as the 
restriction of $D\vp(\nu)$ to $\bmu$-a.e. $\nu\in E$ for any $\vp\in\Phi 
(f)$. This definition includes that $Df$ is only applied to vector fields 
$A\in L^p(E,\bmu;X)$, $p\in [1,\infty]$, with $(i\circ A,\1)=0$ $\bmu$-a.e. 
Let $X\ominus\1:=\{b\in X:(i\circ b,\1)=0\}$. In contrast to $\t f\in {\Bbb 
C}_b^{q,1}(E,F)$ where $D\t f\in L^q(E,\bmu;B(X,Y))$, 
\begin{eqnarray*} 
\mbox{\rm for}\quad f\in {\Bbb C}_b^{q,1}(E)\quad\mbox{\rm we have}\quad 
Df\in L^q(E,\bmu;B(X\ominus\1,Y))\, . 
\end{eqnarray*} 
Let $L_{\rm loc}^\infty(E,\bmu;Y)$ denote the set of all $Y$-valued $\bmu 
$-equivalence classes $f$ on $E$ such that we have the following. For any set 
$K\subseteq G^-$ in Case 1 or $K\subseteq G^+$ in Case 2 which is compact with 
respect to the Prokhorov metric $\pi$ it holds that $\|f(\cdot)\|_Y\le c_K$ 
$\bmu$-a.e. on $K$ for some $c_K>0$. If $Y={\Bbb R}$ then we may use the notation 
$L_{\rm loc}^\infty(E,\bmu)\equiv L_{\rm loc}^\infty (E,\bmu;{\Bbb R})$. We say 
that $f\in {\Bbb C}_b^{q,1} (E)$ belongs to ${\Bbb C}^{q,1}_{b,{\rm loc}}(E)$ if 
\begin{eqnarray*}
Df\in L_{\rm loc}^\infty(E,\bmu;B(X\ominus\1,Y))\, . 
\end{eqnarray*} 

We mention that condition (a4) has to be verified in applications taking into 
consideration the particular choices of $X$ and $G$, see Sections \ref{sec:3} and 
\ref{sec:4}. 
\begin{lemma}\label{Lemma2.2} 
Suppose (a4) and let $q\in [1,\infty]$. If $f,g\in C_b^{q,1}(E)$ then we have 
$fg\in C_b^{q,1}(E)$. Furthermore, $f,g\in C_{b,{\rm loc}}^{q,1}(E)$ implies 
$fg\in C_{b,{\rm loc}}^{q,1}(E)$. In either case it holds that 
\begin{eqnarray*} 
fDg+gDf=D(fg)\, . 
\end{eqnarray*} 
\end{lemma}
Proof. Reviewing the definition of $C_b^{q,1}(E)$ or $C_{b,{\rm loc} 
}^{q,1}(E)$, we get $fg\in C_b^{q,1}(E)$ or, respectively, $fg\in C_{ 
b,{\rm loc}}^{q,1}(E)$. The formula $fDg+gDf=D(fg)$ is a byproduct when 
checking property (ii). 
\qed 
\bigskip 

As already noted in the beginning of Subsection \ref{sec:1:1}, $C_0(D)$ 
is separable. Thus there is a sequence $h_n\in C_0(D)$, $n\in {\Bbb N}$, 
such that the closed linear span of $h_1,h_2,\ldots $ with respect to the 
sup-norm is $C_0(D)$. In the application of Section \ref{sec:4} we will 
work with an explicit choice of this sequence, cf. Remark (5) of Section 
\ref{sec:4}. Introduce $C_b(E)$ as the space of all bounded continuous 
real functions on $E$ and 
\begin{eqnarray}\label{2.1}
\t C_b^1(E):=\left\{f(\nu)=\vp((h_1,\nu),\ldots,(h_r,\nu)),\ \nu\in E: 
\vphantom{\t f}\vp\in C_b^1({\Bbb R}^r),\ r\in {\Bbb N}\right\}.\qquad
\end{eqnarray}

\noindent
{\bf Remarks }(1) Let us note that any function $f(\nu)$, $\nu$ belonging 
to $E$, has a formal representation $f(\nu)=\vp((h_1,\nu),(h_2,\nu),\ldots 
\, )$. For this, recall that the closed linear span with respect to the 
sup-norm of $h_1,h_2,\ldots\, $ is $C_0(D)$. The latter implies also that 
every $\nu\in E$ is uniquely determined by $(h_1,\nu),(h_2,\nu),\ldots\, $. 
In particular the set $\t C_b^1(E)$ separates the points of $E$. 
\medskip 

\nid 
(2) As mentioned in Subsection \ref{sec:1:1}, convergence on $E$ is by the 
Portmanteau theorem in the form of \cite{EK86}, Theorem 3.3.1, equivalent 
to narrow convergence of probability measures. For this recall also that 
$(D,\rho)$ is a separable metric space. Thus the functions in $\t C_b^1 
(E)$ are continuous. Moreover by the Stone-Weierstrass theorem and the 
last sentence of Remark (1), for each compact set $K\subseteq E$, the set 
$\{f|_K:f\in\t C_b^1(E)\}$ is dense in the space $C_b(K)$ of all (bounded 
and) continuous functions on $K$. 
\medskip 

\nid 
(3) Suppose (a4). Obviously, $\t C_b^1(E)\subseteq C_b^{\infty,1}(E)$ and 
for $f$ and $\vp$ as in (\ref{2.1}), the gradient of $f$ in the sense of 
(i)-(iii) is given by  
\begin{eqnarray*}
&&\hspace{-.5cm}Df(\nu)=\sum_{i=1}^r\frac{\partial\vp}{\partial x_i}\left( 
(h_1,\nu),\ldots ,(h_r,\nu)\right)h_i \\ 
&&\hspace{.5cm}\equiv\sum_{i=1}^r\frac{\partial\vp}{\partial x_i}\left((h_1 
,\nu),\ldots ,(h_r,\nu)\right)\, (h_i,i\circ\cdot)\, .  
\end{eqnarray*} 
Here we emphasize that, for the sake of definiteness, we will only work with 
$Df(\nu)(h)$ where $h\in X\ominus\1$, cf. the definitions following Lemma  
\ref{LemmaA.1}. 

Let us use the notation $(X\ominus\1)^\ast\equiv B(X\ominus\1,{\Bbb R})$. We 
recall that $X$ is densely and continuously embedded in $L^v(D)$ which implies 
that $X\ominus\1$ is densely and continuously embedded in $L^v(D)\ominus\1$. 
In particular, there is $C>0$ such that for all $x\in X\ominus\1$ we have 
$\|x\|_{L^v(D)}\le C\|x\|_X$. Denoting for $z^\ast\in B(L^v(D)\ominus\1, 
{\Bbb R})$ by $i^\ast z^\ast$ the restriction of $z^\ast$ to the dense set 
$X\ominus\1$, we recall also that 
\begin{eqnarray*}
&&\hspace{-.5cm}\|i^\ast z^\ast\|_{(X\ominus\1)^\ast}=\sup_{x\in X\ominus\1} 
\frac{|z^\ast(x)|}{\|x\|_X}\le\sup_{z\in L^v(D)\ominus\1}C\, \frac{|z^\ast 
(z)|}{\|z\|_{L^v(D)}}=C\|z^\ast\|_{B(L^v(D)\ominus\1,{\Bbb R})}<\infty\, . 
\end{eqnarray*} 
In other words, $i^\ast z^\ast\in (X\ominus\1)^\ast$. In this sense, $i^\ast$ 
is the continuous inclusion of $B(L^v(D)\ominus\1,{\Bbb R})\equiv L^w(D)\ominus 
\1$ into $(X\ominus\1)^\ast$ where $1/v+1/w=1$ for $1<v<\infty$ and $w=\infty$ 
if $v=1$. We note that $(f,i\circ\cdot)={}_{(X\ominus\1)^\ast}\langle i^\ast 
f,\cdot\, \rangle_{X\ominus\1}$, $f\in L^w(D)\ominus\1$, i. e. $i^\ast:L^w(D) 
\ominus\1\mapsto (X\ominus\1)^\ast$ is the dual operator of the densely defined 
inclusion $i:X\ominus\1\mapsto L^v(D)\ominus\1$. For $f\in\t C_b^1(E)$ and for 
$f$ and $\vp$ as above we obtain  
\begin{eqnarray*}
Df(\nu)=\sum_{i=1}^r\frac{\partial\vp}{\partial x_i}\left((h_1,\nu),\ldots ,
(h_r,\nu)\right)\, {}_{(X\ominus\1)^\ast}\langle i^\ast h_i,\cdot\, \rangle_{ 
X\ominus\1}\, .  
\end{eqnarray*} 

\subsection{General Hypotheses on $\bmu$ and $A^f$ }\label{sec:2:1}

Let $1\le v<\infty$. Let us recall from Paragraph \ref{sec:1:1:5} that 
$G,G^-,G^+,G^-_{-t},G^+_t\in {\cal B}(E)$, $t\in [0,1]$. Introduce 
$\gamma^- G:=G\setminus G^-$ as  well as ${\cal B}(\gamma^- G):=\{B\cap 
\gamma^- G:B\in {\cal B}(E)\}$ in Case 1, and $\gamma^+ G:=G\setminus 
G^+$ as well as ${\cal B}(\gamma^+ G):=\{B\cap\gamma^+ G:B\in {\cal B} 
(E)\}$ in Case 2. We observe that, if Case 1 and Case 2 hold 
simultaneously, $\gamma^- G=\gamma^+ G=\emptyset$. 

We note also that, in Case 1, $\t \nu\in G\setminus G^-_{-1}$ can uniquely 
be identified with $(\nu,t)\in\gamma^-G\times [0,1]$ such that $\t \nu= 
\nu_{-t}\equiv\nu_{-t}(\nu)$. Furthermore, in Case 2, $\t \nu\in G\setminus 
G^+_1$ can uniquely be identified with $(\nu,t)\in\gamma^+ G\times [0,1]$ 
such that $\t \nu=\nu_t\equiv\nu_t(\nu)$. Let us recall that by Paragraph 
\ref{sec:1:1:5}, $[-1,1]\ni t\mapsto\nu_{-t}$ is $\bmu$-a.e. continuous with 
respect to the metric $\pi$.  
\medskip 

The following proposition will motivate condition (jj) below. 
\begin{proposition}\label{Proposition2.25} 
Suppose (a1) and (a3). (a) If, in Case 1, $\bmu(G\setminus G^-_{-1}) 
>0$ then $\chi_{G\setminus G^-_{-1}}\cdot\bmu$, as a measure on the 
$\sigma$-algebra $\{B\cap (G\setminus G^-_{-1}): B\in {\cal B}(E)\}$, 
can be disintegrated into 
\begin{eqnarray*}
\chi_{G\setminus G^-_{-1}}\cdot\bmu=\int_{\gamma^-G}\lambda_{-1} 
(\nu)\, \bmu^-_\gamma(d\nu)
\end{eqnarray*}
such that $\bmu^-_\gamma$ is a probability measure on $(\gamma^-G,{\cal 
B}(\gamma^- G))$. Furthermore, $\lambda_{-1} (\nu)$, $\nu\in\gamma^-G$, 
are measures on $[0,1]$ endowed with its $\sigma$-algebra ${\cal B}([0, 
1])$ such that $\lambda_{-1}(\nu)([0,1])=\bmu(G\setminus G^-_{-1})$ for 
$\bmu^-_\gamma$-a.e. $\nu\in \gamma^-G$. The probability measure 
$\bmu^-_\gamma$ is unique. For $\bmu^-_\gamma$-a.e. $\nu\in \gamma^-G$ 
the measures $\lambda_{-1}(\nu)$ are also unique. 

Replacing $G^-_{-1}$ as well as $\gamma^- G$ by $G^+_1$ as well as 
$\gamma^+ G$, a corresponding assertion holds in Case 2 for some unique 
probability measure $\bmu^+_\gamma$ on $(\gamma^+G,{\cal B}(\gamma^+ G))$ 
and some family $\lambda_1\equiv\lambda_1(\nu)$ of measures on $([0,1], 
{\cal B}([0,1]))$, instead of $\bmu^-_\gamma$ and $\lambda_{-1}\equiv 
\lambda_{-1}(\nu)$. \\ 
(b) Assume Case 1 and $\bmu(G\setminus G^-_{-1})>0$. Suppose that the 
measures $\bmu\circ\nu_{-t}$, $t\in [0,1]$, are equivalent and that the 
right derivative, 
\begin{eqnarray}\label{2.5*}
\rho^-_0:=\left.\frac{d^+}{dt}\right|_{t=0}\frac{d\bmu\circ\nu_{-t}}{d\bmu} 
\quad\mbox{\it exists as a derivative in } L^1(E,\bmu)\, .
\end{eqnarray}
Suppose furthermore, that 
\begin{eqnarray*}
\rho^-_\gamma(\nu):=\lim_{t\downarrow 0}\frac1t\lambda_{-1}(\nu)([0,t]) 
\quad\mbox{\it exists as a limit in } L^1(\gamma^-G,\bmu^-_\gamma)\, .
\end{eqnarray*}
Then, for all $f\in C_b(E)$ for which $\frac{d^+}{dt} f(\nu_t)$ exists as 
a derivative in $L^1(E,\bmu)$, we have 
\begin{eqnarray}\label{2.6*}
\int\left.\frac{d^+}{dt}\right|_{t=0}f(\nu_t)\, d\bmu=\int f\rho^-_0\, 
d\bmu+\int_{\gamma^- G}f\rho^-_\gamma\, d\bmu^-_\gamma\, .  
\end{eqnarray}
If in Case 1, $\bmu(G\setminus G^-_{-1})=0$ then assuming (\ref{2.5*}) we 
have (\ref{2.6*}) for $\rho^-_\gamma=0$. \\ 
(c) Assume Case 2 and replace in the hypotheses of part (b) the symbols 
$\bmu\circ\nu_{-t}$, $\bmu^-_\gamma$, $\lambda_{-1}$, $\rho^-_0$, 
$\rho^-_\gamma$, and $\gamma^- G$ by $\bmu\circ\nu_t$, $\bmu^+_\gamma$, 
$\lambda_1$, $\rho^+_0$, $\rho^+_\gamma$, and $\gamma^+ G$. Then 
\begin{eqnarray*}
\int\left.\frac{d^+}{dt}\right|_{t=0}f(\nu_{-t})\, d\bmu=\int f\rho^+_0\, 
d\bmu+\int_{\gamma^+ G}f\rho^+_\gamma\, d\bmu^+_\gamma  
\end{eqnarray*}
for all $f\in C_b(E)$ for which $\frac{d^+}{dt} f(\nu_{-t})$ exists as a 
derivative in $L^1(E,\bmu)$. 
\end{proposition}
Proof. For part (a) we check with \cite{AGS05}, Subsection 5.3, and the 
references given therein. What we have to verify is that $(G\setminus 
G^-_{-1},\pi)$ as well as $(\gamma^- G,\pi)$ are separable Radon spaces.  

Since we always suppose that $(D,\rho)$ is separable, $(E,\pi)$ is 
separable as well, cf. \cite{EK86}, Theorem 3.1.7. Thus $(G\setminus 
G^-_{-1},\pi)$ as well as $(\gamma^- G,\pi)$ are separable metric spaces. 
Furthermore, ${\cal S}^v=E\cap\{h\cdot\lambda:h\in L^v(D)\}\in {\cal B} 
(E)$ according to Lemma \ref{Lemma1.0} and $({\cal S}^v,\pi)$ is by Lemma 
\ref{Lemma1.2} a Radon space. Since $G\setminus G^-_{-1}\in {\cal B}(E)$ 
as well as $\gamma^- G\in {\cal B}(E)$ and $G\setminus G^-_{-1}\subseteq 
{\cal S}^v$ as well as $\gamma^-G\subseteq {\cal S}^v$ it holds that 
$G\setminus G^-_{-1}\in {\cal B}({\cal S}^v)$ and $\gamma^- G\in {\cal B} 
({\cal S}^v)$. Thus the separable metric spaces $(G\setminus G^-_{-1},\pi)$ 
as well as $(\gamma^-G,\pi)$ are also Radon spaces. 
\medskip 

\nid
We prove part (b). First we demonstrate that for $f\in C_b(E)$ it holds 
that 
\begin{eqnarray}\label{2.50*}
\lim_{t\downarrow 0}\frac1t\int f(\nu_t)\chi_{G\setminus G^-_{-t}}\, d\bmu 
=\int_{\gamma^- G}f\rho^-_\gamma\, d\bmu^-_\gamma\, .  
\end{eqnarray}
We fix $t\in (0,1]$ for a moment and note that the measure $\frac1t\, 
\chi_{G\setminus G^-_{-t}}\cdot\bmu$, as a measure on the 
$\sigma$-algebra $\{B\cap (G\setminus G^-_{-t}):B\in {\cal B}(E)\}$, can be 
represented by 
\begin{eqnarray*}
\frac1t\, \chi_{G\setminus G^-_{-t}}\cdot\bmu=\frac1t\int_{\gamma^-G} 
\lambda_{-1}(\nu)|_{[0,t]}\, \bmu^-_\gamma(d\nu)
\end{eqnarray*}
where $\lambda_{-1}(\nu)|_{[0,t]}$ is the restriction of $\lambda_{-1}(\nu)$ 
to $[0,t]$. For $f\in C_b(E)$, this leads to 
\begin{eqnarray*} 
&&\hspace{-.5cm}\frac1t\int f(\nu_t)\chi_{G\setminus G^-_{-t}}\, d\bmu= 
\frac1t\int\int_{[0,t]}f(\nu_{t-s})\lambda_{-1}(\nu)(ds)\, \bmu^-_\gamma 
(d\nu) \\ 
&&\hspace{.5cm}=\frac1t\int f(\nu_u)\lambda_{-1}([0,t])(\nu)\, \bmu^-_\gamma 
(d\nu) 
\end{eqnarray*}
for some $u\equiv u(\nu,t)\in [0,t]$, $t\in (0,1]$. Now we observe that 
$[0,1]\ni t\mapsto f(\nu_u)\equiv f(\nu,u(t))$ converges in the measure 
$\bmu^-_\gamma$ to the restriction $f|_{\gamma^-G}$ of $f\in C_b(E)$ to 
$\gamma^-G$, boundedly on $\gamma^-G$ as $t\to 0$. Together with the 
hypothesis that $\frac1t\lambda_{-1}([0,t])(\cdot)$ converges in $L^1 
(\gamma^-G,\bmu^-_\gamma)$ to $\rho^-_\gamma$ as $t\to 0$, this implies 
(\ref{2.50*}). 

The claim is now a straightforward conclusion of the hypotheses of part 
(b) of the proposition and (\ref{2.50*}). For all $f\in C_b(E)$ for which 
$\frac{d^+}{dt} f(\nu_t)$ exists as a derivative in $L^1(E,\bmu)$ we 
obtain 
\begin{eqnarray*}
&&\hspace{-.5cm}\int\left.\frac{d^+}{dt}\right|_{t=0}f(\nu_t)\, d\bmu 
=\lim_{t\downarrow 0}\frac1t\int (f(\nu_t)-f(\nu))\, \bmu(d\nu) \\ 
&&\hspace{.5cm}=\lim_{t\downarrow 0}\frac1t\int\left(f(\nu_t)\cdot 
\chi_G(\nu_t)-f(\nu)\right)\, d\bmu+\lim_{t\downarrow 0}\frac1t\int f 
(\nu_t)\left(1-\chi_G(\nu_t)\right)\, d\bmu \\ 
&&\hspace{.5cm}=\lim_{t\downarrow 0}\frac1t\int f\left(\frac{d\bmu\circ 
\nu_{-t}}{d\bmu}-1\right)\, d\bmu +\lim_{t\downarrow 0}\frac1t\int f 
(\nu_t)\chi_{G\setminus G_{-t}}(\nu)\, d\bmu\\ 
&&\hspace{.5cm}=\int f\rho^-_0\, d\bmu+\int_{\gamma^- G}f\rho^-_\gamma 
\, d\bmu^-_\gamma\, .  
\end{eqnarray*}
\qed 
\medskip

We aim to prove an absolute continuity result relative to $\bmu$ under the 
flow $W^+$. For the subsequent assumptions, suppose (a1)-(a4) and note that 
the Radon-Nikodym derivative 
\begin{eqnarray*}
a^f(\cdot):=d( A^f \cdot )/d\lambda\equiv d( A^f\nu)/d\lambda 
\end{eqnarray*}
takes values in $L^v(D)$ by Subsection \ref{sec:1:1} and is a function of 
$\nu\in E$. We stress that by the definition of the set $G$ in Subsection 
\ref{sec:1:1} and relations (\ref{1.3}), (\ref{1.4}) we have $(A^f,\1)= 
(a^f,\1)=0$ $\bmu$-a.e. on $E$. This guarantees well-definiteness of 
hypothesis (jj) below. For simplicity in the notation, for example in the 
subsequent condition (j), we  shall use the symbol $a^f(\cdot)$ for both, 
an element of $X$ as well as the embedded element in $L^v(D)$. 
\begin{itemize}
\item[(j)] There exists $p>1$ such that 
\begin{eqnarray*}
a^f\in L^p(E,\bmu;X)\, .  
\end{eqnarray*}
\item[(jj)] For $p$ given by (j), $1/p+1/q=1$, there exist 
\begin{itemize}
\item[$\bullet$] a unique {\it divergence} $\delta (A^f)\in L^1(E,\bmu)$ of 
the measure $\bmu$ relative to the vector field $A^f\cdot$ and the gradient 
$D$ and 
\item[$\bullet$] a finite measure $\bmu_\gamma$ on $(\gamma^- G,{\cal B} 
(\gamma^- G))$ in Case 1 or $(\gamma^+ G,{\cal B}(\gamma^+ G))$ in Case 2 
\end{itemize} 
such that 
\begin{eqnarray*}
\int Df(a^f)\, d\bmu=-\int f\, \delta (A^f)\, d\bmu+b\int f\, d\bmu_\gamma 
\, , \quad f\in C^{q,1}_b(E),  
\end{eqnarray*}
where $b=1$ in Case 1 and $b=-1$ in Case 2. If Case 1 and Case 2 hold 
simultaneously then $\mu_\gamma$ is the zero measure. 
\item[(jjj)] 
\begin{eqnarray*} 
\delta (A^f)\in L_{\rm loc}^\infty(E,\bmu) 
\end{eqnarray*} 
\end{itemize} 
where we mention that in (jj) we require $\delta (A^f)\in L^1(E,\bmu)$ 
for definiteness. 
\bigskip

The hypotheses (a1),(a2) of Subsection \ref{sec:1:1} can be regarded as 
conditions on maps of elements and subsets of $E$ with respect to the 
trajectory $[-1,1]\ni t\mapsto\nu_t=U^+(t,\nu)$. In addition hypotheses 
(a3),(a4) of Subsection \ref{sec:1:1} and the introduction to Section 
\ref{sec:2}, together with condition (j), are posted in order to 
guarantee well-definiteness of the objects of the differential calculus 
we are going to use. 

On the other hand, (jj) and (jjj) are very compressed conditions on both, 
the measure $\bmu$ and the vector field $A^f\cdot$. In Proposition 
\ref{Proposition2.2} below we formulate conditions sufficient for 
(j)-(jjj) such that the requirements on the measure $\bmu$ are, as good 
as possible, separated from the requirements on the vector field $A^f\cot$. 

\subsection{Specified Hypotheses on $\bmu$ and $A^f$ }\label{sec:2:11}

Throughout the whole Subsection \ref{sec:2:11} we shall suppose that 
Case 1 and Case 2 hold simultaneously, i. e. that $\mu_\gamma$ is the 
zero measure. Without further mentioning, throughout this subsection 
we also will assume that (a1)-(a4) hold. 

In this motivating introductory paragraph, let $D$ be a bounded open 
subset of ${\Bbb R}^d$ with a globally Lipschitz continuous boundary 
or, equivalently, boundary of H\"older class $C^{0,1}$. Such a boundary 
is even a minimally smooth boundary in the sense of \cite{EE87}, Section 
V.4.4. Furthermore, let $\lambda$ be the Lebesgue measure on $(D,{\cal 
B}(D))$, $X:=W^{1,v'}(D)$ and $Y:=L^v(D)$, where $1\le v'<d$ and $v\in 
[v',dv'/(d-v'))$. The canonical embedding $i:X\mapsto Y$ is dense and, 
according to \cite{EE87}, Theorems V.3.7 and V.4.13, compact. 

According to \cite{Ma02} and \cite{Ma03} there exists a Schauder basis 
$b^o_n$, $n\in {\Bbb N}$, in $X$ such that $i\circ b^o_n$, $n\in {\Bbb N}$, 
is a Schauder basis in $Y$. As a standard consequence, $b_1:=\1$, $b_n:= 
b^o_n-(\1,b^o_n)\, \1$, $n\ge 2$, is a Schauder basis in $X$ such that $i 
\circ b_n$, $n\in {\Bbb N}$, is a Schauder basis in $Y$ and 
\begin{eqnarray}\label{A.2}
b_1=\1\quad\mbox{\rm as well as}\quad \int b_n\, d\lambda=0\, , 
\quad n\ge 2.
\end{eqnarray} 

\subsubsection{General Setup of Subsection \ref{sec:2:11}}\label{sec:2:2:1}  
Now we drop the particular choice of the separable locally compact metric 
space $(D,\rho)$, the measure $\lambda$ on $D$ equipped with its Borel 
$\sigma$-algebra ${\cal B}$, and the Banach space $X$, i. e. we leave $D$, 
$\lambda$, and $X$ unspecified. However, we keep $Y=L^v(D)$. Motivated by 
the above we suppose that $X$ is densely and continuously canonically embedded 
in $Y$ which we here also denote by $i:X\mapsto Y$. In particular, 
there is $C>0$ just depending on $X$ such that 
\begin{eqnarray}\label{A.3}
\|i\circ b\|_Y\le C\, \|b\|_X\, ,\quad b\in X. 
\end{eqnarray} 
Furthermore, we assume that $X$ admits a Schauder basis $b_n$, 
$n\in {\Bbb N}$, satisfying (\ref{A.2}) such that $i\circ b_n$, 
$n\in {\Bbb N}$, is a Schauder basis in $Y$. This setup we keep 
in force throughout the whole Subsection \ref{sec:2:11}. 

Let $b_n^\ast\in X^\ast$, $n\in {\Bbb N}$, denote the biorthogonal 
functionals associated with the Schauder basis $b_n$, $n\in 
{\Bbb N}$, in $X$ and let $i^{\ast -1}\circ b_n^\ast\in Y^\ast$, 
$n\in {\Bbb N}$, denote the biorthogonal functionals associated 
with the Schauder basis $i\circ b_n$, $n\in {\Bbb N}$, in $Y$. 
Here $i^{\ast -1}$ can be seen as the dual of the inverse operator 
of $i$. Due to (\ref{A.2}) we have $i^{\ast -1}\circ b_1^\ast=\1
/\lambda(D)$, i. e. for $h\in Y\ominus\1\equiv\{g\in Y:(g,\1)=0 
\}$ it holds that $h=\sum_{n=2}^\infty {}_{Y^\ast}\langle i^{\ast 
-1}\circ b_n^\ast,h\rangle_Y\, i\circ b_n$. This property will 
be crucial in Paragraph \ref{sec:2:2:5}. For example, it will 
characterize the type of integrand we are going to use in the 
stochastic integral in (jjj'). 

\subsubsection{A class of measures $\bmu$}\label{sec:2:2:2} In this 
paragraph we formulate a certain {\it summability condition} 
(\ref{A.4}) guaranteeing the well-definiteness of the class of 
probability measures $\bmu$ introduced by (\ref{A.5}) and (\ref{A.6}) 
below. 

Denote by $f_n:\{\alpha\cdot\lambda:\alpha\in Y\}\mapsto {\Bbb R}$ 
the map $\alpha\cdot\lambda\mapsto {}_{Y^\ast}\langle i^{\ast -1} 
\circ b_n^\ast,\alpha\rangle_Y$, $n\in {\Bbb N}$. Let $\t \mu_1,\t 
\mu_2,\ldots $ be probability measures on $({\Bbb R},{\cal B}({\Bbb 
R}))$ and consider cylindrical sets of the form  
\begin{eqnarray*}
\t A:=\left\{f_{n_1}^{-1}(A_1)\right\}\cap\ldots\cap\left\{f_{n_k}^{-1} 
(A_k)\right\}\, ,\quad A_1,\ldots ,A_k\in {\cal B}({\Bbb R}), 
\end{eqnarray*} 
where $n_1,\ldots ,n_k\in {\Bbb N}$ with $n_i\neq n_j$ for $i\neq j$. 
From (\ref{A.3}), 
\begin{eqnarray*}
\sum_{n=1}^\infty{}_{Y^\ast}\langle i^{\ast -1}\circ b_n^\ast,\alpha 
\rangle_Y\, i\circ b_n =\alpha\ \mbox{\rm in }L^v(D)\, ,\quad\alpha\in Y 
\equiv L^v(D), 
\end{eqnarray*} 
the Borel-Cantelli lemma, and the Kolmogorov extension theorem in the form 
of e. g. \cite{Ta11}, Theorem 2.4.3, or \cite{K06}, Theorem 14.36, we deduce 
the following implication. If there is a sequence $c_n>0$, $n\in {\Bbb N}$, 
of constants with $\sum_{n=1}^\infty c_n\|b_n\|_X<\infty$ such that 
\begin{eqnarray}\label{A.4}
\sum_{n=1}^\infty\t \mu_n\left(\left\{f_n(\alpha\cdot\lambda):|f_n(\alpha 
\cdot\lambda)|>c_n,\ \alpha\in Y\right\}\right)<\infty 
\end{eqnarray} 
then 
\begin{eqnarray}\label{A.5} 
\t \bmu(\t A):=\t \mu_{n_1}(A_1)\cdot\ldots\cdot\t \mu_{n_k}(A_k)\, ,\quad 
A_1,\ldots ,A_k\in {\cal B}({\Bbb R}), 
\end{eqnarray} 
$n_1,\ldots ,n_k\in {\Bbb N}$ with $n_i\neq n_j$ for $i\neq j$ defines a 
probability measure on $\{h\cdot\lambda:h\in L^v(D)\}$ endowed with the 
$\sigma$-algebra ${\cal C}$ generated by the above cylindrical sets $\t A$. 

Let us assume the existence of such a sequence $c_n>0$, $n\in {\Bbb N}$.  
We recall that the set $G\in {\cal B}(E)$ defined in Subsection \ref{sec:1:1} 
is a subset of $\{h\cdot\lambda:h\in L^v(D),\ \|h\|_{L^1(D)}=1,\ h\ge 0\}$ 
with $\bmu(G)=1$. Let $\vp:\left(\{h\cdot\lambda:h\in L^v(D)\},{\cal C} 
\vphantom{\D l^1}\right)\mapsto\left(G,\{B\cap G:B\in {\cal B}(E)\}\vphantom 
{\D l^1}\right)$ be a measurable bijective map. Now we can define a probability 
measure $\bmu$ on $(E,{\cal B}(E))$ by $\bmu(E\setminus G)=0$ and 
\begin{eqnarray}\label{A.6}
\bmu(A):=\t\bmu\circ\vp^{-1}(A)\quad\mbox{\rm for all}\quad A\in\{B\cap G: 
B\in {\cal B}(E)\}. 
\end{eqnarray}

\subsubsection{The derivative of the measure $\bmu$ in direction of 
a basis vector}\label{sec:2:2:3} This paragraph is dedicated to a 
preliminary definition of the term {\it differentiability} of a 
measure over probability measures. 

Here the motivation comes from the derivatives in the sense of {\it 
Fomin} and {\it Skorokhod}, see \cite{Bo10}, particularly Subsections 
3.1, 3.3, 3.6. We say that the measure $\bmu$ is {\it Fomin type 
differentiable} in direction of $i\circ b$, $b\in X$ with $\int b\, 
d\lambda=0$, if  
\begin{eqnarray*} 
d_{i\circ b}\bmu(A):=\lim_{t\to 0}\frac{\bmu(A+t\cdot i\circ b\cdot 
\lambda)-\bmu(A)}{t}\quad\mbox{\rm exists for all}\ A\in {\cal B}(E), 
\end{eqnarray*} 
where $A+t\cdot i\circ b\cdot\lambda:=\{(h+t\cdot i\circ b)\cdot 
\lambda:h\cdot\lambda\in A,\ (h+t\cdot i\circ b)\cdot\lambda\in E 
\}$. The {\it Fomin type derivative} $d_{i\circ b}\bmu$ is by 
\cite{Ry02}, Proposition C.4, a finite signed measure on $(E,{\cal B} 
(E))$. 

Set $\nu+t\cdot i\circ b\cdot\lambda=(h+t\cdot i\circ b)\cdot 
\lambda$ if $\nu=h\cdot\lambda\in E$. We say that the measure $\bmu$ 
is {\it Skorokhod type differentiable} in direction of $i\circ b$, 
$b\in X$ with $\int b\, d\lambda=0$, if for all $f\in C_b(E)$ the 
function 
\begin{eqnarray*}
t\mapsto\int f(\nu-t\cdot i\circ b\cdot\lambda)\, \bmu (d\nu)
\end{eqnarray*} 
is differentiable, where by definition, $f(\nu+t\cdot i\circ b\cdot 
\lambda)=0$ if $\nu+t\cdot i\circ b\cdot\lambda\not\in E$. By 
\cite{Bo07}, Theorem 8.7.1, there is a finite signed measure $\bnu$ 
on $(E,{\cal B}(E))$ such that 
\begin{eqnarray*} 
\int f\, d\bnu=-\lim_{t\to 0}\frac{\int f(\nu+t\cdot i\circ b\cdot 
\lambda)\, \bmu (d\nu)-\int f(\nu)\, \bmu(d\nu)}{t}\, ,\quad  f\in 
C_b(E). 
\end{eqnarray*} 
We call $\bnu$ the {\it Skorokhod type derivative} of $\bmu$ in 
direction of $i\circ b$, $b\in X$ with $\int b\, d\lambda=0$. 

Formally, if $f$ is the indicator $\chi_A$ of some open set $A\in 
{\cal B}(E)$ then $d_{i\circ b}\bmu(f)$ coincides with $d_{i\circ b} 
\bmu(A)$. For the relation between Fomin and Skorokhod derivative 
in the case of a locally convex state space, see \cite{Bo07}, 
Corollaries 3.6.6 and 3.6.7. 
\medskip

Let us specify the objects of Paragraph \ref{sec:2:2:2}. Suppose that 
$\vp:\left(\{h\cdot\lambda:h\in L^v(D)\},{\cal C}\vphantom{\D l^1} 
\right)\ni h\cdot\lambda=\mu\mapsto\vp(\mu)=g\cdot\lambda=\nu\in G$ 
generates a bounded, bijective, and Fr\'echet differentiable 
map $\psi:L^v(D)\ni h\mapsto\psi(h)=g\in\{k\in L^v(D):k\cdot\lambda 
\in G\}=:\gamma$ with bounded derivative $\nabla\psi$ such that for 
each $h\in L^v(D)$, $\nabla\psi(h):L^v(D)\mapsto\{i\circ k:k\in X\}$. 
Assume also that its inverse $\psi^{-1}$ can be extended to some in 
$L^v(D)$ open neighborhood of $\gamma$ such that $\psi^{-1}$ is 
two times Fr\'echet differentiable in every point of $\gamma$ and 
$\psi^{-1}$ as well as its derivatives are bounded on $\gamma$. 
Correspondingly extend $\vp^{-1}$ to some $\Gamma\supseteq G$. 
For some $0<\ve<1$ and $0<|t|<\ve$ set $\vp^{-1}(\nu+t\cdot i\circ 
b\cdot\lambda):=0$ if $\nu+t\cdot i\circ b\cdot\lambda\not\in\Gamma$ 
and define for $n\ge 2$
\begin{eqnarray*}
B_{n,t}\equiv B_{n,t}(\nu):=\frac1t\left(\vp^{-1}\left(\nu+t\cdot i 
\circ b_n\cdot\lambda\right)-\vp^{-1}(\nu)\right)\, ,\quad\nu\in G 
\, . 
\end{eqnarray*}
Let $h_{n,t}\cdot\lambda\equiv\nu_{n,t}:=\vp^{-1}\left(\nu+t\cdot i 
\circ b_n\cdot\lambda\right)$ and $h\cdot\lambda\equiv\nu:=\vp^{-1} 
(\nu)$. By the recent hypotheses, the limit $\lim_{t\to 0}\frac1t( 
h_{n,t}-h)$ is well-defined in $L^v(D)$. Thus the limit $B_{n,0} 
\equiv B_{n,0}(\nu):=\lim_{t\to 0}B_{n,t}(\nu)$, $\nu\in G$, $n\ge 
2$, is also well-defined. 

We pick up the situation of the preceding paragraph. We define a 
measure ${\bf m}$ on the set $\{h\cdot\lambda\in E:h\in L^v(D)\}$ 
endowed with the $\sigma$-algebra generated by all subsets $\{h 
\cdot\lambda\in E:h\in B\}$ where $B\in {\cal B}(L^v(D))$. We set 
${\bf m}(E\setminus G)=0$ and 
\begin{eqnarray*}
{\bf m}(A):=\t\bmu\circ\vp^{-1}(A)\quad\mbox{\rm for all}\quad A= 
\{h\cdot\lambda:h\in B\}\cap G\quad\mbox{\rm and}\quad B\in {\cal 
B}(L^v(D)). 
\end{eqnarray*} 

Let $r:L^v(D)\mapsto {\Bbb R}$ be a Fr\'echet differentiable test 
function and let $f(g\cdot\lambda):=r(g)$, $g\in L^v(D)$. Let us use 
the symbol $\nabla$ to denote the Fr\'echet derivative and also to 
denote $\nabla f(g\cdot\lambda):=\nabla r(g)$, $g\in L^v(D)$, as well 
as $\nabla(f\circ\vp)(h\cdot\lambda):=\nabla(r\circ\psi)(h)$, $h\in 
L^v(D)$. Suppose that $f$ as well as $\nabla f$ are bounded. Recalling 
condition (\ref{A.2}) we obtain for $n\ge 2$ 
\begin{eqnarray}\label{A.60}
&&\hspace{-.5cm}\int (\nabla f,b_n)\, d{\bf m}=\lim_{t\to 0}\frac{ 
\int f(\nu+t\cdot i\circ b_n\cdot\lambda)\, {\bf m}(d\nu)-\int f(\nu) 
\, {\bf m}(d\nu)}{t}\nonumber \\ 
&&\hspace{.5cm}=\lim_{t\to 0}\frac{\int f(\nu+t\cdot i\circ b_n 
\cdot\lambda)\, \t \bmu\circ\vp^{-1}(d\nu)-\int f(\nu)\, \t \bmu 
\circ\vp^{-1}(d\nu)}{t}\nonumber \\ 
&&\hspace{.5cm}=\lim_{t\to 0}\frac{\int f\circ\vp\left(\vp^{-1}(\nu 
+t\cdot i\circ b_n\cdot\lambda)\vphantom{\int}\right)\, \t \bmu\circ 
\vp^{-1}(d\nu)-\int f(\nu)\, \t \bmu\circ\vp^{-1}(d\nu)}{t}\nonumber 
 \\ 
&&\hspace{.5cm}=\lim_{t\to 0}\frac{\int f\circ\vp\left(\vp^{-1}(\nu) 
+t\cdot B_{n,t}(\nu)\vphantom{\int}\right)\, \t \bmu\circ\vp^{-1} 
(d\nu)-\int f(\nu)\, \t \bmu\circ\vp^{-1}(d\nu)}{t}\nonumber \\
&&\hspace{.5cm}=\lim_{t\to 0}\frac{\int f\circ\vp\left(\mu +t\cdot 
B_{n,t}\circ\vp(\mu)\vphantom{\int}\right)\, \t \bmu (d\mu)-\int f 
\circ\vp(\mu)\, \t \bmu (d\mu)}{t}\nonumber \\
&&\hspace{.5cm}=\int\left(\nabla(f\circ\vp),B_{n,0}\circ\vp\vphantom 
{l^1}\right)\, d\t \bmu\nonumber \\
&&\hspace{.5cm}=\sum_{m=1}^\infty\int\left(\nabla(f\circ\vp),i\circ 
b_m\vphantom{l^1}\right)\cdot\Psi_{n,m}\, d\t \bmu   
\end{eqnarray} 
where 
\begin{eqnarray*}
\Psi_{n,m}:={{}_{Y^\ast}\langle i^{\ast -1}\circ b_m^\ast,B_{n,0}\circ 
\vp\rangle_Y}\, ,\quad m\in {\Bbb N}. 
\end{eqnarray*} 

Let denote $L^v_\lambda(D):=\{h\cdot\lambda:h\in L^v(D)\}$. Now we 
specify $\t\bmu$. We assume $\t \bmu$ to be defined by (\ref{A.5}) 
with $\t \mu_m=d_m\cdot {\it leb}$. We assume furthermore $d_m\in 
C^1_b({\Bbb R})$ and $d_m>0$, $m\in {\Bbb N}$, as well as $\lim_{|x| 
\to\infty}d_m(x)=0$ such that 
\begin{eqnarray}\label{A.605} 
\hat{\delta}(b_n)\circ\vp:=\sum_{m=1}^\infty \left(\Psi_{n,m}\cdot 
\left(\frac{d'_m}{d_m}\right)\circ f_m+\frac{\partial\Psi_{n,m}} 
{\partial i\circ b_m}\right)\, ,\quad n\ge 2,
\end{eqnarray} 
converges in $L^1(L^v_\lambda(D),\t\bmu)$. Next we remind of $f_m 
(\alpha\cdot\lambda)={}_{Y^\ast}\langle i^{\ast -1}\circ b_m^\ast, 
\alpha\rangle_Y$, $\alpha\in Y$, $m\in {\Bbb N}$. We obtain for 
cylindrical test functions of the form 
\begin{eqnarray*} 
f\circ\vp\equiv\Phi\left({}_{Y^\ast}\langle i^{\ast -1}\circ 
b_{m_1}^\ast,\, \cdot\, \rangle_Y,\ldots\, ,{}_{Y^\ast}\langle i^{ 
\ast -1}\circ b_{m_k}^\ast,\, \cdot\, \rangle_Y\right) 
\end{eqnarray*} 
and $n\ge 2$
\begin{eqnarray}\label{A.61} 
&&\hspace{-.5cm}\int (\nabla f,b_n)\, d{\bf m}\nonumber \\ 
&&\hspace{.5cm}=-\sum_{j=1}^k\int\Phi\left(x_{m_1},\ldots\, ,x_{m_k} 
\right)\Psi_{n,m_j}\left({\T\sum_{m'=1}^\infty}x_{m'}\ i\circ b_{m'} 
\right)\frac{d'_{m_j}(x_{m_j})}{d_{m_j}(x_{m_j})}\, \t \mu(dx_{m_1}), 
\ldots\, ,\t \mu(dx_{m_k})\nonumber \\ 
&&\hspace{1cm}-\sum_{j=1}^k\int\Phi\left(x_{m_1},\ldots\, ,x_{m_k} 
\right)\frac{\partial\Psi_{n,m_j}\left({\T\sum_{m'=1}^\infty}x_{m'} 
\ i\circ b_{m'} \right)}{\partial x_{m_j}}\, \t \mu(dx_{m_1}),\ldots 
\, ,\t \mu(dx_{m_k})\nonumber \\ 
&&\hspace{.5cm}=-\int f\circ\vp\ \sum_{j=1}^k\left(\Psi_{n,m_j}\cdot 
\left(\frac{d'_{m_j}}{d_{m_j}}\right)\circ f_{m_j}+\frac{\partial 
\Psi_{n,m_j}}{\partial i\circ b_{m_j}}\right)\, d\t \bmu\nonumber \\ 
&&\hspace{.5cm}=-\int f\hat{\delta}(b_n)\, d{\bf m}\, .  
\end{eqnarray} 
The first line of (\ref{A.60}) together with the last line of 
(\ref{A.61}) including the summability of (\ref{A.605}), describe the 
form of differentiability we are implicitly imposing on ${\bf m}$ in 
direction of $i\circ b_n$, $n\ge 2$. One may compare it with the above 
defined {\it Skorokhod type differentiability}. 
\medskip 

Let ${\cal B}_G(E)$ be the $\sigma$-algebra on $E$ generated by $E 
\setminus G$ and $\{A\cap G:A\in \mathcal{B}(E)\}$. Let 
\begin{eqnarray*}
\t \delta(b_n):=E_{\bf m}\left.\left(\hat{\delta}(b_n)\right|{\cal B}_G 
(E)\right) 
\end{eqnarray*} 
denote the conditional expectation of $\hat{\delta}(b_n)$ under the 
$\sigma$-algebra ${\cal B}_G(E)$ with respect to the measure ${\bf 
m}$. Let $\bmu$ be the measure defined on $(E,{\cal B}(E))$ given 
by the restriction of ${\bf m}$ to ${\cal B}_G(E)$. By construction 
this is compatible with (\ref{A.6}).
\begin{proposition}\label{Proposition2.3} 
In the setup of the present paragraph we have 
\begin{eqnarray*}
\int (Df,b_n)\, d\bmu=-\int f\, \t \delta (b_n)\, d\bmu\, ,\quad f\in 
{\Bbb C}^{1,1}_b(E;X,{\Bbb R}),\ n\ge 2. 
\end{eqnarray*} 
\end{proposition} 
Proof. {\it Step 1 } We denote by $\ve_l$ the usual mollifier function 
on ${\Bbb R}^k$ whose support is the closed ball in ${\Bbb R}^l$ about 
the origin with radius $1/l$, $l\in {\Bbb R}$. For a real function $F$ 
on ${\Bbb R}^k$ let us write 
\begin{eqnarray*} 
&&\hspace{-.5cm}F\left({}_{Y^\ast}\langle i^{\ast -1}\circ b_1^\ast,\, 
\cdot\, \rangle_Y,\ldots\, ,{}_{Y^\ast}\langle i^{\ast -1}\circ b_k^\ast, 
\, \cdot\, \rangle_Y\right)\ast\delta_l \\ 
&&\hspace{.5cm}:=F\ast\ve_l\left({}_{Y^\ast}\langle i^{\ast -1}\circ 
b_1^\ast,\, \cdot\, \rangle_Y,\ldots\, ,{}_{Y^\ast}\langle i^{\ast -1} 
\circ b_k^\ast,\, \cdot\, \rangle_Y\right)\, . 
\end{eqnarray*} 
Furthermore let ${\cal C}_k$ denote the $\sigma$-algebra on 
$L^v_\lambda(D)\equiv\{h\cdot\lambda:h\in L^v(D)\}$ generated by all 
$\{\{h\cdot\lambda:h\in B\}: B\in {\cal C}\}$-measurable functions on 
$L^v_\lambda(D)$ which are independent of $\sum_{m=k+1}^\infty 
{}_{Y^\ast}\langle i^{\ast -1}\circ b_m^\ast,\, \cdot\, \rangle_Y\, i 
\circ b_m$. Set $\|h\cdot\lambda\|_{L^v_\lambda(D)}:=\|h\|_{L^v(D)}$. 

Let $f\in {\Bbb C}^{1,1}_b(E;X,{\Bbb R})$ and let $E_{\t \bmu}\left. 
\left(\cdot\right|{\cal C}_k\right)$ denote the conditional expectation 
under the $\sigma$-algebra ${\cal C}_k$ with respect to the measure 
${\t \bmu}$. Using, among other things, Doob's martingale convergence 
theorem it is standard to verify that there is a strictly increasing 
sequence $l(k)\equiv l(k,f)$, $k\in {\Bbb N}$, such that 
\begin{eqnarray*}
\Phi^{(k)}:=\left(E_{\t \bmu}\left.\left(f\circ\vp\right|{\cal C}_k 
\right)\right)\ast\delta_{l(k)}\stack{k\to\infty}{\lra}f\circ\vp 
\end{eqnarray*} 
$\t \bmu$-a.e., uniformly boundedly in $L^\infty(L^v_\lambda(D),\t 
\bmu)$, and hence boundedly in $L^1(L^v_\lambda(D),\t \bmu)$. By the 
definition of $B_{n,0}$ and the properties of $\vp$, for $n\ge 2$ we 
obtain 
\begin{eqnarray*}
\left(\nabla\Phi^{(k)},B_{n,0}\circ\vp\vphantom{l^1}\right)=\left( 
E_{\t \bmu}\left.\left(\lim_{t\to 0}\frac1t\left(f\circ\vp\left(\cdot 
+t\cdot B_{n,t}\circ\vp\right)-f\circ\vp\vphantom{l^1}\right)\right| 
{\cal C}_k\right)\right)\ast\delta_{l(k)}\, . 
\end{eqnarray*} 
Using again Doob's martingale convergence theorem, it follows that 
\begin{eqnarray}\label{A.63}
\left(\nabla\Phi^{(k)},B_{n,0}\circ\vp\vphantom{l^1}\right)\stack 
{k\to\infty}{\lra}\lim_{t\to 0}\frac1t\left(f\circ\vp\left(\cdot 
+t\cdot B_{n,t}\circ\vp\right)-f\circ\vp\vphantom{l^1}\right)\quad 
\mbox{\rm in } L^1(L^v_\lambda(D),\t\bmu)\, . 
\end{eqnarray} 

\noindent
{\it Step 2 } Reviewing (\ref{A.60}) as well as (\ref{A.61}) and 
applying (\ref{A.63}), we obtain 
\begin{eqnarray*}
&&\hspace{-.5cm}\int (Df,b_n)\, d\bmu=\int (Df,b_n)\, d{\bf m}= 
\int\lim_{t\to 0}\frac1t\left(f(\nu+t\cdot i\circ b_n\cdot\lambda) 
-f(\nu)\vphantom{l^1}\right)\, {\bf m}(d\nu) \\ 
&&\hspace{.5cm}=\int\lim_{t\to 0}\frac1t\left(f\circ\vp\left(\mu 
+t\cdot B_{n,t}\circ\vp(\mu)\right)-f\circ\vp(\mu)\vphantom{l^1} 
\right)\, \t\bmu(d\mu) \\ 
&&\hspace{.5cm}=\lim_{k\to\infty}\int\left(\nabla\Phi^{(k)},B_{n,0} 
\circ\vp\vphantom{l^1}\right)\, d\t \bmu\nonumber \\ 
&&\hspace{.5cm}=-\lim_{k\to\infty}\int\Phi^{(k)}\ \sum_{m=1}^k\left( 
\Psi_{n,m}\cdot\left(\frac{d'_m}{d_m}\right)\circ f_m+\frac{\partial 
\Psi_{n,m}}{\partial i\circ b_m}\right)\, d\t \bmu\nonumber \\ 
&&\hspace{.5cm}=-\int f\circ\vp\ \sum_{m=1}^\infty\left(\Psi_{n,m} 
\cdot\left(\frac{d'_m}{d_m}\right)\circ f_m+\frac{\partial\Psi_{n,m}} 
{\partial i\circ b_m}\right)\, d\t \bmu\nonumber \\ 
&&\hspace{.5cm}=-\int f\, \hat{\delta}(b_n)\, d{\bf m}=-\int f\, \t 
\delta (b_n)\, d\bmu\, ,\quad n\ge 2\, ; 
\end{eqnarray*} 
for the second last equality sign use uniform boundedness of $\Phi^{ 
(k)}$, $k\in {\Bbb N}$, in $L^\infty(L^v_\lambda(D),\t \bmu)$, and the 
convergence in (\ref{A.605}) with respect to $L^1(L^v_\lambda(D),\t\bmu)$. 
\qed

\subsubsection{The trace of the gradient}\label{sec:2:2:4} For $g\in 
{\Bbb C}^{p,1}_b(E)$ let $g_n:={}_{Y^\ast}\langle i^{\ast -1}\circ 
b_n^\ast,g\rangle_Y$, $n\in {\Bbb N}$. Provided that the sum 
$\sum_{n=2}^\infty (Dg_n,b_n)$ is finite, let us call this sum the 
{\it trace} of the operator $Dg(\cdot)\equiv (Dg,\cdot)\in B(X\ominus 
\1,Y)$ where $X\ominus\1:=\{b\in X:(i\circ b,\1)=0\}$. Note that 
below we use both notations $(\cdot\, ,\, \cdot)$ and ${}_{Y^\ast} 
\langle\cdot\, ,\, \cdot\rangle_{Y}$ in order to better distinguish 
between the representation of an operator applied to some (test) 
function from certain coordinate representations using the Schauder 
basis. 

As before, let $\nu+t\cdot i\circ b\cdot\lambda=(h+t\cdot i\circ b)\cdot 
\lambda$ if $\nu=h\cdot\lambda\in F$ and $b\in X$. In the present paragraph, 
by (\ref{A.2}) and (i)-(iii) as well as Lemma \ref{LemmaA.1}, for $g\in 
{\Bbb C}^{p,1}_b(E)$ the gradient $Dg$ and the directional derivatives 
$\lim_{t\to 0}\frac1t\left(g(\nu+t\cdot i\circ b_n\cdot\lambda)-g(\nu) 
\right)$ are $\bmu$-a.e. well-defined for all $n\ge 2$ but not necessarily 
for $n=1$. 
\medskip 

The following calculation illustrates the choice of the term {\it trace} 
of the operator $Dg(\cdot)\in B(X\ominus\1,Y)$. The independence 
relative to the Schauder basis $b_n\in X$, $n\in {\Bbb N}$, under certain 
conditions will be a side result of Proposition \ref{Proposition2.2}, see 
Remark (5) below. Let $g\in {\Bbb C}^{p,1}_b(E)$ with $\sum_{n=2}^\infty 
(Dg_n,b_n)<\infty$ $\bmu$-a.e. We have $\bmu$-a.e. 
\begin{eqnarray}\label{A.65}
&&\hspace{-.5cm}\sum_{n=2}^\infty (Dg_n(\nu),b_n)=\sum_{n=2}^\infty\frac{ 
\partial g_n}{\partial b_n}(\nu)=\sum_{n=2}^\infty\lim_{t\to 0}\frac1t 
\left(g_n(\nu+t\cdot i\circ b_n\cdot\lambda)-g_n(\nu)\right)\nonumber \\ 
&&\hspace{.5cm}=\sum_{n=2}^\infty\lim_{t\to 0}\frac1t\ {}_{Y^\ast}\left 
\langle i^{\ast -1}\circ b_n^\ast\, ,\, g(\nu+t\cdot i\circ b_n\cdot\lambda)- 
g(\nu)\right\rangle_{Y}\nonumber \\ 
&&\hspace{.5cm}=\sum_{n=2}^\infty\ {\vphantom{\lim_{t\to 0}}}_{Y^\ast} 
\left\langle i^{\ast -1}\circ b_n^\ast\, ,\, \lim_{t\to 0}\frac1t\left(g(\nu+ 
t\cdot i\circ b_n\cdot\lambda)-g(\nu)\right)\right\rangle_{Y}\nonumber \\ 
&&\hspace{.5cm}=\sum_{n=2}^\infty\ {}_{Y^\ast}\left\langle i^{\ast -1}\circ 
b_n^\ast\, ,\, \left(Dg(\nu),b_n\right)\right\rangle_{Y}\, . 
\end{eqnarray}

\subsubsection{Conditions sufficient for (j)-(jjj) separating requirements 
on $\bmu$ from the requirements on $A^f$}\label{sec:2:2:5} We are no 
longer focusing on the particular form (\ref{A.5}), (\ref{A.6}) of the 
measure $\bmu$. In the remainder of this subsection we shall work with 
the following set of hypotheses motivated by Paragraphs 
\ref{sec:2:2:2}-\ref{sec:2:2:4}. We remind of the definition of $a^f$ 
in the beginning of Subsection \ref{sec:2:1}. 
\begin{itemize}
\item[(j')] There exists $p>1$ such that 
\begin{eqnarray*}
a^f\in {\Bbb C}^{p,1}_{b,{\rm loc}}(E)\, .  
\end{eqnarray*}
Furthermore, the trace of the operator $Da^f(\cdot)\equiv (Da^f,\cdot) 
\in B(X\ominus\1,Y)$ converges in $L^1(E,\bmu)$ to an element belonging 
to $L^1(E,\bmu)\cap L_{\rm loc}^\infty(E,\bmu)$. 
\item[(jj')] For $n\ge 2$, there is a unique {\it derivative } 
$\t \delta (b_n)\in L^1(E,\bmu)$ {\it of the measure $\bmu$ in direction 
of $\ i\circ b_n\cdot\lambda$} such that 
\begin{eqnarray*}
\int (Df,b_n)\, d\bmu=-\int f\, \t \delta (b_n)\, d\bmu\, ,\quad f\in 
{\Bbb C}^{1,1}_b(E;X,{\Bbb R})\equiv C^{1,1}_b(E). 
\end{eqnarray*}
\item[(jjj')] For any {\it integrand} $\vp\in {\Bbb C}^{1,1}_b(E)$, with 
$(\vp,\1)=0$ $\bmu$-a.e., i. e. $\vp=\sum_{n=2}^\infty {}_{Y^\ast}\langle 
i^{\ast -1} \circ b_n^\ast,\vp \rangle_Y\, i\circ b_n$, the {\it Ogawa 
type stochastic integral} 
\begin{eqnarray*} 
-\sum_{n=2}^\infty\ {}_{Y^\ast}\langle i^{\ast -1}\circ b_n^\ast\, ,\, 
\vp\rangle_{Y}\, \t \delta (b_n) 
\end{eqnarray*} 
converges in $L^1(E,\bmu)$ to an element belonging to $L^1(E,\bmu)\cap 
L_{\rm loc}^\infty(E,\bmu)$.
\end{itemize} 
Let us remind of (\ref{A.2}). As in Subsection \ref{sec:2:1}, we shall 
use the symbol $a^f\equiv a^f(\cdot)$ for both, an element of $X$ as well 
as the embedded element in $Y$. Also recall $(a^f,\1)=0$ $\bmu$-a.e. on 
$E$ as stated in Subsection \ref{sec:2:1}. For $\nu\in E$, let us use the 
notations $a^f\equiv a^f(\nu)=\sum_{n=2}^\infty a_n\cdot b_n\in X$ and 
$a^f\equiv a^f(\nu)=\sum_{n=2}^\infty a_n\cdot i\circ b_n\in Y$ where the 
\begin{eqnarray*}
a_n\equiv a_n(\nu)={}_{X^\ast}\langle b_n^\ast,a^f(\nu)\rangle_X={}_{ 
Y^\ast}\langle i^{\ast -1}\circ b_n^\ast,a^f(\nu)\rangle_Y\in {\Bbb R} 
\, ,\quad n\ge 2, 
\end{eqnarray*}
are the coefficients of $a^f$ with respect to the Schauder bases $b_n$, 
$n\in {\Bbb N}$, in $X$ and $i\circ b_n$, $n\in {\Bbb N}$, in $Y$. 
\begin{proposition}\label{Proposition2.2} {\rm (Skorokhod type stochastic 
integral equals Ogawa type integral minus trace of the gradient of the 
integrand.)} Suppose (j')-(jjj'). We have (j)-(jjj) and 
\begin{eqnarray*}
\delta(A^f)=\sum_{n=2}^\infty\left((Da_n,b_n)+a_n\t \delta(b_n)\right) 
\end{eqnarray*} 
where the infinite sum converges in $L^1(E,\bmu)$. 
\end{proposition} 
Proof. By means of condition 
(j') we verify $a_n={}_{Y^\ast}\langle i^{\ast -1}\circ b_n^\ast,a^f 
\rangle_Y\in {\Bbb C}^{p,1}_{b,{\rm loc}}(E;X,{\Bbb R})$, $n\in {\Bbb 
N}$. The sum $a^f\equiv a^f (\nu)=\sum_{n=2}^\infty a_n\circ b_n$ 
converges $\bmu$-a.e. in $X$-norm. Keeping in mind that with the basis 
constant $K$ of $b_n$, $n\in {\Bbb N}$, we have 
\begin{eqnarray*}
\left\|\sum_{n=2}^N a_n(\nu)\cdot b_n\right\|_X\le K\|a^f(\nu)\|_X
\end{eqnarray*}
for all $N\ge 2$ and $\bmu$-a.e. $\nu\in E$, cf. \cite{AK06}, 
Proposition 1.1.4 and Definition 1.1.5, the sum $\sum_{n=2}^\infty a_n 
\cdot b_n$ converges even in $L^1(E,\bmu;X)$. Thus, for $f\in C^{q,1}_b 
(E)$ the sum ${\T\sum_{n=2}^\infty}(Df,b_n)\cdot a_n$ converges in 
$L^1(E,\bmu)$ and 
\begin{eqnarray*}
&&\hspace{-.5cm}\int (Df,a^f)\, d\bmu=\int\left(Df,\sum_{n=2}^\infty 
a_n\cdot b_n\right)\, d\bmu=\int\sum_{n=2}^\infty(Df,b_n)\cdot a_n\, 
d\bmu\, . 
\end{eqnarray*}
Lemma \ref{Lemma2.2} for $f\in C^{q,1}_b(E)\equiv {\Bbb C}^{q,1}_b(E;X,
{\Bbb R})$ as well as $a_n\in {\Bbb C}^{p,1}_{b,{\rm loc}}(E;X,{\Bbb R}) 
\subseteq {\Bbb C}^{p,1}_b(E;X,{\Bbb R})$ implies $f\cdot a_n\in {\Bbb C 
}^{p\wedge q,1}_b(E;X,{\Bbb R})$, $n\ge 2$. Condition (j') says that 
$\sum_{n=2}^\infty(Da_n,b_n)$ converges in $L^1(E,\bmu)$. Thus Lemma 
\ref{Lemma2.2} yields the convergence of $\sum_{n=2}^\infty(D(f\cdot a_n) 
,b_n)$ in $L^1(E,\bmu)$ and 
\begin{eqnarray*} 
&&\hspace{-.5cm}\int (Df,a^f)\, d\bmu=-\int f\cdot\sum_{n=2}^\infty 
(Da_n,b_n)\, d\bmu+\int\sum_{n=2}^\infty(D(f\cdot a_n),b_n)\, d\bmu \\ 
&&\hspace{.5cm}=-\int f\cdot\sum_{n=2}^\infty(Da_n,b_n)\, d\bmu+\sum_{n 
=2}^\infty\int(D(f\cdot a_n),b_n)\, d\bmu \\ 
&&\hspace{.5cm}=-\int f\cdot\sum_{n=2}^\infty(Da_n,b_n)\, d\bmu-\sum_{n 
=2}^\infty\int f\cdot a_n\, \t \delta(b_n)\vphantom{l^1}\, d\bmu\, , 
\quad f\in C^{q,1}_b(E), 
\end{eqnarray*}
the last line by (jj'). Condition (jjj') says that $\sum_{n=2}^\infty 
a_n\t \delta(b_n)$ converges in $L^1(E,\bmu)$ and thus  
\begin{eqnarray*} 
\int (Df,a^f)\, d\bmu=-\int f\cdot\sum_{n=2}^\infty\left((D a_n,b_n)+ 
a_n\, \t \delta(b_n)\right)\, d\bmu\, ,\quad f\in C^{q,1}_b(E). 
\end{eqnarray*}
This also verifies (jj). Together with (j') and (jjj') we even get (jjj). 
Condition (j) is a direct consequence of (j'). 
\qed 
\bigskip 

\nid
{\bf Remarks }(4) {\it ($\bmu$-a.e. existence of the Ogawa type integral 
not assuming (jjj')) } Suppose for this remark that there is a constant 
$0<C<\infty$ such that $\|\t \delta (b_n)\|_{L^1(E,\sbmu)}\le C\|b_n 
\|_X$ for all $n\ge 2$. This condition corresponds to a basic property 
in Wiener space analysis; the interval $[0,1]$, the space $\{f\in C([0, 
1]):f(0)=0\}$, and the Cameron-Martin space play there, respectively, 
the role of $D$, $Y$, and $X$. 

Let $\vp\in {\Bbb C}^{1,1}_b(E)$. Assume that there is a sequence $c_n 
\equiv c_n(\vp)>0$, $n\ge 2$, of constants with $\sum_{n=2}^\infty 
c_n\|b_n\|_X<\infty$ such that for $E_n:=\left\{\nu\in E:|{}_{Y^\ast} 
\langle i^{\ast -1}\circ b_n^\ast\, ,\, \vp(\nu)\rangle_{Y}|>c_n\right 
\}$ it holds that 
\begin{eqnarray}\label{A.7} 
\sum_{n=2}^\infty\bmu(E_n)<\infty\, . 
\end{eqnarray} 
Let $E^{(n)}:=\left\{\nu\in E:|{}_{Y^\ast}\langle i^{\ast -1}\circ 
b_k^\ast\, ,\, \vp (\nu)\rangle_Y|\le c_k\ \mbox{\rm for all}\ k\ge 
n\right\}$, $n\ge 2$, and let $\bmu_n$ denote the restriction of the 
measure $\bmu$ to $E^{(n)}$. By the Borel-Cantelli lemma and (\ref{A.7}) 
it holds that $\bmu(\bigcup_{n\ge 2}E^{(n)})=\bmu(\bigcap_{n\ge 2} 
\bigcup_{k\ge n}E_k)=1$. Furthermore by $\|\t \delta (b_n)\|_{L^1(E, 
\sbmu)}\le C\|b_n\|_X$, $n\ge 2$, the sum   
\begin{eqnarray*}
-\sum_{k=2}^\infty\ {}_{Y^\ast}\langle i^{\ast -1}\circ b_k^\ast\, ,\, 
\vp\rangle_{Y}\, \t \delta (b_k)\quad\mbox{\rm converges in}\ L^1(E^{(n)} 
,\bmu_n) 
\end{eqnarray*} 
for all $n\ge 2$. As a consequence, for each $n\ge 2$ there is a subsequence 
$N_k(n)$, $k\in {\Bbb N}$, of $N_k(n-1)$, $k\in {\Bbb N}$, such that 
\begin{eqnarray*}
-\sum_{j=2}^{N_k(n)}\ {}_{Y^\ast}\langle i^{\ast -1}\circ b_j^\ast\, ,\, 
\vp\rangle_{Y}\, \t \delta(b_j)\quad\mbox{\rm converges}\ \bmu_n\mbox 
{\rm -a.e. on}\ E^{(n)} 
\end{eqnarray*} 
as $k\to\infty$. Thus, 
\begin{eqnarray*}
-\sum_{n=2}^{N_n(n)}\ {}_{Y^\ast}\langle i^{\ast -1}\circ b_n^\ast\, ,\, 
\vp\rangle_{Y}\, \t \delta (b_n)\quad\mbox{\rm converges}\ \bmu\mbox 
{\rm -a.e. on}\ E 
\end{eqnarray*} 
as $n\to\infty$ to the Ogawa type integral $-\sum_{n=2}^\infty\ {}_{Y^\ast} 
\langle i^{\ast -1}\circ b_n^\ast\, ,\, \vp\rangle_{Y}\, \t \delta (b_n)$.
\bigskip

\noindent
(5) For this remark, let $0<v<\infty$, i. e. $Y\equiv L^v(D)$ is 
reflexive. Suppose (jj'), (jjj') and, to be compatible with Remark (4),
that there is a constant $0<C<\infty$ such that $\|\t \delta (b_n)\|_{L^1 
(E,\sbmu)}\le C\|b_n\|_X$ for all $n\ge 2$. Suppose moreover that $\gamma 
\equiv\gamma(\nu):=\sum_{n=2}^\infty\t \delta(b_n)(\nu)\ i\circ b_n$ 
converges $\bmu$-a.e. in $Y=L^v(D)$. In the same way as in Remark (4), 
this condition relates also to a basic property in Wiener space analysis. 
The objective of this remark is to demonstrate that, in the present 
situation, {\it the Ogawa type integral and the trace of the gradient of 
the integrand are independent of any Schauder basis $b_n$, $n\in {\Bbb N}$, 
of $X$ satisfying (\ref{A.2}).} 

We have $\bmu$-a.e. 
\begin{eqnarray*}
\t \delta(b_n)={}_{Y^\ast}\langle i^{\ast-1}\circ b_n^\ast,\gamma 
(\cdot)\rangle_Y\, ,\quad n\ge 2.
\end{eqnarray*} 
Conversely, the relation $\int (Df,b_n)\, d\bmu=-\int f\, {}_{Y^\ast 
}\langle i^{\ast-1}\circ b_n^\ast,\gamma\rangle_Y\, d\bmu$ for all 
$n\ge 2$ and $f\in {\Bbb C}^{1,1}_b(E;X,{\Bbb R})$ for one fixed 
basis of $X$ with (\ref{A.2}) provides this equality  with the same 
$\gamma$ for any basis of $X$ with (\ref{A.2}). In other words, the 
$\bmu$-a.e. convergence of $\gamma\equiv\gamma(\nu):=\sum_{n=2}^\infty 
\t \delta(b_n)(\nu)\ i\circ b_n$ in $Y=L^v(D)$ is independent of the 
chosen basis of $X$ with (\ref{A.2}).

Now let us fix the class of integrands $\vp$ for this remark. Let  
$\vp\in {\Bbb C}^{1,1}_b(E)$ with $(\vp,\1)=0$ $\bmu$-a.e. Suppose 
that with $\vp_n\equiv\vp_n(\nu):={}_{Y^\ast}\langle i^{\ast -1} 
\circ b_n^\ast\, ,\, \vp\rangle_{Y}$, $n\ge 2$, the sum 
\begin{eqnarray*}
\vp^\ast\equiv\vp^\ast(\nu):=\sum_{n=2}^\infty\vp_n\ i^{\ast -1}\circ 
b_n^\ast\quad\mbox{\rm converges}\quad\bmu\mbox{\rm -a.e. in } Y^\ast. 
\end{eqnarray*} 

For this class of integrands $\vp$ we have $\bmu$-a.e. 
\begin{eqnarray*}
\vp_n={}_{Y^\ast}\langle \vp^\ast\, ,\, i\circ b_n\rangle_Y\, , 
\quad n\ge 2.
\end{eqnarray*} 
Since $Y^\ast\simeq L^w(D)$ where $1/v+1/w=1$ and $Y$ is reflexive, 
by R. C. James' theorem, the functionals $i^{\ast -1}\circ b_n^\ast 
\in Y^\ast$, $n\in {\Bbb N}$, form a Schauder basis of $Y^\ast$ 
which can be represented by a Schauder basis $\beta_n^\ast$, $n\in 
{\Bbb N}$, of $L^w(D)$, cf. \cite{AK06}, Theorem 3.2.19 together 
with Proposition 3.2.8. Its biorthogonal functionals can be 
represented by $i\circ b_n$, $n\in {\Bbb N}$, and $\vp^\ast\in 
Y^\ast$ can be identified with some $f^\ast\in L^w(D)$. The $\bmu 
$-a.e. convergence of $\vp^\ast\equiv\vp^\ast(\nu)=\sum_{n=2}^\infty 
\vp_n\ i^{\ast -1}\circ b_n^\ast\quad\mbox{\rm in}\quad Y^\ast$ is 
therefore equivalent with the $\bmu$-a.e. convergence of 
\begin{eqnarray*}
f^\ast\equiv f^\ast(\nu)=\sum_{n=2}^\infty(f^\ast,i\circ b_n)\,  
\beta_n^\ast\quad\mbox{\rm in }L^w(D). 
\end{eqnarray*} 

The latter convergence is independent of the Schauder basis 
$\beta_n^\ast$, $n\in {\Bbb N}$, of $L^w(D)$ and of its 
biorthogonal functionals, the basis $i\circ b_n$, $n\in {\Bbb N}$, 
of $Y\equiv L^v(D)$. Thus the convergence of 
\begin{eqnarray*}
\vp^\ast=\sum_{n=2}^\infty {}_{Y^\ast}\langle\vp^\ast\, ,\, i\circ b_n 
\rangle_Y\ i^{\ast -1} \circ b_n^\ast\quad\mbox{\rm in}\quad Y^\ast 
\end{eqnarray*} 
is independent of the Schauder basis $i\circ b_n$, $n\in {\Bbb N}$, of 
$Y$. Noting that the above identity ${}_{Y^\ast}\langle\vp^\ast\, ,\, 
i\circ b_n\rangle_Y={}_{Y^\ast}\langle i^{\ast-1}\circ b_n^\ast\, ,\, 
\vp\rangle_{Y}$ holds independent of the basis $i\circ b_n$, $n\in 
{\Bbb N}$, of $Y$, we can finally state that the 
convergence of 
\begin{eqnarray*}
\vp^\ast=\sum_{n=2}^\infty {}_{Y^\ast}\langle i^{\ast -1}\circ b_n^\ast 
\, ,\, \vp\rangle_{Y}\ i^{\ast -1} \circ b_n^\ast\quad\mbox{\rm in}\quad 
Y^\ast 
\end{eqnarray*} 
is independent of the Schauder basis $i\circ b_n$, $n\in {\Bbb N}$, 
of $Y$. 

As in Proposition \ref{Proposition2.2} we get $\bmu$-a.e. the 
representation 
\begin{eqnarray*}
-\delta(\vp\cdot\lambda)=-\sum_{n=2}^\infty((D\vp_n,b_n)-\vp_n\t \delta 
(b_n))  
\end{eqnarray*} 
of the Skorokhod type integral. For the Ogawa type integral we obtain 
$\bmu$-a.e. 
\begin{eqnarray*}
-\sum_{n=2}^\infty\vp_n\t \delta(b_n)=-\sum_{n=2}^\infty {}_{Y^\ast} 
\langle \vp^\ast\, ,\, i\circ b_n\rangle_Y\ {}_{Y^\ast}\langle i^{\ast 
-1}\circ b_n^\ast,\gamma\rangle_Y=-{}_{Y^\ast}\langle \vp^\ast,\gamma 
\rangle_Y\, . 
\end{eqnarray*} 
The two last relations show that under the conditions formulated for 
this remark both sums, $-\sum_{n=2}^\infty\vp_n\t \delta(b_n)$ as well 
as $\sum_{n=2}^\infty(D\vp_n,b_n)$, are $\bmu$-a.e. independent of any 
Schauder basis $b_n$, $n\in {\Bbb N}$, satisfying (\ref{A.2}). 

\subsection{Change of Measure and Quasi-invariance}\label{sec:2:12}

Let us return to the general setup. Hypotheses (j)-(jjj) are the general 
assumptions on $\bmu$ in the paper in order to establish absolute 
continuity under the flow $W^+$. The subsequent conditions (jv) and (v) 
focus mainly on the {\it sensitivity} of the trajectory $\nu_t\equiv 
U^+(t,\nu)$, $t\in [-1,1]$, with respect to the initial value. Conditions 
(jv) and (v) can either be verified directly for the particular PDE as 
for example in Section \ref{sec:4}, or by using a certain Fr\'echet 
differentiability, cf. Proposition \ref{Proposition2.4} and Section 
\ref{sec:3}. 

Throughout the remainder of this section, let $p>1$ be the number specified 
in condition (j) and $q$ be given by 
\begin{eqnarray*}
\frac1p+\frac1q=1\, . 
\end{eqnarray*} 

For definiteness in the conditions (jv),(v),(jv'), and (v') below we 
introduce the following notion. We call a set $U\subseteq X$ {\it 
pre-neighborhood} of $0\in X$ in the topology of $X$ if for each $h 
\in X$ there is an $a\equiv a(h)>0$ such that $bh\in U$ for all $b\in 
(-a,a)$. Clearly, a neighborhood $U$ of $0\in X$ in the topology of $X$ 
is a pre-neighborhood. 
\begin{itemize}
\item[(jv)] For $\bmu$-a.e. $\nu\equiv\t h\cdot\lambda\in G$ there is a 
pre-neighborhood $U_{\t h}(0)$ of $0\in X$ in the topology of $X$ such 
that in Case 1 $\{\t h+i\circ h:h\in U_{\t h}(0)\}\subseteq D(V_1^+)$, 
and in Case 2 $\{\t h+i\circ h:h\in U_{\t h}(0)\}\subseteq\{V_1^+(h):h 
\in D (V_1^+)\}$. 

For $t\in [0,1]$ in Case 1 and $t\in [-1,0]$ in Case 2, $k\in X$, $g\in C_0 
(D)$, and $\bmu$-a.e. $\nu\in G$ we have   
\begin{eqnarray}\label{2.2} 
\left.\frac{d}{du}\right|_{u=0}\left((\nu+u\cdot i\circ k\cdot\lambda)_t,g 
\right)={}_{X^\ast}\langle r(\nu ,t,g),k\rangle_{X}
\end{eqnarray} 
for some $r(\cdot,t,g)\in L^q(E,\bmu;X^\ast)$. For fixed $g\in C_0(D)$, the 
map $t\mapsto r(\cdot,t,g)\in L^q(E,\bmu;X^\ast)$ is continuous on $t\in [0, 
1]$ in Case 1 and on $t\in [-1,0]$ in Case 2. 
\item[(jv')] For $\bmu$-a.e. $\nu\equiv\t h\cdot\lambda\in G$ there is a 
pre-neighborhood $U_{\t h}(0)$ of $0\in X$ in the topology of $X$ such that 
$\{\t h+i\circ h:h\in U_{\t h}(0)\}\subseteq D(V_1^+)\cap \{V_1^+(h):h\in D 
(V_1^+)\}$. Furthermore, all requirements of condition (jv) hold for all 
$t\in [-1,1]$, regardless of whether we have Case 1 or Case 2. 
\item[(v)] For $t\in [0,1]$ in Case 1 and $t\in [-1,0]$ in Case 2, $g\in 
C_0(D)$, and $\bmu$-a.e. $\nu\in G$ we have $a^f\nu\equiv dA^f\nu/d\lambda 
\in X$ and 
\begin{eqnarray*} 
\left.\frac{d}{du}\right|_{u=0}\left((\nu+u\cdot A^f\nu)_t,g\vphantom{l^1} 
\right)=\left(A^f\nu_t,g\right)\, . 
\end{eqnarray*} 
\item[(v')] All requirements of condition (v) hold for all $t\in [-1,1]$, 
regardless of whether we have Case 1 or Case 2. 
\end{itemize} 
{\bf Remark} (6) Let us look at the mathematical situation in the following 
way. We consider probability measure valued trajectories $[-1,1]\ni t\mapsto 
\nu_t:=U^+(t,\nu)$, $\nu\equiv\nu_0\in G$, defined in Paragraph \ref{sec:1:1:3}. 
Recalling (\ref{1.4}) and looking ahead to Sections \ref{sec:3} and \ref{sec:4} 
of this paper, this gives rise to think of the solution $[-1,1]\ni t\mapsto 
\nu_t$ to a non-linear PDE. At time $t=0$ we observe or control the distribution 
$\bmu$ over $\nu\equiv\nu_0\in G$. We are now interested in the distribution of 
$\nu_t=U^+(t,\nu)$, $t\in [0,1]$, in Case 1 of Paragraph \ref{sec:1:1:5}. In 
Case 2 we are interested in the distribution of $\nu_{-t}=U^+(-t,\nu)$, $t\in 
[0,1]$. 
\medskip

For the next theorem introduce the following. In Case 1, let $\tau\equiv\tau 
(\nu)$, $\nu\in G$, denote the first exit time of the trajectory $\nu_t$, $t\in 
[0,1]$, from the set $G$, i. e. $\tau:=\inf\{t\in (0,1]:\nu_t\not\in G\}$. If 
$\{t\in (0,1]:\nu_t\not\in G\}=\emptyset$ we set $\tau:=\infty$. In Case 2, let 
$\tau\equiv\tau(\nu):=\inf\{t\in (0,1]: \nu_{-t}\not\in G\}$, $\nu\in G$, and 
set $\tau:=\infty$ if $\{t\in (0,1]:\nu_{-t}\not\in G\}=\emptyset$. 
\begin{theorem}\label{Theorem2.3} 
Suppose (a1)-(a4) and (j)-(v). (a) In Case 1, all measures $\bmu\circ\nu_{ 
-t}$, $t\in [0,1]$, are equivalent. In Case 2, all measures $\bmu\circ\nu_t$, 
$t\in [0,1]$, are equivalent. Note that $\bmu\circ\nu_0=\bmu$. The Radon-Nikodym 
derivatives have versions 
\begin{eqnarray*} 
r_{-t}:=\frac{d\bmu\circ\nu_{-t}}{d\bmu}\quad\mbox{\rm (Case 1),}\quad 
r_t:=\frac{d\bmu\circ\nu_t}{d\bmu}\quad\mbox{\rm (Case 2)}
\end{eqnarray*}
such that for all $t\in [0,1]$ we have $\bmu$-a.e. 
\begin{eqnarray}\label{2.3} 
r_{-t}(\nu)&&\hspace{-.5cm}=\exp\left\{-\int_{s=0}^t\delta (A^f)(\nu_{-s}) 
\, ds\right\}\quad\mbox{\rm (Case 1),}\nonumber \\ 
&& \\ 
r_t(\nu)&&\hspace{-.5cm}=\exp\left\{\int_{s=0}^t\delta (A^f)(\nu_s)\, ds 
\right\}\quad\mbox{\rm (Case 2)}\, .\nonumber 
\end{eqnarray} 
If we have (jv') as well as (v') and $\bmu_\gamma$ is the zero measure, for 
$\bmu$-a.e. $\nu\in E$ with $\tau(\nu)<\infty$ and all $u\in (0,\tau]$ it  
holds that 
\begin{eqnarray}\label{2.3*} 
\int_0^u \delta (A^f)(\nu_{\tau -s})\, ds&&\hspace{-.5cm}=\infty\quad\mbox 
{\rm (Case 1),}\nonumber \\ 
&& \\ 
\int_0^u \delta (A^f)(\nu_{-\tau +s})\, ds&&\hspace{-.5cm}=-\infty\quad\mbox 
{\rm (Case 2)}\, . \nonumber 
\end{eqnarray} 
(b) For $f\in L^\infty(E,\bmu)$ and $t\in [0,1]$, we have 
\begin{eqnarray*} 
\frac{d}{dt}\int f\, d\bmu\circ\nu_{-t}=-\int f\delta (A^f)(\nu_{-t}) 
\, d\bmu\circ\nu_{-t}\quad\mbox{\rm (Case 1)} 
\end{eqnarray*}
\begin{eqnarray*} 
\frac{d}{dt}\int f\, d\bmu\circ\nu_t=\int f\delta (A^f)(\nu_t)\, d\bmu 
\circ\nu_t\quad\mbox{\rm (Case 2)}\, . 
\end{eqnarray*}
(c) In Case 1, for any $\ve>0$ there exists a compact set $K_\ve\subseteq 
G^-$ with $\ \sup_{t\in [0,1]}\bmu\circ\nu_{-t}(E\setminus K_\ve)<\ve$ such 
that the derivative $\frac{d}{dt}\bmu\circ\nu_{-t}(B)$ exists for all $t 
\in [0,1]$ and $B\in {\cal B}(E)$ with $B\subseteq K_\ve$. In Case 1, this 
defines a finite signed measure 
\begin{eqnarray*}
\frac{d}{dt}\bmu\circ\nu_{-t}\equiv\left(\frac{d}{dt}\bmu\circ\nu_{-t} 
\right)_\ve\quad\mbox{\rm on}\quad\left(K_\ve,\{A\cap K_\ve:A\in {\cal 
B}(E)\}\right) 
\end{eqnarray*}
which is absolutely continuous with respect to $\bmu$ restricted to $K_\ve$. 
The analogue holds in Case 2, replacing $G^-$ by $G^+$ and $\bmu\circ\nu_{-t}$ 
by $\bmu\circ\nu_t$. \\ 
(d) Identifying for $f\in\t C^1_b(E)$ and $t\in [-1,1]$ the composition 
$f(\nu_t)$ with 
\begin{eqnarray*} 
\left\{
\begin{array}{ll} 
f(\nu_t) & {\rm if }\ \nu\in G \\ 
0 & {\rm otherwise}
\end{array}
\right.  
\end{eqnarray*} 
the following holds. In Case 1, the derivative $\frac{d}{dt} f(\nu_t)$ 
exists for $t\in [0,1]$ as a derivative in $L^1(E,\bmu)$ and we have 
$f(\nu_t)\in C_b^{q,1}(E)$ whenever $f\in\t C^1_b(E)$. In Case 2, $\frac 
{d}{dt}f(\nu_t)$ exists for $t\in [-1,0]$ in $L^1(E,\bmu)$ and we have 
$f(\nu_t)\in C_b^{q,1}(E)$ for $f\in\t C^1_b(E)$. 
\end{theorem}
\begin{theorem}\label{Theorem2.7} 
Under the hypotheses of Theorem \ref{Theorem2.3} we have the following. 
(a) For $f\in C_b(E)$ and $t\in [-1,1]$ identify the composition $f(\nu_t)$ 
with 
\begin{eqnarray*} 
\left\{
\begin{array}{ll} 
f(\nu_t) & {\rm if }\ \nu\in G  \\ 
0 & {\rm otherwise}
\end{array}
\right. \, .
\end{eqnarray*} 
For $f\in C_b(E)$ and $t\in [0,1]$, we have 
\begin{eqnarray*} 
\frac{d}{dt}\int f\, d\bmu(d\nu_{-t})=-\int f(\nu_t)\delta (A^f)(\nu)\, \bmu 
(d\nu)+\int f(\nu_t)\, d\bmu_\gamma\quad\mbox{\rm (Case 1)} 
\end{eqnarray*}
\begin{eqnarray*} 
\frac{d}{dt}\int f\, d\bmu(d\nu_t)=\int f(\nu_{-t})\delta (A^f)(\nu)\, \bmu 
(d\nu)+\int f(\nu_{-t})\, d\bmu_\gamma\quad\mbox{\rm (Case 2).} 
\end{eqnarray*}
(b) Denoting by $\l_{[0,1]}$ the Lebesgue measure on $([0,1],{\cal B} 
([0,1]))$ the following are equivalent. 
\begin{itemize} 
\item[(vj)] It holds that 
\begin{eqnarray*} 
\exp\left\{-\int_0^t\delta (A^f)(\nu_{-s})\, ds\right\}\in L^1(\gamma^- 
G\times [0,1]\, ,\, \bmu_\gamma\times\l_{[0,1]})\quad\mbox{\rm (Case 1)}  
\end{eqnarray*}
in the sense of $(\nu,t)\in \gamma^- G\times [0,1]$, or
\begin{eqnarray*} 
\exp\left\{\int_0^t\delta (A^f)(\nu_{s})\, ds\right\}\in L^1(\gamma^+ G 
\times [0,1]\, ,\, \bmu_\gamma\times\l_{[0,1]})\quad\mbox{\rm (Case 2)}  
\end{eqnarray*}
in the sense of $(\nu,t)\in \gamma^+G\times [0,1]$. 
\item[(vj')] For all $B\in {\cal B}(E)$ and $t\in [0,1]$ we have in Case 1  
\begin{eqnarray*}
\bmu(\nu_{-t}:\nu\in B)-\bmu\circ\nu_{-t}(B)=\int_{s=0}^t\int_B\exp\left\{- 
\int_{u=s}^t\delta (A^f)(\nu_{-u})\, du\right\}\, \bmu_\gamma(d\nu_{-s})\, 
ds\, .  
\end{eqnarray*} 
In Case 2 it holds for all $B\in {\cal B}(E)$ and $t\in [0,1]$ that 
\begin{eqnarray*}
\bmu(\nu_t:\nu\in B)-\bmu\circ\nu_t(B)=\int_{s=0}^t\int_B\exp\left\{\int_{u= 
s}^t\delta (A^f)(\nu_u)\, du\right\}\, \bmu_\gamma(d\nu_s)\, ds\, .  
\end{eqnarray*} 
\end{itemize}
(c) Suppose (vj). The measure $\bmu_\gamma$ is the zero measure if and only 
if 
\begin{eqnarray*} 
\lim_{t\to 0}\frac1t\left(\int f\, d\bmu(\nu_{-t})-\int f\, d\bmu\circ 
\nu_{-t}\right)=0\quad\mbox{\rm (Case 1)}\, , 
\end{eqnarray*} 
\begin{eqnarray*} 
\lim_{t\to 0}\frac1t\left(\int f\, d\bmu(\nu_t)-\int f\, d\bmu\circ 
\nu_t\right)=0\quad\mbox{\rm (Case 2)}\, , 
\end{eqnarray*} 
for all $f\in C_b(E)$. 
\end{theorem}

Choosing in Theorem \ref{Theorem2.7} (b) $B=E$ an immediate consequence is the 
following. 
\begin{corollary}\label{Corollary2.8} 
{\it (Loss of mass distribution in terms of $\delta(A^f)$ and $\bmu_\gamma$) } 
Under the hypotheses of Theorem \ref{Theorem2.7} we have 
\begin{eqnarray*} 
1-\bmu\circ\nu_{-t}(E)=\int_{s=0}^t\int\exp\left\{-\int_{u=0}^s\delta (A^f) 
(\nu_{-u})\, du\right\}\, \bmu_\gamma(d\nu)\, ds\, ,\quad t\in [0,1], 
\quad\mbox{\rm (Case 1)} 
\end{eqnarray*}
and 
\begin{eqnarray*} 
1-\bmu\circ\nu_t(E)=\int_{s=0}^t\int\exp\left\{\int_{u=0}^s\delta (A^f)(\nu_u) 
\, du\right\}\, \bmu_\gamma(d\nu)\, ds\, ,\quad t\in [0,1],\quad\mbox{\rm 
(Case 2).}  
\end{eqnarray*}
\end{corollary} 
{\bf Remark} (7) Condition (jj) refers to the form of differentiability of 
$\bmu$ we apply to prove Theorem \ref{Theorem2.3} and Theorem \ref{Theorem2.7}. 
Conversely, condition (jj) follows directly from conditions (j), (v), and 
the statement of Theorem \ref{Theorem2.7} (a) assuming that for $f\in 
C_b^{q,1}(E)$ the limit $\left.\frac{d}{dt}\right|_{t=0} f(\nu_t)$ exists in 
$L^1(E,\bmu)$. In this sense we understand the necessity of condition (jj). 
For a formal motivation of condition (jj) we refer to Proposition 
\ref{Proposition2.25}. 
\bigskip

The following proposition shows how to replace the technical condition (v) 
and (\ref{2.2}) in (jv) by Fr\'echet differentiability. In the application 
of Section \ref{sec:3}, we will choose $v=1$ and verify the Fr\'echet 
differentiability of Proposition \ref{Proposition2.4} in order to get 
Theorem \ref{Theorem2.3}. In contrast, for the model considered in Section 
\ref{sec:4}, the choice of $v=2$ and the direct verification of conditions 
(jv) and (v) give Theorem \ref{Theorem2.3} in an efficient way. For the 
remainder of this subsection assume that the hypotheses (a1)-(a4) of 
Subsection \ref{sec:1:1} hold as far as they are needed. Let $\nabla$ 
denote the Fr\'echet derivative. Furthermore, for $\ve>0$ introduce 
$G^\ve:=\{g\in L^v(D):\|g-h\|_{L^v(D)}<\ve$ for some $h\cdot\lambda\in G\}$. 
\begin{proposition}\label{Proposition2.4} 
Let $\ve>0$ and let $\t {\cal V}$ be an open set in the topology of $L^v 
(D)$ with $G^\ve\subseteq\t {\cal V}\subseteq {\cal V}$. Suppose (j) and 
assume $G\subseteq E\cap\{V_t^+(h)\cdot\lambda:h\in\t {\cal V}\}$ for all 
$t\in [0,1]$. Suppose also that for all $\t h\in\t {\cal V}$ there is a 
pre-neighborhood $U_{\t h}(0)$ of $0\in X$ in the topology of $X$ such that 
$\{\t h+i\circ h:h\in U_{\t h}(0)\}\subseteq D(V_1^+)\cap\{V_1^+(h):h\in 
D(V_1^+)\}$. 
\medskip

Let the map $\Phi_t:\t {\cal V}\ni\t h\mapsto\Phi_t(\t h):=V^+_t\t h 
\in L^v(D)$ be Fr\'echet differentiable at every $\t h\in \t {\cal V}$, 
for all $t\in [0,1]$. In addition to the in Paragraph \ref{sec:1:1:2} 
hypothesized injectivity of the map $\Phi_t$, assume that for all $t 
\in [0,1]$ the inverse $\Phi_{-t}\equiv (\Phi_t)^{-1}$ is a continuous 
map from $\{V_t^+(h):h\in \t {\cal V}\}$ onto $\t {\cal V}$, both sets 
endowed with the topology of $L^v(D)$. 

Furthermore, suppose that for all $\t h\in \t {\cal V}$ and $t\in [0,1]$ 
the Fr\'echet derivative 
\begin{eqnarray*} 
\nabla\Phi_t(\t h):L^v(D)\mapsto L^v(D)\ \mbox{\it is injective} 
\end{eqnarray*} 
and its inverse 
\begin{eqnarray*}
\left(\nabla\Phi_t(\t h)\right)^{-1}:\left\{\nabla\Phi_t(\t h)(h):h\in 
L^v(D)\right\}\mapsto L^v(D)\ \mbox{\it is bounded.}  
\end{eqnarray*} 

Then we have (v') and the derivative (\ref{2.2}) in condition (jv) exists 
for $\bmu$-a.e. $\nu\in G$ and all $\, t\in [-1,1]$, $g\in C_0(D)$, as 
well as $k\in X$, for some $r(\nu,t,g)\in X^\ast$. 
\end{proposition} 

In contrast to other papers focusing on the formula (\ref{2.3}), cf. 
Subsection \ref{sec:1:2}, we have no direct assumption on the 
smoothness of the map $\nu\mapsto\nu_t=U^+(t,\nu)$, $t\in [-1,1]$. 
Furthermore, there is no assumption on some class of functions defined 
$\bmu$-a.e. on $E$, to be invariant under the map $\nu\mapsto\nu_t$. 
Conditions (jv) and (v) will solely be used in Step 1 of the proof of 
Theorem \ref{Theorem2.3}. There we demonstrate that cylindrical 
functions $f\in\t C_b^1(E)$ composed with $\nu\mapsto\nu_t$, $t\in 
[-1,1]$, belong to $C_b^{q,1}(E)$ and have gradients that are 
continuous on $t\in [-1,1]$. The following proposition demonstrates 
that precisely this property implies necessity of the hypotheses (jv) 
and (v). 
\begin{proposition}\label{Proposition2.5} 
Let $t\in [-1,1]$ and $g\in C_0(D)$. For any $g\in C_0(D)$ identify 
the composition $(\nu_t,g)$ with 
\begin{eqnarray*} 
\Psi_{t,g}:=\left\{
\begin{array}{ll} 
(\nu_t,g) & {\rm if }\ \nu\in G  \\ 
0 & {\rm otherwise}
\end{array}
\right. \, .
\end{eqnarray*}  
Assume that for any $g\in C_0(D)$ 
\begin{eqnarray*} 
\Psi_{t,g}\in C_b^{q,1}(E)\quad\mbox{\rm such that }\ [-1,1]\ni t\mapsto 
D\Psi_{t,g}\in L^q(E,\bmu,X^\ast)\ \mbox{\rm is continuous}. 
\end{eqnarray*} 
(a) Then we have (jv'). \\ 
(b) Suppose in addition (a4), (j), (jj), and the statements of Theorems 
\ref{Theorem2.3} (c) and \ref{Theorem2.7} (a). Then we have (v). 
\end{proposition} 
Let us mention that for the proof of part (b) we apply Theorems 
\ref{Theorem2.3} (c) and \ref{Theorem2.7} (a) in a way that we have to 
distinguish between Case 1 and Case 2. Therefore, part (b) of Proposition 
\ref{Proposition2.5} refers just to condition (v), not to (v').

\subsection{Proofs of Subsection \ref{sec:2:12}}\label{sec:2:2}

{\bf Proof of Theorem \ref{Theorem2.3}} (a) Let us prove the formula 
$r_{-t}(\nu)=\exp\{-\int_{s=0}^t\delta (A^f)(\nu_{-s})\, ds\}$, $t\in 
[0,1]$, under the hypotheses of Case 1. En passant below (\ref{2.85}), 
we shall demonstrate that if we have (jv') as well as (v') and 
$\bmu_\gamma$ is the zero measure then for $\bmu$-a.e. $\nu\in E$ with 
$\tau(\nu)<\infty$ it holds that $\int_0^u\delta (A^f)(\nu_{\tau -s}) 
\, ds=\infty$ for all $u\in (0,\tau]$.  
\medskip 

Let $t\in [-1,1]$. Moreover, let $f\in\t C_b^1(E)$ and $\vp$ be 
related to $f$ as in (\ref{2.1}) and define 
\begin{eqnarray*} 
f^F(\nu):=\vp\left((h_1,\nu),\ldots ,(h_r,\nu)\right)\, ,\quad\nu\in F. 
\end{eqnarray*} 
Let us introduce 
\begin{eqnarray*} 
f^F\circ\nu_t:=\left\{
\begin{array}{ll} 
f^F(\nu_t) & {\rm if }\ \nu_t=h_t\cdot\lambda\ \ {\rm and }\ h_t\in L^v(D) \\ 
0 & {\rm otherwise}
\end{array}
\right. \, .
\end{eqnarray*} 
Furthermore, let $f\circ\nu_t$ be the restriction of $f^F\circ\nu_t$ to $E$. 
Throughout this proof, we will always consider these functions just $\bmu$-a.e. 
Recalling the definition of $G\in {\cal B}(E)$ in Paragraph \ref{sec:1:1:5}, the 
property $\bmu(G)=1$ tells us that there is no conflict in the interpretation of 
the symbols $f\circ\nu_t$ or $f(\nu_t)$ as the composition of $f$ with $\nu_t$. 
However, as we shall see in the subsequent Step 1, the above definition of $f\circ 
\nu_t\equiv f(\nu_t)$, $\nu\in G$, allows us to meaningfully introduce directional 
derivative and gradient. Recall also condition (a1) of Paragraph \ref{sec:1:1:5} 
together with the definition of $\bmu\circ\nu_{-t}$ in Paragraph \ref{sec:1:1:4} 
and note that $\bmu\circ\nu_{-t}(d\nu)\not\equiv\bmu(d\nu_{-t})$ in general. 
\medskip 

\nid
{\it Step 1 } Suppose Case 1. In this step we verify that for $f\in\t C_b^1 
(E)$ and $t\in [0,1]$ we have $f\circ\nu_t\in C^{q,1}_b(E)$ and we determine 
the gradient of $f\circ\nu_t$. We show also that $\frac{d}{dt} f(\nu_t)$ 
exists as a derivative in $L^1(E,\bmu)$. Recall that by (j), $a^f\nu\equiv 
dA^f\nu/d\lambda\in X$ for $\bmu$-a.e. $\nu\in G$. 
\medskip

According to condition (v), for $\bmu$-a.e. $\nu\in G$ we find  
\begin{eqnarray}\label{2.19**} 
\frac{d}{dt}f(\nu_t)&&\hspace{-.5cm}=\sum_{i=1}^n\frac{\partial\vp} 
{\partial x_i}\left((h_1,\nu_t),\ldots ,(h_r,\nu_t)\right)\cdot\frac{d} 
{dt}(\nu_t,h_i) \nonumber \\ 
&&\hspace{-.5cm}=\sum_{i=1}^n\frac{\partial\vp}{\partial x_i}\left((h_1, 
\nu_t),\ldots ,(h_r,\nu_t)\right)\cdot (A^f\nu_t,h_i) \nonumber \\ 
&&\hspace{-.5cm}=\sum_{i=1}^n\frac{\partial\vp}{\partial x_i}\left((h_1, 
\nu_t),\ldots ,(h_r,\nu_t)\right)\cdot\left.\frac{d}{du}\right|_{u=0} 
\left((\nu+u\cdot A^f\nu)_t,h_i\right) \nonumber \\ 
&&\hspace{-.5cm}=\left.\frac{d}{du}\right|_{u=0}f^F\left((\nu+u\cdot A^f 
\nu)_t\right)\, . 
\end{eqnarray} 
With 
\begin{eqnarray}\label{2.20**} 
&&\hspace{-.5cm}\left.\frac{d}{du}\right|_{u=0}f^F((\nu+u\cdot i\circ 
h\cdot\lambda)_t)\nonumber \\ 
&&\hspace{.5cm}=\sum_{i=1}^r\frac{\partial\vp}{\partial x_i}\left((h_1, 
\nu_t),\ldots ,(h_r,\nu_t)\right)\cdot{}_{X^\ast}\langle r(\nu ,t,h_i), 
h\rangle_X
\end{eqnarray} 
for all $h\in X$ and $\bmu$-a.e. $\nu\in G$ by condition (jv), we have  
\begin{eqnarray}\label{2.21**}  
D(f^F\circ\nu_t)=\sum_{i=1}^r\frac{\partial\vp}{\partial x_i}\left((h_1, 
\nu_t),\ldots ,(h_r,\nu_t)\right)\cdot {}_{X^\ast}\langle r(\nu ,t,h_i), 
\cdot\rangle_X\, . 
\end{eqnarray} 
Condition (jv) says also that $D(f^F\circ\nu_t)\in L^q(E,\bmu;X^\ast)$.  
Together with (\ref{2.20**}) and (\ref{2.21**}) this yields $f^F\circ 
\nu_t\equiv f^F\circ\nu_t(\cdot)\in C^{q,1}_b(F,E)$ and therefore 
\begin{eqnarray*} 
f\circ\nu_t\equiv f\circ\nu_t(\cdot)\in C^{q,1}_b(E)\, ,\quad f\in\t 
C_b^1(E),\ t\in [0,1], 
\end{eqnarray*} 
with $q>1$ given in condition (jj). Furthermore, (\ref{2.21**}) implies 
also  
\begin{eqnarray}\label{2.5}  
D(f\circ\nu_t)=\sum_{i=1}^r\frac{\partial\vp}{\partial x_i}\left((h_1, 
\nu_t),\ldots ,(h_r,\nu_t)\right)\cdot {}_{(X\ominus\1)^\ast}\langle 
r(\nu ,t,h_i),\cdot\rangle_{X\ominus\1}\, .
\end{eqnarray} 
With (\ref{2.19**}) and (\ref{2.20**}) this gives 
\begin{eqnarray}\label{2.4} 
\frac{d}{dt}f(\nu_t)=D(f\circ\nu_t)(a^f\nu)\, ,\quad t\in [0,1],  
\end{eqnarray} 
for $\bmu$-a.e. $\nu\in G$. From (\ref{2.5}) and condition (jv) we obtain 
the continuity of $[0,1]\ni t\mapsto D(f\circ\nu_t)$ in $L^q(E,\bmu; 
(X\ominus\1)^\ast)$. Relation (\ref{2.4}) and condition (j) imply now that 
the derivative $\frac{d}{dt}f(\nu_t)$ exists as a derivative in $L^1(E,\bmu)$ 
in the sense of part (c).

We emphasize that in case of (jv') and (v') everything in Step 1 remains true 
for all $t\in [-1,1]$. 
\medskip 

\nid  
{\it Step 2 } We derive an equation for the measure valued function $[0,1] 
\ni t\to\bmu\circ\nu_{-t}$ with initial value $\bmu$. For this, we follow 
an idea outlined in the proof of Theorem 9.2.4 in \cite{Bo97}. We also 
demonstrate in this step that the measures $\bmu\circ\nu_{-t}$ and $\bmu$ 
are equivalent. 
\medskip 

Together with conditions (j) and (jj) relation (\ref{2.4}) gives 
\begin{eqnarray*} 
&&\hspace{-.5cm}\frac{d}{dt}\int f(\nu_t)\, \bmu(d\nu )=\int\frac{d}{dt} 
f(\nu_t)\, \bmu(d\nu ) \\ 
&&\hspace{.5cm}=\int D(f\circ\nu_t)(a^f\nu)\, \bmu(d\nu) \\ 
&&\hspace{.5cm}=-\int f(\nu_t)\delta (A^f)(\nu)\, \bmu(d\nu)+\int f(\nu_t) 
\, d\bmu_\gamma\, , \quad t\in [0,1],\ f\in \t C^1_b(E); 
\end{eqnarray*}
where we stress that the interchange of $\int$ and $\frac{d}{dt}$ in 
the first line is justified by the fact that $\frac{d}{dt} f(\nu_t)$ 
exists by Step 1 as a derivative in $L^1(E,\bmu)$. 
\bigskip 

For $t\in [0,1]$, $f\in \t C^1_b(E)$, we have 
\begin{eqnarray*} 
\int f(\nu_t)\, \bmu(d\nu )-\int f\, d\bmu=-\int_0^t\int f(\nu_s) 
\delta (A^f)\, d\bmu\, ds+\int_0^t\int f(\nu_s)\, d\bmu_\gamma\, ds 
\end{eqnarray*}
and therefore 
\begin{eqnarray*}
&&\hspace{-.5cm}\int f(\nu)\, \bmu(d\nu_{-t})-\int f\, d\bmu \\ 
&&\hspace{.5cm}=-\int_0^t\int f(\nu)\delta (A^f)(\nu_{-s})\, \bmu 
(d\nu_{-s})\, ds+\int_0^t\int f(\nu)\, \bmu_\gamma(d\nu_{-s})\, ds 
\, .  
\end{eqnarray*}

Let $\ve>0$. According to Lemma \ref{Lemma1.1} (a) and condition (jj) 
there is a compact subset ${\cal K}_\ve\equiv{\cal S}^v_m$ of $(E,\pi)$ 
such that 
\begin{eqnarray*}
(2+\|\delta(A^f)\|_{L^1(E,\sbmu)})\cdot\sup_{t\in [-1,1]}\bmu\circ 
\nu_{-t}(E\setminus{\cal K}_\ve)+\sup_{t\in [-1,1]}\bmu_\gamma\circ\nu_{-t} 
(\gamma^- G\setminus{\cal K}_\ve)\le\ve\, . 
\end{eqnarray*}
As above, we use the symbol $\chi$ to denote the indicator function. 
Recalling Remark (2) of the present section we obtain for all $f\in 
C_b(E)$ and all $t\in [0,1]$ 
\begin{eqnarray}\label{2.19*}
&&\hspace{-.5cm}\left|\int f\chi_{{\cal K}_\ve}\, \bmu(d\nu_{-t})-\int 
f\chi_{{\cal K}_\ve}\, d\bmu\right.\nonumber \\ 
&&\hspace{.5cm}\left.+\int_0^t\int f\chi_{{\cal K}_\ve}\delta (A^f) 
(\nu_{-s})\, d\bmu(\nu_{-s})\, ds-\int_0^t\int f\chi_{{\cal K}_\ve}\, 
\bmu_\gamma(d\nu_{-s})\, ds\right|\le\ve\cdot\|f\|\, .\qquad
\end{eqnarray} 
For any finite intersection $\t B$ of the open balls $B_\nu(r)$ with 
center $\nu\in E$ and radius $r>0$, the indicator function $\chi_{\t 
B}$ can be boundedly point-wise approximated by functions belonging 
to $C_b(E)$. For example use the approximation $f_n(\nu):=(n\cdot 
\inf_{\nu'\in\t B^c}\pi(\nu,\nu'))\wedge 1$ as $n\to\infty$. To show 
continuity of $f_n$ use $\pi(\nu_1,\nu_2)\ge|\inf_{\nu'\in\t B^c}\pi 
(\nu_1,\nu')-\inf_{\nu'\in\t B^c}\pi(\nu_2,\nu')|$. As a consequence 
we obtain (\ref{2.19*}) for all such $f=\chi_{\t B}$ and by using the 
monotone class theorem, for all $f=\chi_B$, $B\in {\cal B}(E)$. 

If we have (jv') as well as (v') then everything of Step 2 up to this  
point holds for all $t\in [-1,1]$. For this recall also the last 
sentence of Step 1. 
\medskip

In addition, for any $\nu\equiv\nu_0\in G^-$ and $s\in [0,1]$ it holds 
that $\nu_{-s}\in G^-$, cf. Paragraph \ref{sec:1:1:5}. Therefore, from 
(\ref{2.19*}) we obtain 
\begin{eqnarray}\label{2.75} 
\bmu(\nu_{-t}:\nu\in B)-\bmu(B)=-\int_{s=0}^t\int_B\delta (A^f)(\nu_{ 
-s})\, \bmu(d\nu_{-s})\,ds\, ,\quad t\in [0,1], 
\end{eqnarray}
for all $B\in {\cal B}(E)$ with $B\subseteq G^-$. If $\bmu_\gamma$ is 
the zero measure and we have (jv') as well as (v') then (\ref{2.75}) 
holds for all $t\in [-1,1]$. Set  
\begin{eqnarray}\label{2.8} 
a(t,B):=\left\{
\begin{array}{cl}
\frac{\displaystyle-\int_B\delta (A^f)(\nu_{-t})\, \bmu(d\nu_{-t})} 
{\displaystyle\bmu(\nu_{-t}:\nu\in B)\vphantom{\displaystyle\int}} 
&\ \mbox{\rm if }\ \bmu(\nu_{-t}:\nu\in B)>0 \\ 
0 &\ \mbox{\rm otherwise }\vphantom{\displaystyle\int}
\end{array}\right.\, , 
\end{eqnarray}
for all $B\in {\cal B}(E)$ with $B\subseteq G^-$ and $t\in [0,1]$. 
From the last two relations we obtain 
\begin{eqnarray}\label{2.85} 
\bmu(\nu_{-t}:\nu\in B)-\bmu(B)=\int_{s=0}^ta(s,B)\cdot\bmu(\nu_{-s}: 
\nu\in B)\, ds     
\end{eqnarray}
for all $B\in {\cal B}(E)$ with $B\subseteq G^-$ and $t\in [0,1]$. 

If $\bmu_\gamma$ is the zero measure and we have (jv') as well as 
(v') then (\ref{2.8}) can be defined for all $t\in [-1,1]$ and 
(\ref{2.85}) holds for all $t\in [-1,1]$. For $B\in {\cal B}(E)$ with 
$B\subseteq G^-$ and $\bmu(B)>0$ such that $\tau (\nu)\le 1$ for 
$\bmu$-a.e. $\nu\in B$, we observe $\lim_{t\to -1}\bmu(\nu_{-t}:\nu\in B) 
=0$ and therefore $\inf\{\int_{s=0}^ta(s,B)\, ds:t\in [-1,0]\}=-\infty$. 
Assuming that (\ref{2.3*}) was not true we could deduce the existence of 
such a set $B$ in a way that $\int_0^u\delta (A^f)(\nu_{\tau -s})\, ds$ 
is bounded from above uniformly for all $u\in (0,\tau]$, for $\bmu$-a.e. 
$\nu\in B$. By (\ref{2.8}) this contradicts $\inf\{\int_{s=0}^ta(s,B)\, 
ds:t\in [-1,0]\}=-\infty$. Thus (\ref{2.3*}) holds. 
\medskip

Since for $f\in L^\infty(E)$ with $f=0$ on $E\setminus G^-$ it holds 
that $\int f(\nu)\, \bmu(d\nu_{-t})=\int f\, d\bmu\circ\nu_{-t}$, cf. 
Paragraph \ref{sec:1:1:4}, relation (\ref{2.85}) yields 
\begin{eqnarray*} 
\bmu\circ\nu_{-t}(B)-\bmu(B)=\int_{s=0}^ta(s,B)\cdot\bmu\circ 
\nu_{-s}(B)\, ds     
\end{eqnarray*}
for all $B\in {\cal B}(E)$ with $B\subseteq G^-$ and $t\in [0,1]$.    
Note here that $\bmu\circ\nu_{-t}(G\setminus G^-)=0$ for all $t\in 
[0,1]$ according to Paragraph \ref{sec:1:1:4}. 
\medskip

Now recall the definition of the space $L^\infty_{\rm loc}$ in 
the introductory part of Section \ref{sec:2}, condition (jjj) 
in Subsection \ref{sec:2:1}, Lemma \ref{Lemma1.1} (b) and its 
proof, Remark (3) of Section \ref{sec:1}, and condition (a2). 

For any $\ve>0$ there exists a compact set $K_\ve\equiv\hat{\cal 
S}^-\subseteq G^-$ such that $\sup_{t\in [0,1]}\bmu\circ\nu_{-t} 
(E\setminus K_\ve)<\ve$, see Lemma \ref{Lemma1.1} (b). 
According to Lemma \ref{Lemma1.1} (b), $\{\nu_{-s}:\nu\in K_\ve,\ s\in 
[0,1]\}$ is contained in some compact subset $K'_\ve\equiv {\cal S}^-$ 
of $G^-$ in Case 1. We note that $\sup_{t\in [0,1]}\bmu\circ\nu_{-t}(E 
\setminus K'_\ve)<\ve$. Thus by condition (jjj) and (\ref{2.8}) 
\begin{eqnarray}\label{2.10} 
|a(s,B)|\le\|\delta (A^f)\cdot\chi_{K'_\ve}\|_{L^\infty(E,\sbmu)}< 
\infty\, ,\quad s\in [0,1],\ B\in {\cal B}(E),\ B\subseteq K_\ve,
\end{eqnarray}
in Case 1. Therefore,  
\begin{eqnarray}\label{2.9} 
\bmu\circ\nu_{-t}(B)=\exp\left\{\int_{s=0}^t a(s,B)\, ds\right\}\cdot 
\bmu(B)\, , \quad t\in [0,1],\ B\in {\cal B}(E),\ B\subseteq K_\ve,  
\end{eqnarray}
in Case 1. Recalling Remark (3) of Section \ref{sec:1}, (\ref{2.10}) 
and (\ref{2.9}) say that the measures $\bmu\circ\nu_{-t}$ and $\bmu$ 
are equivalent under the hypotheses of Case 1. 
\medskip 

\nid
{\it Step 3 } Until the end of this proof suppose the hypotheses 
of Case 1. The objective of this step is to show well-definiteness 
of $\frac{d}{dt}\bmu\circ\nu_{-t}$, $t\in [0,1]$, as a finite signed 
measure on $\left(K_\ve,\{B\cap K_\ve:B\in {\cal B}(E)\}\right)$, i. e. 
we show (c). Here $\ve>0$ is arbitrary and $K_\ve\subseteq E$ is the 
compact set defined in Step 2. 

Let $\ve>0$, $B\in {\cal B}(E)$ such that $B\subseteq K_\ve$. We 
observe that for $t\in [0,1]$  
\begin{eqnarray*}
a(t,B)=-(\bmu\circ\nu_{-t}(B))^{-1}\cdot\int_B\delta (A^f)(\nu_{-t}) 
\, d\bmu\circ\nu_{-t} 
\end{eqnarray*} 
if $\bmu\circ\nu_{-t}(B)=\bmu(\nu_{-t}:\nu\in B)>0$. If $\bmu(B)=0$ 
then by the result of Step 2 we have $a(t,B)=0$ for all $t\in [0,1]$, 
cf. (\ref{2.8}). Assuming $\bmu(B)>0$, by (\ref{2.10}) and (\ref{2.9}), 
\begin{eqnarray}\label{2.11}
\bmu\circ\nu_{-t}(B)\quad \mbox{\rm is positive and continuous on}\quad 
t\in [0,1]. 
\end{eqnarray} 

What we now aim to verify is that $a(t,B)$ is continuous on $t\in [0,1]$. 
We approximate $\chi_{\t B}$ by $f\in C_b(E)$, boundedly point-wise, as 
described below (\ref{2.19*}) where $\t B$ is a finite intersection of 
open balls $B_\nu(r)$ with center $\nu\in E$ and radius $r>0$. Recalling 
$\delta (A^f)\in L^1(E,\bmu)\cap L_{\rm loc}^\infty(E,\bmu)$ and that 
$K_\ve\equiv \hat{\cal S}^-\subseteq G^-$ was chosen in Step 2 such that 
$\sup_{t\in [0,1]}\bmu\circ\nu_{-t}(E\setminus K_\ve)<\ve$, we obtain for 
$s,t\in [0,1]$ 
\begin{eqnarray*} 
&&\hspace{-.5cm}\left|\int\chi_{\t B}(\nu)\cdot\delta (A^f)(\nu_{-s})\, 
\bmu\circ\nu_{-s}(d\nu)-\int\chi_{\t B}(\nu)\cdot\delta (A^f)(\nu_{-t}) 
\, \bmu\circ\nu_{-t}(d\nu)\right|\nonumber \\ 
&&\hspace{.5cm}\le\left|\int(\chi_{\t B}\circ\nu_s-f\circ\nu_s)\cdot 
\delta(A^f)\, d\bmu-\int(\chi_{\t B}\circ\nu_t-f\circ\nu_t)\cdot\delta 
(A^f)\, d\bmu\right|\nonumber \\ 
&&\hspace{1.0cm}+\left|\int(f\circ\nu_s-f\circ\nu_t)\cdot\delta(A^f)\, 
d\bmu\right|\nonumber \\ 
&&\hspace{.5cm}\le\|\delta (A^f)\cdot\chi_{K_\ve}\|_{L^\infty(E,\sbmu)} 
\left(\int\left|\chi_{{\t B}}\circ\nu_s-f\circ\nu_s\right|\, d\bmu+\int 
\left|f\circ\nu_t-\chi_{{\t B}}\circ\nu_t\right|\, d\bmu\right)+\kappa 
\nonumber \\ 
&&\hspace{1.0cm}+\left|\int(f\circ\nu_s-f\circ\nu_t)\cdot\delta(A^f)\, 
d\bmu\right| 
\end{eqnarray*} 
for some $\kappa\equiv\kappa(\ve)>0$ with $\kappa\to 0$ as $\ve\to 0$. 
In other words, 
\begin{eqnarray*}
\int\chi_{\t B}(\nu)\cdot\delta (A^f)(\nu_{-t})\, \bmu\circ\nu_{-t}(d\nu) 
\quad \mbox{\rm is continuous on}\quad t\in [0,1]. 
\end{eqnarray*} 
We mention that the collection of all finite intersections $\t B$ of open 
balls generate ${\cal B}(E)$. 
In order to apply the monotone class theorem introduce the set ${\cal H}$ 
of all bounded measurable functions $g$ for which $\int g(\nu)\delta (A^f) 
(\nu_{-t})\, \bmu\circ\nu_{-t}(d\nu)$ is continuous on $t\in [0,1]$. It 
remains to show that for a sequence of non-negative functions $g_n\in {\cal 
H}$ that increase to a bounded function $g$ we have $g\in {\cal H}$. However 
this follows from a calculation similar to the previous chain of inequalities 
replacing there $\chi_{\t B}$ by $g$ and $f$ by $g_n$. We have verified that 
for all $g\in L^\infty(E;\bmu)$ 
\begin{eqnarray}\label{2.12} 
\int g(\nu)\delta (A^f)(\nu_{-t})\, \bmu\circ\nu_{-t}(d\nu)\ \mbox{\rm is 
continuous on}\ t\in [0,1]. 
\end{eqnarray}
This and (\ref{2.11}) yield continuity of $a(t,B)$ on $t\in [0,1]$ for 
$B\in {\cal B}(E)$ such that $B\subseteq K_\ve$. 
\medskip

From (\ref{2.9}) it follows now that, for any $\ve>0$, $\frac{d}{dt} 
\bmu\circ\nu_{-t}(B)\equiv\left(\frac{d}{dt}\bmu\circ\nu_{-t}\right)_\ve(B)$ 
exists for $B\in {\cal B}(E)$ with $B\subseteq K_\ve$ as the finite limit of 
\begin{eqnarray*} 
\lim_{h\to 0}\frac{1}{h}\left(\bmu\circ\nu_{-t-h}(B)-\bmu \circ\nu_{-t} 
(B)\right)\, ,\quad t\in [0,1]. 
\end{eqnarray*}
According to the Nikodym Convergence Theorem, see e. g. \cite{Ry02}, 
Proposition C.4, $\frac{d}{dt}\bmu\circ\nu_{-t}\equiv\left(\frac{d}{dt} 
\bmu\circ\nu_{-t}\right)_\ve$ is therefore for each $t\in [0,1]$ and 
$\ve>0$ a finite signed measure on $\left(K_\ve,\{B\cap K_\ve:B\in {\cal B} 
(E)\}\right)$. From (\ref{2.10}) and (\ref{2.9}) we derive the estimate 
\begin{eqnarray}\label{2.13}
&&\hspace{-.5cm}\frac{1}{h}\left|\bmu\circ\nu_{-t-h}(B)-\bmu\circ\nu_{-t} 
(B)\right|\nonumber \\ 
&&\hspace{.5cm}=\frac1h\left|\exp\left\{\int_{s=0}^{t+h}a(s,B)\, ds\right 
\}-\exp\left\{\int_{s=0}^t a(s,B)\, ds\right\}\right|\cdot\bmu(B) 
\nonumber \\ 
&&\hspace{.5cm}\le\left\|\delta (A^f)\cdot\chi_{K'_\ve}\right\|_{L^\infty 
(E,\sbmu)}\cdot\exp\left\{\left\|\delta(A^f)\cdot\chi_{K'_\ve}\right\|_{ 
L^\infty(E,\sbmu)}\right\}\cdot\bmu(B)  \vphantom{\int_{s=0}^{t+h}} 
\end{eqnarray}
where $\ve>0$ is arbitrary and $B\in {\cal B}(E)$ with $B\subseteq K_\ve$. 
Among other things this estimate says that the measures $\frac{d}{dt}\bmu 
\circ\nu_{-t}\equiv\left(\frac{d}{dt}\bmu\circ\nu_{-t}\right)_\ve$, $t\in 
[0,1]$, are absolutely continuous with respect to $\bmu$ restricted to 
$K_\ve$. 
\medskip 

\nid 
{\it Step 4 } In this step we construct the versions $r_{-t}$ of the 
Radon-Nikodym derivatives of $\bmu\circ\nu_{-t}$ with respect to $\bmu$ 
such that 
\begin{eqnarray*} 
[0,1]\ni t\to r_{-t}=\frac{d\bmu\circ\nu_{-t}}{d\bmu} 
\end{eqnarray*}
is $\bmu$-a.e. absolutely continuous on $([0,1],{\cal B}([0,1]))$. 
Let $\ve$ and $K_\ve$ as in Steps 2 and 3. 
\medskip 

Let $\bmu\circ\nu_{-\cdot}\equiv (\bmu\circ\nu_{-\cdot})_\ve$ be the 
finite signed measure on the product space $([0,1]\times K_\ve,{\cal B} 
([0,1]\times K_\ve))$ given by 
\begin{eqnarray*} 
\bmu\circ\nu_{-\cdot}((s,t]\times B):=\int_s^t\frac{d}{du}\bmu\circ 
\nu_{-u}(B)\, du\equiv\bmu\circ\nu_{-t}(B)-\bmu\circ\nu_{-s}(B)\, , 
\end{eqnarray*}
where $0\le s<t\le 1$ and $B\in {\cal B}(E)$ with $B\subseteq K_\ve$. 
For the finiteness of $(\bmu\circ\nu_{-\cdot})_\ve$, see (\ref{2.13}). 
From (\ref{2.9}) we deduce also the absolute continuity of $\bmu\circ 
\nu_{-\cdot}$ with respect to the product measure $dt\times\bmu$, where 
$dt$ stands for the Lebesgue measure on $([0,1],{\cal B}([0,1]))$. 

For the next conclusion we keep Theorem 58 of Chapter V in \cite{DM82} 
in mind. It implies that for any $\ve>0$ there is a version $\rho_{ 
\ve,-t}(\nu)$, $(t,\nu)\in [0,1]\times K_\ve$, of the Radon-Nikodym 
derivative $d\bmu\circ\nu_{-\cdot}/d(dt\times\bmu)$ restricted to 
$[0,1]\times K_\ve$ such that for all $t\in [0,1]$ the function 
\begin{eqnarray}\label{2.14} 
\mbox{\rm $\rho_{\ve,-t}$ is a version of the Radon-Nikodym derivative }\ 
\frac{d\left(\frac{d}{dt}\bmu\circ\nu_{-t}\right)_\ve}{d\bmu}\, . 
\end{eqnarray}
Looking at (\ref{2.13}) we may even assume that for $\bmu$-a.e. $\nu\in 
K_\ve$
\begin{eqnarray}\label{2.15} 
\rho_{\ve,-t}(\nu)\le\left\|\delta(A^f)\cdot\chi_{K'_\ve}\right\|_{L^\infty 
(E,\sbmu)}\cdot\exp\left\{\left\|\delta (A^f)\cdot\chi_{K'_\ve}\right 
\|_{L^\infty(E,\sbmu)}\right\}\, ,\quad t\in [0,1]. 
\end{eqnarray}
With this version we can define 
\begin{eqnarray}\label{2.16}
r_{\ve,-t}:=\int_0^t\rho_{\ve,-s}\, ds+1\, ,\quad t\in [0,1].   
\end{eqnarray}
We obtain from (\ref{2.14})-(\ref{2.16})
\begin{eqnarray}\label{2.17}
\int_B r_{\ve,-t}\, d\bmu&&\hspace{-.5cm}=\int_B\int_0^t\rho_{\ve,-s}\, ds 
\, d\bmu +\bmu(B)\nonumber \\ 
&&\hspace{-.5cm}=\int_0^t\int_B\rho_{\ve,-s}\, d\bmu\, ds+\bmu(B)=\bmu\circ 
\nu_{-t}(B)\, ,  
\end{eqnarray}
$t\in [0,1]$, $B\in {\cal B}(E)$ with $B\subseteq K_\ve$ for some $\ve>0$. 
This means that, for all $0\le t\le 1$, there is a function $r_{-t}\in L^1 
(E,\bmu)$ which is a version of the Radon-Nikodym derivative of $\bmu\circ 
\nu_{-t}$ with respect to $\bmu$ such that $r_{-t}=r_{\ve,-t}$ $\bmu$-a.e. 
on $K_\ve$, $\ve>0$, $t\in [0,1]$. 
\medskip 

\nid 
{\it Step 5 } Let $\ve>0$. We recall that $\delta (A^f)\in L^1(E,\bmu)$ 
according to (jjj) and that by (\ref{2.15}), (\ref{2.16}) we  have 
$r_{-t}\cdot\chi_{K_\ve}\in L^\infty(E,\bmu)$, $t\in [0,1]$. Using that 
for $B\in {\cal B}(E)$ with $B\subseteq K_\ve\subseteq G^-$ and $t\in 
[0,1]$ it holds that $\bmu(\nu_{-t}:\nu\in B)=\bmu\circ\nu_{-t}(B)$ cf. 
Paragraph \ref{sec:1:1:4}, we get from (\ref{2.75}) and (\ref{2.17}) 
\begin{eqnarray}\label{2.33}
&&\hspace{-.5cm}\int_B r_{-t}(\nu)\, \bmu(d\nu)-\bmu(B)=\bmu(\nu_{-t}: 
\nu\in B)-\bmu(B)\nonumber \\ 
&&\hspace{.5cm}=-\int_{s=0}^t\int_B\delta (A^f)(\nu_{-s})r_{-s}(\nu)\, 
\bmu(d\nu)\, ds\, . 
\end{eqnarray}
Recalling from Step 2, that $K'_\ve:=\{\nu_{-s}:\nu\in K_\ve,\ s\in 
[0,1]\}$ is also a compact subset of $G^-$, by means of Fubini's theorem 
and condition (jjj), we may conclude that 
\begin{eqnarray*}
&&\hspace{-.5cm}\int_B r_{-t}(\nu)\, \bmu(d\nu)-\bmu(B)=-\int_B\int_{s 
=0}^t\delta (A^f)(\nu_{-s})r_{-s}(\nu)\, ds\, \bmu(d\nu) 
\end{eqnarray*}
for all $t\in [0,1]$ and $B\in {\cal B}(E)$ with $B\subseteq K_\ve$. We 
obtain  
\begin{eqnarray}\label{2.19}
r_{-t}(\nu)-1=-\int_{s=0}^t\delta (A^f)(\nu_{-s})r_{-s}(\nu)\, ds 
\end{eqnarray}
first for all rational $t$ belonging to $[0,1]$ and $\bmu$-a.e. $\nu\in 
K_\ve$. Then, noting that (\ref{2.15}) and (\ref{2.16}) show that $t\to 
r_{-t}$ is $\bmu$-a.e. continuous on $[0,1]$, we get (\ref{2.19}) 
$\bmu$-a.e. for all $t\in [0,1]$, i. e. for $\bmu$-a.e. $\nu\in G^-$ and 
all $t\in [0,1]$. An immediate consequence of this and (\ref{2.17}) is 
(\ref{2.3}). 
\medskip

Using the monotone class theorem, the identity 
\begin{eqnarray}\label{2.37}
\int f(\nu)\delta (A^f)(\nu_{-t}) r_{-t}(\nu)\, \bmu (d\nu)=\int f\delta 
(A^f)(\nu_{-t})\, d\bmu\circ\nu_{-t}\, ,\quad t\in [0,1],
\end{eqnarray}
and $\delta (A^f)\in L^1(E,\bmu)$, it follows from (\ref{2.33}) that 
\begin{eqnarray*}
\int f(\nu)r_{-t}(\nu)\, \bmu  (d\nu)-\int f\, d\bmu=-\int_{s=0}^t\int f 
\delta (A^f)(\nu_{-s})\, d\bmu\circ\nu_{-s}\, ds 
\end{eqnarray*}
for all $t\in [0,1]$ and all $f\in L^\infty(E,\bmu)$. Part (b) is now a 
consequence of this equation for all $f\in L^\infty(E,\bmu)$ and the 
continuity of the expression (\ref{2.37}) in $t\in [0,1]$, cf. (\ref{2.12}). 
Part (d) has been proved in Step 1. 
\qed
\bigskip 

\nid 
{\bf Proof of Theorem \ref{Theorem2.7} } We prove the theorem for Case 1.
\medskip 

\nid
{\it Step 1 } In this step we derive an equation for $\bmu(\nu_{-t}:\nu 
\in\cdot\setminus G^-)$, $t\in [0,1]$. We pick up the proof of Theorem 
\ref{Theorem2.3} at (\ref{2.19*}) and conclude 
\begin{eqnarray*} 
\int f\, \bmu(d\nu_{-t})-\int f\, d\bmu =-\int_{s=0}^t\int f\delta (A^f) 
(\nu_{-s})\, \bmu(d\nu_{-s})\, ds+\int_0^t\int f\, \bmu_\gamma(d\nu_{-s}) 
\, ds 
\end{eqnarray*}
for all $f\in C_b(E)$ and all $t\in [0,1]$. This is part (a). The remainder 
of this step and Steps 2 through 4 below are devoted to the proof of the 
equivalence stated in part (b). 

We follow the proof of Theorem \ref{Theorem2.3} until (\ref{2.75}) and obtain    
\begin{eqnarray*} 
&&\hspace{-.5cm}\bmu(\nu_{-t}:\nu\in B)-\bmu(B)\vphantom{\int} \\ 
&&\hspace{.5cm}=-\int_{s=0}^t\int_B\delta (A^f)(\nu_{-s})\, \bmu(d 
\nu_{-s})\, ds+\int_0^t\bmu_\gamma(\nu_{-s}:\nu\in B)\, ds    
\end{eqnarray*}
for all $B\in {\cal B}(E)$ and all $t\in [0,1]$. Let us replace the set 
$B$ by $B\setminus G^-$. We obtain 
\begin{eqnarray}\label{2.45}
&&\hspace{-.5cm}\bmu(\nu_{-t}:\nu\in B\setminus G^-)\vphantom{\int} 
\nonumber \\ 
&&\hspace{.5cm}=-\int_{s=0}^t\int_{B\setminus G^-}\delta (A^f)(\nu_{-s}) 
\, \bmu(d\nu_{-s})\, ds+\int_0^t\bmu_\gamma(\nu_{-s}:\nu\in B\setminus 
G^-)\, ds    
\end{eqnarray} 
for all $B\in {\cal B}(E)$ and all $t\in [0,1]$. Given $\delta(A^f)$ 
and $\bmu_\gamma$, this relation can be considered as an equation for 
$\bmu(\nu_{-t}:\nu\in\cdot\setminus G^-)$, $t\in [0,1]$. 
\medskip 

\nid 
{\it Step 2 } Let us deal with uniqueness relative to this equation. 
Throughout this step, let $\bmu_1$ and $\bmu_2$ be two measures satisfying 
the hypotheses of Theorem \ref{Theorem2.7} such that $\bmu_1(\nu_{-t}:\nu 
\in\cdot\setminus G^-)$, $t\in [0,1]$, as well as $\bmu_2(\nu_{-t}:\nu\in 
\cdot\setminus G^-)$, $t\in [0,1]$, solve (\ref{2.45}). It follows that 
\begin{eqnarray}\label{2.46}
|\bmu_1-\bmu_2|(\nu_{-t}:\nu\in B\setminus G^-)\le\int_{s=0}^t\int_{B 
\setminus G^-}|\delta (A^f)(\nu_{-s})|\, |\bmu_1-\bmu_2|(d\nu_{-s})\, ds
\end{eqnarray} 
for all $B\in {\cal B}(E)$ and all $t\in [0,1]$. Here $|\bmu_1-\bmu_2|$ is, 
according to the Hahn decomposition of the signed measure $\bmu_1-\bmu_2$, 
defined by $|\bmu_1-\bmu_2|:=(\bmu_1-\bmu_2)^++(\bmu_1-\bmu_2)^-$. With 
\begin{eqnarray*} 
A(t,B):=\left\{
\begin{array}{cl}
\frac{\displaystyle\int_{B}|\delta (A^f)(\nu_{-t})|\, |\bmu_1-\bmu_2|(d\nu_{ 
-t})}{\displaystyle|\bmu_1-\bmu_2|(\nu_{-t}:\nu\in B)\vphantom{\displaystyle 
\int}} &\ \mbox{\rm if }\ |\bmu_1-\bmu_2|(\nu_{-t}:\nu\in B)>0 \\ 
0 &\ \mbox{\rm otherwise }\vphantom{\displaystyle\int}
\end{array}\right. 
\end{eqnarray*}
we have 
\begin{eqnarray*}
&&\hspace{-.5cm}\exp\left\{-\int_0^t A(s,B\setminus G^-)\, ds\right\}\cdot 
\int_0^t\int_{B\setminus G^-}|\delta (A^f)(\nu_{-s})|\, |\bmu_1-\bmu_2|(d 
\nu_{-s})\, ds \\ 
&&\hspace{.5cm}=\int_0^t\left(|\bmu_1-\bmu_2|(\nu_{-s}:\nu\in B)-\int_0^s 
\int_{B\setminus G^-}|\delta (A^f)(\nu_{-u})|\, |\bmu_1-\bmu_2|(d\nu_{-u}) 
\, du\right)\times \\ 
&&\hspace{1.5cm}\times\, A(s,B\setminus G^-)\cdot\exp\left\{-\int_0^s A(u,B 
\setminus G^-)\, du\right\}\, dt\le 0
\end{eqnarray*} 
since the expression in parentheses is non-positive by (\ref{2.46}). Using 
again (\ref{2.46}), it follows that for all $B\in {\cal B}(E)$ and all $t\in 
[0,1]$ we have $|\bmu_1-\bmu_2|(\nu_{-t}:\nu\in B\setminus G^-)=0$. In other 
words, the problem formulated below (\ref{2.45}) has at most one solution 
$\bmu(\nu_{-t}:\nu\in B\setminus G^-)$ for all $B\in {\cal B}(E)$ and all $t 
\in [0,1]$. 
\medskip 

\nid 
{\it Step 3 } Existence of a solution $\bmu(\nu_{-t}:\nu\in\cdot\setminus 
G^-)$, $t\in [0,1]$, to (\ref{2.45}) is obvious. It is given by the measure 
$\bmu$ introduced in the hypotheses of this theorem. In this step we are 
concerned with the representation of $\bmu(\nu_{-t}:\nu\in\cdot\setminus 
G^-)$, $t\in [0,1]$, in terms of $\delta(A^f)$ and $\bmu_\gamma$. Throughout 
this step we suppose (vj). Let us firstly verify well-definiteness of the 
expression 
\begin{eqnarray}\label{2.47}
I(t,B):=\int_{s=0}^t\int_{B\setminus G^-}\exp\left\{-\int_s^t\delta(A^f) 
(\nu_{-u})\, du\right\}\, \bmu_\gamma(d\nu_{-s})\, ds\, , 
\end{eqnarray}
$B\in {\cal B}(E)$, $t\in [0,1]$. Recalling that $\bmu_\gamma$ is 
concentrated on $\gamma^-G\equiv G\setminus G^-$ it turns out that 
\begin{eqnarray*}
I(t,B)\le\int_{s=0}^t\int_{\gamma^-G}\exp\left\{-\int_0^{t-s}\delta(A^f) 
(\nu_{-u})\, du\right\}\, \bmu_\gamma(d\nu)\, ds\, . 
\end{eqnarray*}
Condition (vj) says now that $I(t,B)<\infty$ for all $B\in {\cal B}(E)$ 
and all $t\in [0,1]$. Secondly, let us verify that $\bmu(\nu_{-t}:\nu\in 
\cdot\setminus G^-)=I(t,B)$, $t\in [0,1]$, is a representation of the 
unique solution to (\ref{2.45}) whenever $B\in {\cal B}(E)$. From 
(\ref{2.47}) and condition (vj) we deduce the continuous differentiability 
of $[0,1]\ni t\mapsto I(t,B)$. For $B\in {\cal B}(E)$ and $t\in [0,1]$ we 
get from (\ref{2.47})
\begin{eqnarray*} 
&&\hspace{-.5cm}I(t,B)=\int_{s=0}^t\frac{d}{ds}I(s,B)\, ds \\ 
&&\hspace{.5cm}=\int_{s=0}^t\frac{d}{ds}\left[\int_{u=0}^s\int_{B\setminus 
G^-}\exp\left\{-\int_u^s\delta(A^f)(\nu_{-v})\, dv\right\}\, \bmu_\gamma(d 
\nu_{-u})\, du\right]\, ds \\ 
&&\hspace{.5cm}=-\int_{s=0}^t\int_{u=0}^s\int_{B\setminus G^-}\delta (A^f) 
(\nu_{-s})\, \exp\left\{-\int_u^s\delta(A^f)(\nu_{-v})\, dv\right\}\, 
\bmu_\gamma(d\nu_{-u})\, du\, ds \\ 
&&\hspace{1.0cm}+\int_{s=0}^t\bmu_\gamma(\nu_{-s}:\nu\in B\setminus G^-)\, 
ds \\ 
&&\hspace{.5cm}=-\int_{s=0}^t\int_{B\setminus G^-}\delta (A^f)(\nu_{-s})\, 
I(s,d\nu)\, ds+\int_0^t\bmu_\gamma(\nu_{-s}:\nu\in B\setminus G^-)\, ds\, . 
\end{eqnarray*} 
In other words, we have verified that, for $B\in {\cal B}(E)$, $[0,1]\mapsto 
I(t,B)$ is a representation of the unique solution to (\ref{2.45}). 
\medskip 

\nid 
{\it Step 4 } We are now concerned with the equivalence of (vj) and (vj'). 
Let us first suppose (vj). From Step 3 we obtain  
\begin{eqnarray*} 
&&\hspace{-.5cm}\int_{s=0}^t\int_{B\setminus G^-}\exp\left\{-\int_s^t 
\delta (A^f)(\nu_{-u})\, du\right\}\, \bmu_\gamma(d\nu_{-s})\, ds=\bmu 
(\nu_{-t}:\nu\in B\setminus G^-) \\ 
&&\hspace{.5cm}=\bmu(\nu_{-t}:\nu\in B)-\bmu(\nu_{-t}:\nu\in B\cap G^-) 
\, ,\quad B\in {\cal B}(E),\ t\in [0,1].\vphantom{\int}
\end{eqnarray*} 
Recalling $\bmu(\nu_{-t}:\nu\in B\cap G^-)=\bmu\circ\nu_{-t}(B\cap G^-)= 
\bmu\circ\nu_{-t}(B)$ from the definition and discussion in Paragraph 
\ref{sec:1:1:4}, this yields 
\begin{eqnarray*}
&&\hspace{-.5cm}\bmu(\nu_{-t}:\nu\in B)-\bmu\circ\nu_{-t}(B)=\int_{s=0}^t 
\int_{B\setminus G^-}\exp\left\{-\int_s^t\delta(A^f)(\nu_{-u})\, du\right\}
\, \bmu_\gamma(d\nu_{-s})\, ds \\ 
&&\hspace{.5cm}=\int_{s=0}^t\int_B\exp\left\{-\int_s^t\delta (A^f)(\nu_{-u}) 
\, du\right\}\, \bmu_\gamma(d\nu_{-s})\, ds\, ,\quad B\in {\cal B}(E),\ t 
\in [0,1]. 
\end{eqnarray*} 
For the last line, we have taken into consideration that $\bmu_\gamma$ 
is concentrated on $\gamma^-G\equiv G\setminus G^-$ and therefore $\bmu_\gamma 
(\nu_{-s}:\nu\in B\cap G^-)=0$, $s\in [0,1]$. We have shown that (vj) implies 
(vj'). 

Now suppose (vj'). Choosing $B=E$ and $t=1$ in (vj'), and substituting 
$\nu_{-s}\lra\nu$ we get
\begin{eqnarray*}
1\ge\int_{s=0}^1\int\exp\left\{-\int_0^{1-s}\delta (A^f)(\nu_{-u})\, du 
\right\}\, \bmu_\gamma(d\nu)\, ds\, . 
\end{eqnarray*} 
The right-hand side is the norm of $\exp\left\{-\int_0^t\delta (A^f) 
(\nu_{-s})\, ds\right\}$ in $L^1(\gamma^-G\times [0,1]\, ,\, \bmu_\gamma 
\times\l_{[0,1]})$ in the sense of $(\nu,t)\in \gamma^- G\times [0,1]$, 
i. e. we have (vj'). 
\medskip 

\nid 
{\it Step 5 } Part (c) follows from (a) and Theorem \ref{Theorem2.3} (b).
\qed
\bigskip 

\nid 
{\bf Proof of Proposition \ref{Proposition2.4} } {\it Step 1 } For 
$t\in [0,1]$ and $h\in\t {\cal V}$, let us use the notation $h_t\equiv 
\Phi_t(h):=V^+_th$. We recall from the formulation of this proposition 
that $\Phi_{-t}$ stands for the inverse of $\Phi_t$, $t\in [0,1]$, 
i. e. $\Phi_{-t}$ maps $\{h_t:h\in \t {\cal V}\}$ bijectively onto 
$\t {\cal V}$. Let us also use the notation $h_{-t}\equiv\Phi_{-t}h$, 
$h\in\{\t h_t:\t h\in \t {\cal V}\}$, $t\in [0,1]$. In this step we 
prove the existence of certain directional derivatives of $\Phi_{-t}$ 
and calculate them. 

Let $t\in [0,1]$. For well-definiteness in the calculation below, let 
us keep in mind the following. In this Proposition we suppose that for 
all $\t h\in \t {\cal V}$ there is a pre-neighborhood $U_{\t h}(0)$ of 
$0\in X$ in the topology of $X$ such that $\{\t h+i\circ h:h\in U_{\t h} 
(0)\}\subseteq\{V_1^+(h):h\in D(V_1^+)\}\subseteq\{V_t^+(h):h\in D 
(V_t^+)\}$, recall also (\ref{1.2}). According to the hypothesis $G 
\subseteq E\cap\{V_t^+(h)\cdot\lambda:h\in\t {\cal V}\}$ and $\bmu(G)=1$, 
for $\bmu$-a.e. $\nu=h\cdot\lambda\in G$ and $u\in {\Bbb R}$ such that 
$|u|$ is sufficiently small, it holds that 
\begin{eqnarray}\label{2.48}
&&\hspace{-.5cm}x=\frac1u\left(\Phi_t\left((h+u\cdot x)_{-t}\right)- 
\Phi_t\left(h_{-t}\right)\right)\nonumber \\ 
&&\hspace{.5cm}=\nabla\Phi_t(h_{-t})\left(\frac{(h+u\cdot x)_{-t}- 
h_{-t}}{u}\right)+\frac1u R\left((h+u\cdot x)_{-t}-h_{-t}\right)\, ,
\end{eqnarray} 
where $x=i\circ\t x$ with $\t x\in X$. Moreover,  
\begin{eqnarray*}
\lim_{\|f\|_{L^v(D)}\to 0}R(f)/\|f\|_{L^v(D)}=0\, . 
\end{eqnarray*}
It follows now from (\ref{2.48}) that 
\begin{eqnarray*}
&&\hspace{-.5cm}\left(\nabla\Phi_t(h_{-t})\right)^{-1}(x)-\frac{R\left( 
(h+u\cdot x)_{-t}-h_{-t}\right)}{\left\|(h+u\cdot x)_{-t}-h_{-t}\right 
\|_{L^v(D)}}\cdot\frac{\left\|(h+u\cdot x)_{-t}-h_{-t}\right\|_{L^v(D)}} 
{u} \\ 
&&\hspace{.5cm}=\frac{(h+u\cdot x)_{-t}-h_{-t}}{u}\, ,\quad x=i\circ\t x, 
\  \t x\in X. 
\end{eqnarray*} 
Letting $u\to 0$ it follows from the hypothesized continuity of 
$\Phi_{-t}$ that $(h+u\cdot x)_{-t}-h_{-t}\to 0$ which with ({2.49}) 
implies 
\begin{eqnarray*}
\left(\nabla\Phi_t(h_{-t})\right)^{-1}(x)=\lim_{u\to 0}\frac{(h+u\cdot x 
)_{-t}-h_{-t}}{u}\equiv\lim_{u\to 0}\frac{\Phi_{-t}(h+u\cdot x)-\Phi_{-t} 
(h)}{u} 
\end{eqnarray*} 
whenever $x=i\circ\t x$ and $\t x\in X$. Using the symbol $\partial_x$ 
for the directional derivative in direction $x$ we may refer to this 
result as 
\begin{eqnarray}\label{2.50}
\partial_x\Phi_{-t}(h)=\left(\nabla\Phi_t(h_{-t})\right)^{-1}(x)\, , 
\quad x=i\circ\t x,\ \t x\in X,  
\end{eqnarray} 
for $\bmu$-a.e. $\nu\in G$ and $\nu=h\cdot\lambda$, $t\in [0,1]$. 
\medskip

\noindent
{\it Step 2 } For $h\cdot\lambda=\nu$ and $h\in \t {\cal V}$ recall the 
notation $A^f\nu=a^fh\cdot\lambda$ of Subsection \ref{sec:2:1}. On the 
one hand, for $h\in \t {\cal V}$ and $t\in [0,1]$ we have 
\begin{eqnarray*}
\left.\frac{d}{du}\right|_{u=0}(h+u\cdot a^fh)_t=\nabla\Phi_t(h)(a^fh) 
\end{eqnarray*}
as an ordinary directional derivative in $L^v(D)$. For the well-definiteness 
of the left-hand side, recall that $\t {\cal V}$ is an open set in $L^v(D)$. 
On the other hand, it holds that 
\begin{eqnarray}\label{2.51}
a^fh_t=\left.\frac{d}{du}\right|_{u=0}h_{t+u}=\left.\frac{d}{du}\right 
|_{u=0}\Phi_t\circ h_u=\nabla\Phi_t(h)(a^fh) 
\end{eqnarray}
for all $t\in [0,1]$ by the chain rule for Fr\'echet derivatives. By the 
hypothesis $G\subseteq E\cap\{h\cdot\lambda:h\in\t {\cal V}\}$ and $\bmu(G) 
=1$, we get for $\bmu$-a.e. $\nu\in G$ and $t\in [0,1]$ 
\begin{eqnarray*}
\left.\frac{d}{du}\right|_{u=0}\left((\nu+u\cdot A^f\nu)_t,g\vphantom{l^1} 
\right)=\left(A^f\nu_t,g\vphantom{l^1}\right)\, .  
\end{eqnarray*}
Together with (j) this is (v') for $t\in [0,1]$. 
\medskip

\noindent
{\it Step 3 } It follows from (\ref{2.51}) that for $t\in [0,1]$ and $h\in 
\{\t h_t:\t h\in \t {\cal V}\}$  
\begin{eqnarray*}
a^fh\in\left\{\nabla\Phi_t(h_{-t})(f):f\in L^v(D)\right\}\quad\mbox{\rm as 
well as}\quad\left(\nabla\Phi_t(h_{-t})\right)^{-1}(a^fh)=a^fh_{-t}\, .
\end{eqnarray*}
With (j) and (\ref{2.50}) we obtain 
\begin{eqnarray*}
\frac{(h+u\cdot a^fh)_{-t}-h_{-t}}{u}\stack{n\to\infty}{\lra}a^fh_{-t}\quad 
\mbox{\rm in } L^v(D)\vphantom{\frac11}
\end{eqnarray*}
whenever $t\in [0,1]$ and $h\in\{\t h_t:\t h\in \t {\cal V}\}$. An immediate 
consequence of the last limit and hypothesis $G\subseteq E\cap\{V_t^+(h) 
\cdot\lambda:h\in \t {\cal V}\}$ together with $\bmu(G)=1$ is, that for 
$\bmu$-a.e. $\nu\in G$ we get 
\begin{eqnarray*}
\left.\frac{d}{du}\right|_{u=0}\left((\nu+u\cdot A^f\nu)_{-t},g\vphantom 
{l^1}\right)=\left(A^f\nu_{-t},g\vphantom{l^1}\right)\, ,\quad t\in [-1,0].    
\end{eqnarray*}
Switching $t$ to $-t$, this is (v') for $t\in [-1,0]$. 
\medskip

\noindent
{\it Step 4 } Let $g\in C_0(D)$. According to the setup of the proposition, 
for $t\in [0,1]$ and $\nu=h\cdot\lambda\in E\cap\{f\cdot\lambda:f\in 
\t {\cal V}\}$, i. e. for $\bmu$-a.e. $\nu=h\cdot\lambda\in G$, we have 
\begin{eqnarray*}
&&\hspace{-.5cm}\left|\left.\frac{d}{du}\right|_{u=0}\left((\nu+u\cdot 
i\circ k\cdot\lambda)_t,g\vphantom{l^1}\right)\right|=\left|\left( 
\partial_{i\circ k}\Phi_t(h),g\vphantom{l^1}\right)\right|=\left|\left( 
\nabla\Phi_t(h)(i\circ k),g\right)\right| \\ 
&&\hspace{.5cm}\le\left\|\left(\nabla\Phi_t(h)(\cdot),g\right)\right 
\|_{L^w(D)}\, \|i\circ k\|_{L^v(D)}\vphantom{\frac11} \\ 
&&\hspace{.5cm}\le\left\|\left(\nabla\Phi_t(h)(\cdot),g\right)\right 
\|_{L^w(D)}\, \|i\|_{B(X,L^v(D))}\, \|k\|_{X}\, ,\quad k\in X,\vphantom 
{\frac11} 
\end{eqnarray*} 
where $1/v+1/w=1$ for $1<v<\infty$ and $w=\infty$ if $v=1$.
Here, the equality of $\left.\frac{d}{du}\right|_{u=0}\int_{y\in D} 
(h+u\cdot i\circ k)_t(y)\ldots\, dy$ with $\int_{y\in D}\left.\frac{d} 
{du}\right|_{u=0}(h+u\cdot i\circ k)_t(y)\ldots\, dy$, representing the 
first equality sign, is correct because the derivative $\left.\frac{d} 
{du}\right|_{u=0}(h+u\cdot i\circ k)_t$ is a limit in $L^v(D)$ and 
$\left.\frac{d}{du}\right|_{u=0}\int_{y\in D}(h+u\cdot i\circ k)_t(y) 
g(y)\, dy$ is the duality pairing of this limit and $g\in C_0(D)\subset 
L^w(D)$. 

Similarly, from (\ref{2.50}) we obtain for $t\in [-1,0]$ and $\bmu$-a.e. 
$\nu=h\cdot\lambda\in G$ 
\begin{eqnarray*}
&&\hspace{-.5cm}\left|\left.\frac{d}{du}\right|_{u=0}\left((\nu+u\cdot 
i\circ k\cdot\lambda)_t,g\vphantom{l^1}\right)\right|=\left|\left( 
\partial_{i\circ k}\Phi_t(h),g\vphantom{l^1}\right)\right|=\left|\left( 
\left(\nabla\Phi_t(h_t)\right)^{-1}(i\circ k),g\right)\right| \\ 
&&\hspace{.5cm}\le\left\|\left(\left(\nabla\Phi_t(h_t)\right)^{-1}(\cdot) 
,g\right)\right\|_{L^w(D)}\, \|i\circ k\|_{L^v(D)}\vphantom{\frac11} \\ 
&&\hspace{.5cm}\le\left\|\left(\left(\nabla\Phi_t(h_t)\right)^{-1}(\cdot)
,g\right)\right\|_{L^w(D)}\, \|i\|_{B(X,L^v(D))}\, \|k\|_{X}\, ,\quad k 
\in X\, .\vphantom{\frac11} 
\end{eqnarray*} 
Thus, for all $t\in [-1,1]$, $X\ni k\mapsto\left.\frac{d}{du}\right|_{u=0} 
\left((\nu+u\cdot i\circ k\cdot\lambda)_t,g\vphantom{l^1}\right)$ is a 
bounded linear functional on $X$ which we denote by $r(\nu,t,g)\in X^\ast$.  
In other words, we have shown that the derivative (\ref{2.2}) in condition 
(jv) exists for $\bmu$-a.e. $\nu\in G$ and all $\, t\in [-1,1]$, $g\in C_0 
(D)$, as well as $k\in X$, for some $r(\nu,t,g)\in X^\ast$. 
\qed 
\bigskip 

\nid
{\bf Proof of Proposition \ref{Proposition2.5} } Part (a), i. e. the 
verification of (jv'), is an immediate consequence of the hypothesis 
formulated in the proposition. Recall from Paragraphs \ref{sec:1:1:5} 
and \ref{sec:1:1:4} the definition of $G$ and $\bmu(G)=1$. 
\medskip 

\nid
(b) We prove the claim under the hypotheses of Case 1. Let $t\in [0,1]$ 
and $f\in\t C_b^1(E)$. Using Lemma \ref{Lemma2.2}, for which we have 
supposed (a4), and (j) as well as (jj) we get
\begin{eqnarray*} 
&&\hspace{-.5cm}\int\left.\frac{d}{du}\right|_{u=0}\left((\nu+u\cdot 
A^f\nu)_t,g\vphantom{l^1}\right)\cdot f\, d\bmu=\int\left(D\Psi_{t,g}, 
a^f\right)\cdot f\, d\bmu \\ 
&&\hspace{.5cm}=\int\left(D(\Psi_{t,g}\cdot f),a^f\right)\, d\bmu - 
\int\left(Df,a^f\right)\cdot\Psi_{t,g}\, d\bmu \\ 
&&\hspace{.5cm}=-\int\Psi_{t,g}\cdot f\cdot\delta(A^f)\, d\bmu+\int 
\Psi_{t,g}\cdot f\, d\bmu_\gamma -\int\left(Df,a^f\right)\cdot\Psi_{t,g} 
\, d\bmu\, .    
\end{eqnarray*} 
Theorem \ref{Theorem2.3} (d) states the existence of all derivatives 
$\frac{d}{dt} f(\nu_t)$, $f\in\t C^1_b(E)$, and hence also of $\frac{d}{dt} 
(\nu_t,g)$ identified with $\frac{d}{dt}\Psi_{t,g}$, as derivatives 
in $L^1(E,\bmu)$. It follows now from Theorem \ref{Theorem2.7} (a) that 
\begin{eqnarray*} 
&&\hspace{-.5cm}\int\left.\frac{d}{du}\right|_{u=0}\left((\nu+u\cdot 
A^f\nu)_t,g\vphantom{l^1}\right)\cdot f\, d\bmu \\ 
&&\hspace{.5cm}=\left.\frac{d}{ds}\right|_{s=0}\int \Psi_{t,g}\cdot 
f\, d\bmu(d\nu_{-s})-\int\left(Df,a^f\right)\cdot\Psi_{t,g}\, d\bmu 
 \\ 
&&\hspace{.5cm}=\left.\frac{d}{ds}\right|_{s=0}\int\Psi_{t,g}(\nu_s) 
\cdot f\, d\bmu \\  
&&\hspace{1cm}+\left[\left.\frac{d}{ds}\right|_{s=0}\int\Psi_{t,g}\cdot 
f(\nu_s)\, d\bmu-\int\left(Df,a^f\right)\cdot\Psi_{t,g}\, d\bmu\right] 
 \\ 
&&\hspace{.5cm}=\int\left.\frac{d}{ds}\right|_{s=0}\Psi_{t,g}(\nu_s) 
\cdot f\, d\bmu\, .
\end{eqnarray*} 
Recalling the definition of $\Psi_{t,g}$ we obtain 
\begin{eqnarray*} 
&&\hspace{-.5cm}\int\left.\frac{d}{du}\right|_{u=0}\left((\nu+u\cdot 
A^f\nu)_t,g\vphantom{l^1}\right)\cdot f\, d\bmu=\int\left(A^f\nu_t,g 
\vphantom{l^1}\right)\cdot f\, d\bmu\, ,\quad t\in [0,1]. 
\end{eqnarray*} 
In other words, we have (v) in Case 1. The verification of (v) in Case 2 
is similar. 
\qed

\section{Boltzmann Type Equation}
\label{sec:3}
\setcounter{equation}{0} 

In this section, we shall focus on the mathematical description of a 
rarefied gas in a vessel with diffusive and specular reflective boundary 
conditions. 

Let $d\in\{2,3\}$. For a bounded set $\Gamma\subseteq\mathbb{R}^d$ let 
$\overline{\Gamma}$ denote its closure in $\mathbb{R}^d$. Let $\Omega 
\subset\mathbb{R}^d$ be a bounded domain, called the {\it physical 
space}. Let us suppose that $\partial\Omega=\bigcup_{i=1}^{n_\partial} 
\overline{\Gamma_i}$ where the $\Gamma_i$ are $(d-1)$-dimensional 
manifolds which are smooth up to the boundary and, at the same time, 
open sets in the topology of $\partial\Omega$, $i\in \{1,\ldots,n_ 
\partial\}$. Furthermore, assume that for $i\neq j$ the intersection 
$\overline{\Gamma_i}\cap\overline{\Gamma_j}$ either is the empty set or 
a smooth closed $(d-2)$-dimensional manifold which, by definition, for 
$d=2$ is just a single point in $\mathbb{R}^2$. 

In case of $d=3$, for any two neighboring $\Gamma_i$ and $\Gamma_j$ 
let us call $x\in\overline{\Gamma_i}\cap\overline{\Gamma_j}$ an {\it 
end point} if for any sufficiently small open ball $B_x$ in $\mathbb 
{R}^d$ with center $x$ the set $\partial B_x\cap (\overline{\Gamma_i} 
\cap\overline{\Gamma_j})$ consists of exactly one point. For the 
boundary $\partial\Omega$ let us also assume that there is $\xi\in 
(0,\pi)$ such that for any two neighboring $\Gamma_i$ and $\Gamma_j$ 
and any $x\in\overline{\Gamma_i}\cap\overline{\Gamma_j}$, not being 
an end point if $d=3$, we have the following. For the angle $\xi(i,j; 
x)$ between the rays $R_i$ and $R_j$ from $x$ tangential to $\Gamma_i$ 
and $\Gamma_j$ respectively, and orthogonal to $\overline{\Gamma_i} 
\cap\overline{\Gamma_j}$ if $d=3$, it holds that 
\begin{eqnarray*} 
\xi<\xi(i,j;x)<2\pi-\xi\, . 
\end{eqnarray*}
The latter guarantees compatibility with Assumption A in the proof 
of Theorem 2.1 in \cite{CPW98}, as we have discussed in Remark 3 
of \cite{Lo20}. Theorem 2.1 of \cite{CPW98} contributes to the proof 
of Proposition \ref{Proposition3.3**} (a) below, cf. Remarks 3 and 4 
in \cite{Lo20}. For illustrations of the boundary and the boundary 
conditions below we also refer to \cite{Lo20}, Sections 1 and 2. 

Let $V:=\{v\in \mathbb{R}^d:0<v_{min}<|v|<v_{max}<\infty\}$ be the 
{\it velocity space} and let $\lambda>0$. Denote by $n(r)$ the outer 
normal at 
\begin{eqnarray*} 
r\in\partial^{(1)}\Omega:=\bigcup_{i=1}^{n_{\partial}}\Gamma_i\, ,  
\end{eqnarray*} 
indicate the inner product in $\mathbb{R}^d$ by ``$\circ$", and let 
$R_r(v):=v-2v\circ n(r)\cdot n(r)$ for $(r,v)\in\partial^{(1)}\Omega 
\times V$. 
\medskip 

For $(r,v,t)\in\Omega\times V\times [0,\infty)$, consider the {\it 
Boltzmann type equation} 
\begin{eqnarray}\label{3.1} 
\frac{d}{dt}\, p(r,v,t)=-v\circ\nabla_rp(r,v,t)+\lambda Q(p,p)\, (r,v,t)
\end{eqnarray} 
with boundary conditions 
\begin{eqnarray}\label{3.2} 
p(r,v,t)=\omega\, p(r,R_r(v),t) +(1-\omega)J(r,t)(p)M(r,v)\, ,
\end{eqnarray} 
$r\in\partial^{(1)}\Omega$, $v\circ n(r)\le 0$, 
for some $\omega\in [0,1)$ and initial probability density $p(0,\cdot,\cdot) 
:=p_0$ on $\Omega\times V$. Consider also its {\it integrated (mild) version} 
\begin{eqnarray}\label{3.3} 
p(r,v,t)=S(t)\, p_0(r,v)+\lambda\int_0^t S(t-s)\, Q(p,p)\, (r,v,s) 
\, ds\, . 
\end{eqnarray} 

The following global conditions (k)-(wkk) on the terms in 
(\ref{3.1})-(\ref{3.3}) and the shape of $\Omega$ will be in force 
throughout this section.
\begin{itemize} 
\item[(k)] For all $t\ge 0$ and all $r\in\partial^{(1)}\Omega$, the 
function $J$ is given by 
\begin{eqnarray*} 
J(r,t)(p)=\int_{v\circ n(r)\ge 0}v\circ n(r)\, p(r,v,t)\, dv\, .  
\end{eqnarray*} 
\item[(kk)] The real function $M$ on $\{(r,v):r\in\partial^{(1)}\Omega,
\ v\in V,\ v\circ n(r)\le 0\}$ has positive lower and upper bounds 
$M_{\rm min}$ and $M_{\rm max}$, is continuous on every $\Gamma_i$, $i\in 
\{1,\ldots ,n_\partial\}$, and satisfies
\begin{eqnarray*} 
\int_{v\circ n(r)\le 0}|v\circ n(r)|\, M(r,v)\, dv=1\, ,\quad r\in 
\partial^{(1)}\Omega. 
\end{eqnarray*} 
\item[(kkk)] $S(t)$, $t\ge 0$, is given by the solution to the initial 
boundary value problem 
\begin{eqnarray*}
\left(\frac{d}{dt}+v\circ\nabla_r\right)(S(t)p_0)(r,v)=0\, ,\quad (r,v,t) 
\in \Omega\times V\times [0,\infty)\, ,  
\end{eqnarray*} 
$S(0)p_0=p_0$, and 
\begin{eqnarray*} 
(S(t)p_0)(r,v)=\omega\, (S(t)p_0)(r,R_r(v))+(1-\omega)J(r,t)(S(\cdot)p_0) 
M(r,v)\, , 
\end{eqnarray*} 
$t>0$, for all $(r,v)\in\partial^{(1)}\Omega\times V$ with $v\circ n(r)\le 
0$. We call $S(t)$, $t\ge 0$, {\it Knudsen type semigroup}. If $\omega=0$ 
we indicate this in the notation by $S_0(t)$, $t\ge 0$. We  follow 
\cite{CPW98} and call $S_0(t)$, $t\ge 0$, just {\it Knudsen semigroup}. 
\item[(kw)] Denoting by $\chi$ the indicator function and setting $p:=0$ 
as well as $q:=0 $ on $\Omega\times (\mathbb{R}^d\setminus V)\times 
[0,\infty)$, the {\it collision operator} $Q$ is given by 
\begin{eqnarray*}
Q(p,q)(r,v,t)&&\hspace{-.5cm}=\frac12\int_{\Omega}\int_V\int_{S_+^{d-1}} 
B(v,v_1,e)h_\gamma(r,y)\chi_{\{(v^\ast ,v_1^\ast)\in V\times V\}}\times \\ 
&&\hspace{-1.0cm}\times\left(\left(p(r,v^\ast ,t)q(y,v_1^\ast ,t)-p(r,v,t) 
q(y,v_1,t)\vphantom{l^1}\right)\vphantom{\dot{f}}\right.\vphantom{\int} \\ 
&&\hspace{-.7cm}\left.+\left(q(r,v^\ast ,t)p(y,v_1^\ast ,t)-q(r,v,t)p(y, 
v_1,t)\vphantom{l^1}\right)\vphantom{\dot{f}}\right)\, de\, dv_1\, dy\, .  
\vphantom{\int} 
\end{eqnarray*} 
Here $S^{d-1}$ is the unit sphere. Moreover, $S^{d-1}_+\equiv S^{d-1}_+(v- 
v_1):=\{e\in S^{d-1}:e\circ(v-v_1)>0\}$, $v^\ast :=v-e\circ (v-v_1)\, e$, 
$v_1^\ast:=v_1+e\circ (v-v_1)\, e$ for $e\in S^{d-1}_+$ as well as $v,v_1 
\in V$, and $de$ refers to the normalized Lebesgue measure on $S^{d-1}_+$.  
\item[(w)] The {\it collision kernel} $B$ is assumed to be non-negative 
and continuous on $V\times V\times S^{d-1}$, symmetric in $v$ and $v_1$, 
and satisfying $B(v^\ast,v_1^\ast,e)=B(v,v_1,e)$ for all $v,v_1\in V$ and 
$e\in S^{d-1}_+$ for which $(v^\ast ,v_1^\ast)\in V\times V$. Furthermore, 
there is a constant $b\in (0,\infty)$ such that with $\varphi(v_1) \equiv 
\varphi(v;v_1,e) :=-v+2v_1+(|v-v_1|^2/e\circ(v-v_1))\, e$, defined for 
$(v,v_1,e)\in \mathbb{R}^d \times \mathbb{R}^d \times S^{d-1}$, we have 
\begin{eqnarray*}
B(v,\varphi(v_1),e)\le 2^{1-d}b\cdot\left(\frac{v-v_1}{|v-v_1|}\circ e 
\right)^2
\end{eqnarray*} 
for all $(v,v_1,e)\in V \times V \times S^{d-1}_+$ such that $\varphi(v_1) 
\in V$.
\end{itemize} 
For a discussion of the last hypothesis we refer to \cite{Lo20}, Section 2. 
\begin{itemize}
\item[(wk)] $h_\gamma$ is a continuous function on $\overline{\Omega 
\times\Omega}$ which is non-negative and symmetric, and vanishes for $|r-y| 
\ge\gamma>0$.  
\end{itemize} 

Let us turn to the last of the global conditions, namely condition (wkk) 
below. For this, recall the structure of $\partial\Omega$ from the beginning 
of this section. For $(y,v)\in\overline{\Omega}\times V$ introduce 
\begin{eqnarray*}
T_\Omega\equiv T_\Omega(y,v):=\inf\{s>0:y-sv\not\in\Omega\}  
\end{eqnarray*} 
and 
\begin{eqnarray*} 
y^-\equiv y^-(y,v):=y-T_\Omega(y,v)v\, . 
\end{eqnarray*} 
For all $(y,v)\in\bigcup_{i=1}^{n_\partial}\Gamma_i\times V=\partial^{ 
(1)}\Omega\times V$ with $v\circ n(y)\le 0$ as well as $y^-(y,R_y(v))\in 
\partial^{(1)}\Omega$, let 
\begin{eqnarray*}
\sigma(y,v):=\left(y^-(y,R_y(v)),R_y(v)\right)\, . 
\end{eqnarray*} 
We note that there exist $m_\partial\equiv m_\partial(\sigma)\in\mathbb{N}$ 
and mutually disjoint sets $G_i\subseteq\{(r,v):r\in\partial^{(1)}\Omega, 
\ v\in V,\ v\circ n(r)\le 0\}$, $i\in\{1,\ldots ,m_\partial\}$, satisfying 
\begin{eqnarray*} 
\bigcup_{i=1}^{m_\partial}\overline{G_i}=\overline{\{(r,v):r\in 
\partial^{(1)}\Omega,\ v\in V,\ v\circ n(r)\le 0\}} 
\end{eqnarray*} 
such that the following holds. {\it For each $i\in\{1,\ldots ,m_\partial 
\}$ there exists $j\in\{1,\ldots ,m_\partial\}$ such that $\sigma$ maps 
$G_i$ bijectively and continuously to $G_j$.}

In this way we understand $\sigma$ as a map defined for $(y,v)\in 
\partial^{(1)}\Omega\times V$ with $v\circ n(y)\le 0$ and $y^-(y,R_y(v)) 
\in\partial^{(1)}\Omega$. For its $k$-fold composition $\sigma^k(y,v)= 
(y^{(k)}(y,v),v^{(k)}(y,v))\equiv (y^{(k)},v^{(k)})$, $k\in\mathbb{N}$, 
we suppose also that by iteration $y^{(j)}\in\partial^{(1)}\Omega$, $j\in 
\{1,\ldots ,k\}$, and we indicate this by the phrase ``\,for a.e. $y\in 
\partial^{(1)}\Omega\, $". Set $(y^{(0)},v^{(0)}):=(y,v)$. An important 
role will play the following condition on the shape of $\Omega$. 
\begin{itemize} 
\item[(wkk)] There exist $k^{(0)}\in\mathbb{N}$ and $\sigma_{min}>0$ 
such that for a.e. $(y,v)\in\partial^{(1)}\Omega\times V$ with $v\circ 
n(y)\le 0$ and all $j\in\mathbb{Z}_+$ we have 
\end{itemize} 
\vspace{-.5cm}
\begin{eqnarray*}
\sigma_{min}\le\left|y^{(j)}-y^{(j+1)}\right|+\left|y^{(j+1)}-y^{(j+2)} 
\right|+\ldots +\left|y^{(j+k^{(0)}-1)}-y^{(j+k^{(0)})}\right|\, .
\end{eqnarray*} 
Condition (wkk) is primarily used in the proofs of the subsequent 
Proposition \ref{Proposition3.3**} (a) and Theorem \ref{Theorem5.8}, 
which are Lemma 3.1, Corollary 4.6 (a), and Theorem 5.8 in \cite{Lo20}. 
These assertions are fundamental for the whole section. The verifiability 
of condition (wkk) is discussed in Remark 10 of \cite{Lo20}. 

\subsection{Global Solutions}\label{sec:3:1} 

In this subsection we summarize the results of \cite{Lo20} we need for the 
further analysis of the present paper. Throughout the whole subsection 
suppose (k)-(wkk). 

In order to formulate conditions under which the Knudsen-type semigroup 
$S(t)$, $t\ge 0$, can be extended to a group we need some preparations, 
for more details we refer to Remarks 8-10 in \cite{Lo20}. Set $\mathbf{M}
:=\{(r,w)\in\partial^{(1)}\Omega\times V$ with $w\circ n(r)\ge 0\}$. For  
$(r,w)\in\mathbf{M}$ let us introduce $\mathbb{T}^0(r,w):=(r,w)$, 
\begin{eqnarray*}
\mathbb{T}(r,w):=\left(r^-(r,w),R_{r^-}(w)\right)\, , 
\end{eqnarray*} 
and the abbreviations $l_k\equiv l_k(r,w):=T_\Omega\left(\mathbb{T}^k 
(r,w)\right)$, $k\in\mathbb{Z}_+$. Noting that 
\begin{eqnarray*}
\mathbb{T}^{-1}(r,w)=\left(r^-(r,-R_r(w)),R_r(w)\right) 
\end{eqnarray*} 
for all $(r,w)\in\mathbf{M}$ we also abbreviate $l_k\equiv l_k(r,w):= 
T_\Omega \left(\mathbb{T}^k(r,w)\right)$,  $-k\in\mathbb{N}$. Let 
diam$(\Omega)=\sup\{|r_1-r_2|:r_1,r_2\in\Omega\}$ denote the diameter 
of $\Omega$. 

The following well known facts from the ergodic theory of $d$-dimensional 
mathematical billiards, $d=2,3$, can be taken from \cite{CM03}. First 
recall the Ergodic Theorem of Birkhoff-Khinchin, Theorem II.1.1 together 
with the Corollary II.1.4. Then recall from \cite{CM03}, (IV.2.3) and 
the formula above (IV.2.3), that the billiard map is measure preserving 
with respect to some probability measure $\nu$ on $\{(r,v)\in\partial^{ 
(1)}\Omega\times S^{d-1}:v\circ n(r)\ge 0\}$. In fact, the Radon-Nikodym 
derivative of $\nu$ with respect to the Lebesgue measure $\lambda$ is 
\begin{eqnarray*}
\frac{d\nu}{d\lambda}=c_\nu\cdot n(r)\circ v 
\end{eqnarray*} 
where $c_\nu>0$ is the normalizing constant. It follows that for a.e. 
$(r,w)\in\mathbf{M}$ there is a $\tau(r,w)\le\, $diam$(\Omega)/v_{min}$ 
such that 
\begin{eqnarray}\label{3.72*}
\lim_{n\to\infty}\frac1n\sum_{k=0}^{n-1}l_k=\lim_{n\to\infty}\frac1n 
\sum_{k=1}^nl_{-k}=\tau(r,w)\, . 
\end{eqnarray} 
Relation (\ref{3.72*}) holds even a.e. on $\{(r,v)\in\partial^{(1)} 
\Omega\times S^{d-1}:v\circ n(r)\ge 0\}$, where $S^{d-1}$ denotes 
the unit sphere. 

Introduce 
\begin{eqnarray*}
c_\tau:=\mathop{\mathrm{ess~inf}}\limits_{(r,w)\in\mathbf{M}}\tau(r,w) 
\cdot |w| 
\end{eqnarray*} 
as well as 
\begin{eqnarray*}
C_\tau\equiv C_{\tau,k_0}:=\mathop{\mathrm{ess~sup}}\limits_{(r,w)\in 
\mathbf{M},\, k\ge k_0}\frac{k\cdot\tau(r,w)-\sum_{i=1}^kl_i(r,w)}{k 
\cdot\tau (r,w)}\, ,\quad k_0\in\mathbb{N}. 
\end{eqnarray*} 
\smallskip

\begin{theorem}\label{Theorem3.19} (Theorem 4.5 in \cite{Lo20}) 
Suppose 
\begin{eqnarray}\label{3.8*}
(1-C_{\tau,k_0})\, c_\tau>0\quad\mbox{\rm for some }k_0\in\mathbb{N}   
\end{eqnarray} 
and 
\begin{eqnarray*}
2^{-\frac{1}{k_0}}<\omega<1\, . 
\end{eqnarray*} 
Then $S(t)$, $t\ge 0$, extends to a strongly continuous group in 
$L^1(\Omega\times V)$ which we will denote by $S(t)$, $t\in\mathbb{R}$. 
\end{theorem}

Verifiability of (\ref{3.8*}), also in relation to hypothesis (wkk), is 
discussed in \cite{Lo20}, Remark 10. Let $(A,D(A))$ denote the 
infinitesimal generator of the strongly continuous group $S(t)$, $t\in 
\mathbb{R}$, in $L^1(\Omega\times V)$. We continue with a spectral property. 
\begin{proposition}\label{Proposition3.3*} (Corollary 4.4 (b) in 
\cite{Lo20}). There exists $m_1<0$ such that for all $\mu\in\mathbb{C}$ 
with $\mathfrak{Re}\mu<m_1$ and $g\in L^1(\Omega\times V)$ the equation 
$\mu f-Af=g$ has a unique solution $f\in D(A)$.
\end{proposition}

In addition we have the following boundedness properties. 
\begin{proposition}\label{Proposition3.3**} (a) (Lemma 3.1 and Corollary 
4.6 (a) in \cite{Lo20}) Let the conditions of Theorem \ref{Theorem3.19} 
be satisfied. Let $p_0\in L^1(\Omega\times V)$ and $t'\in\mathbb{R}$ such 
that $S(t')p_0\in L^\infty(\Omega\times V)$. Then there are finite real 
numbers $p_{0,{\rm min}}$ and $p_{0,{\rm max}}$ such that 
\begin{eqnarray}\label{3.11}
p_{0,{\rm min}}\le S(t)\, p_0\le p_{0,{\rm max}}\quad\mbox{\rm a.e. on } 
\Omega\times V
\end{eqnarray} 
for all $t\ge t'$. In particular, if $S(t')p_0\ge 0$ a.e. and $\|1/ 
S(t')p_0\|_{L^\infty(\Omega\times V)}<\infty$ then there exists $p_{ 
0,{\rm min}}>0$. \\ 
(b) (Part of Corollary 4.6 (c) in \cite{Lo20}) There exists $M\ge 1$ 
such that with $m_1$ given by Proposition \ref{Proposition3.3*} we have 
for all $p_0\in L^1(\Omega\times V)$ 
\begin{eqnarray}\label{4.52}
\|S(t)\, p_0\|_{L^1(\Omega\times V)}\le Me^{m_1 t}\|p_0\|_{L^1(\Omega 
\times V)}\, ,\quad t\le 0.  
\end{eqnarray} 
\end{proposition}

Let $\1$ denote the function constant to one on $\Omega\times V$. Furthermore, 
for $t\in [0,\infty]$ introduce 
\begin{eqnarray}\label{3.181}
c^{\1}_{t,\rm max}:=\sup\|S(\tau)\, \1\|_{L^\infty(\Omega\times V)} 
\end{eqnarray} 
where, for $0\le t<\infty$ the supremum is taken over $\tau\in [0,t]$ and 
for $t=\infty$ the supremum is taken over $\tau\in\mathbb{R}_+$. 
The previous lemma yields $0<c^{\1}_{t,\rm max}<\infty$, $t\in [0,\infty]$. 
\medskip

Before formulating the global existence and uniqueness result of 
\cite{Lo20} let us recall the meaning of the parameter $\lambda$ 
in (\ref{3.3}) and (\ref{5.29}) below. In \cite{Lo20} we proceed as 
follows. First we verify existence and uniqueness of the solution 
to (\ref{3.3}) on some interval $t\in [0,T]$ where $T\cdot\lambda>0$ 
is constant. The major tool for this is Banach's fixed point theorem. 
Then we use this local result to prove existence and uniqueness of 
the solution to (\ref{3.3}) on $t\in [0,\infty)$. Here, we apply the 
theory of first order differential equations in Banach space and, as 
an important technical feature, Proposition \ref{Proposition3.3**} (a) 
for $t\in [0,\infty)$. The key to extend existence and uniqueness of 
the solution to (\ref{5.29}) on $t\in[\tau_0,\infty)$ for a certain 
$\tau_0\equiv\tau_0(\lambda)$ with $\lim_{\lambda\to 0}\tau_0(\lambda) 
=-\infty$ is Theorem \ref{Theorem3.19}.
\medskip

Recall the constant $p_{0,{\rm max}}$ from Proposition 
\ref{Proposition3.3**} (a) and the function $h_\gamma$ from hypothesis 
(wk). Let $b$ be the constant introduced in hypothesis (w), and let 
$c^{\1}_{\infty,\rm max}$ be the constant given by (\ref{3.181}). 

Introduce $a_1:=8\|h_\gamma\|\|B\|$, $-\tau_0\equiv -\tau_0(\lambda):= 
\max_{x\ge 1}\log(x)/(-m_1+a_1\lambda x)$ and note that this Definition 
implies the just mentioned $\lim_{\lambda\to 0}\tau_0(\lambda)=-\infty$. 
Let $K>1$ be the unique number that satisfies $-\tau_0=\log(K)/(-m_1+a_1 
\lambda K)$. The numbering (2')-(4') and (2'') is taken over from 
\cite{Lo20}. 
\begin{theorem}\label{Theorem5.8} (Theorem 5.8 in \cite{Lo20})
Suppose that the conditions of Theorem \ref{Theorem3.19} are satisfied.\\ 
(a) Let $p_0\in L^1(\Omega\times V)$ with $\|p_0\|_{L^1(\Omega\times V)}=1$ 
and suppose that there are constants $0<c\le C<\infty$ with $c\le p_0\le C$ 
a.e. on $\Omega\times V$. Then there is a unique map $[\tau_0,\infty)\ni 
t\mapsto p_t(p_0)\equiv p(\cdot,\cdot ,t)\in L^1(\Omega\times V)$ which is 
continuous with respect to the topology of $L^1(\Omega\times V)$ such that 
for every $\tau_0\le\tau\le 0$ 
\begin{eqnarray}\label{5.29} 
p(r,v,t)=S(t-\tau)\, p(r,v,\tau)+\lambda\int_\tau^t S(t-s)\, Q(p,p)\, (r,v,s) 
\, ds 
\end{eqnarray} 
a.e.~on $(r,v,t)\in\Omega\times V\times [\tau,\infty)$ and $p(\cdot ,\cdot , 
0)=p_0$. We have the following.
\begin{itemize} 
\item[(2')] We have $\|p(\cdot ,\cdot ,t)\|_{L^1(\Omega\times V)}\le K$, 
$t\ge\tau_0$ and $\|p(\cdot ,\cdot ,t)\|_{L^1(\Omega\times V)}=1$, $t\ge 0$. 
\item[(3')] If $p(\cdot,\cdot ,s)\ge c'$ a.e. on $\Omega\times V$ for some 
$c'>0$ and some $s\in [\tau_0,0]$ then there exists a strictly decreasing 
positive function $[s,\infty)\ni t\mapsto c_t\equiv c_t(c',s)$ such that 
$p(\cdot ,\cdot ,t)\ge c_t$ a.e. on $\Omega\times V$. 
\item[(4')] For $s$ introduced in (3') and $t\ge s$ we have 
\begin{eqnarray*}  
\|p(\cdot ,\cdot ,t)\|_{L^\infty(\Omega\times V)}\le p_{0,{\rm max}} 
\cdot\exp\left\{\lambda\|h_\gamma\|\, b\, c^{\1}_{\infty,\rm max}\cdot 
|t|\vphantom{l^1}\right\}\, .  
\end{eqnarray*}
\end{itemize} 
Conversely, for all $\tau_0<0$ there is a $\lambda>0$ such that 
(\ref{5.29}). If (\ref{5.29}) is modified by replacing $p$ with $p/(1 
\vee\|p/K\|_{L^1(\Omega\times V)})$, then the modified equation has a 
unique continuous solution $\mathbb{R}\ni t\mapsto \t p_t(p_0)\equiv 
\t p(\cdot,\cdot ,t)\in L^1(\Omega\times V)$ with $\t p_t(p_0)=p_t(p_0)$ 
for all $t\ge\tau_0$. The following holds.  
\begin{itemize} 
\item[(2'')] For all $t\in\mathbb{R}$ we have $\int_\Omega\int_V \t 
p(\cdot,\cdot ,t)\, dv\, dr =1$ and for all $t\le 0$ it holds that 
\begin{eqnarray*}
\left\|\t p_t(p_0)\right\|_{L^1(\Omega\times V)}\le\exp\left\{-t\cdot 
\left(-m_1+a_1\lambda K\vphantom{l^1}\right)\right\}\, . 
\end{eqnarray*}
\end{itemize} 
(b) Let $p_0\in D(A)$ with $\|p_0\|_{L^1(\Omega\times V)}=1$ and suppose 
that there are constants $0<c\le C<\infty$ with $c\le p_0\le C$ a.e. on 
$\Omega\times V$. Then there is a unique map $[\tau_0,\infty)\ni t\mapsto 
p_t(p_0)\equiv p(\cdot,\cdot ,t)\in D(A)$ which is continuous on $[\tau_0, 
\infty)$ and continuously differentiable on $(\tau_0,\infty)$ with respect 
to the topology of $L^1(\Omega\times V)$ such that 
\begin{eqnarray*}
\frac{d}{dt}\, p(\cdot ,\cdot ,t)=Ap(\cdot ,\cdot ,t)+\lambda\, Q(p,p)\, 
(\cdot ,\cdot ,t)\, ,\quad t\ge\tau_0, 
\end{eqnarray*} 
and $p(\cdot ,\cdot, 0)=p_0$. Here, $d/dt$ denotes the derivative in $L^1 
(\Omega\times V)$. At $t=\tau_0$ it is the right derivative. This map 
coincides with the solution to the equation (\ref{5.29}) if there $p(\cdot 
,\cdot, 0)=p_0\in D(A)$. We have properties (2')-(4') and (2'') of part (a).
\end{theorem} 

\subsection{Differentiability of the Solution Map to the Boltzmann Type 
Equation}\label{sec:3:2}

Suppose (k)-(wkk) throughout this subsection. By hypothesis (kk), $S(t)$, 
$t\ge 0$, given in (kkk) maps $L^1(\Omega\times V)$ linearly to $L^1(\Omega 
\times V)$ with operator norm one. We observe furthermore that by the 
definitions in (kw)-(wk), for fixed $t\ge 0$, $Q$ maps $(p(\cdot,\cdot,t), 
q(\cdot,\cdot,t))\in L^1(\Omega\times V)\times L^1(\Omega\times V)$ to 
$L^1(\Omega\times V)$ such 
that 
\begin{eqnarray}\label{3.5}
&&\hspace{-.5cm}\|Q(p(\cdot,\cdot,t),q(\cdot,\cdot,t))\|_{L^1(\Omega\times 
V)}\vphantom{\dot{f}}\nonumber \\ 
&&\hspace{.5cm}\le 2\|h_\gamma\|\|B\|\|p(\cdot,\cdot,t)\|_{L^1(\Omega\times 
V)}\|q(\cdot,\cdot,t)\|_{L^1(\Omega\times V)}
\end{eqnarray} 
where $\|\cdot\|$ denotes the sup-norm. 
\medskip 

Let $T>0$. Let $(L^1(\Omega\times V))^{[0,T]}$ be the space of all 
measurable real functions $f(r,v,t)$, $(r,v)\in\Omega\times V$, $t\in 
[0,T]$, such that $F_f(t):=f(\cdot,\cdot,t)\in L^1(\Omega\times V)$ for 
all $t\in [0,T]$ and $F_f\in C([0,T];L^1(\Omega\times V))$. With the norm 
\begin{eqnarray*}
\|f\|_{1,T}:=\sup_{t\in [0,T]}\|f(\cdot,\cdot,t)\|_{L^1(\Omega\times V)}
\end{eqnarray*} 
$(L^1(\Omega\times V))^{[0,T]}$ is a Banach space. Let $p_0\in L^1(\Omega 
\times V)$ and let $p,q\in (L^1(\Omega\times V))^{[0,T]}$. For $(r,v,t)\in 
\Omega\times V\times [0,T]$ set 
\begin{eqnarray*} 
\Psi(p_0,p)(r,v,t):=S(t)\, p_0(r,v)+\lambda\int_0^t S(t-s)\, Q(p,p)\, 
(r,v,s)\, ds\, .
\end{eqnarray*} 
where the integral converges in $L^1(\Omega\times V)$. We have $\Psi 
(p_0,p)\in (L^1(\Omega\times V))^{[0,T]}$ and according to (\ref{3.5})
\begin{eqnarray}\label{3.7} 
&&\hspace{-.5cm}\|\Psi(p_0,p)\|_{1,T}\le\|p_0\|_{L^1(\Omega\times V)}+ 
\lambda T\|Q(p,p)\|_{1,T}\nonumber \\ 
&&\hspace{.5cm}\le\|p_0\|_{L^1(\Omega\times V)}+2\lambda T\|h_\gamma\| 
\|B\|\|p\|_{1,T}^2\, . 
\end{eqnarray} 
Moreover, if $\lambda\le 1/(16 T\|h_\gamma\|\|B\|)$ and $\|p_0\|_{L^1 
(\Omega\times V)}<3/2$ then (\ref{3.7}) leads to 
\begin{eqnarray}\label{3.8}
\|p\|_{1,T}<2\quad {\rm implies}\quad \|\Psi (p_0,p)\|_{1,T}<2\, .  
\end{eqnarray} 
\begin{proposition}\label{Proposition3.2} 
(a) (Proposition 5.2 (a) in  \cite{Lo20}) Let $T>0$. Fix the parameter 
in the Boltzmann type equation $\lambda\le 1/(16 T\|h_\gamma\|\|B\|)$. 
Let $p_0\in L^1(\Omega\times V)$ with $\|p_0\|_{L^1(\Omega\times V)}< 
3/2$. \\ 
There is a unique element $p\equiv p(p_0)\in (L^1(\Omega\times V))^{[0, 
T]}$ such that  
\begin{eqnarray*}
p=\Psi\left(p_0,p\right)\, . 
\end{eqnarray*} 
The map $\{p_0'\in L^1(\Omega\times V):\|p_0'\|_{L^1(\Omega\times V)} 
<3/2\}\ni p_0\mapsto p(p_0)\in (L^1(\Omega\times V))^{[0,T]}$ is 
continuous. We have $\|p(p_0)\|_{1,T}<2$ if $\|p_0\|_{L^1(\Omega\times 
V)}<3/2$. \\ 
(b) Let $0<t\le T$ and choose $\lambda>0$ such that 
\begin{eqnarray*}
-\tau_1\equiv -\tau_1(\lambda):=\max_{x\ge\frac32}\left(\log(x)-\log 
\frac32\right)/(-m_1+a_1\lambda x)\ge T 
\end{eqnarray*} 
as well as $\lambda< 1/(4KT\|h_\gamma\|\|B\|)$ where $K\equiv K(\lambda) 
>\frac32$ is the unique number satisfying $-\tau_1=(\log(K)-\log 
\frac32)/(-m_1+a_1\lambda K)$; for the notation recall also Proposition 
\ref{Proposition3.3*} and Theorem \ref{Theorem5.8}. Moreover, let $q_t 
\in L^1(\Omega\times V)$ with $\|q_t\|_{L^1(\Omega\times V)}<\frac32$. 
There is a unique element $q\equiv q(q_t)\in (L^1(\Omega\times V))^{[0, 
t]}$ such that with $p_0:=q(\cdot,\cdot,0)$ it holds that 
\begin{eqnarray*}
q=\Psi|_{[0,t]}\left(p_0,q\right)\quad\mbox{\rm and}\quad q_t=q(\cdot, 
\cdot,t)  
\end{eqnarray*} 
where $\Psi|_{[0,t]}$ denotes the restriction of $\Psi$ to $\Omega\times 
V\times [0,t]$. The map $\{q_t'\in L^1(\Omega\times V):\|q_t'\|_{L^1 
(\Omega \times V)}<\frac32\}\ni q_t\mapsto q(q_t)\in (L^1(\Omega\times V) 
)^{[0,t]}$ is continuous. We have $\|q(q_t)\|_{1,t}<K$ if $\|q_t\|_{L^1 
(\Omega\times V)}<\frac32$.
\end{proposition} 
Proof. It remains to verify part (b). Relations (5.32) and (5.33) of the 
proof of Theorem 5.8 in \cite{Lo20} have been proved for $f\in L^1(\Omega 
\times V)$ while (5.34) and (5.35) hold in the stated form only for $f\in 
L^1(\Omega\times V)$ with $\|f\|_{L^1(\Omega\times V)}=1$. Recall that 
$K>0$ in (5.32) is non-specified. Let us modify (5.34) and (5.35) for $f 
\in L^1(\Omega\times V)$ with $\|f\|_{L^1(\Omega\times V)}\le\frac32$. In 
order to have for such $f$ and the solution $\t T(-t)f$ of (5.33) on $-t 
\in [-t_1,0]$ the relation $\|\t T(-t)f\|_{L^1(\Omega\times V)}\le K$, let 
us specify $K$ and $t_1$ as follows. For given $\lambda>0$, let $t_1=\max_{ 
x\ge\frac32}(\log(x)-\log\frac32)/(-m_1+a_1\lambda x)$ and let $K>\frac32$ 
be the unique number satisfying $t_1=(\log(K)-\log \frac32)/(-m_1+a_1\lambda 
K)$. As in the proof of Theorem 5.8 in \cite{Lo20} it follows now that 
equation (5.34) holds on $t\in [-t_1,0]$ for all $f\in L^1(\Omega\times V)$ 
with $\|f\|_{L^1(\Omega\times V)}\le\frac32$. In addition, equation (5.35) 
of \cite{Lo20} holds for all $u\in [0,-\tau_1]:=[0,t_1]$ and all $p_0\equiv 
p(\cdot,\cdot,t)$ with $\|p_0\|_{L^1(\Omega\times V)}\le\frac32$. Existence 
and uniqueness of $q\equiv q(q_t)\in (L^1(\Omega\times V))^{[0,t]}$ as stated 
follows whenever $\|q_t\|_{L^1(\Omega\times V)}\le\frac32$. 
%

The estimate in part (b) is now a consequence of the estimate (5.32) in \cite
{Lo20}. The continuity statement follows for $\lambda<1/(4KT\|h_\gamma\|\|B 
\|)$ from the following adjustment of (5.6) in \cite{Lo20} based on the just 
modified backward equation (5.35) of \cite{Lo20} and (\ref{4.52}) of the 
present paper, 
\begin{eqnarray*}
&&\hspace{-.5cm}\|q(q_{1,t})-q(q_{2,t})\|_{1,T}\nonumber \\ 
&&\hspace{.5cm}\le\max_{t\in [0,T]}\|S(-t)\, (q_{1,t}-q_{2,t})\|_{L^1 
(\Omega\times V)}+\lambda T\left\|Q\left(q(q_{1,t})+q(q_{2,t}),q(q_{1,t}) 
-q(q_{2,t})\vphantom{l^1}\right)\right\|_{1,T} \nonumber \\ 
&&\hspace{.5cm}\le Me^{-m_1T}\|q_{1,t}-q_{2,t}\|_{L^1(\Omega\times V)}+2 
\lambda T\|h_\gamma\|\|B\|\|q(q_{1,t})+q(q_{2,t})\|_{1,T}\|q(q_{1,t})-q 
(q_{2,t})\|_{1,T}\, ,   
\end{eqnarray*} 
where $q_{1,t},q_{2,t}\in\{q_t'\in L^1(\Omega\times V):\|q_t'\|_{L^1(\Omega 
\times V)}<\frac32\}$. 
\qed
\medskip

Part (b) of the preceding proposition will be applied in Corollary 
\ref{Corollary3.5} below. Part (a) is used in order to prepare part 
(c) of the following proposition. Let ${\cal N}$ denote the set of all 
non-negative $p_0\in L^1(\Omega\times V)$ with $\|1/p_0\|_{L^\infty 
(\Omega\times V)}<\infty$, $\|p_0\|_{L^\infty(\Omega\times V)}<\infty$, 
and $\|p_0\|_{L^1(\Omega\times V)}=1$. For all $p_0\in {\cal N}$, we may 
and do assume (\ref{3.11}) with $p_{0,{\rm min}}>0$ and $p_{0,{\rm max}} 
<\infty$. For $\ve>0$ define 
\begin{eqnarray*}
{\cal N}^\ve:=\left\{p_0\in L^1(\Omega\times V):\|p_0-p_0'\|_{L^1(\Omega 
\times V)}<\ve\ {\rm for\ some\ }p_0'\in {\cal N}\right\}\, . 
\end{eqnarray*} 
\begin{proposition}\label{Proposition3.3} 
Let $T>0$. \\ 
(a) Let $q\in (L^1(\Omega\times V))^{[0,T]}$. The map $L^1(\Omega\times V) 
\ni q_0\mapsto\Psi(q_0,q)\in (L^1(\Omega\times V))^{[0,T]}$ is Fr\'echet 
differentiable. The Fr\'echet derivative has the representation 
\begin{eqnarray*}
\left(\vphantom{l^1}\nabla_1\Psi(q_0,q)(\cdot)\right)(r,v,t)=(S(t)(\cdot )) 
(r,v)\, ,\quad (r,v)\in \Omega\times V,\ t\in [0,T]. 
\end{eqnarray*} 
(b) Let $p_0\in L^1(\Omega\times V)$. The map $(L^1(\Omega\times V))^{[0, 
T]}\ni q\mapsto\Psi (p_0,q)\in (L^1(\Omega\times V))^{[0,T]}$ is Fr\'echet 
differentiable. The Fr\'echet derivative has the representation 
\begin{eqnarray*}
\quad\left(\vphantom{l^1}\nabla_2\Psi(p_0,q)(\cdot)\right)(r,v,t)=2\lambda 
\int_0^tS(t-s)\, Q(\cdot,q)\, (r,v,s)\, ds\, ,
\end{eqnarray*} 
$(r,v)\in \Omega\times V$, $t\in [0,T]$. \\ 
(c) Let $0<\ve\le 1/2$ and, for the parameter in the Boltzmann type equation,
let $\lambda\le 1/(16 T\|h_\gamma\|\|B\|)$. The map ${\cal N}^{\ve}\ni p_0 
\mapsto p(p_0)\in (L^1(\Omega\times V))^{[0,T]}$ is continuously Fr\'echet 
differentiable. The Fr\'echet derivative has the representation 
\begin{eqnarray*}
\nabla p(p_0)(\cdot )=\left(\sum_{k=0}^\infty\left(\nabla_2\Psi\left(p_0,p( 
p_0)\vphantom{l^2}\right)\right)^k\circ\nabla_1\Psi\left(p_0,p(p_0)\vphantom 
{l^2}\right)\right)(\cdot ) 
\end{eqnarray*} 
where the exponent $k$ refers to $k$-fold composition. The right-hand side 
converges absolutely in the operator norm. 
\end{proposition}
Proof. (a) The derivative with respect to the first variable of $\Psi(\cdot, 
\cdot)$ in direction of $h_0\in L^1(\Omega\times V)$ is $\partial_{1,h_0} 
\Psi (q_0,q)(r,v,t)=(S(t)h_0)(r,v)$, $(r,v)\in\Omega\times V$, $t\in [0,T]$, 
independent of $q_0\in {\Bbb N}$. Thus, $\Psi(q_0+h_0,q)-\Psi(q_0,q)- 
\partial_{1,h_0}\Psi(q_0,q)=0$. 
\medskip 

\nid 
(b) We obtain the following derivative with respect to the second variable of 
$\Psi(\cdot,\cdot)$ in direction of $h\in (L^1(\Omega\times V))^{[0,T]}$, 
\begin{eqnarray*}
\partial_{2,h}\Psi(p_0,q)(r,v,t)=2\lambda\int_0^t S(t-s)\, Q(h,q)(r,v,s)\, ds 
\, , 
\end{eqnarray*} 
$(r,v)\in\Omega\times V$, $t\in [0,T]$. Therefore, 
\begin{eqnarray*}
&&\hspace{-.5cm}\frac{\left\|\Psi(p_0,q+h)-\Psi(p_0,q)-\partial_{2,h}\Psi(p_0, 
q)\right\|_{1,T}}{\|h\|_{1,T}} \\ 
&&\hspace{.5cm}=\lambda\left\|\int_0^t S(t-s)\left(Q(h,h)(\cdot,\cdot,s)/\|h 
\|_{1,T}\vphantom{l^1}\right)\, ds\right\|_{1,T}  
\end{eqnarray*} 
which converges to zero as $\|h\|_{1,T}\to 0$ by (\ref{3.5}). 
\medskip 

\nid 
(c) Let us regard for a moment $p_0$ as a fixed element of ${\cal N}^\ve$. 
Note that the choice $\ve\le 1/2$ corresponds to (\ref{3.8}) and thus to the 
hypotheses of Proposition \ref{Proposition3.2}. Furthermore, let $p\equiv 
p(p_0)$ be as above. We are interested in linear continuous solutions $f:L^1 
(\Omega\times V)\mapsto (L^1(\Omega\times V))^{[0,T]}$ to 
\begin{eqnarray}\label{3.13} 
f=\nabla_1\Psi(p_0,p)+\nabla_2\Psi(p_0,p)\circ f\, .  
\end{eqnarray} 
By (\ref{3.5}) and (\ref{3.8}), the particular choice of $\lambda$, and 
parts (a) and (b) of the present proposition, the following expression is a 
well-defined linear continuous solution $f$ to (\ref{3.13}), 
\begin{eqnarray}\label{3.14}
f(\cdot )=\left(\sum_{k=0}^\infty\left(\nabla_2\Psi(p_0,p)\right)^k\circ 
\nabla_1\Psi (p_0,p)\right)(\cdot )\, . 
\end{eqnarray} 
In particular, the right-hand side converges absolutely in the operator 
norm by the choice of $\lambda$ and Proposition \ref{Proposition3.2} (a). 
Uniqueness of the solution to (\ref{3.13}) is obtained by assuming the 
existence of another solution $\t f$ to (\ref{3.13}). Then $g:=f-\t f$ is a 
non-trivial solution to $g=\nabla_2\Psi(p_0,p)\circ g$. In particular, there 
is a test function $\vp\in L^1(\Omega\times V)$ with $g(\vp)\neq 0$ in $(L^1 
(\Omega\times V))^{[0,T]}$. From part (b) we obtain 
\begin{eqnarray*}
g(\vp)(r,v,s)=2\lambda\int_0^t S(t-s)\, Q(g(\vp),p)\, (r,v,s)\, ds\, ,
\end{eqnarray*} 
$(r,v)\in \Omega\times V$, $t\in [0,T]$. But taking into consideration the 
special choice of $\lambda$ this contradicts (\ref{3.5}) combined with 
(\ref{3.8}). 
\medskip 

Since (\ref{3.13}) is an equation for the G\^ateaux derivative 
\begin{eqnarray*}
\nabla^{G}p(p_0)=\nabla^{G}\Psi(p_0,p(p_0))=\nabla_1\Psi(p_0,p)+\nabla_2 
\Psi(p_0,p)\circ\nabla^{G}p(p_0)
\end{eqnarray*}
we conclude that $p(p_0)$ is G\^ateaux differentiable 
on $p_0\in {\cal N}^\ve$ with representation of the derivative $\nabla^{G}p 
(p_0)$ on the right-hand side of (\ref{3.14}). For the particular chain rule, 
cf. for example \cite{Ze86}, Proposition 4.10. Continuity of the right-hand 
side of (\ref{3.14}) with respect to $p_0\in {\cal N}^\ve$ in the operator 
norm follows from parts (a) and (b) of the present proposition, Proposition 
\ref{Proposition3.2} (a), and (\ref{3.5}). This implies the statement. 
\qed 

\subsection{Conditions (jv) and (v) via Fr\'echet Differentiability} 
\label{sec:3:3} 

In this subsection we prepare the proof of quasi-invariance under the 
flow generated by the solution map of the Boltzmann type equation 
constructed in Subsections \ref{sec:3:1} and \ref{sec:3:2}. For this 
we shall suppose (j)-(jjj) of Section \ref{sec:2} and (k)-(wkk) of the 
present section. 

Looking at Theorem \ref{Theorem2.3} and Theorem \ref{Theorem2.7} we have 
to show that (jv) and (v) of Section \ref{sec:2} are satisfied. We will 
use Propositions \ref{Proposition2.4} and \ref{Proposition3.3} for this. 
In order to verify the conditions of Theorem \ref{Theorem2.3} we need to 
arrange compatibility with the general framework of Subsection \ref{sec:1:1}. 
For this, we put $D:=\Omega\times V$ and choose the measure $\lambda\equiv
\l$ to be the Lebesgue measure on $(\Omega\times V,{\cal B}(\Omega\times 
V))$. Formally we set $v:=1$, $r:=3/2$, and ${\cal V}:=B^1(2)$. Choose 
$S\ge 1$ and $T=S+1$. For $p_0\in {\cal V}$ let  
\begin{eqnarray*}
V^+_tp_0=\t V^+_tp_0:=p(p_0)(\cdot,\cdot,t)\equiv p_t(p_0)\, ,\quad t\in 
[0,T]. 
\end{eqnarray*} 
By the definitions of Subsection \ref{sec:1:1} the flow $W^+$ and trajectory 
$U^+$ can now be introduced. Moreover, the generator $A^f$ of $U^+$ (and 
$W^+$) can directly be derived from (\ref{3.1}) with boundary conditions 
(\ref{3.2}). This setup is particularly used in order to apply Propositions 
\ref{Proposition3.3} and implicitly also Proposition \ref{Proposition3.2} 
to prove Corollary \ref{Corollary3.5} below. 

Therefore, we shall suppose in Corollary \ref{Corollary3.5} that the parameter 
of the Boltzmann type equation $\lambda$ satisfies the conditions of Proposition 
\ref{Proposition3.2} (a) and (b) which include the condition on $\lambda$ of 
Proposition \ref{Proposition3.3}. However, in the following two lemmas, Lemma 
\ref{Lemma3.7} and Lemma \ref{Lemma3.4}, the parameter of the Boltzmann type 
equation $\lambda$ can be chosen arbitrarily.  
\medskip

In addition, we will use Theorem \ref{Theorem5.8} in order to verify the 
general conditions (a1)-(a3) of Subsection \ref{sec:1:1}. Depending on 
the framework we refer to, we shall also use the notation  
\begin{eqnarray*}
p_t(p_0)\cdot\l\equiv p(\cdot,\cdot,t)\cdot\l\equiv h_t\cdot\l\equiv\nu_t 
\, ,\quad t\in\mathbb{R}, 
\end{eqnarray*} 
if $\nu\equiv h\cdot\l\equiv h_0\cdot\l\equiv p_0\cdot\l\in G$. As above, 
let $(A,D(A))$ denote the infinitesimal generator of the strongly continuous 
group $S(t)$, $t\in\mathbb{R}$, in $L^1(\Omega\times V)$. For $0<m<M<\infty$ 
introduce 
\begin{eqnarray*}
{\cal N}_{m,M}:=\{p_0\in {\cal N}:m\le p_0\le M\ \mbox{\rm a.e.}\} 
\end{eqnarray*}
and 
\begin{eqnarray*}
{\cal N}_{m,M}^D:=\{p_0\in {\cal N}:p_0\in D(A),\ m\le p_0\le M,\ |Ap_0|\le 
M\ \mbox{\rm a.e.}\}\, .
\end{eqnarray*}
\begin{lemma}\label{Lemma3.7} 
(a) For $0<m<M<\infty$, the set $\hat{\cal N}_{m,M}:=\{p_0\cdot\l :p_0\in 
\mathcal{N}_{m,M}\}$ is a compact subset of $(E,\pi)$. In particular, $\hat 
{\cal N}_{m,M}\in\mathcal{B}(E)$. \\ 
(b) For $0<m<M<\infty$, the set $\hat{\cal N}_{m,M}^D:=\{p_0\cdot\l :p_0 
\in\mathcal{N}_{m,M}^D\}$ is a compact subset of $(E,\pi)$. In particular, 
$\hat{\cal N}_{m,M}^D\in\mathcal{B}(E)$.  
\end{lemma}
Proof. (a) We pick up the idea from the proof of Lemma \ref{Lemma1.2} to show 
that $\hat{\cal N}_{m,M}$ is compact in $(E,\pi)$. Let $h_n\cdot\l\in\hat 
{\mathcal{N}}_{m,M}$, $n\in {\Bbb N}$, be an arbitrary sequence. We shall 
demonstrate that it has a subsequence converging relative to the metric $\pi$ 
with a limit belonging to $\hat{\cal N}_{m,M}$.

The sequence $h_n$, $n\in {\Bbb N}$, is uniformly integrable by the uniform 
upper bound $M$. Furthermore, $h_n$, $n\in {\Bbb N}$, converges weakly in 
$L^1(\Omega\times V)$ on some subsequence $(n_k)$ of indices by the just 
mentioned uniform integrability, $\|h_n\|_{L^1(\Omega\times V)}=1$, $n\in 
{\Bbb N}$, and the Dunford-Pettis theorem. In other words, there exists $h\in 
L^1(\Omega\times V)$ such that $h_{n_k}\stack{k\to\infty}{\lra}h$ weakly in 
$L^1(\Omega\times V)$. 

The bounds $m\le h$ and $h\le M$ on the limiting element can be verified by 
checking $0\le\int\vp (h_{n_k}-m)\, d\, \l\stack {k\to\infty}{\lra}\int\vp(h-m) 
\, d\, \l$ and $0\le\int\vp (M-h_{n_k})\, d\, \l\stack {k\to\infty}{\lra}\int 
\vp(M-h)\, d\, \l$ for all $\vp\in L^\infty(\Omega\times V)$ with $\vp\ge 0$. 
The relation $\|h\|_{L^1(\Omega\times V)}=1$ follows from $1=\int h_{n_k}\, d 
\, \l\stack{k\to\infty}{\lra}\int h\, d\, \l$. Thus $h\cdot\l\in\hat{\cal N}_{ 
m,M}$. We have in particular shown that $h_{n_k}\cdot\l\stack{k\to\infty} 
{\lra}h\cdot\l$ narrowly. The Portmanteau theorem in the form of [12], Theorem 
3.3.1, says that this is equivalent to the convergence relative to the metric 
$\pi$. Thus $\hat{\cal N}_{m,M}$ is compact. 
\medskip

\nid
(b) Let $h_n\cdot\l\in\hat{\cal N}_{m,M}^D$, $n\in {\Bbb N}$, be an arbitrary 
sequence. Since $\Omega\times V$ is bounded we have $\|Ah_n\|_{L^1(\Omega 
\times V)}\le\l(\Omega\times V)\cdot M$, $n\in {\Bbb N}$. Recall also $|Ah_n 
|\le M$ a.e., $n\in {\Bbb N}$. 

Reviewing the proof of part (a), the only thing we have to show is that $h\in 
D(A)$. As we shall see, this is a consequence of part (a) of this lemma and 
the closedness of $(A,D(A))$. Indeed by part (a), $h_n$, $n\in {\Bbb N}$, and 
$Ah_n$, $n\in {\Bbb N}$, converge both weakly in $L^1(\Omega\times V)$ on some 
subsequence $(n_k)$ of indices to $h\in L^1(\Omega\times V)$ and, respectively, 
to some $g\in L^1(\Omega\times V)$. The closedness of $(A,D(A))$ in $L^1(\Omega 
\times V)$ implies now $h\in D(A)$ (and $Ah=g$), see for example \cite{EN06}, 
Definition A.5 (a) (d), or \cite{Ka95}, Problem 5.12 of Chapter 2. 
\qed 
\medskip

Set
\begin{eqnarray*}
\hat{\cal N}:=\textstyle{\bigcup}_{N\in\mathbb{N}}\hat{\cal N}_{1/N,N}=
\{p_0\cdot\l\in E:m\le p_0\le M\ \mbox{\rm a.e. for some } 0<m<M<\infty\}  
\end{eqnarray*} 
and 
\begin{eqnarray*}
&&\hspace{-.5cm}\hat{\cal N}^D:=\textstyle{\bigcup}_{N\in\mathbb{N}}\hat
{\cal N}_{1/N,N}^D=\{p_0\cdot\l\in E:p_0\in D(A),\ m\le p_0\le M, \\ 
&&\hspace{.5cm}\ |Ap_0|\le M\ \mbox{\rm a.e. for some } 0<m<M<\infty\}\, . 
\vphantom{\sum}
\end{eqnarray*} 
Recall that for the $\sigma$-algebra ${\cal B}(\hat{\cal N})$ on $\hat 
{\cal N}$ generated by the open subsets of $\hat{\cal N}$ with respect to 
$\pi$ we have ${\cal B}(\hat{\cal N})=\{B\cap\hat{\cal N}:B\in {\cal B}(E) 
\}$. Choose $t_0>0$. This $t_0$ will be used below in order to verify 
hypothesis (a2) of Subsection \ref{sec:1:1}. It appears only in relation 
to the set $G$ which we define in the next lemma. 
\begin{lemma}\label{Lemma3.4} 
Suppose that the hypotheses of Theorem \ref{Theorem3.19} are satisfied. 
 \\ 
(a) For all $t\ge 0$, the map $p_0\cdot\l\mapsto p_t(p_0)\cdot\l$ is a 
measurable injection $\hat{\mathcal{N}}\mapsto\hat{\mathcal{N}}$ with 
measurable inverse. \\ 
(b) With 
\begin{eqnarray*}
G:=\left\{p_{1+t_0}(h)\cdot\l:h\cdot\l\in\hat{\mathcal{N}}^D\right\}
\end{eqnarray*}
we have $G\in {\cal B}(E)$ and (a1)-(a3) of Subsection \ref{sec:1:1} and 
(a4) of Lemma \ref{LemmaA.1}. \\ 
(c) We have Case 2.  
\end{lemma} 
Proof. {\it Step 1 } In Steps 1 and 2 we verify (a). Let $0<m<M<\infty$. 
As a consequence of Theorem \ref{Theorem5.8} (a), $\hat{\cal N}_{m,M}\ni 
p_0\cdot\l\mapsto p_t(p_0)\cdot\l$ is injective, note in particular 
properties (3') and (4') and the paragraph below (4'). First we show that 
on a certain interval $t\in [0,\t T]$ the map $\hat{\cal N}_{m,M}\ni p_0 
\cdot\l\mapsto p_t(p_0)\cdot\l$ is continuous relative to the metric $\pi$. 
It is then a consequence of \cite{Bo95}, Corollary 2 in \S~9.4 of Chapter 
I, that $\{p_t(p_0)\cdot\l:p_0\in B\}\in {\cal B}(E)$ for all $B\in {\cal 
B}(\hat{\mathcal N}_{m,M})$. 

Let $h_n\cdot\l\in\hat{\cal N}_{m,M}$, $n\in {\Bbb N}$, be converge to $h 
\cdot\l\in\hat{\cal N}_{m,M}$ with respect to $\pi$. The map $p_t$ maps 
$\hat{\cal N}_{m,M}$ into some compact space $\hat{\cal N}_{m',M'}$, cf. 
Theorem \ref{Theorem5.8} (a) and Lemma \ref{Lemma3.7} (a). Therefore, it 
is reasonable to proceed as follows. First we show that any two subsequences 
of $p_t(h_n)$, $n\in {\Bbb N}$, converging with respect to $\pi$, have the 
same limit which we denote by $\hat{p}_t$, $t\in [0,\t T]$. Then we verify 
$\hat{p}_t=p_t(h)$, $t\in [0,\t T]$. The present Step 1 is devoted to 
preparing the just mentioned procedure. 
\medskip

Until the end of Step 2, choose 
\begin{eqnarray*}
0<\t T\le 1/(6\lambda\|h_\gamma\|\|B\|) 
\end{eqnarray*} 
where $\lambda$ is the parameter in the Boltzmann type equation. Let $t\in 
[0,\t T]$. For the next calculation we mention that for the map 
\begin{eqnarray*}
{\Bbb R}^{2d}\ni (v,v_1)\mapsto (v^\ast,v_1^\ast):=\left(v-e\circ (v-v_1) 
\, e\, ,\, v_1+e\circ (v-v_1)\, e\vphantom{l^1}\right) 
\end{eqnarray*} 
the absolute value of the Jacobian determinant is one. It follows now 
from (kw) that for $p,q\in (L^1(\Omega\times V))^{[0,\t T]}$ with 
$\|p\|_{1,\t T}=\|q\|_{1,\t T}=1$ and $g\in L^\infty(\Omega\times V)$ 
with $\|g\|_{L^\infty(\Omega\times V)}=1$ we have 
\begin{eqnarray}\label{3.15} 
&&\hspace{-.5cm}\left|\int_\Omega\int_V g(r,v)\lambda\int_0^t S(t-s)Q(p-q, 
p+q)(r,v,s)\, ds\, dv\, dr\right|\nonumber \\ 
&&\hspace{.5cm}\le\lambda\int_\Omega\int_V\int_0^t\left|S(t-s)Q(p-q,p+q) 
(r,v,s)\right|\, ds\, dv\, dr\nonumber \\ 
&&\hspace{.5cm}\le\lambda\int_0^t\left\|S(t-s)Q(p-q,p+q)(\cdot,\cdot,s) 
\right\|_{L^1(\Omega\times V)}\, ds\nonumber \\ 
&&\hspace{.5cm}\le\lambda\int_0^t\left\|Q(p-q,p+q)(\cdot,\cdot,s)\right 
\|_{L^1(\Omega\times V)}\, ds\nonumber \\ 
&&\hspace{.5cm}\le\frac{\lambda}{2}\|B\|\|h_\gamma\|\cdot\int_0^t 
\int_\Omega\int_V\int_{S_+^{d-1}}\int_{\Omega}\int_V\chi_{\{(v^\ast , 
v_1^\ast)\in V\times V\}}\times\nonumber \\ 
&&\hspace{1.0cm}\times\left(|p-q|(r,v^\ast ,s)\cdot|p+q|(y,v_1^\ast ,s) 
\vphantom{\dot{f}}\right.\vphantom{\int}\nonumber \\ 
&&\hspace{1.5cm}+|p-q|(r,v,s)\cdot|p+q|(y,v_1,s)\vphantom{l^1}\vphantom 
{\int}\nonumber \\ 
&&\hspace{1.5cm}+|p+q|q(r,v^\ast ,s)\cdot|p-q|(y,v_1^\ast ,s)\vphantom 
{\int}\nonumber \\ 
&&\hspace{1.5cm}\left.+|p+q|(r,v,s)\cdot|p-q|(y,v_1,s)\vphantom{\dot{f}} 
\right)\, dv_1\, dy\, de\, dv\, dr\, ds\vphantom{\int}\nonumber \\
&&\hspace{.5cm}\le 4\lambda\|B\|\|h_\gamma\|\cdot\int_0^t\int_\Omega\int_V 
|p-q|(r,v,s)\, dv\, dr\, ds\, .   
\end{eqnarray} 
We may now conclude 
\begin{eqnarray}\label{3.151} 
&&\hspace{-.5cm}\left|\int_\Omega\int_V g(r,v)\lambda\int_0^t S(t-s)Q(p-q, 
p+q)(r,v,s)\, ds\, dv\, dr\right|\nonumber \\ 
&&\hspace{.5cm}\le 4\lambda \t T\|B\|\|h_\gamma\|\cdot\|p-q\|_{1,\t T}\, . 
\vphantom{\int}
\end{eqnarray} 

\noindent
{\it Step 2 } We recall the idea from the beginning of Step 1. In particular 
$\|h_n\|_{L^1(\Omega\times V)}=1$ and $h_n\le M$ on $\Omega\times V$, $n\in 
\mathbb{N}$, yield uniform integrability of $h_n$, $n\in\mathbb{N}$. The 
Dunford-Pettis theorem says now that any subsequence $h_{n_k}$, $k\in\mathbb 
{N}$, contains another subsequence converging weakly in $L^1(\Omega\times V)$. 

By the narrow convergence of $h_n\cdot\l$, $n\in\mathbb{N}$, to $h\cdot\l$ 
according to the supposed convergence with respect to $\pi$, this even 
yields 
\begin{eqnarray*}
\int_\Omega\int_V gh_n\, dv\, dr\stack{n\to\infty}{\lra}\int_\Omega\int_V 
gh\, dv\, dr\quad\mbox{\rm for all}\quad g\in L^\infty(\Omega\times V)\, . 
\end{eqnarray*} 
The dual space of $L^1(\Omega\times V)$ is isomorphic to $L^\infty(\Omega 
\times V)$. As a consequence of Theorem \ref{Theorem3.19} and \cite{Pa83} 
there is an adjoint group $S^\ast (t)$, $t\in\mathbb{R}$, on $L^\infty( 
\Omega\times V)$. For $g\in L^\infty(\Omega\times V)$ we obtain
\begin{eqnarray}\label{3.152}
&&\hspace{-.5cm}\lim_{n\to\infty}\int_\Omega\int_V g\, S(t)h_n\, dv\, dr 
=\lim_{n\to\infty}\int_\Omega\int_V h_n\, S^\ast(t)g\, dv\, dr\nonumber \\ 
&&\hspace{.5cm}=\int_\Omega\int_V h\, S^\ast(t)g\, dv\, dr =\int_\Omega 
\int_V g\, S(t)h\, dv\, dr\, ,\quad t\in\mathbb{R}.              
\end{eqnarray} 

Simultaneously for all $t\in [0,\t T]$ let us consider two convergent 
subsequences with respect to $\pi$, $p_t(h_{n_k})$, $k\in {\Bbb N}$, 
and $p_t(h_{n_l})$, $l\in {\Bbb N}$, with limits $p_t^1$ and $p_t^2$, 
respectively. We aim to show that $p_t^1=p_t^2$, $t\in [0,\t T]$. For 
this we estimate $\|p^1-p^2\|_{1,\t T}$ using 
\begin{eqnarray}\label{3.19**}
\|p^1-p^2\|_{1,\t T}=\sup_{ t\in [0,\t T]}\sup\left\{\int_\Omega\int_V 
g(p_t^1-p_t^2)\, dr\, dv: \|g\|_{L^\infty (\Omega\times V)}=1\right\}\, . 
\end{eqnarray} 
Replacing in (\ref{3.152}) $h_n$ first by $h_{n_k}$ and then by $h_{n_l}$ 
and looking at the difference, the right-hand side is zero. Furthermore, 
we obtain from plugging first (\ref{3.3}) and then (\ref{3.151}) into the 
right-hand side of (\ref{3.19**})
\begin{eqnarray*}
\|p^1-p^2\|_{1,\t T}\le 4\lambda \t T\|B\|\|h_\gamma\|\cdot\|p^1-p^2\|_{1,
\t T}\le\frac23\cdot\|p^1-p^2\|_{1,\t T}
\end{eqnarray*} 
where, for the second $\le$, we recall the choice of $\t T$ in the beginning 
of Step 1. In other words $\|p^1-p^2\|_{1,\t T}=0$. As discussed in Step 1, 
this shows the unique existence of $\hat{p}_t\in \hat{\cal N}_{m',M'}$ with 
$\lim_{n\to\infty}\pi(p_t(h_n),\hat{p}_t)=0$, $t\in [0,\t T]$, i. e. $p_t(h_n) 
\stack{n\to\infty}{\lra}\hat{p}_t$ in the sense of narrow convergence, cf. 
Subsection \ref{sec:1:1}. 
\medskip

Next we show $\hat{p}_t=p_t(h)$, $t\in [0,\t T]$. For this we recall $B(v^\ast, 
v_1^\ast,e)=B(v,v_1,e)$ from hypothesis (w) and that the absolute value of 
the Jacobian of the map $(v,v_1)\mapsto (v^\ast,v_1^\ast)$ is one. This map 
has for fixed $e\in S^{d-1}_+$ an inverse which is given by $(v^{-\ast}, 
v_1^{-\ast}):=\left(v^\ast+\frac12 e\circ (v_1^\ast-v^\ast)\, e\, ,\, v_1^\ast 
+\frac12 e\circ (v^\ast-v_1^\ast)\, e\right)$. 

By the reasoning of the first two paragraphs of Step 2, $p_t(h_n)\stack{n\to 
\infty} {\lra}\hat{p}_t$ holds even weakly in $L^1(\Omega\times V)$ for all $t\in 
[0,\t T]$. As a consequence $p_t(h_n)\cdot p_t(h_n)$ converges weakly in $L^1$ 
to $\hat{p}_t\cdot\hat{p}_t$ on the product space $((\Omega\times V)^2,{\cal B} 
((\Omega\times V)^2))$ as $n\lra\infty$, $t\in [0,\t T]$. Thus, for $g\in 
L^\infty(\Omega\times V)$ and $g_{t-s}:=S^\ast(t-s)g$, $0\le s\le t\le\t T$, 
\begin{eqnarray*}
&&\hspace{-.5cm}\int_\Omega\int_V g(r,v)\lambda\int_0^t S(t-s)Q(p(h_n),p 
(h_n))(r,v,s)\, ds\, dv\, dr \\ 
&&\hspace{.5cm}=\int_0^t\lambda\int_\Omega\int_V g_{t-s}(r,v)Q(p(h_n),p(h_n)) 
(r,v,s)\, dv\, dr\, ds \\ 
&&\hspace{0.5cm}=\int_0^t\lambda\int_\Omega\int_V\int_{\Omega}\int_V\left( 
\int_{S_+^{d-1}}g_{t-s}(r,v^{-\ast})B(v,v_1,e)h_\gamma(r,y)\, 
de\right)\times \\ 
&&\hspace{1.0cm}\times\, p(h_n)(y,v_1,s)p(h_n)(r,v,s)\, dv_1\, dy\, dv\, dr\ 
ds\vphantom{\int_0^t} \\ 
&&\hspace{1.0cm}-\int_0^t\lambda\int_\Omega\int_V\int_{\Omega}\int_V\left(
 g_{t-s}(r,v)\cdot\int_{S_+^{d-1}}B(v,v_1,e)h_\gamma 
(r,y)\chi_{\{(v^\ast ,v_1^\ast)\in V\times V\}}\, de\right)\times \\ 
&&\hspace{1.0cm}\times\, p(h_n)(y,v_1,s)p(h_n)(r,v,s)\, dv_1\, dy\, dv\, dr\ 
ds\vphantom{\int_0^t} \\ 
&&\hspace{.2cm}\stack{n\to\infty}{\lra}\int_0^t\lambda\int_\Omega\int_V g_{t-s} 
(r,v)Q(\hat{p}_s,\hat{p}_s)(r,v)\, dv\, dr\, ds \\ 
&&\hspace{0.5cm}=\int_\Omega\int_V g(r,v)\lambda\int_0^t S(t-s)Q(\hat{p}_s, 
\hat{p}_s)(r,v)\, ds\, dv\, dr 
\end{eqnarray*} 
where the interchange of $\lim_{n\to\infty}$ and $\int_0^t$ is due to dominated 
convergence. Together with (\ref{3.152}) and $\hat{p}_0=h=p_0(h)$ this says 
$\hat{p}_t=p_t(h)$, $t\in [0,\t T]$. In other words, we have reached the goal 
formulated in the first two paragraphs of Step 1, i. e. we have proved that the map 
\begin{eqnarray}\label{3.20**}
\hat{\cal N}_{m,M}\ni p_0\cdot\l\mapsto p_t(p_0)\cdot\l\quad\mbox{\rm is 
continuous for } t\in [0,\t T].
\end{eqnarray} 
We have verified part (a) for $p_0\cdot\l\in\hat{\mathcal{N}}_{m,M}$ and 
$t\in [0,\t T]$. Recalling that $\hat{\cal N}_{m,M}$ is compact with respect 
to $\pi$ by Lemma \ref{Lemma3.7} (a) and taking into consideration \cite{Bo95}, 
Corollary 2 in \S~9.4 of Chapter I, the inverse of (\ref{3.20**}) 
\begin{eqnarray}\label{3.21**}
&&\hspace{-.5cm}\left\{p_t(p_0)\cdot\l:p_0\cdot\l\in\hat{\cal N}_{m,M} 
\right\}\ni\hat{p}_0\cdot\l\nonumber \\ 
&&\hspace{.5cm}\mapsto p_{-t}(\hat{p}_0)\cdot\l \quad\mbox{\rm is also 
continuous for } t\in [0,\t T].
\end{eqnarray} 

Now recall the notation introduced after Lemma \ref{Lemma3.7}. It 
follows from Theorem \ref{Theorem5.8} (a), properties (3') and (4') and 
the paragraph below (4'), that 
for all $t\in [0,\t T]$ the map $p_0\cdot\l\mapsto p_t(p_0)\cdot\l$ 
is a measurable injection $\hat{\mathcal{N}}\mapsto\hat{\mathcal{N}}$ with 
measurable inverse, here measurability with respect to $\mathcal{B}(\hat 
{\mathcal{N}})$. By iteration we obtain that $p_0\cdot\l\mapsto p_t(p_0) 
\cdot\l$ is a measurable injection $\hat{\mathcal{N}}\mapsto\hat{\mathcal 
{N}}$ with measurable inverse for all $t\ge 0$. We have completed the proof 
of part (a). 
\medskip 

\nid 
{\it Step 3 } In this step we are concerned with part (b). As a byproduct 
we shall verify part (c) of the lemma. We observe that $G=\{p_{t_0+1}(h) 
\cdot\l:h\cdot\l\in\hat{\mathcal{N}}^D\}\in {\cal B}(E)$ follows from Lemma 
\ref{Lemma3.7} and the just proven part (a) of this lemma. Furthermore, 
condition (a1) is also a consequence of part (a) of the present lemma  
and of the continuity of the map $t\mapsto p_t(p_0)$ mentioned in Theorem 
\ref{Theorem5.8} (a) together with the definition of $G$. Condition (a4) 
is trivial under the choice 
\begin{eqnarray}\label{3.23**}
X=L^1(\Omega\times V)\quad\mbox{\rm and }\quad E\ni\nu=k\cdot\l\mapsto k 
\cdot\l=\mu\, ,  
\end{eqnarray} 
cf. Lemma \ref{LemmaA.1}. Regarding condition (a3), let us keep in mind 
that $\nu=h\cdot\l\in G$ satisfies $h\in D(A)$ and that $\bmu$ is a 
probability measure on $(E,{\cal B}(E))$ such that the closure of $G$ in 
$(E,\pi)$ is supp $\bmu$ and we have $\bmu(G)=1$. Therefore Theorem 
\ref{Theorem5.8} (b) and the definition of $G$ imply condition (a3). 
\medskip

In order to verify condition (a2), let us apply properties (3') and (4') 
of Theorem \ref{Theorem5.8} (a). For $N_1\in (0,1)$, it follows that there 
is $0<\t N_1\le N_1$ such that 
\begin{eqnarray}\label{3.20*}
\nu_t\in\hat{\mathcal{N}}_{\t N_1,1/\t N_1}\ \mbox{\rm whenever }\nu_{ 
-(1+t_0)}\in\hat{\mathcal{N}}_{N_1,1/N_1}\ \mbox{\rm and }t\in [-(1+t_0) 
,1+t_0]. 
\end{eqnarray} 

Now let $\nu_{-(1+t_0)}\in\hat {\mathcal{N}}_{N_1,1/N_1}^D\subset\hat 
{\mathcal{N}}_{N_1,1/N_1}$. For $\nu_t=h_t\cdot\l$, $t\in [-(1+t_0),1 
+t_0]$, we have $h_t\in D(A)$ by Theorem \ref{Theorem5.8} (b). 
Furthermore, the map $[-(1+t_0),1+t_0]\ni t\mapsto Q(h_t,h_t)$ is 
continuous in $L^1(\Omega\times V)$ by the continuity of $t\mapsto 
h_t$ according to Theorem \ref{Theorem5.8} (a), relation (\ref{3.20*}), 
and calculations similar to those for the last two $``\le"$ signs in 
(\ref{3.15}). Thus for any $-(1+t_0)\le\tau\le 0$ and $\tau\le t$
\begin{eqnarray*}
&&\hspace{-.5cm}Ah_t(r,v)=AS(t-\tau)\, h_\tau(r,v)+\lambda\lim_{\ve\downarrow 
0}\frac{S(\ve)-id}{\ve}\int_0^{t-\tau} S(s)\, Q(h_{t-s},h_{t-s})\, (r,v)\, ds
\nonumber \\ 
&&\hspace{.5cm}=AS(t-\tau)\, h_\tau(r,v)+\lambda \lim_{\ve\downarrow 0}\frac{1} 
{\ve}\left(\int_{\ve}^{t-\tau+\ve} S(s)\, Q(h_{t-s},h_{t-s})\, (r,v)\, ds 
\right.\nonumber \\ 
&&\hspace{1.0cm}\left.-\int_0^{t-\tau}S(s)\, Q(h_{t-s},h_{t-s})\, (r,v)\, ds 
\right)\nonumber \\ 
&&\hspace{1.0cm}+\lambda \lim_{\ve\downarrow 0}\frac{1}{\ve}\left(\int_{\ve}^{ 
t-\tau+\ve} S(s)\, Q(h_{t-s+\ve},h_{t-s+\ve})\, (r,v)\, ds\right.\nonumber \\ 
&&\hspace{1.0cm}\left.-\int_\ve^{t-\tau+\ve}S(s)\, Q(h_{t-s},h_{t-s})\, (r,v) 
\, ds\right)\nonumber \\ 
&&\hspace{.5cm}=S(t-\tau)\, Ah_\tau(r,v)+\lambda S(t-\tau)Q(h_\tau,h_\tau)\, 
(r,v)-\lambda Q(h_t,h_t)\, (r,v)\vphantom{\int_0^\ve}\nonumber \\ 
&&\hspace{1.0cm}+\int_0^{t-\tau}S(s)\left(\nabla Q(h_{t-s},h_{t-s})\left(\frac 
{d}{dt}h_{t-s}\right)\right)\, (r,v)\, ds 
\end{eqnarray*} 
where $\nabla Q(h_{t-s},h_{t-s})$ denotes a Fr\'echet derivative, which 
is determined by calculations similar to those in the proof of Proposition 
\ref{Proposition3.3} (b). Noting that the test function $d/dt\, h_{t-s}$ 
corresponds to equation (\ref{3.1}) and recalling $h_{t-s}\in D(A)$ from 
the beginning of this paragraph we have shown that 
\begin{eqnarray}\label{3.23}
&&\hspace{-.5cm}Ah_t(r,v)=S(t-\tau)\, Ah_\tau(r,v)+\lambda S(t-\tau)Q(h_\tau, 
h_\tau)\, (r,v)-\lambda Q(h_t,h_t)\, (r,v)\vphantom{\int_0^\ve}\nonumber \\ 
&&\hspace{1.0cm}+2\int_0^{t-\tau}S(s)\left(Q\left(h_{t-s},Ah_{t-s}+ 
\lambda Q(h_{t-s},h_{t-s})\vphantom{l^1}\right)\vphantom{\dot{f}}\right) 
\, (r,v)\, ds 
\end{eqnarray} 
for $-(1+t_0)\le\tau\le 0$ and $\tau\le t$. Fixing now $t=0$ and letting 
$\tau\in [-(1+t_0),0]$ be variable we get 
\begin{eqnarray*} 
&&\hspace{-.5cm}Ah_\tau(r,v)=S(\tau)\, Ah_0(r,v)-\lambda Q(h_\tau, 
h_\tau)\, (r,v)+\lambda S(\tau)\, Q(h_0,h_0)\, (r,v)\vphantom 
{\int_0^\ve} \\
&&\hspace{1.0cm}-2\int_{\tau}^0S(s-\tau)\left(Q\left(h_{\tau-s},Ah_{ 
\tau-s}+\lambda Q(h_{\tau-s},h_{\tau-s})\vphantom{l^1}\right)\vphantom 
{\dot{f}}\right)\, (r,v)\, ds\, . 
\end{eqnarray*} 
For the next calculation set $C:=2\lambda\|B\|\|h_\gamma\|\l(\Omega 
\times V)$ where $B$ and $h_\gamma$ are the terms defined in (kw), 
(k), and (wk) of the introductory part to Section 3. For the 
definition of $Q$ recall also that the modulus of the Jacobian of 
the map $(v,v_1)\mapsto (v^\ast,v_1^\ast)$ is one. We obtain  
\begin{eqnarray*} 
&&\hspace{-.5cm}\|Ah_\tau\|_{L^\infty}\le\|S(\tau)Ah_0\|_{L^\infty}+ 
\|\lambda Q(h_\tau,h_\tau)\, \|_{L^\infty}+\lambda\|S(\tau)Q(h_0,h_0) 
\|_{L^\infty}\vphantom{\int_0^\ve} \\
&&\hspace{1.0cm}+2\int_{\tau}^0S(s-\tau)\left(\left\|Q\left(h_{\tau 
-s},Ah_{\tau-s}+\lambda Q(h_{\tau-s},h_{\tau-s})\vphantom{l^1}\right)
\right\|_{L^\infty}\cdot\1\right)\, ds \\ 
&&\hspace{.5cm}\le c^{\1}_{\infty,\rm max}\|Ah_0\|_{L^\infty}+C\left( 
\|h_\tau\|_{L^\infty}^2+c^{\1}_{\infty,\rm max}\cdot\|h_0\|_{L^\infty 
}^2\right)\vphantom{\int_0^\ve} \\
&&\hspace{1.0cm}+2Cc^{\1}_{\infty,\rm max}\int_{\tau}^0\|h_{\tau-s} 
\|_{L^\infty}\left(\|Ah_{\tau-s}\|_{L^\infty}+\lambda C\|h_{\tau-s} 
\|_{L^\infty}^2\right)\, ds\, . 
\end{eqnarray*} 
By $\nu\equiv\nu_0=h_0\cdot\l$ where $\nu_{-(1+t_0)}\in\hat{\mathcal 
{N}}_{N_1,1/N_1}^D$ and property (4') of Theorem \ref{Theorem5.8} there 
exist $C_0>0$ such that $\|h_{\tau-s}\|_{L^\infty}\le C_0$ for $-(1+ 
t_0)\le\tau\le s\le 0$. Therefore, with $C_1:=c^{\1}_{\infty,\rm max} 
/N_1+CC_0^2(1+c^{\1}_{\infty,\rm max})+2\lambda C^2C_0^3 c^{\1}_{\infty 
,\rm max}$ and $C_2:=2CC_0c^{\1}_{\infty,\rm max}$ it holds that 
\begin{eqnarray*} 
&&\hspace{-.5cm}\|Ah_\tau\|_{L^\infty}\le C_1+C_2\int_{\tau}^0\|Ah_{ 
\tau-s}\|_{L^\infty}\, ds\, ,\quad -(1+t_0)\le\tau\le 0, 
\end{eqnarray*} 
which by Gr\"onwall's inequality implies $\|Ah_{-(1+t_0)}\|_{L^\infty} 
<\infty$. Now we use this for fixed $\tau=-(1+t_0)$ and variable $t\in 
[-(1+t_0),(1+t_0)]$ in (\ref{3.23}). Adjusting the calculations from 
(\ref{3.23}) to here, Gr\"onwall's inequality results in the existence 
of $C_3\in (0,\infty)$ with $\|Ah_t\|_{L^\infty}<C_3$ for all $t\in 
[-(1+t_0),(1+t_0)]$. Together with (\ref{3.20*}) this implies that there 
is some $0<N^{(1)}\le N_1$ such that 
\begin{eqnarray}\label{3.22*}
\nu_t\in\hat{\mathcal{N}}_{N^{(1)},1/N^{(1)}}^D\ \mbox{\rm whenever } 
\nu_{-(1+t_0)}\in\hat{\mathcal{N}}_{N_1,1/N_1}^D\ \mbox{\rm and }t\in 
[-(1+t_0),1+t_0]. 
\end{eqnarray} 
Since for all those $\nu_t=h_t\cdot\l$ we have $N^{(1)}\le h_t\le 
1/N^{(1)}$, the set 
\begin{eqnarray*}
\left\{h_t:\nu_t=h_t\cdot\l\ \mbox{\rm for some }\nu_{-(1+t_0)}\in\hat 
{\mathcal{N}}_{N_1,1/N_1}^D\ \mbox{\rm and }t\in [-(1+t_0),1+t_0]\right 
\}
\end{eqnarray*} 
is uniformly integrable. Recalling the concrete definition of ${\cal 
S}^1_m$ in Step 1 of the proof of Lemma \ref{Lemma1.2} we conclude 
that there is an index $m$ with $m_1\le N_1^{(1)}$ such that 
\begin{eqnarray}\label{3.23*} 
\nu_t\in\mathcal{S}^1_m\ \mbox{\rm whenever}\ \nu_{-(1+t_0)}\in 
\hat{\mathcal{N}}_{N_1,1/N_1}^D\ \mbox{\rm and}\ t\in [-(1+t_0), 
1+t_0]\, . 
\end{eqnarray} 

For an index $N\equiv (N_1,N_2,\ldots\, )$ as in condition (a2) 
introduce 
\begin{eqnarray*}
G_N:=\left\{p_{1+t_0}(h)\cdot\l:h\cdot\l\in\hat{\mathcal{N}}_{N_1,
1/N_1}^D\right\}\, . 
\end{eqnarray*}  
We observe $G_N\subseteq G_{\t N}$ if $\t N_n\le N_n$ for all 
$n\in {\Bbb N}$. In addition, we recall that the time range $[0, 
\t T]$ in (\ref{3.20**}) is by definition independent of $m$ and 
$M$. Since $\{\nu_{\t T}:\nu\in G_N\}$ is by properties (3') and 
(4') of Theorem \ref{Theorem5.8} a subset of some $\hat{\mathcal{N}}_{m',M'}$ we obtain by iteration that (\ref{3.20**}) (and 
hence also (\ref{3.21**})) holds for all $t\ge 0$, i. e. 
\begin{eqnarray}\label{3.27*}
\hat{\cal N}_{m,M}\ni p_0\cdot\l\mapsto p_t(p_0)\cdot\l\quad\mbox 
{\rm is continuous for } t\ge 0.
\end{eqnarray} 
It follows now from Lemma \ref{Lemma3.7} (b) that $G_N$, as a 
continuous image of a compact set in $(E,\pi)$, is closed in $(E, 
\pi)$. Furthermore, we have (1) of condition (a2) by definition and 
(2) of condition (a2) is a consequence of (\ref{3.23*}). In addition, 
(\ref{3.22*}) yields the existence of an index $\hat{N}$ such that 
\begin{eqnarray}\label{3.27}
\nu_t\in G_{\hat{N}}\ \mbox{\rm whenever}\ \nu\in G_N\ \mbox{\rm 
and}\ t\in [-(1+t_0),1+t_0]\, ,
\end{eqnarray} 
i. e. we have shown (3) of condition (a2). In particular, relation 
(\ref{3.27}) for $t\in [0,1]$ is needed for the next conclusion.
\medskip 

We have $\{\nu_t:\nu\in G\}\subseteq G$, $t\in [0,1]$, by (\ref{3.27}). 
In addition, by the result of Step 2, $G\ni\nu\mapsto\nu_t\in G$ is 
injective for all $t\in [0,1]$. En passant we have proved (c) of the 
lemma. 
\medskip

Part (4) of condition (a2) follows now from the definition of $G_N$ 
together with (\ref{3.27*}). We have completed the verification of (a2) 
and therefore completed the proof of part (b) of the lemma. 
\qed 
\medskip

Let us keep in mind the text prior to Proposition \ref{Proposition2.4} 
and recall the definition of ${\cal N}^\ve$ above Proposition 
\ref{Proposition3.3}. Moreover, let us recall the notation $B^1(r)=\{f 
\in L^1(\Omega\times V):\|f\|_{L^1(\Omega\times V)}<r\}$ and that we 
have chosen $r=3/2$, ${\cal V}=B^1(2)$, and $T=S+1$ for some $S\ge 1$. 
\begin{corollary}\label{Corollary3.5} 
(to Propositions \ref{Proposition2.4} and \ref{Proposition3.3}.) 
Suppose that the conditions of Theorem \ref{Theorem3.19} are satisfied.
For the parameter in the Boltzmann type equation, let $\lambda$ satisfy 
the conditions of Proposition \ref{Proposition3.2} (a) and (b) for some 
$T>1$. In addition, let $q\in (1,\infty)$ be arbitrary. Then we have 
conditions (jv') and (v') of Section \ref{sec:2}. 
\end{corollary} 
Proof. {\it Step 1 } In this step we verify the preliminary conditions 
of Proposition \ref{Proposition2.4} and calculate the Fr\'echet 
derivative $\nabla (p(p_0))(\cdot,\cdot,t)$ of the map $p_0\mapsto p_t 
(p_0)$ at fixed time $t\in [0,T]$. 

Choose $0<\ve<1/2$ and set $\t {\cal V}:={\cal 
N}^\ve$. In order to apply Proposition \ref{Proposition2.4} let us 
recall that $G=\{p_{1+t_0}(h)\cdot\l:h\cdot\l\in\hat{\mathcal{N}}^D 
\}$ from Lemma \ref{Lemma3.4}. This implies 
\begin{eqnarray}\label{3.30}
&&\hspace{-.5cm}G=\{p_{t}(q_0)\cdot\l:q_0\in\{(q_0')_{t_0}:q_0'\cdot 
\l\in\hat{\mathcal{N}}^D\}\}\nonumber \\ 
&&\hspace{.5cm}\subset\{p_{t}(q_0)\cdot\l:q_0\in\mathcal{N}\}\subset
\{p_{t}(q_0)\cdot\l:q_0\in\t {\cal V}\}\nonumber \\ 
&&\hspace{.5cm}\subseteq E\cap\{V_t^+(h)\cdot\l:h\in\t {\cal V}\}\, , 
\quad t\in [0,1],
\end{eqnarray}
where the first inclusion holds by Theorem \ref{Theorem5.8} (3') and 
(4'). By the same reference, we have $G\subset\hat{\mathcal{N}}$. This 
yields $G^\ve\subseteq {\cal N}^\ve=\t {\cal V}$, cf. the formulation 
of Proposition \ref{Proposition2.4}. From Paragraph \ref{sec:1:1:2} 
we take $B^1(3/2)\subset{\cal V}\subseteq D(V^+_1)$. Furthermore, 
Proposition \ref{Proposition3.2} (b) asserts that $B^1(3/2)\subseteq 
\{V_1^+(h):h\in D(V_1^+)\}$. 

In accordance with (\ref{3.23**}) the map $i$ is just the identity. As 
pre-neighborhoods $U(0)$ of $0\in X=L^1(\Omega\times V)$ we just take 
open balls in $L^1(\Omega\times V)$ with center 0. In this setup it is 
obvious that for all $\hat{h}\in\t {\cal V}\subseteq B^1(3/2)$ there is a 
(pre-)neighborhood $U_{\hat{h}}(0)$ of $0\in X$ in the topology of $X$ such 
that $\{\hat{h}+i\circ h:h\in U_{\hat{h}}(0)\}=\{\hat{h}+h:h\in U_{\hat{h}} 
(0)\}\subseteq B^1(3/2)\subseteq D(V_1^+)\cap\{V_1^+(h):h\in D(V_1^+)\}$. 

Furthermore, we note that the Fr\'echet differentiability of the solution
map to the Boltzmann type equation of Proposition \ref{Proposition3.3} 
(c) is even stronger than the Fr\'echet differentiability required in 
Proposition \ref{Proposition2.4}. Indeed, by  $0<\ve<1/2$ and the definition 
of the norm $\|\cdot\|_{1,T}$ we have for $p_0\in {\cal N}^\ve=\t {\cal V} 
\subseteq B^1(3/2)$,  $t\in [0,S+1]=[0,T]$, and $h\in L^1(\Omega\times V)$,   
\begin{eqnarray}\label{3.16}
&&\hspace{-.5cm}\frac{\left\|p_t(p_0+h)-p_t(p_0)-\nabla p_t(p_0)(h)\right 
\|_{L^1(\Omega\times V)}}{\|h\|_{L^1(\Omega\times V)}}\nonumber \\ 
&&\hspace{.5cm}\le\frac{\left\|p(p_0+h)-p(p_0)-\nabla p(p_0)(h)\right\|_{ 
1,T}}{\|h\|_{L^1(\Omega\times V)}}\to 0\quad\mbox{\rm as}\ \|h\|_{L^1( 
\Omega\times V)}\to 0\, . \qquad
\end{eqnarray} 
We observe that $\nabla p_t(p_0)(h)=\nabla p(p_0)(h)(\cdot ,\cdot,t)$, 
$t\in [0,T]$, where, on the left-hand side, the symbol $\nabla$ refers 
to the Fr\'echet derivative as a map ${\cal N}^\ve\to B\left(L^1(\Omega 
\times V),L^1(\Omega\times V)\right)$ whereas, on the right-hand side, 
the symbol $\nabla$ refers to the Fr\'echet derivative as a map ${\cal 
N}^\ve\to B\left(L^1(\Omega\times V),(L^1(\Omega\times V))^{[0,T]} 
\right)$. 
\medskip 

\nid
{\it Step 2 } Here we are concerned with the injectivity of the Fr\'echet 
derivative $\nabla p_t(p_0)=\nabla p(p_0)(\cdot ,\cdot,t)$. We also verify 
the remaining conditions of Proposition \ref{Proposition2.4}. According to 
the last paragraph of the proof of Proposition \ref{Proposition3.3}  
$f(h)=\nabla p(p_0)(h)$ satisfies for $h\in L^1(\Omega\times V)$ the 
equation (\ref{3.13}), which takes by Proposition \ref{Proposition3.3} 
(a) and (b) the form
\begin{eqnarray}\label{3.32*} 
f(t)(h)=S(t)h+2\lambda\int_0^tS(t-s)\, Q(f(s)(h),p_s(p_0))\, ds\, ,
\quad t\in [0,T].   
\end{eqnarray} 
As a consequence, it holds that $h=f(0)(h)$ and $f(t_1+t_2)(h)=f(t_1) 
\circ f(t_2)(h)$ for all $t_1,t_2\ge 0$ with $t_1+t_2\le T$ and $h\in 
L^1(\Omega\times V)$. We substitute first $t$ by $u$ and then $p_0$ by 
$p_{t-u}(p_0)$, and correspondingly $h$ by $f(t-u)(h)$. After that we 
apply $S(-u)$ to both sides. We obtain 
\begin{eqnarray}\label{3.32} 
&&\hspace{-.5cm}f(t-u)(h)=S(-u)f(t)(h)-2\lambda\int_0^{u}S(-s)\, Q(f 
(t-u+s)(h),p_{t-u+s}(p_0))\, ds\nonumber \\ 
&&\hspace{.5cm}=S(-u)f(t)(h)-2\lambda\int_0^{u}S(-u+s)\, Q(f(t-s)(h), 
p_{t-s}(p_0))\, ds\, ,\quad u\in [0,t],  \qquad 
\end{eqnarray} 
$h\in L^1(\Omega\times V)$. Given $f(t)(h)$ and looking at (\ref{3.32}) 
as an equation for $f(t-u)(h)$, $u\in [0,t]$, uniqueness of the solution 
would imply injectivity of the map $L^1(\Omega\times V)\ni h\mapsto f(t) 
(h)=\nabla p_t(p_0)(h)\in L^1(\Omega\times V)$  where $p_0\in {\cal 
N}^\ve=\t {\cal V}$ and $t\in [0,T]$. Furthermore, $p_0\cdot\l\in G 
\subseteq\{p_{1}(q_0)\cdot\l:q_0\in\t {\cal V}\}$ by (\ref{3.30}). Thus, 
the uniqueness of the solution to (\ref{3.32}) would also yield 
\begin{eqnarray}\label{3.33} 
\nabla p_{-t}(p_0)=\left(\nabla p_t(p_0)\right)^{-1}\, ,\quad t\in [-1,1], 
\ p_0\cdot\l\in G; 
\end{eqnarray} 
for this recall Theorem \ref{Theorem5.8} (a) and note that more background 
material is contained in Step 2 of the proof of Theorem 5.8 in \cite{Lo20}, 
especially (5.33)-(5.35). 
Assume now that the solution to (\ref{3.32}), as an equation for $f(t-u) 
(h)$, $u\in [0,t]$, $t\in [0,T]$, with given initial value $f(t)(h)=\nabla 
p_t(p_0)(h)$, is not unique. Then there would be a non-trivial map $[0,t]\ni 
u\mapsto\vp(-u)\in L^1(\Omega\times V)$ with $\vp(-u)=-2\lambda\int_{0}^u 
S(-u+s)\, Q(\vp(-s),p_{t-s}(p_0))\, ds$ which by (\ref{3.5}) and the choice 
of $\lambda$ would imply 
\begin{eqnarray*} 
\|\vp(-u)\|_{L^1(\Omega\times V)}\le 4\lambda t\|B\|\|h_\gamma\|\|\vp(- 
\cdot)\|_{1,t}\le\frac14\|\vp(-\cdot)\|_{1,t}\, ,\quad u\in [0,t],  
\end{eqnarray*} 
and the contradiction $0<\|\vp(-\cdot)\|_{1,t}\le \frac14\|\vp(-\cdot)\|_{ 
1,t}$. Therefore the map $L^1(\Omega\times V)\ni h\mapsto f(t)(h)=\nabla 
p_t(p_0)(h)\in L^1(\Omega\times V)$, $p_0\in \t {\cal V}$, is injective and 
we have (\ref{3.33}). 

From equation (\ref{3.32}) it follows that  
\begin{eqnarray*} 
&&\hspace{-.5cm}\|f(t-u)(h)\|_{L^1(\Omega\times V)}\le C_t\|f(t)(h)\|_{L^1 
(\Omega\times V)} \\ 
&&\hspace{.5cm}+4\lambda t\|B\|\|h_\gamma\|C_t\cdot\int_0^u\|f(t-s)(h)\|_{ 
L^1(\Omega\times V)}\, ds\, ,\quad h\in L^1(\Omega\times V),\ u\in [0,t], 
\end{eqnarray*} 
where the constant $C_t$ can be specified as $C_t=Me^{-m_1t}$, see Proposition 
\ref{Proposition3.3**} (b). Now Gr\"onwall's inequality gives with $C:=\exp 
\left\{4\lambda t\|B\|\|h_\gamma\|C_t\cdot t \right\}$ 
\begin{eqnarray}\label{3.35} 
\|h\|_{L^1(\Omega\times V)}\le CC_t\|f(t)(h)\|_{L^1(\Omega\times V)}=CC_t 
\|\nabla p_t(p_0)(h)\|_{L^1(\Omega\times V)}\, .
\end{eqnarray} 
In other words, $\left(\nabla p_t(p_0)\right)^{-1}:\left\{\nabla p_t(p_0) 
(h):h\in L^1(\Omega\times V)\right\}\mapsto L^1(\Omega\times V)$ is bounded 
for all $p_0\in\t {\cal V}$ and $t\in [0,1]$.  

In addition, as a consequence of Proposition \ref{Proposition3.2} (b), for 
all $t\in [0,1]$ the map $\{p_t(p_0'):p_0'\in\t {\cal V}\}\subset\{q\in 
L^1(\Omega\times V):\|q\|_{L^1(\Omega\times V)}\le 3/2\}\ni p_t(p_0)\mapsto 
p_0\in L^1(\Omega\times V)$ is continuous in the topology of $L^1(\Omega 
\times V)$.

Summing up Steps 1 and 2, we have verified the conditions of Proposition 
\ref{Proposition2.4}. In other words, we have (v') and the derivative 
(\ref{2.2}) in condition (jv) exists for $\bmu$-a.e. $\nu\in G$ and all 
$\, t\in [-1,1]$, $g\in C_0(\Omega\times V)$, as well as $k\in X=L^1(\Omega 
\times V)$, for some $r(\nu,t,g)\in X^\ast\cong L^\infty(\Omega\times V)$. 
\medskip

\noindent
{\it Step 3 } We still have to complete the verification of (jv'). For this, 
we keep in mind that the map $i$ is the identity in $L^1(\Omega\times V)
=X$, i. e. $X^\ast\cong L^\infty(\Omega\times V)$. We aim to verify that 
for $g\in C_0(\Omega\times V)$ the term $r(\cdot,t,g)$ given by (\ref{2.2}) 
has a finite upper bound in the norm of $L^\infty(E,\bmu;L^\infty(\Omega 
\times V))$ which is uniform in $t\in [-1,1]$. Moreover we shall show that 
he map 
\begin{eqnarray}\label{3.18} 
t\mapsto r(\cdot,t,g)\ \mbox{\rm is continuous for } t\in [-1,1],\ \mbox 
{\rm in }L^q(E,\bmu;L^\infty(\Omega\times V)).
\end{eqnarray} 
for any $q\in (1,\infty)$. On the right-hand side of (\ref{2.2}) we obtain 
for $g\in C_0(\Omega\times V)$, $t\in [-1,1]$, and $h\in L^1(\Omega\times 
V)$ 
\begin{eqnarray*}
\left.\frac{d}{du}\right|_{u=0}\left((\nu+u\cdot h\cdot\lambda)_t,g\right)
=\left(\left.\frac{d}{du}\right|_{u=0}(\nu+u\cdot h\cdot\lambda)_t,g\right)
=\left(\nabla\nu_t,g\right)(h)\, .  
\end{eqnarray*} 
This means that $r(\nu,t,g)$ represents the bounded linear functional 
$\left(\nabla\nu_t,g\right)$ on $L^1(\Omega\times V)$, $t\in [-1,1]$, for 
well-definiteness recall (\ref{3.33}). From (\ref{3.32*}) and (\ref{3.32}) 
we obtain continuity of $[-1,1]\ni t\mapsto r(p_0\cdot\l,t,g)\in L^\infty  
(\Omega\times V)$ for all $p_0\cdot\l\in G$. Taking into consideration 
the concrete values of the constants in (\ref{3.35}) it follows that 
$r(\cdot,t,g)$ has a finite upper bound in the norm of $L^\infty(E,\bmu; 
L^\infty(\Omega\times V))$ which is uniform in $t\in [-1,0]$. Similarly 
one may obtain this relation for $t\in [0,1]$ by replacing in the 
derivation of (\ref{3.35}) equation (\ref{3.32}) by (\ref{3.32*}). By 
dominated convergence we get (\ref{3.18}). Together with the results of 
Steps 1 and 2, we have verified the conditions (jv') and (v') of Section 
\ref{sec:2} for arbitrary $q\in (1,\infty)$.
\qed 
\medskip 

\noindent 
{\bf Remarks }(1) As a byproduct, from Proposition \ref{Proposition3.3} 
(a) as well as (c) and (\ref{3.16}) we get the subsequent representation 
of the Fr\'echet derivative $\nabla p_t(p_0)$ for all $p_0\in {\cal 
N}^\ve=\t {\cal V}$, and $t\in [0,T]$. Let $\Psi|_{[0,t]}$ denote the 
restriction of $\Psi$ to $\Omega\times V\times [0,t]$. Using a test 
function $h\in L^1(\Omega\times V)$ we have 
\begin{eqnarray*}
&&\hspace{-.5cm}\nabla p_t(p_0)(h)=\left(\sum_{k=0}^\infty\left(\nabla_2 
\Psi|_{[0,t]}\left(p_0,p(p_0)\vphantom{l^2}\right)\right)^k\circ\nabla_1 
\Psi|_{[0,t]}\left(p_0,p(p_0)\vphantom{l^2}\right)\right)(h)(\cdot,\cdot 
,t) \\ 
&&\hspace{.5cm}=\sum_{k=0}^\infty\left(\left(\nabla_2\Psi|_{[0,t]}\left( 
p_0,p(p_0)\vphantom{l^2}\right)\right)^k\circ S(\cdot)h\right)(t)\, , 
\end{eqnarray*} 
where the absolute convergence of the infinite sums in $L^1(\Omega\times V)$ 
is obtained from the choice of $\lambda$, Proposition \ref{Proposition3.3} 
(b), and Proposition \ref{Proposition3.2} (a). 
\medskip

\nid
(2) Assume hypotheses (j)-(jjj). We have verified Theorem \ref{Theorem2.3} 
and Theorem \ref{Theorem2.7} (a). In case of condition (vj), the following 
holds. Condition (vj') of Theorem \ref{Theorem2.7} (b), Theorem \ref{Theorem2.7} 
(c), and Corollary \ref{Corollary2.8} are valid. 
 
\section{Integration by Parts for a Fleming-Viot Type PDE}\label{sec:4}
\setcounter{equation}{0} 

We consider a system $\{X_1,\ldots ,X_n\}$ of $n$ particles in a bounded 
$d$-dimensional domain $D$, $d\ge 2$, with smooth boundary $\partial D$. 
During periods none of the particles $X_1,\ldots ,X_n$ hit the boundary 
$\partial D$, the system behaves like $n$ independent $d$-dimensional 
Brownian motions. When one of the particles hits the boundary $\partial 
D$, then it instantaneously jumps back to $D$ and relocates according to 
a probability distribution $\eta\, (dx)$ on $(D,{\cal B}(D))$. The 
probability distribution $\eta\, (dx)$ depends on the location of the 
remaining $n-1$ particles in a way that relocation is more likely close 
to one of those particles than elsewhere. Such models have their 
background in the sciences as explained, for example, in \cite{CBH05}. 
The asymptotic behavior of such a system $\{X_1,\ldots ,X_n\}$ as $n\to 
\infty$ has been investigated in case of $\eta\, (dx)=\frac{1}{n-1} 
\sum_{i\neq j}^n\delta_{X_i}(dx)$ in \cite{BHM00},\cite{GK04},\cite{Lo051}; 
here $X_j$ being the particle that jumps and $\delta_{X_i}$ denoting the 
Dirac measure concentrated at $X_i$. The nonlinear PDE (\ref{4.2}) describing 
the limiting behavior is taken from \cite{GK04}. 

\subsection{The Asymptotic PDE and Its Solution}\label{sec:4:1} 

We recall that $F$ denotes the set of all finite signed measures on 
$(D,{\cal B})$. In particular let $\l$ denote the Lebesgue measure on 
$(D,{\cal B})$. Let $h_1,h_2,\ldots\, $ be the eigenfunctions of the 
Dirichlet Laplacian on $D$ corresponding to the eigenvalues $0>2 
\lambda_1>2\lambda_2\ge\ldots\, $, normalized in $L^2(D)$ such that 
$h_1>0$. For $t\in [0,1]$ define the space 
\begin{eqnarray*} 
H(t)=\left\{h\cdot\l\in F:\sum_{n=1}^\infty\lambda_n^2e^{-2t\lambda_n}\cdot 
(h_n,h)^2<\infty\right\}
\end{eqnarray*}
which becomes a Hilbert space with the norm $\|h\cdot\l\|_{H(t)}=(\sum_{n 
=1}^\infty\lambda_n^2e^{-2t\lambda_n}\cdot (h_n,h)^2)^{1/2}$, $h\cdot\l\in 
H(t)$. Abbreviate $H\equiv H(0)$. 
\medskip 

Denote by $p(t,x,y)$, $t\ge 0$, $x,y\in D$, the transition density 
function of a Brownian motion on $D$ killed when hitting $\partial D$. 
For $h\in L^2(D)\setminus\{0\}$, $h\cdot\l\equiv\nu$ and $t\ge 0$, let 
$|\nu|:=|h|\cdot\l $, i. e. $|\nu|(D)=\int_D |h|\, dx$, and set 
\begin{eqnarray*}
u(t,y):=\int_{x\in D}p(t,x,y)\, \nu(dx) 
\end{eqnarray*}
as well as 
\begin{eqnarray*}
z(t)\equiv z(\nu ,t):=\frac{1}{|\nu|(D)}\int_{y\in D}\left|\int_{x\in D} 
p(t,x,y)\, \nu(dx)\right|\, dy\, . 
\end{eqnarray*}
For $h\in L^2(D)\setminus\{0\}$, $h\cdot\l\equiv\nu$, and $t\ge 0$ define  
\begin{eqnarray}\label{4.1}
v(t,y):=\frac{1}{z(\nu,t)}u(t,y)=\frac{1}{z(\nu,t)}\int_{x\in D}p(t,x,y) 
\, \nu(dx)\, ,\quad y\in D. 
\end{eqnarray}
For definiteness, if $h=0$ we set $u(t,y)=v(t,y)=0$ and $z(\nu,t)=1$, 
$y\in D$, $t\ge 0$. Now let $\nu\in E\cap H$, $t\ge 0$, and $x\in D$. 
We observe that $v(t,x)$ satisfies 
\begin{eqnarray}\label{4.2} 
\left\{
\begin{array}{l}
\D\frac{\partial}{\partial t}v(t,x)=\frac12\Delta v(t,x)-\frac{z'(\nu, 
t)}{z(\nu,t)}v(t,x)  \\
\D v(t,x)|_{x\in \partial D}=0,\ v(t,\cdot)\cdot\l\stack {t\to 0}{\Ra}\nu 
\vphantom{\int}
\end{array} 
\right.
\end{eqnarray}
where $\partial/\partial t$ is a derivative in $L^2(D)$, $\Delta=2\sum_{ 
i=1}^\infty\lambda_i(h_i,v)_{L^2(D)}h_i$, $v\cdot\l\in H$, is the Laplace 
operator on $D$ and $\Ra$ indicates narrow convergence of finite signed 
measures. In particular we mention that, for $t=0$, $z'(0)$ is the right 
derivative. Below, $\1$ will denote the constant function on $D$ taking 
the value 1.
\medskip

\nid 
{\bf Remarks }(1) Let $t\ge 0$ and $\nu=h\cdot\l$ with $h\in L^2(D)\setminus 
\{0\}$. It holds that 
\begin{eqnarray*}
v(t,\cdot)=(z(\nu,t))^{-1}\sum_{j=1}^\infty e^{\lambda_j t}(h_j,\nu)h_j= 
\frac{\sum_{j=1}^\infty e^{\lambda_j t}(h_j,\nu)h_j}{\int_{y\in D}\left| 
\sum_{j=1}^\infty e^{\lambda_j t}(h_j,\nu)h_j(y)\right|\, dy}\cdot |\nu|(D) 
\end{eqnarray*} 
and, in particular for $\mu:=v(t,\cdot)\cdot\l$, we have $|\mu|(D)=|\nu| 
(D)$. Among other things this shows that $\nu\mapsto v(t,\cdot)\cdot\l$ maps 
$E\cap H$ into $E\cap H$. Furthermore, 
\begin{eqnarray*} 
v(t,\cdot)\cdot\l&&\hspace{-.5cm}=\frac{\sum_{j=1}^\infty e^{\lambda_j t} 
(h_j,\nu)h_j\cdot\l}{\sum_{j=1}^\infty e^{\lambda_j t}(h_j,\nu)(h_j,\1)} 
 \\ 
&&\hspace{-.5cm}=\frac{(h_1,\nu)h_1\cdot\l+\sum_{j=2}^\infty e^{(\lambda_j 
-\lambda_1)t}(h_j,\nu)h_j\cdot\l}{(h_1,\nu)(h_1,\1)+\sum_{j=2}^\infty 
e^{(\lambda_j-\lambda_1)t}(h_j,\nu)(h_j,\1)} \\ 
&&\hspace{-0.9cm}\stack{t\to\infty}{\Ra}\frac{h_1\cdot\l}{(h_1,\1)}\, , 
\quad\nu\in E\cap H\, .\vphantom{\frac{\sum}{\sum}}
\end{eqnarray*} 
The evolution of the infinite particle system is hence asymptotically 
directed toward the point $h_1\cdot\l/(h_1,\1)$ of $E$ in the narrow 
topology. 
\medskip 

\nid
(2) For all $\nu=h\cdot\l\in E\cap H$, the expression $\frac12\Delta h 
=\sum_{i=1}^\infty\lambda_i(h_i,h)_{L^2(D)}h_i$ is well-defined in $L^2 
(D)$. This is one motivation to take special consideration to the space 
$H$. We will also write $\Delta\nu$ for $(\Delta h)\cdot\l$. 
\medskip 

Next, we focus on compatibility with the framework of Sections 
\ref{sec:1} and \ref{sec:2}. For the terms of Subsection \ref{sec:1:1}, 
let the exponent $v=2$ and choose arbitrary $S\ge 1$, $r>1$. Let 
$\lambda\equiv \l$ be the Lebesgue measure on $(D,{\cal B})$. Moreover 
set ${\cal V}:=L^2(D)$ and $\t V_t^+h:=v(t,\cdot)$, $t\in [0,S+1]$, 
defined by (\ref{4.1}) for $\nu=h\cdot\l$ and $h\in {\cal V}$. Note 
also that in this application we have $(\t V_t^+,{\cal V})=(V_t^+, 
D(V_t^+))$, $t\in [0,S+1]$. 

We obtain for $t\in [0,S+1]$, $\mu=\t V_t^+h\cdot \l$ for some $h\in 
{\cal V}$, and $\nu=h\cdot\l$
\begin{eqnarray}\label{4.21}
\sum_{j=1}^\infty e^{-\lambda_j t}(h_j,\mu)h_j\cdot\l=\sum_{j=1}^\infty 
(z(\nu,t))^{-1}\cdot(h_j,\nu)h_j\cdot\l=(z(\nu,t))^{-1}\cdot\nu\, , 
\end{eqnarray} 
recall also Remark (1) of this section. This says that, for given $\mu$, 
the objects $h\in {\cal V}$, and $\nu=h\cdot\l$ are uniquely determined by 
\begin{eqnarray}\label{4.22}
\nu=h\cdot\l=\frac{\sum_{j=1}^\infty e^{-\lambda_j t}(h_j,\mu)h_j\cdot 
\l}{\int_{y\in D}\left|\sum_{j=1}^\infty e^{-\lambda_j t}(h_j,\mu)h_j
(y)\right|\, dy}\cdot|\mu|(D)\, . 
\end{eqnarray} 
Furthermore, 
\begin{eqnarray}\label{4.23}
(z(\nu,t))^{-1}=\frac{1}{|\mu|(D)}\int_{y\in D}\left|\sum_{j=1}^\infty 
e^{-\lambda_j t}(h_j,\mu)h_j(y)\right|\, dy\, . 
\end{eqnarray} 
In other words, the map ${\cal V}\ni h\mapsto\t V_t^+h \in {L^2(D)}$ is for 
all $t\in [0,S+1]$ injective. In particular, $H\ni\nu=h\cdot\l\mapsto\t V_t^+ 
h\cdot\l\in H(t)$ is for all $t\in [0,1]$ a bijection. 
\medskip

By the definitions of Subsection \ref{sec:1:1} we are now able to 
introduce the flow $(W^+(t,\cdot),$ $D(W^+(t,\cdot)))$, $t\in[0,S+1]$, 
the trajectory $(U^+(s,\cdot),D(U^+(s,\cdot)))$, $s\in[-1,1]$, and 
\begin{eqnarray*}
\nu_s\equiv U^+(s,\nu)\, ,\quad \nu\in D(U^+(s,\cdot)),\ s\in [-1,1]. 
\end{eqnarray*} 
Recall that all $D(U^+(s,\cdot))$, $s\in[-1,1]$, are identical by 
definition. Recall also from Remark (1) of this section that the map 
${\cal V}\ni h\mapsto\t V_t^+h\in L^2(D)$ preserves for all $t\in 
[0,S+1]$ the norm in $L^1(D)$. Thus, the flow property $W^+(t,\cdot) 
\circ W^+(s,\nu)=W^+(s+t,\nu)$ for $\nu\in D(W^+(s+t,\cdot))$ and $s,t 
\in [0,S+1]$ with $s+t\in [0,S+1]$ is an immediate consequence of 
(\ref{4.1}). In order to introduce the set $G$ for this application in 
the following lemma, choose $t_0>0$. This will be the $t_0$ appearing 
in hypothesis (a2) of Subsection \ref{sec:1:1}. 
\begin{lemma}\label{Lemma4.0} 
(a) There is a compact set of the form ${\cal S}^2_m$ introduced in 
Lemma \ref{Lemma1.2} such that $\left\{W^+(1,h\cdot\l):h\in L^2(D),\ h 
\ge 0,\ \|h\|_{L^1(D)}=1\right\}\subseteq {\cal S}^2_m$. \\ 
(b) It holds that 
\begin{eqnarray*}
G&&\hspace{-.5cm}:=\left\{\nu\in E\cap H(1+t_0):\nu_{-1-t_0}\ge 0 
\right\} \\ 
&&\hspace{-.4cm}\equiv\{W^+(1+t_0,\nu):\nu\in E\cap H\}\in {\cal B} 
(E)\, .  
\end{eqnarray*} 
(c) For all $t\ge 0$ the map ${\cal S}^2=E\cap\{h\cdot\l:h\in L^2(D)\} 
\ni\nu\mapsto\nu_t$ is injective. It holds that $\{\nu_t:\nu\in G\} 
\subseteq G$, $t\ge 0$, i. e. we have Case 2. \\ 
(d) For all $t\ge 0$ the map ${\cal S}^2=E\cap\{h\cdot\l:h\in L^2(D)\} 
\ni\nu\mapsto\nu_t$ is continuous in $(E,\pi)$. We have conditions 
(a1)-(a3) of Subsection \ref{sec:1:1} and (a4) of Lemma \ref{LemmaA.1}. 
\end{lemma} 
Proof. {\it Step 1 } We verify (a). According to \cite{SV09}, Subsection 
5.5.1, and \cite{DS84}, Theorem 3.2, there exist constants $0<\t c_t<c_t 
<\infty$ such that $\t c_th_1(x)h_1(y)\le p(t,x,y)\le c_th_1 (x)h_1(y)$, 
$x,y\in D$, $t>0$, where here $h_1$ denotes the eigenfunction of the 
Dirichlet Laplacian on $D$ corresponding to the eigenvalue $2\lambda_1>$.  
For $\nu=h\cdot\l\in E$ with $h\in L^2(D)$ this has the consequence 
\begin{eqnarray*}
v(1,y)=\frac{\int_{x\in D} p(1,x,y)\, \nu(dx)}{\int_{y\in D}\int_{x\in D} 
p(1,x,y)\, \nu(dx)\, dy\vphantom{\dot{f}}}\le\frac{c_1}{\t c_1}\cdot\frac 
{h_1(y)}{(h_1,\1)}\, ,\quad y\in D. 
\end{eqnarray*}
We obtain $\|v(1,\cdot)\|_{L^2}\le c_1/(\t c_1\cdot (h_1,\1))$. In the 
context of Lemma \ref{Lemma1.2} (b) and its proof we have therefore $W^+ 
(1,h\cdot\l)\in {\cal S}^2_m$ where $m$ is the greatest integer smaller 
than $c_1/(\t c_1\cdot (h_1,\1))+1$. We recall that ${\cal S}^2_m$ is 
compact in $(E,\pi)$, cf. Lemma \ref{Lemma1.2} (b) and its proof.  
\medskip 

\nid
{\it Step 2 } We prove part (b). For the following conclusions let us 
keep in mind that $H\ni\nu_{-1-t_0}\equiv h\cdot\l\mapsto\nu\equiv W^+ 
(1+t_0,h\cdot\l)\in H(1+t_0)$ is a bijection, cf. (\ref{4.21})-(\ref 
{4.23}), and the paragraph below these references. For $N\in {\Bbb N}$ 
set 
\begin{eqnarray*} 
&&\hspace{-.5cm}G_{N}:=\left\{\nu\in E\cap H(1+t_0):\nu_{-1-t_0}\equiv 
h\cdot\l,\ \right. \vphantom{W^+(t_0)}\\ 
&&\hspace{1.5cm}\left. h\ge 0,\ \|h\|_{L^2(D)}\vee\|h\cdot\l\|_H\le N 
\right\} \\ 
&&\hspace{.5cm}=\left\{W^+(1+t_0,h\cdot\l):h\cdot\l\in H,\ \right. \\ 
&&\hspace{1.5cm}\left. h\ge 0,\ \|h\|_{L^1(D)}=1,\ \|h\|_{L^2(D)}\vee 
\|h\cdot\l\|_H\le N\right\}\, . 
\end{eqnarray*} 
Since $\{W^+(1+t_0,h\cdot\l):h\cdot\l\in H,\ h\ge 0\}\subseteq\{W^+(1, 
h\cdot\l):h\cdot\l\in H,\ h\ge 0\}$ the set $G_{N}$ is by part (a) of 
this lemma a subset of some compact set in $(E,\pi)$ of the form ${\cal 
S}^2_m$. Consequently, for any Cauchy sequence with respect to $(E,\pi)$ 
of the form $\nu^{(n)}\equiv W^+(1+t_0,h^{(n)}\cdot\l)\in G_{N}$, $n\in 
{\Bbb N}$, there is a limit $\nu\in{\cal S}^2_m$ and we have the 
following. Since $\|h^{(n)}\|_{L^2(D)}\vee\|h^{(n)}\cdot\l\|_H\le N$, 
$n\in {\Bbb N}$, there is a subsequence $h^{(n_k)}$, $k\in {\Bbb N}$, 
such that $h^{(n_k)}$, $k\in {\Bbb N}$, converges weakly in $L^2(D)$ to 
some $h\in L^2(D)$ and $h^{(n_k)}\cdot\l$, $k\in {\Bbb N}$, converges 
weakly in $H$ to some $\t h\cdot\l\in H$ as $k\to\infty$. Moreover we 
have $\|h\|_{L^2 (D)}\vee\|\t h\cdot\l\|_H\le N$ and $h\ge 0$. Choose 
\begin{eqnarray*} 
\t g\cdot\l\in\t H:=\left\{\t f\cdot\l\in H:\sum_{j=1}^\infty\lambda_j^4 
(\t f,h_j)^2<\infty\right\} 
\end{eqnarray*} 
and let $g:=\sum_{j=1}^\infty\lambda_j^2(\t g,h_j)h_j$. Noting that 
$\t H$ is dense in $H$ and that $h^{(n_k)}\cdot\l$, $k\in {\Bbb N}$, is 
bounded in $H$ from 
\begin{eqnarray*}
&&\hspace{-.5cm}(\t g\cdot\l,\t h\cdot\l)_H=\lim_{k\to\infty}(\t g\cdot 
\l,h^{(n_k)}\cdot\l)_H \\  
&&\hspace{.5cm}=\lim_{k\to\infty}\left(\sum_{j=1}^\infty\lambda_j^2(\t g, 
h_j)h_j ,\sum_{j=1}^\infty(h^{(n_k)},h_j)h_j\right) \\  
&&\hspace{.5cm}=\lim_{k\to\infty}\left(g,h^{(n_k)}\right)=(g,h)=(\t g 
\cdot\l,h\cdot\l)_H
\end{eqnarray*} 
it follows that $h\cdot\l=\t h\cdot\l\in H$. Furthermore, for any $f 
\in C_b(D)$ it holds that 
\begin{eqnarray}\label{4.25}
(f,\nu^{(n_k)})=\frac{\sum_{j=1}^\infty e^{\lambda_j(1+t_0)}(h_j,h^{ 
(n_k)})(f,h_j)}{\sum_{j=1}^\infty e^{\lambda_j(1+t_0)}(h_j,h^{(n_k)}) 
(\1,h_j)}\stack{k\to\infty}{\lra}\frac{\sum_{j=1}^\infty e^{\lambda_j 
(1+t_0)}(h_j,h)(f,h_j)}{\sum_{j=1}^\infty e^{\lambda_j(1+t_0)}(h_j,h) 
(\1,h_j)}\, . 
\end{eqnarray} 
For the convergences of the infinite sums recall $\|h^{(n_k)}\|_{L^2 
(D)}\le N$, $k\in {\Bbb N}$, as well as $\|h\|_{L^2(D)}\le N$ and 
$\liminf_{j\to\infty}-\lambda_j/j^{2/d}>0$ cf. \cite{Va92}. Since by 
the Portmanteau theorem in the form of \cite{EK86}, Theorem 3.3.1, the 
right-hand side of (\ref{4.25}) coincides with $(f,\nu)$, this says 
that $\nu=W^+(1+t_0,h\cdot\l)\in G_{N}$. We have verified that $G_{N}$ 
is closed in $(E,\pi)$. This yields 
\begin{eqnarray*}
G=\bigcup_{N=1}^\infty G_{N}\in {\cal B}(E)\, .  
\end{eqnarray*} 

\noindent
{\it Step 3 } Let us demonstrate (c). As already demonstrated in 
(\ref{4.21})-(\ref{4.23}), and the paragraph below these relations, 
the map ${\cal S}^2=E\cap\{h\cdot\l:h\in L^2(D)\}\ni\nu\mapsto\nu_t$ 
is injective for all $t\in [0,S+1]$ and by iteration for all $t\ge 0$. 
We have $G=\{W^+(1+t_0,\nu):\nu\in E\cap H\}$. If $\nu=W^+(1+t_0,\t \nu) 
\in G$ for some $\t \nu\in E\cap H$ then $\nu_t=W^+(1+t_0,\cdot)\circ 
W^+(t,\t \nu)$, $t\ge 0$. Since $W^+(t,\t \nu)\in E\cap H$, we have 
$\nu_t\in G$, i. e. $\{\nu_t:\nu\in G\}\subseteq G$, $t\ge 0$. 
\medskip

\nid
{\it Step 4 } In this step, we show (a1) and (a2). Let us begin with 
(a2). Properties (1) and (2) of (a2) are a byproduct of Step 2 of the 
present proof. Concerning property (3) of (a2) we observe that for $N 
\in\mathbb{N}$ and $W^+(1+t_0,h_{-(1+t_0)}\cdot\l)\in G_N$ for some 
$h_{-(1+t_0)}\cdot\l\in E\cap H$ with $ h_{-(1+t_0)}\ge 0,\ \|h_{-(1+ 
t_0)}\|_{L^2(D)}\vee\|h_{-(1+t_0)}\cdot\l\|_H\le N$ the term 
\begin{eqnarray*} 
h_t\cdot\l:=W^+(t,\cdot)\circ W^+(1+t_0,h_{-(1+t_0)}\cdot\l)=W^+(1+t_0 
+t,h_{-(1+t_0)}\cdot\l)
\end{eqnarray*} 
satisfies  
\begin{eqnarray*} 
\left(h_{1+t_0+t}\right)_{-(1+t_0)}=h_t=\frac{\sum_{j=1}^\infty e^{ 
\lambda_j(1+t_0+t)}(h_j,h_{-(1+t_0)})h_j\cdot\l}{\sum_{j=1}^\infty 
e^{\lambda_j(1+t_0+t)}(h_j,h_{-(1+t_0)})(h_j,\1)}\, , 
\end{eqnarray*} 
$t\in [-(1+t_0),0]$. Recalling the inequality from the beginning 
of Step 1 follows that 
\begin{eqnarray*} 
\|h_t\|_{L^2(D)}^2\le\frac{\sum_{j=1}^\infty (h_j,h_{-(1+t_0)})^2} 
{\left(\int\int p(1+t_0,x,y)h_{-(1+t_0)}(y)\, dy\, dx\right)^2}\le 
\frac{\|h_{-(1+t_0)}\|_{L^2(D)}^2}{\t c^2_{1+t_0}(h_1,\1)^2} 
\end{eqnarray*} 
and similarly 
\begin{eqnarray*} 
\|h_t\|_H^2\le\frac{\|h_{-(1+t_0)}\|_H^2}{\t c^2_{1+t_0}(h_1,\1)^2} 
\, ,\quad t\in [-(1+t_0),0]. 
\end{eqnarray*} 
Thus, $\|h_t\|_{L^2(D)}\vee\|h_t\cdot\l\|_H\le N/(\t c_{1+t_0}(h_1, 
\1))$ for all $t\in [-(1+t_0),0]$, i. e. we have property (3) of (a2). 

For the verification of hypothesis (a2) it remains to show property 
(4). For this, let us demonstrate a stronger property, namely that 
for $t\ge 0$, the map ${\cal S}^2\ni\nu\mapsto\nu_t$ is continuous 
with respect to $\pi$. This can be obtained from the Portmanteau 
theorem as follows. Because of $\liminf_{j\to\infty}-\lambda_j/j^{ 
2/d}>0$ and $\|h_j\|\le (-2\lambda_j)^{d/4}$, $j\in {\Bbb N}$, see 
\cite{Va92}, \cite{BB99}, and \cite{Da73}, for $t>0$ and $f\in C_b(D)$ 
the sum $\sum_{j=1}^\infty e^{\lambda_jt}|(f,h_j)|\|h_j\|$ converges, 
$\sum_{j=1}^\infty e^{\lambda_j t}(f,h_j)h_j$ converges in $C_b(D)$, 
and we have 
\begin{eqnarray*} 
&&\hspace{-.5cm}\left(f,(\nu^{(n)})_t\right)=\frac{\sum_{j=1}^\infty 
e^{\lambda_j t}(h_j,\nu^{(n)})(f,h_j)}{\sum_{j=1}^\infty e^{\lambda_j 
t}(h_j,\nu^{(n)})(h_j,\1)}\stack{n\to\infty}{\lra}\frac{\sum_{j=1 
}^\infty e^{\lambda_j t}(h_j,\nu)(f,h_j)}{\sum_{j=1}^\infty e^{\lambda_j  
t}(h_j,\nu)(h_j,\1)}=\left(f,\nu_t\right)
\end{eqnarray*} 
whenever $\nu^{(n)}\in {\cal S}^2$ converges to $\nu\in {\cal S}^2$ 
in the Prokhorov metric $\pi$ as $n\to\infty$. 

Let us verify (a1) of Subsection \ref{sec:1:1} in the remainder of this 
step. Let $t\in [0,1]$. We recall that by parts (a) and (b) of this 
lemma, $G\in {\cal B}(E)$ is a subset of the compact metric space 
$({\cal S}^2_m,\pi)$ where $m\in {\Bbb N}$ has been fixed in part (a). 
According to (c) the map ${\cal S}^2_m\ni\nu\mapsto\nu_t$ is continuous 
in $(E,\pi)$. Theorem 6.9.7 of \cite{Bo07} states that there is a set 
$B\in {\cal B}({\cal S}^2_m)$ with $\{\nu_t:\nu\in B\}=\{\nu_t:\nu\in 
{\cal S}^2_m\}$ such that the map $B\ni\nu\mapsto\nu_t$ is injective 
and its inverse is Borel measurable. However the injectivity of $G\ni 
\nu\mapsto\nu_t$ implies that we may assume $G\subseteq B$ and hence 
the map $\{\t \nu_t:\t \nu\in G\}\ni\nu_t\mapsto\nu\in G$ is Borel 
measurable. As already mentioned, $G\ni\nu\mapsto\nu_t$ is continuous 
and hence Borel measurable.

Let us recall that ${\cal S}^2=E\cap\{h\cdot\l:h\in L^2(D)\}$ and $G= 
\{W^+(1+t_0,\nu):\nu\in E\cap H\}$ yield ${\cal S}^2\supseteq\{\t\nu_{ 
-t}:\t \nu\in G\}$. Furthermore, ${\cal S}^2_m$, $m\in {\Bbb N}$, are 
compact subsets of $(E,\pi)$. Part (c) and the already proved first 
part of (d) say that for all $t\in [0,1]$ the map ${\cal S}^2\supseteq 
\{\t\nu_{-t}:\t \nu\in G\}\ni\nu_{-t}\mapsto\nu\in G$ is continuous and 
injective. It follows now that for all $m\in {\Bbb N}$ the map $\{\t 
\nu_{-t}:\t \nu\in G\}\cap {\cal S}^2_m\ni\nu_{-t}\mapsto\nu\in G$ is 
continuous and hence Borel measurable. As above, it turns out that $G 
\cap\{\t \nu_t:\t \nu\in {\cal S}^2_m\}\ni\nu\mapsto\nu_{-t}\in\{\t 
\nu_{-t}:\t \nu\in G\}\cap {\cal S}^2_m$ is Borel measurable for all 
$m\in {\Bbb N}$. With Lemma \ref{Lemma1.2} (b) we find Borel 
measurability of $\{\t\nu_{-t}:\t \nu\in G\}\ni\nu_{-t}\mapsto\nu\in G$ 
and $G\ni\nu\mapsto\nu_{-t}$. 

Furthermore, for $h\cdot\l\in G\subseteq H(1+t_0)$ and $t\in [-1,1]$, 
the sum $\sum_{j=1}^\infty e^{\lambda_j(1-t)}(h_j,h_{-1})h_j$ converges 
in $L^2(D)$. As a consequence, $[-1,1]\ni t\mapsto h_{-t}$ is continuous 
in $L^2(D)$. We have completed the verification of (a1). 
\medskip

\nid  
{\it Step 5 } First we verify condition (a3) of Subsection \ref{sec:1:1}.  
Let us begin with some preparations. For $\nu=h\cdot\l\in E\cap H$ it 
follows that 
\begin{eqnarray}\label{4.26}
z'(0)\equiv z'(\nu,0)&&\hspace{-.5cm}=\left.\frac{d}{dt}\right|_{t=0} 
\int_{y\in D}\int_{x\in D}p(t,x,y)\, \nu(dx)\, dy\nonumber \\ 
&&\hspace{-.5cm}=\int_{y\in D}\left.\frac{d}{dt}\right|_{t=0}\int_{x 
\in D}p(t,x,y)\, h(x)\, dx\, dy\nonumber \\ 
&&\hspace{-.5cm}=\int_{y\in D}\frac12(\Delta h)(y)\, dy  
\end{eqnarray}
where the interchange of $\left.\frac{d}{dt}\right|_{t=0}$ and $\int_{ 
y\in D}$ in the second line applies because the derivative $\left.\frac 
{d}{dt}\right|_{t=0}\int_{x\in D}p(t,x,\, \cdot\, )\, h(x)\, dx$ is a 
limit in $L^2(D)$ and $\int_{y\in D}\left.\frac{d}{dt}\right|_{t=0} 
\int_{x\in D}p(t,x,y)\, h(x)\, dx\, dy$ is the inner product of this 
limit with $\1\in L^2(D)$. 

It follows from (\ref{4.26}) that we have for $\nu=h\cdot\l\in E\cap H$   
\begin{eqnarray*}
|z'(\nu,0)|\le\frac12\, (\l(D))^{\frac12}\left\|\Delta h\right\|_{L^2(D)}< 
\infty\, . 
\end{eqnarray*}
Keeping Remark (1) of this section in mind we obtain the following 
representations of the generator $A^f:E\cap H\to L^2(D)$ of $U^+$. We 
have 
\begin{eqnarray}\label{4.3}
A^f\nu&&\hspace{-.5cm}=\left.\frac{d}{dt}\right|_{t=0}\nu_t=\frac12 
\Delta\nu -z'(0)\, \nu \nonumber \\ 
&&\hspace{-.5cm}=\sum_{i=1}^\infty\left(\lambda_i-\int\sum_{j=1}^\infty 
\lambda_j(h_j,\nu)h_j\, dx\right)(h_i,\nu)\cdot h_i\cdot\l 
\end{eqnarray}
which by (\ref{4.26}) exists for every $\nu\in E\cap H$ in $L^2(D)$. 
For $\nu\equiv h\cdot\l\in E$ and $t\in [0,1]$ it follows similar to 
(\ref{4.3}) from (\ref{4.21}) and (\ref{4.22}) that 
\begin{eqnarray*}
&&\hspace{-.5cm}\frac{dA^f\nu_{-t}(\cdot)}{d\, \l} \\ 
&&\hspace{.5cm}=z(\nu_{-t},t)\cdot\sum_{n=1}^\infty\left(\lambda_n- 
z(\nu_{-t},t)\int\sum_{j=1}^\infty\lambda_je^{-t\lambda_j}(h_j,\nu)h_j 
\, dx\right)e^{-t\lambda_n}(h_n,h)\cdot h_n
\end{eqnarray*}
if the right-hand side exists in $L^2(D)$. Thus the limit $dA^f\nu_{-t} 
(\cdot)/d\,\l=\lim_{s\to 0}\frac1s\left(h_{s-t}-h_{-t}\right)$ exists 
in $L^2(D)$ for a fixed $t\in [0,1]$ if and only if $\sum_{n=1}^\infty 
\lambda_n^2e^{-2t\lambda_j}(h_n,\nu)^2<\infty$ that is, if and only if 
$\nu\in E\cap H(t)$. The definition of $H(t)$, $t\in [0,1+t_0]$, yields 
\begin{eqnarray*}
&&\hspace{-.5cm}G\equiv\{\nu\in E\cap H(1+t_0):\nu_{-1-t_0}\ge 0\} 
\subset\left\{\nu=h\cdot\l\in E:\vphantom{\frac{dA^f\nu_{-t}(\cdot)} 
{d\, \l}\lim_{s\to 0}}\right. \\ 
&&\hspace{.5cm}\left.\frac{dA^f\nu_{-t}(\cdot)}{d\, \l}=\lim_{s\to 0} 
\frac1s(h_{s-t}-h_{-t})\ \mbox{\rm exists in}\ L^2(D)\ \mbox{\rm for all} 
\ t\in [0,1]\right\}\, .  
\end{eqnarray*}
In addition, for $\nu\in E\cap H(1+t_0)$ and $t\in [-1,0)$ relation  
(\ref{4.3}) implies the existence of 
\begin{eqnarray*}
&&\hspace{-.5cm}\frac{dA^f\nu_{-t}(\cdot)}{d\, \l} \\ 
&&\hspace{.5cm}=\sum_{n=1}^\infty\frac{1}{z(\nu,-t)}\left(\lambda_n- 
\frac{1}{z(\nu,-t)}\int\sum_{j=1}^\infty\lambda_je^{-t\lambda_j}(h_j, 
\nu)h_j\, dx\right)e^{-t\lambda_n}(h_n,\nu)\cdot h_n 
\end{eqnarray*}
in $L^2(D)$ and, as a consequence, yields the inclusion  
\begin{eqnarray*}
&&\hspace{-.5cm}G\subset\left\{\nu=h\cdot\l\in E:\vphantom{\frac{dA^f 
\nu_{-t}(\cdot)}{d\, \l}=\lim_{s\to 0}}\right. \\ 
&&\hspace{.5cm}\left.\frac{dA^f\nu_{-t}(\cdot)}{d\, \l}=\lim_{s\to 0} 
\frac1s(h_{s-t}-h_{-t})\ \mbox{\rm exists in}\ L^2(D)\ \mbox{\rm for 
all}\ t\in [-1,1]\right\}\, .  
\end{eqnarray*}
Let $\bnu$ be a probability measure on $(E,{\cal B}(E))$ as introduced 
in Paragraph \ref{sec:1:1:4}. In particular $\bnu(G)=1$ implies that 
the limit in (\ref{1.3}), (\ref{1.4}), namely $dA^f\nu_{-t}(\cdot)/d 
\,\l=\lim_{s\to 0}\frac1s(h_{s-t}-h_{-t})$, exists in $L^2(D)$ for 
all $t\in [-1,1]$ for $\bnu$-a.e. $\nu\equiv h\cdot\l\in G$. In other 
words, it guarantees (a3) of Subsection \ref{sec:1:1}. 

As for condition (a4) we remark that $X=\left\{h\in L^2(D):\sum_{n=1 
}^\infty e^{-2\lambda_n}\cdot (h_n,h)^2<\infty\right\}$ is going to be 
the space we use for the verification of conditions (jv) and (v). It is 
also compatible with condition (j). For all this, see Lemma \ref{Lemma4.2} 
and its proof below. From part (b) of the present lemma, we obtain 
$G\subset H(1)\subset X$ which implies (a4). 
\qed
\medskip

\nid 
{\bf Remark} (3) Let us continue Remark (2) of the present section. We 
recall that, for $\nu=h\cdot\l$ and $U^+(t,\nu)=h_t\cdot\l$, $t\in [-1, 
1]$, the right derivative $\left.\frac{d}{dt}\right|_{t=0}h_t$ exists in 
$L^2(D)$ if and only if $\nu=h\cdot\l\in E\cap H$, cf. (\ref{4.26}) and 
(\ref{4.3}) and the text below these references. In this sense, $A^f\nu= 
\left.\frac{d}{dt}\right|_{t=0}\nu_t$ is well-defined. Together with 
Remark (2), this is the actual motivation to take special consideration 
to the space $H$. Under certain circumstances, we also will follow a 
trajectory backward in time, cf. Theorem \ref{Theorem2.3} or hypothesis 
(a2), the latter being related to the choice of $G$. For this, we may 
need to restrict the space $H$ to $H(t)$, $t\in [0,1+t_0]$. 
\medskip 

Let $\bmu$ be an arbitrary probability measure on $(E,{\cal B}(E))$ 
satisfying the hypotheses of Paragraph \ref{sec:1:1:4}. Next we provide 
a technical detail which we will use several times throughout the 
remainder of this section. 
\begin{lemma}\label{Lemma4.11} 
We have $(z(\cdot ,1))^{-1}\in L^\infty(E,\bmu)$ which by $\bmu$-a.e. 
monotonicity implies $$(z(\cdot ,t))^{-1}\in L^\infty(E,\bmu)\, ,\quad 
t\in [0,1]\, .$$
\end{lemma}
Proof. We recall $\bmu(G)=1$ from the basic hypotheses of Paragraph 
\ref{sec:1:1:4}. Let $0<\t c_1<c_1<\infty$ be the constants introduced 
in Step 1 of the proof of Lemma \ref{Lemma4.0}. For $\nu=\t h\cdot\l\in 
G$ with $\nu=W^+(1,\t h_{-1}\cdot\l)$ and $\t h_{-1}\cdot\l\in E\cap H 
(t_0)\subset E\cap H$ it holds that 
\begin{eqnarray*} 
&&\hspace{-.5cm}(h_1,\nu)=\frac{\int\int p(1,x,y)\t h_{-1}(x)h_1(y) 
\, dx\, dy}{\int\int p(1,x,y)\t h_{-1}(x)\, dx\, dy}\ge\frac{\t c_1(h_1, 
\t h_{-1})(h_1,h_1)}{c_1(h_1,\t h_{-1})(h_1,\1)}=\frac{\t c_1}{c_1(h_1, 
\1)}\, .
\end{eqnarray*} 
and hence 
\begin{eqnarray*} 
(z(\nu ,1))^{-1}=\left(\int\int p(1,x,y)\, \nu(dx)\, dy\right)^{-1}\le 
\left(\t c_1(h_1,\nu)(h_1,\1)\right)^{-1}\vphantom{\t f}\le\frac{c_1} 
{\t c_1^2}\, ,\quad \nu\in G. 
\end{eqnarray*} 
\qed
\medskip

We continue with some preparations for the next lemma. Because of $\nu\in 
G^+\equiv G_0^+$ implies $\nu\in G_t$ for some $t>0$ we have $\nu_{-t} 
\in G$, and therefore $\nu_{-1-t}:=(\nu_{-t})_{-1}\in E\cap H(t_0)\subset 
E\cap H$. We obtain $\nu_{-1-t}=(\nu_{-1})_{-t}$ and 
\begin{eqnarray*}
\frac{\nu_{-1}(dy)}{dy}=\frac{1}{z(\nu_{-1-t},t)}\int_D p(t,x,y)\, \nu_{ 
-1-t}(dx)>0\, ,\quad y\in D,  
\end{eqnarray*}
which we abbreviate with $\nu_{-1}>0$. Consequently $G^+\subseteq\{\nu 
\in G:\nu_{-1}>0\}$. From the hypothesis $\bmu(G\setminus G^+)=0$ in 
Paragraph \ref{sec:1:1:4} it follows that 
\begin{eqnarray}\label{4.115}
\bmu(\{\nu\in G:\nu_{-1}>0\})=1\, . 
\end{eqnarray}
A suitable choice of the space $X$ appearing in conditions (jv) as well 
as (v) and (jv') as well as (v') is 
\begin{eqnarray*} 
X=\left\{h\in L^2(D):\sum_{n=1}^\infty e^{-2\lambda_n}\cdot (h_n,h)^2< 
\infty\right\}
\end{eqnarray*}
which becomes a Hilbert space with the norm $\|h\|_X=\left(\sum_{n=1 
}^\infty e^{-2\lambda_n}\cdot (h_n,h)^2\right)^{1/2}$, $h\in X$. 
\begin{lemma}\label{Lemma4.2} 
Let $q=2$. Then conditions (jv') and (v') of Section \ref{sec:2} are 
satisfied. 
\end{lemma}
Proof. {\it Step 1 } For the sake of definiteness, let us first verify 
the first part of condition (jv'). We recall that by the initial choice of 
$v=2$, ${\cal V}=L^2(D)$, and $V_1^+(h)=\sum_{j=1}^\infty e^{\lambda_j}(h, 
h_j)h_j$ for any $h\in {\cal V}$, we have $D (V_1^+)=L^2(D)$. Let $\nu 
\equiv\t h\cdot\l\in G$ and note that, because of $G\subset H(1)$, it 
holds that $\t h_{-1}=\sum_{j=1}^\infty e^{-\lambda_j}(\t h,h_j)h_j\in L^2 
(D)$. Choose 
\begin{eqnarray*}
U_{\t h}(0):=\left\{f\in X:\|f\|_X<1\right\}=\left\{\sum_{j=1}^\infty e^{ 
\lambda_j}(g,h_j)h_j:\|g\|_{L^2(D)}<1\right\} 
\end{eqnarray*} 
as a (pre-)neighborhood of $0\in X$ in the topology of $X$. Then  
\begin{eqnarray*} 
&&\hspace{-.5cm}\left\{\t h+i\circ h:h\in U_{\t h}(0)\right\}=\left\{ 
\sum_{j=1}^\infty e^{\lambda_j}\left(\frac{\t h_{-1}}{z(\t h_{-1}\cdot 
\l,1)}+g,h_j\right)h_j:\|g\|_{L^2(D)}<1\right\} \\ 
&&\hspace{.5cm}\subseteq\left\{\sum_{j=1}^\infty e^{\lambda_j}(h,h_j) 
h_j:h\in L^2(D)\right\}=\left\{V_1^+(h):h\in D (V_1^+)\right\}\subseteq 
L^2(D)=D(V_1^+). 
\end{eqnarray*} 
{\it Step 2 } We verify (v'). By definition, $G\subseteq H(1)$ 
which includes that $\nu\in G$ implies $a^f\nu\equiv dA^f\nu/d\, \l  
\in X$. According to (\ref{4.115}) for $\bmu$-a.e. $\nu\in G$ we have 
$\nu_{-1}>0$. By (\ref{4.1}) this provides us with $\nu_s>0$, $s\in 
[-1,1]$, which we will keep in mind throughout this proof. Let $t\in 
[-1,1]$. We obtain for $\nu\in G$ with $\nu_{-1}>0$, i. e. for $\bmu$-a.e. 
$\nu\in G$, 
\begin{eqnarray}\label{4.11*}  
&&\hspace{-.5cm}\left.\frac{d}{du}\right|_{u=0}|\nu+u\cdot A^f\nu|(D)= 
\left.\frac{d}{du}\right|_{u=0}\int_{y\in D}\left|\frac{\nu(dy)}{dy}+ 
u\cdot\frac{A^f\nu(dy)}{dy}\right|\, dy\nonumber \\ 
&&\hspace{.5cm}=\int_{y\in D}\frac{A^f\nu(dy)}{dy}\, dy=0  
\end{eqnarray} 
as well as 
\begin{eqnarray}\label{4.12*} 
&&\hspace{-.5cm}\left.\frac{d}{du}\right|_{u=0}\int_{y\in D}\left| 
\sum_{j=1}^\infty e^{\lambda_jt}(h_j,\nu+u\cdot A^f\nu)h_j(y)\right| 
\, dy\nonumber \\ 
&&\hspace{.5cm}=\left.\frac{d}{du}\right|_{u=0}\int_{y\in D}\left| 
\frac{\nu_t(dy)}{dy}+u\cdot\sum_{j=1}^\infty e^{\lambda_jt}(h_j,A^f 
\nu)h_j(y)\right|\, dy\nonumber \\ 
&&\hspace{.5cm}=\sum_{j=1}^\infty e^{\lambda_jt}(h_j,A^f\nu)(h_j,\1)
\, . 
\end{eqnarray} 

The following equality we get from (\ref{4.22}) if $t\in [-1,0)$ and 
from Remark (1) of this section if $t\in [0,1]$. We obtain     
\begin{eqnarray}\label{4.13*} 
(\nu+u\cdot A^f\nu)_t=\frac{\sum_{j=1}^\infty e^{\lambda_jt}(h_j,\nu 
+u\cdot A^f\nu)h_j\cdot |\nu+u\cdot A^f\nu|(D)}{\int_{y\in D}\left| 
\sum_{j=1}^\infty e^{\lambda_jt}(h_j,\nu+u\cdot A^f\nu)h_j(y)\right| 
\, dy}\, .
\end{eqnarray} 
Keeping in mind (\ref{4.11*})-(\ref{4.13*}) as well as (\ref{4.2}) 
and (\ref{4.3}) the following calculations are straightforward. 
For $g\in C_0(D)$ and $\bmu$-a.e. $\nu\in G$ we have 
\begin{eqnarray*} 
&&\hspace{-.5cm}\left.\frac{d}{du}\right|_{u=0}\left((\nu+u\cdot A^f 
\nu)_t,g\right) \\ 
&&\hspace{.5cm}=\left.\frac{d}{du}\right|_{u=0}\frac{\sum_{j=1 
}^\infty e^{\lambda_jt}(h_j,\nu+u\cdot A^f\nu)(h_j,g)}{\int_{y\in D} 
\left|\sum_{j=1}^\infty e^{\lambda_jt}(h_j,\nu+u\cdot A^f\nu)h_j(y) 
\right|\, dy} \\ 
&&\hspace{.5cm}=\frac{\sum_{j=1}^\infty e^{\lambda_jt}(h_j,A^f\nu)(h_j, 
g)}{\sum_{j=1}^\infty e^{\lambda_jt}(h_j,\nu)(h_j,\1)}-\frac{\sum_{j=1 
}^\infty e^{\lambda_jt}(h_j,\nu)(h_j,g)\cdot\sum_{j=1}^\infty 
e^{\lambda_jt}(h_j,A^f\nu)(h_j,\1)}{\left(\sum_{j=1}^\infty e^{\lambda_jt} 
(h_j,\nu)(h_j,\1)\right)^2} \\ 
&&\hspace{.5cm}=\frac{\sum_{j=1}^\infty (\lambda_j-z'(\nu,0))\, 
e^{\lambda_jt}(h_j,\nu)(h_j,g)}{\sum_{j=1}^\infty e^{\lambda_jt}(h_j, 
\nu)(h_j,\1)} \\ 
&&\hspace{1.0cm}-\frac{\sum_{j=1}^\infty e^{\lambda_jt}(h_j,\nu)(h_j,g) 
\cdot\sum_{j=1}^\infty (\lambda_j-z'(\nu,0))\, e^{\lambda_jt}(h_j,\nu) 
(h_j,\1)}{\left(\sum_{j=1}^\infty e^{\lambda_jt}(h_j,\nu)(h_j,\1)\right 
)^2} \\ 
&&\hspace{.5cm}=\left(\textstyle{\frac12}\Delta \nu_t,g\right)-\left(\nu_t, 
g\right)\cdot\frac{\sum_{j=1}^\infty\lambda_j e^{\lambda_jt}(h_j,\nu)(h_j, 
\1)}{\sum_{j=1}^\infty e^{\lambda_jt}(h_j,\nu)(h_j,\1)} \\ 
&&\hspace{.5cm}=\left(\textstyle{\frac12}\Delta (\nu_{-1})_{1+t},g\right) 
-\left((\nu_{-1})_{1+t},g\right)\cdot\frac{\sum_{j=1}^\infty \lambda_j 
e^{\lambda_j(1+t)}(h_j,\nu_{-1})(h_j,\1)}{\sum_{j=1}^\infty e^{\lambda_j 
(1+t)}(h_j,\nu_{-1})(h_j,\1)} \\ 
&&\hspace{.5cm}=\left(\textstyle{\frac12}\Delta (\nu_{-1})_{1+t},g\right) 
-\left(\frac{z'(\nu_{-1},t+1)}{z(\nu_{-1},t+1)}\cdot(\nu_{-1})_{1+t},g 
\right) \\
&&\hspace{.5cm}=\left(A^f(\nu_{-1})_{1+t},g\right)=\left(A^f\nu_t,g\right) 
\vphantom{\left(\frac12\right)}\, , 
\end{eqnarray*} 
for the interchange of $\left.\frac{d}{du}\right|_{u=0}$ with $\int_{ 
y\in D}$ or $\sum_{j=1}^\infty$ recall the justification after (\ref{4.26}). 
\medskip 

\noindent
{\it Step 3 } We verify the second part of (jv'). Recalling $\nu_s>0$, $s\in 
[-1,1]$, for $\bmu$-a.e. $\nu\in G$ from the beginning of Step 2, for $h\in X$, 
and $g\in C_0(D)$ we obtain  
\begin{eqnarray*}  
&&\hspace{-.5cm}\left.\frac{d}{du}\right|_{u=0}|\nu+u\cdot h\cdot\l|(D)=\left. 
\frac{d}{du}\right|_{u=0}\int_{y\in D}\left|\frac{\nu(dy)}{dy}+u\cdot h\right| 
\, dy=(h,\1)\, . 
\end{eqnarray*} 
Introducing 
\begin{eqnarray*}  
Z(\nu,t):=\left\{ 
\begin{array}{ccc} 
(z(\nu_t,-t))^{-1} &\ {\rm if}\ & t\in [-1,0) \\ 
z(\nu,t) &\ {\rm if}\ & t\in [0,1]
\end{array}
\right.\, ,\quad \nu\in G, 
\end{eqnarray*} 
and keeping in mind (\ref{4.21}) for $t\in [-1,1]$, this leads to 
\begin{eqnarray*} 
&&\hspace{-.5cm}\left.\frac{d}{du}\right|_{u=0}\int_{y\in D}\left| 
\sum_{j=1}^\infty e^{\lambda_jt}(h_j,\nu+u\cdot h\cdot\l)h_j(y)\right| 
\, dy \\ 
&&\hspace{.5cm}=\left.\frac{d}{du}\right|_{u=0}\int_{y\in D}\left| 
Z(\nu,t)\cdot\frac{\nu_t(dy)}{dy}-u\cdot\sum_{j=1}^\infty e^{\lambda_jt} 
(h_j,h)h_j(y)\right|\, dy  \\ 
&&\hspace{.5cm}=\sum_{j=1}^\infty e^{\lambda_jt}(h_j,h)(h_j,\1)\, . 
\end{eqnarray*} 
According to (\ref{4.22}) if $t\in [-1,0)$ and Remark (1) of this section 
if $t\in [0,1]$, this leads to  
\begin{eqnarray*} 
&&\hspace{-.5cm}\left.\frac{d}{du}\right|_{u=0}\left((\nu+u\cdot h\cdot 
\l)_t,g\right) \\ 
&&\hspace{.5cm}=\left.\frac{d}{du}\right|_{u=0}\frac{\sum_{j=1}^\infty 
e^{\lambda_jt}(h_j,\nu+u\cdot h\cdot\l)(h_j,g)\cdot|\nu+u\cdot h\cdot\l| 
(D)}{\int_{y\in D}\left|\sum_{j=1}^\infty e^{\lambda_jt}(h_j,\nu+u\cdot 
h\cdot\l)h_j(y)\right|\, dy} \\ 
&&\hspace{.5cm}=\frac{\sum_{j=1}^\infty e^{\lambda_jt}(h_j,h)(h_j,g)} 
{\sum_{j=1}^\infty e^{\lambda_jt}(h_j,\nu)(h_j,\1)}-(\nu_t,g)\cdot 
\frac{\sum_{j=1}^\infty e^{\lambda_jt}(h_j,h)(h_j,\1)}{\sum_{j=1}^\infty 
e^{\lambda_jt}(h_j,\nu)(h_j,\1)}+(\nu_t,g)\cdot(h,\1)  
\end{eqnarray*} 
for all $t\in [-1,1]$. Furthermore, from (\ref{4.23}) and $|\nu_s|(D)=\nu_s 
(D)=1$, $s\in [-1,1]$, we deduce 
\begin{eqnarray*} 
\sum_{j=1}^\infty e^{\lambda_jt}(h_j,\nu)(h_j,\1)=Z(\nu,t)\, ,\quad t\in 
[-1,1]. 
\end{eqnarray*} 
The last two equalities result in 
\begin{eqnarray*} 
&&\hspace{-.5cm}\left.\frac{d}{du}\right|_{u=0}\left((\nu+u\cdot h\cdot 
\l)_t,g\right) \\ 
&&\hspace{.5cm}=\frac{1}{Z(\nu,t)}\cdot\sum_{j=1}^\infty e^{\lambda_jt}(h_j,h) 
(h_j,g) \\ 
&&\hspace{1.0cm}-\frac{(\nu_t,g)}{Z(\nu,t)}\cdot\sum_{j=1}^\infty e^{\lambda_jt} 
(h_j,h)(h_j,\1)+(\nu_t,g)\cdot(h,\1) \\
&&\hspace{.5cm}\equiv{}_{X^\ast}\langle r(\nu ,t,g),h\rangle_{X}\, ,\quad 
t\in [-1,1].\vphantom{\frac12}
\end{eqnarray*} 
Taking into consideration $X=X^\ast$ it turns out that 
\begin{eqnarray*}  
&&\hspace{-.5cm}r(\nu ,t,g)=\frac{1}{Z(\nu,t)}\cdot\sum_{j=1}^\infty e^{ 
\lambda_j(t+2)}(h_j,g)h_j \\ 
&&\hspace{1.0cm}-\frac{(\nu_t,g)}{Z(\nu,t)}\cdot\sum_{j=1}^\infty e^{\lambda_j 
(t+2)}(h_j,\1)h_j+(\nu_t,g)\cdot\sum_{j=1}^\infty e^{2\lambda_j}(h_j,\1)h_j 
\, ,\quad t\in [-1,1].  
\end{eqnarray*} 
Furthermore $z(\nu_t,-t)\le 1$ $\mu$-a.e. for $t\in [-1,0)$ by definition, and 
$\|(z(\cdot,t))^{-1}\|_{L^\infty(E,\mu)}\le \|(z(\cdot,1))^{-1}\|_{L^\infty(E, 
\mu)}$ for $t\in [0,1]$ according to Lemma \ref{Lemma4.11}. Thus $\|(Z(\cdot, 
t) )^{-1}\|_{L^\infty(E,\mu)}$ is bounded on $t\in [-1,1]$. In addition $|(\nu_t 
,g)|\le\|g\|$ $\mu$-a.e. for all $t\in [-1,1]$ . As a consequence of these facts, 
the map $[-1,1]\ni t\mapsto r(\cdot,t,g)\in L^2(E,\bmu;X^\ast)$ is bounded. 
Continuity of this map follows now from the continuity of 
\begin{eqnarray*} 
[-1,1]\ni t\mapsto (Z(\nu,t))^{-1}=\left(\int\int p(1+t,x,y)\, \nu_{-1}(dx) 
\, dy\right)^{-1} 
\end{eqnarray*} 
as well as 
\begin{eqnarray*}
[-1,1]\ni t\mapsto (\nu_t,g)=\frac{\int\int p(1+t,x,y)\, \nu_{-1}(dx)g(y) 
\, dy}{\int\int p(1+t,x,y)\, \nu_{-1}(dx)\, dy} 
\end{eqnarray*} 
for all $\nu\in G$. 
\qed 
\medskip 

\nid 
{\bf Remark }(4) In the remaining part of the paper we will suppose  
(j)-(jjj) where we, according to the conventions of Subsection 
\ref{sec:4:1}, implicitly assume that $p=q=2$. By Lemma \ref{Lemma4.2}, 
we have the quasi-invariance results of Theorem \ref{Theorem2.3} and 
Theorem \ref{Theorem2.7} (a). If condition (vj) is satisfied we have in 
addition condition (vj') of Theorem \ref{Theorem2.7} (b), Theorem 
\ref{Theorem2.7} (c), and Corollary \ref{Corollary2.8}. 

\subsection{Application: Integration by Parts}\label{sec:4:2} 

Let $p=q=2$ in conditions (j) and (jj) and keep the last remark in mind. 
As an application of Theorem \ref{Theorem2.3}, we establish an integration 
by parts formula relative to the infinitesimal generator of the semigroup 
associated with (\ref{4.1}), cf. Theorem \ref{Theorem4.6} below. 
\begin{lemma}\label{Lemma4.01} 
There exists $c>0$ such that $\left\|\nu\right\|_{H}<c$ for all $\nu 
\in G$ and thus $\left\|\nu\right\|_{H}\in L^\infty(E,\bmu)$. 
\end{lemma}
Proof. We recall $\bmu(G)=1$ from the basic hypotheses of Paragraph 
\ref{sec:1:1:4}. For $\nu=\t h\cdot\l\in G$ with $\nu=W^+(1,\t h_{-1} 
\cdot\l)$ and $\t h_{-1}\cdot\l\in E\cap H(t_0)\subset E\cap H$ it 
holds that 
\begin{eqnarray*} 
&&\hspace{-.5cm}\left\|\nu\right\|^2_{H}=\frac{\sum_{j=1}^\infty 
\lambda_j^2 e^{2\lambda_j}(h_j,\t h_{-1})^2}{\left(z(\t h_{-1}\cdot 
\l,1)\right)^2}\le C\, \frac{\sum_{j=1}^\infty e^{\lambda_j}(h_j, 
\t h_{-1})^2}{\left(z(\t h_{-1}\cdot\l,1)\right)^2} \\ 
&&\hspace{.5cm}=C\, \frac{\int\int p(1,x,y)\t h_{-1}(x)\t h_{-1}(y) 
\, dx\, dy}{\left(\int\int p(1,x,y)\t h_{-1}(x)\, dx\, dy\right)^2} 
\le\frac{Cc_1}{\t c_1^2(h_1,\1)^2} \, ,\quad \nu\in G,
\end{eqnarray*} 
where the constants $0<\t c_1<c_1<\infty$ have been introduced in Step 
1 of the proof of Lemma \ref{Lemma4.0} and $C>0$ is chosen such that 
$\lambda_j^2e^{\lambda_j}\le C$ for all $j\in {\Bbb N}$. 
\qed
\medskip

Let us recall definition (\ref{2.1}). Let $C_c^1({\Bbb R})$ denote the set 
of all $f\in C^1({\Bbb R})$ which have compact support. Furthermore define $K$ 
to be the set of all non-negative $k\in C(D)\cap L^2(D)$ such that $\lim_{D\ni 
x\to y}k(x)=\infty$, $y\in\partial D$. Introduce 
\begin{eqnarray}\label{4.4}
\t C_{b,c}^1(E):=\left\{f=f_1\cdot\vp_0((k,\cdot)):f_1\in\t C_b^1 (E),\ 
\vp_0\in C_c^1({\Bbb R}),\ k\in K\right\}\, .  
\end{eqnarray} 
\begin{lemma}\label{Lemma4.1} 
We have $z'(\cdot,t)\in L^\infty(E,\bmu)$ and $\|z'(\cdot,t)\|_{L^\infty 
(E,\sbmu)}$ is uniformly bounded on $t\in[0,1]$. Furthermore, we have $(h_n, 
\cdot)\in L^\infty(E,\bmu)$, $n\in {\Bbb N}$. In addition, for $k\in K$ it 
holds that $\left(k, {\T\frac12}\Delta\, \cdot\, \right)\in L^\infty(E,\bmu)$ 
and $(k,\cdot)\in L^\infty(E,\bmu)$. 
\end{lemma}
Proof. According to Lemma \ref{Lemma4.01} we have 
\begin{eqnarray}\label{4.5} 
\sum_{n=1}^\infty\lambda_n^2(h_n,\nu)^2=\|\nu\|^2_H<c\, ,\quad\nu 
\in G.  
\end{eqnarray}
Furthermore,  
\begin{eqnarray*} 
\sum_{n=1}^\infty e^{2\lambda_nt}(h_n,\1)^2=\left\|\int p(t,\cdot ,y) 
\, dy\right\|^2_{L^2(D)}\le\left\|\1\right\|^2_{L^2(D)}=\l(D)\, .  
\end{eqnarray*}
For $\nu\in G$ and $t\in [0,1]$ from the last two relations we obtain 
\begin{eqnarray*}
&&\hspace{-.5cm}\left(z'(\nu,t)\right)^2=\left(\sum_{n=1}^\infty\lambda_n 
e^{\lambda_nt}(h_n,\nu)(h_n,\1)\right)^2\le\sum_{n=1}^\infty\lambda_n^2 
(h_n,\nu)^2\cdot\sum_{n=1}^\infty e^{2\lambda_nt}(h_n,\1)^2\le c\, \l(D)\, . 
\end{eqnarray*}
This says that $\|z'(\cdot,t))\|_{L^\infty(E,\sbmu)}$ is uniformly bounded 
for $t\in[0,1]$. Moreover, it holds that $\left|(h_n,\nu)\right|\le\|h_n\|$, 
$\nu\in G$, which immediately yields $(h_n,\cdot)\in L^\infty(E,\bmu)$, 
$n\in {\Bbb N}$. For $k\in K$ and $\nu=h\cdot\l\in G$ we find 
\begin{eqnarray*}
\left(k,{\T\frac12}\Delta h\right)^2\le\|k\|_{L^2(D)}^2\left\|{\T\frac12} 
\Delta h\right\|_{L^2(D)}^2=\|k\|^2_{L^2(D)}\sum_{n=1}^\infty\lambda_n^2 
(h_n,h)^2 
\end{eqnarray*}
as well as 
\begin{eqnarray*}
(k,\nu)^2\le\|k\|^2_{L^2(D)}\|h\|^2_{L^2(D)}=\|k\|^2_{L^2(D)}\sum_{n=1
}^\infty (h_n,h)^2\, . 
\end{eqnarray*}
Again (\ref{4.5}) shows $\left(k,{\T\frac12}\Delta\, \cdot\, \right) 
\in L^\infty(E,\bmu)$ as well as $(k,\cdot)\in L^\infty(E,\bmu)$. 
\qed 
\medskip 

\nid 
{\bf Remarks} (5) In order to take over the analysis of Section 
\ref{sec:2} we mention that the closed linear span with respect to the 
sup-norm of $h_1,h_2,\ldots\, $ is $C_0(D)$. This follows from 
\begin{eqnarray*}
\sum_{n=1}^\infty e^{\lambda_nt}(h_n,f)h_n=\int_{y\in D} p(t,\cdot,y)f(y) 
\, dy\stack {t\to 0}{\lra}f\, ,\quad f\in C_0(D), 
\end{eqnarray*} 
with respect to the $\sup$-norm, $\sum_{n=1}^\infty(h_n,f)^2<\infty$, 
$\liminf_{n\to\infty}-\lambda_n/n^{2/d}>0$, $n\in {\Bbb N}$, see 
\cite{Va92}, and $\|h_n\|\le (-2\lambda_n)^{d/4}$, $n\in {\Bbb N}$, 
see \cite{Da73} and \cite{BB99}.
\medskip 

\nid 
(6) In this remark, we will use the notation of (\ref{4.4}). Both sets, 
$\t C_b^1(E)$ and $\t C_{b,c}^1(E)$, are dense in $L^2(E,\bmu)$. To see 
this we introduce the set $E^c$ of all probability measures on $\overline 
{D}$ and endow $E^c$ with the Prokhorov metric $\pi$. We note that $(E^c, 
\pi)$ is a compact space. It follows from the Stone-Weierstrass theorem 
and Remark (1) of Section \ref{sec:2} that 
\begin{eqnarray*}
\t C_b^1(E^c)&&\hspace{-.5cm}:=\left\{f(\nu)=\vp\left((h_1,\nu),\ldots, 
(h_r,\nu),{\T\int_{y\in\partial D}g(y)\, \nu(dy)}\right),\ \nu\in E^c: 
\right.\\ 
&&\hspace{0.7cm}\left.\vphantom{\t f}g\in C(\partial D),\ \vp\in C_b^1 
({\Bbb R}^{r+1}),\ r\in {\Bbb N}\right\} 
\end{eqnarray*}
is dense in the space $C(E^c)$ of all continuous functions on $E^c$.  
Recalling the definition of $K$ above (\ref{4.4}), we define 
\begin{eqnarray*}
\t C_{b,c}^1(E^c)&&\hspace{-.5cm}:=\left\{f(\nu)=f_1(\nu)\cdot\vp_0 
\left({\T\int_D k\, d\nu}\right),\ \nu\in E^c:\vphantom{\t f}\right.\\ 
&&\hspace{0.7cm}\left.\vphantom{\t f}f_1\in \t C_b^1(E^c),\ \vp_0\in 
C_c^1({\Bbb R}),\ k\in K\right\}\, . 
\end{eqnarray*}
According to $\bmu(G)=1$ and Lemma \ref{Lemma4.0} (b) we have $\nu=h 
\cdot\l$ for some $h\in L^2(D)$. This gives $\int_D k\, d\nu<\infty$ 
$\bmu$-a.e. We extend the measure $\bmu$ to $(E^c,{\cal B}(E^c))$ by 
setting $\bmu(E^c\setminus E)=0$ and obtain that $\t C_b^1(E^c)$ is 
dense in $L^2(E^c,\bmu)$, which immediately implies that $\t C_{b,c}^1 
(E^c)$ is also dense in $L^2(E^c,\bmu)$. What we have claimed in the 
beginning of this remark follows now by restriction to $E$. 
\medskip 

We recall that the superscript $f$ in $A^f$ indicates the generator of the 
flow associated with (\ref{4.1}). Below, in contrast, $A$ will denote the 
infinitesimal generator of the semigroup associated with (\ref{4.1}). In the 
remainder of this section, let us use $\langle\cdot\, ,\cdot\, \rangle_{L^2}$ 
for the inner product in $L^2(E,\bmu)$ and $\|\cdot\|_{L^p}$ for $\|\cdot 
\|_{L^p(E,\sbmu)}$, $1\le p\le\infty$. 
\medskip 

\nid
{\bf Remark} (7) In (\ref{4.16*}) and several other places below for $t\in 
[0,1]$ we need to set a bound in the norm of $L^\infty(E,\bmu)$ on 
\begin{eqnarray*} 
\left(\frac{d\bmu\circ\nu_t}{d\bmu}\right)^{-1}=r_t^{-1}=\exp\left\{-\int_{ 
s=0}^t\delta (A^f)(\nu_s)\, ds\right\}\, . 
\end{eqnarray*} 
Even more sophisticated, in order to prove Lemma \ref{Lemma4.4} below, we 
need a bound in the norm of $L^\infty(E,\bmu)$ on $\delta (A^f)$. Instead of 
(jjj), for the rest of this section, we suppose 
\begin{itemize}
\item[(jjj'')] 
\begin{eqnarray*} 
\delta (A^f)\in L^\infty(E,\bmu)\, .  
\end{eqnarray*} 
\end{itemize} 
This implies 
\begin{eqnarray}\label{4.24} 
\|r_t^{-1}\|_{L^\infty}\le e^{t\|\delta (A^f)\|_{L^\infty}}<\infty\, . 
\end{eqnarray} 

For the next assertions, Proposition \ref{Proposition4.3} and Lemma 
\ref{Lemma4.4} below, let us recall $\bmu(G)=1$ by Paragraph \ref{sec:1:1:5} 
and that $\nu\in G$ by Lemma \ref{Lemma4.0} (b) implies $\nu\in E\cap H(1+ 
t_0)\subset E\cap H(1)$. The latter yields $\bmu(E\cap H(1))=1$. Moreover,  
let us recall the formal definition of $\nu_t$, $t\ge 0$, in Paragraphs 
\ref{sec:1:1:3} and \ref{sec:1:1:5}. The representation of $\nu_t$ is then 
determined by $\nu_t=v(t,\cdot)\cdot\l$, $t\ge 0$, where $v(t,\cdot)$ is 
given by (\ref{4.1}). 
\begin{proposition}\label{Proposition4.3} 
Assume (j), (jj), and (jjj''). (a) Then $\nu_t$, $\nu\in E\cap H(1)$, $t\ge 
0$ is associated with a strongly continuous semigroup $(T_t)_{t\ge 0}$ of 
bounded linear operators on $L^2(E,\bmu)$ given by 
\begin{eqnarray*}
T_tf(\nu)=f(\nu_t)\, ,\quad \nu\in E\cap H,\ f\in L^2(E,\bmu). 
\end{eqnarray*}
(b) Let $A$ denote its infinitesimal generator. We have $\t C_b^1(E)\cup\t 
C_{b,c}^1(E)\subseteq D(A)$. If $f\in\t C_b^1(E)$ then with $f$ and $\vp$ 
related as in (\ref{2.1}), $\bmu$-a.e. 
\begin{eqnarray*}
Af=\sum_{i=1}^r\frac{\partial\vp}{\partial x_i}\cdot (h_i,\cdot)\left( 
\lambda_i-z'(\cdot ,0)\right)\, .  
\end{eqnarray*}
If $f\in\t C_{b,c}^1(E)$ then with $f$ and $\vp$ as well as $\vp_0$ related 
as in (\ref{4.4}), $\bmu$-a.e. 
\begin{eqnarray*}
Af &&\hspace{-.5cm}=\sum_{i=1}^r\frac{\partial\vp}{\partial x_i}\vp_0\cdot 
\lambda_i(h_i,\cdot)+\vp\vp'_0\cdot\left(k,{\T\frac12}\Delta\, \cdot\, 
\right) \\ 
&&\hspace{4.0cm}-z'(\cdot ,0)\cdot\left(\sum_{i=1}^r\frac{\partial\vp} 
{\partial x_i}\vp_0\cdot (h_i,\cdot)+\vp\vp'_0\cdot(k,\cdot)\right)\, . 
\end{eqnarray*}
\end{proposition} 
Proof. {\it Step 1 } (a) We show that $(T_t)_{t\ge 0}$ forms a strongly 
continuous semigroup of bounded linear operators in $L^2(E,\bmu)$. The 
semigroup property of $(T_t)_{t\ge 0}$, that is $T_sT_tf(\nu)=T_sf(\nu_t) 
=f(\nu_{s+t})=T_{s+t}f(\nu)$, $f\in L^2(E,\bmu)$, $\nu\in E\cap H(1)$, 
$s,t\ge 0$, is a consequence of (\ref{4.1}). Linearity of $(T_t)_{t\ge 
0}$ is obvious. 

According to Theorem \ref{Theorem2.3} (a) and condition (jjj'') we have  
\begin{eqnarray}\label{4.16*}
&&\hspace{-.5cm}\|T_tf\|_{L^2}^2=\int(f(\nu_t))^2\, d\bmu=\int_G(f(\nu_t))^2 
\left(\frac{d\bmu\circ\nu_t}{d\bmu}(\nu)\right)^{-1}\, \bmu\circ\nu_t(d\nu) 
\nonumber \\  
&&\hspace{.5cm}=\int_G(f(\nu_t))^2(r_t(\nu))^{-1}\, d\bmu(\nu_t)\le\|r_t^{-1} 
\|_{L^\infty}\int_G(f(\nu_t))^2\, d\bmu(\nu_t)\nonumber \\ 
&&\hspace{.5cm}=\|r_t^{-1}\|_{L^\infty}\int_{G_t}(f(\nu))^2\, \bmu(d\nu) 
\le e^{t\|\delta (A^f)\|_{L^\infty}}\|f\|^2_{L^2}<\infty, 
\end{eqnarray} 
i. e. $T_tf\in L^2(E,\bmu)$ whenever $f\in L^2(E,\bmu)$, $t\ge 0$. In addition 
we get boundedness of $T_t$ in the operator norm for all $t\ge 0$. In other words, 
$(T_t)_{t\ge 0}$ forms a semigroup of bounded linear operators in $L^2(E,\bmu)$. 

Let us verify that $(T_t)_{t\ge 0}$ is strongly continuous. For this we recall 
(\ref{4.24}). By Remark (6) of this section we can choose $g\in\t C_b^1(E)$ such 
that 
\begin{eqnarray*} 
\left\|f-g\right\|_{L^2}<\frac{\ve}{2}\cdot\left(\sup_{t\in [0,1]}\left 
\|r_t^{-1}\right\|^{\frac12}_{L^\infty}+1\right)^{-1}\, . 
\end{eqnarray*} 
As a consequence of hypothesis (a1), cf. Lemma \ref{Lemma4.0} (d), $[0,1] 
\ni t\mapsto\nu_t$ is $\bmu$-a.e. continuous in $(E,\pi)$. Thus, by $T_tg 
(\nu)=g(\nu_t)$ we may also choose $t\in (0,1]$ such that $\left\|T_tg- 
g\right\|_{L^2}<\ve/2$. Now, as in (\ref{4.16*}), 
\begin{eqnarray*} 
&&\hspace{-.5cm}\left\|T_tf-f\right\|_{L^2}\le\left(\int\left(f(\nu_t)- 
g(\nu_t)\right)^2(r_t(\nu))^{-1}\, \bmu(d\nu_t)\right)^{\frac12}+\left\| 
T_tg-g\right\|_{L^2}+\left\|g-f\right\|_{L^2} \\
&&\hspace{.5cm}\le\left\|f-g\right\|_{L^2}\cdot\left(\left\|r_t^{-1}\right 
\|^{\frac12}_{L^\infty}+1\right)+\left\|T_tg-g\right\|_{L^2}<\ve\, . 
\end{eqnarray*} 
We have shown that $(T_t)_{t\ge 0}$ is strongly continuous.
\medskip

\nid 
{\it Step 2 } (b) For $\t C_b^1(E)\subseteq D(A)$ and the claimed 
representation of $Af$, $f\in\t C_b^1(E)$, it is sufficient to demonstrate 
that the  derivative $\frac{d}{dt}f(\nu_t)$ exists in $L^2(E,\bmu)$, 
$t\in [0,1]$, and that for the pointwise derivative
\begin{eqnarray}\label{4.6}
\left.\frac{d}{dt}\right|_{t=0}f(\nu_t)=\sum_{i=1}^r\frac{\partial\vp 
(\nu)}{\partial x_i}\cdot (h_i,\nu)\left(\lambda_i-z'(\nu,0)\right)\, , 
\quad\nu\in E\cap H(1),\ \mbox{\rm hence }\bmu\mbox{\rm -a.e. } 
\end{eqnarray}

The representation of $f\in\t C_b^1(E)$ as given in (\ref{2.1}) implies 
for the pointwise derivative
\begin{eqnarray}\label{4.8} 
\frac{d}{dt}f(\nu_t)=\sum_{i=1}^r\frac{\partial\vp(\nu)}{\partial x_i} 
\cdot\frac{d}{dt}(h_i,\nu_t)\, ,\quad\nu\in E\cap H(1),\ \mbox{\rm hence } 
\bmu\mbox{\rm -a.e.},\ t\in [0,1]. 
\end{eqnarray}
Recalling Remarks (1), (2) of the present section and (\ref{4.26}), 
(\ref{4.3}) we obtain relation (\ref{4.6}). By the same references we get 
for the pointwise derivative
\begin{eqnarray}\label{4.10}
\frac{d}{dt}(h_j,\nu_t)&&\hspace{-.5cm}=\left(h_j,\frac{d\nu_t}{dt} 
\right)\nonumber \\ 
&&\hspace{-.5cm}=(h_j,\nu_t)\left(\lambda_j-\sum_{i=1}^\infty\lambda_i 
(h_i,\nu_t)(h_i,\1)\right)\nonumber \\ 
&&\hspace{-.5cm}=(h_j,\nu_t)\left(\lambda_j-\frac{z'(\nu,t)}{z(\nu,t)} 
\right)\, ,\quad\nu\in E\cap H(1),\ \mbox{\rm hence }\bmu\mbox{\rm -a.e.}, 
\ t\in [0,1]. \qquad
\end{eqnarray}
Here we take into consideration that $(h_j,\nu_t)\in L^\infty(E,\bmu)$, 
$z'(\cdot,t)\in L^\infty(E,\bmu)$, and $(z(\cdot,t))^{-1}\in L^\infty 
(E,\bmu)$, all three uniformly in $t\in [0,1]$ by Lemma \ref{Lemma4.11} 
and Lemma \ref{Lemma4.1}. Therefore, the derivative (\ref{4.10}) exists 
in $L^2(E,\bmu)$ which implies that the derivative (\ref{4.8}) also exists 
in $L^2(E,\bmu)$. 
\medskip 

\nid 
{\it Step 3 } For the representation of $Af$, $f\in\t C_{b,c}^1(E)$, 
we proceed similar to Step 2. We use again Lemma \ref{Lemma4.11} and 
Lemma \ref{Lemma4.1} to show that 
\begin{eqnarray*}
\frac{d}{dt}(k,\nu_t)&&\hspace{-.5cm}=\sum_{j=1}^\infty(h_j,k)\left( 
h_j,\frac{d\nu_t}{dt}\right) \\ 
&&\hspace{-.5cm}=\sum_{j=1}^\infty(h_j,k)(h_j,\nu_t)\left(\lambda_j- 
\frac{z'(\nu,t)}{z(\nu,t)}\right) \\ 
&&\hspace{-.5cm}=\left(k,{\T\frac12}\Delta\nu_t\right)-(k,\nu_t)\cdot 
\frac{z'(\nu,t)}{z(\nu,t)}
\end{eqnarray*}
in $L^2(E,\bmu)$. 
\qed 
\medskip 

\nid
{\bf Remark} (8) Part (a) of Proposition \ref{Proposition4.3} is used 
in the analysis below. In addition, the particular representation of 
$Af$ in part (b) is of independent interest since it contains information 
on the evolution $T_t f(\nu)=\nu_t$ for probability measures $\nu=h\cdot 
\l$ near the boundary $\partial D$. 
\medskip 

Let $(T^\ast_t)_{t\ge 0}$ denote the adjoint semigroup to $(T_t)_{t\ge 0}$ in 
$L^2(E,\bmu)$. Introduce 
\begin{eqnarray*} 
D(\t A):=\left\{f\in L^2(E,\bmu):\lim_{t\downarrow 0}\frac{T^\ast_tf-f\cdot 
\chi_G(\nu_{-t})}{t}\ \ \mbox{\rm exists in } L^2(E,\bmu)\right\} 
\end{eqnarray*} 
and 
\begin{eqnarray*} 
\t Af:=\lim_{t\downarrow 0}\frac{T^\ast_tf-f\cdot\chi_G(\nu_{-t})}{t}\, ,\quad 
f\in D(\t A). 
\end{eqnarray*} 
For the remainder of this section we need one more regularity assumption 
on the measure $\bmu$.  
\begin{itemize}
\item[(vjj)] For $\bmu$-a.e. $\nu\in E$ the function 
\begin{eqnarray*} 
[0,1]\ni t\mapsto\delta(A^f)(\nu_t)\quad\mbox{\rm is right continuous.}
\end{eqnarray*} 
\end{itemize} 
\begin{lemma}\label{Lemma4.4}
Assume (j),(jj), (jjj''), and (vjj). We have $D(A)=D(\t A)$ and $\t Af=-Af+ 
\t A\1\cdot f$, $f\in\t D(A)$. In particular, we have $\t A\1=-\delta(A^f)$. 
\end{lemma}
Proof. {\it Step 1 } Let us examine the adjoint semigroup. For this we recall 
relation (\ref{2.3}) and the calculations of (\ref{4.16*}). We obtain for 
$f,g\in L^2(E,\bmu)$, and $t\in [0,1]$ 
\begin{eqnarray*}  
&&\hspace{-.5cm}\langle g,T^\ast_tf\rangle_{L^2}=\langle T_tg,f 
\rangle_{L^2} \\ 
&&\hspace{.5cm}=\int g(\nu_t)f(\nu)\, \bmu (d\nu)=\int g(\nu_t)f(\nu) 
(r_t(\nu))^{-1}\, \bmu (d\nu_t) \\ 
&&\hspace{.5cm}=\int_{G_t} g(\nu)f(\nu_{-t})\exp\left\{-\int_{s=0}^t 
\delta(A^f)(\nu_s)\, ds\right\}(\nu_{-t})\, \bmu (d\nu)\, . 
\end{eqnarray*} 
Taking into consideration (\ref{4.24}) this yields 
\begin{eqnarray}\label{4.20} 
T^\ast_t\1=\exp\left\{-\int_{s=0}^t\delta (A^f)(\nu_{s})\, ds\right\} 
(\nu_{-t})\cdot\chi_G(\nu_{-t})\in L^\infty(E,\bmu)   
\end{eqnarray} 
and
\begin{eqnarray*}
T^\ast_tf(\nu)=T^\ast_t\1(\nu)\cdot f(\nu_{-t})\quad{\rm for }\quad\bmu 
\mbox{\rm -a.e. }\nu\in E. 
\end{eqnarray*} 
Recalling again Remark (7) of this section as well as (\ref{4.24}), 
(jjj''), and (vjj) we get  
\begin{eqnarray}\label{4.21*}
&&\hspace{-.5cm}\int\left(\frac{T^\ast_t\1(\nu)-\1\cdot\chi_G(\nu_{-t})} 
{t}+\delta (A^f)\right)^2\, d\bmu\nonumber \\ 
&&\hspace{.5cm}=\int_{G_{-t}}\frac1{t^2}\left(\exp\left\{-\int_{s= 
0}^t\delta (A^f)(\nu_s)\, ds\right\}-1\right)^2\cdot\chi_G\, \bmu(d 
\nu_t)\nonumber \\ 
&&\hspace{1.0cm}+2\int_{G_{-t}}\frac1t\left(\exp\left\{-\int_{s=0}^t 
\delta (A^f)(\nu_s)\, ds\right\}-1\right)\cdot\chi_G\, \delta (A^f) 
(\nu_t)\, \bmu(d\nu_t)\nonumber \\ 
&&\hspace{1.0cm}+\left\|\delta (A^f)\right\|^2_{L^2}\vphantom{\int} 
\nonumber \\ 
&&\hspace{.5cm}=\int_G\frac1{t}\left(\exp\left\{-\int_{s=0}^t\delta 
(A^f)(\nu_s)\, ds\right\}-1\right)\times\nonumber \\ 
&&\hspace{1.5cm}\times\frac1{t}\left(1-\exp\left\{\int_{s=0}^t\delta 
(A^f)(\nu_s)\, ds\right\}\right)\, \bmu(d\nu)\nonumber \\ 
&&\hspace{1.0cm}+2\int_G\frac1t\left(1-\exp\left\{\int_{s=0}^t\delta 
(A^f)(\nu_s)\, ds\right\}\right)\, \delta (A^f)(\nu_t)\, \bmu(d\nu)+ 
\left\|\delta (A^f)\right\|^2_{L^2}\nonumber \\ 
&&\hspace{-.0cm}\stack{t\downarrow 0}{\lra}0\, .\vphantom{\dot{f}}
\end{eqnarray}
We obtain $\t A\1=-\delta(A^f)$. 
\medskip 

\nid
{\it Step 2 } In this step we prove a weak version of the claim.  
Fix an arbitrary $\t a\in L^2(E,\bmu)$ and let $f\in L^2(E,\bmu)$. 
It holds that 
\begin{eqnarray}\label{4.22*} 
&&\hspace{-.5cm}\frac{T^\ast_tf-f\cdot\chi_G(\nu_{-t})}{t}-\t a+\delta (A^f) 
\cdot f=\frac{T^\ast_t\1\cdot f\circ\nu_{-t}-f\cdot\chi_G(\nu_{-t})}{t}-\t a 
+\delta (A^f)\cdot f \nonumber \\ 
&&\hspace{.5cm}=\left(\frac{f\circ\nu_{-t}-f}{t}\cdot\chi_G(\nu_{-t})-\t a 
\right)+\left((f\circ\nu_{-t}-f)\cdot\frac{T^\ast_t\1-\1\cdot\chi_G(\nu_{-t})} 
{t}\right)\nonumber \\ 
&&\hspace{1.0cm}+\left(\frac{T^\ast_t\1-\1\cdot\chi_G(\nu_{-t})}{t}+\delta 
(A^f)\right)\cdot f\, .
\end{eqnarray} 
Let $\vp\in L^\infty(E,\bmu)$ be a test function. We have  
\begin{eqnarray*}
&&\hspace{-.5cm}\left\langle\left(\frac{T^\ast_t\1-\1\cdot\chi_G(\nu_{-t})} 
{t}+\delta (A^f)\right)\cdot f\, ,\, \vp\right\rangle_{L^2}=\left\langle 
\frac{T^\ast_t\1-\1\cdot\chi_G(\nu_{-t})}{t}+\delta (A^f)\, ,\, \vp\cdot f 
\right\rangle_{L^2} \\ 
&&\hspace{-.0cm}\stack{t\downarrow 0}{\lra}0\vphantom{\dot{f}}
\end{eqnarray*} 
according to (\ref{4.21*}). By 
\begin{eqnarray*}
\left\|\frac{T^\ast_t\1-\1\cdot\chi_G(\nu_{-t})}{t}+\delta (A^f)\right\|_{ 
L^\infty}\le \|\delta(A^f)\|_{L^\infty}\left(e^{\|\delta(A^f)\|_{L^\infty}} 
+1\right)\, ,\quad t\in [0,1], 
\end{eqnarray*} 
the last convergence holds for all $\vp\in L^2(E,\bmu)$, i. e. we have weak 
convergence 
\begin{eqnarray}\label{4.23*}
\left(\frac{T^\ast_t\1-\1\cdot\chi_G(\nu_{-t})}{t}+\delta (A^f)\right)\cdot f 
\rightharpoonup 0\quad \mbox{\rm in}\ L^2(E,\bmu). 
\end{eqnarray} 

Let us use the abbreviation 
\begin{eqnarray*}
\Delta_tr_t^{-1}:=\frac{\exp\left\{-\int_{s=0}^t\delta (A^f)(\nu_{s})\, ds 
\right\}-\1}{t}\, . 
\end{eqnarray*} 
Recalling that $T^\ast_t\1\cdot\chi_G(\nu_{-t})=T^\ast_t\1$ by 
(\ref{4.20}), we obtain from Proposition \ref{Proposition4.3} (a), 
Theorem \ref{Theorem2.3} (a), and (jjj'') for $t\in [0,1]$
\begin{eqnarray*} 
&&\hspace{-.5cm}\left|\left\langle(f\circ\nu_{-t}\cdot\chi_G(\nu_{-t}) 
-f)\cdot\frac{T^\ast_t\1-\1\cdot\chi_G(\nu_{-t})}{t}\, ,\, \vp\right 
\rangle_{L^2}\right|\nonumber \\ 
&&\hspace{0.5cm}=\left|\int_G\left(f-f(\nu_t)\vphantom{l^1}\right)\cdot 
\Delta_tr_t^{-1}\cdot\vp(\nu_t)\, \bmu(d\nu_t)\right|\nonumber \\ 
&&\hspace{.5cm}\le\left\|f-f(\nu_t)\right\|_{L^2}\cdot\left\|\Delta_t 
r_t^{-1}\cdot\vp(\nu_t)\cdot r_t\right\|_{L^2}\vphantom{\int}\nonumber 
 \\ 
&&\hspace{.5cm}\le\left\|f-f(\nu_t)\right\|_{L^2}\cdot\|\delta(A^f)\|_{ 
L^\infty}e^{\|\delta(A^f)\|_{L^\infty}}\cdot\|\vp\|_{L^\infty}\cdot e^{ 
\|\delta(A^f)\|_{L^\infty}}\vphantom{\int}\nonumber \\ 
&&\hspace{-.0cm}\stack{t\downarrow 0}{\lra}0\, ,\quad\vp\in L^\infty(E, 
\bmu). \vphantom{\dot{f}}
\end{eqnarray*} 
Recalling again (jjj''), Theorem \ref{Theorem2.3} (a), and (\ref{4.16*}), 
we can estimate 
\begin{eqnarray*}
&&\hspace{-.5cm}\left\|(f\circ\nu_{-t}\cdot\chi_G(\nu_{-t})-f)\cdot\frac 
{T^\ast_t\1-\1\cdot\chi_G(\nu_{-t})}{t}\right\|_{L^2} \\ 
&&\hspace{.5cm}\le\|\delta(A^f)\|_{L^\infty}e^{\|\delta(A^f)\|_{L^\infty}} 
\cdot\left\|f\circ\nu_{-t}\cdot\chi_G(\nu_{-t})-f\right\|_{L^2}\vphantom 
{\int} \\ 
&&\hspace{.5cm}\le\|\delta(A^f)\|_{L^\infty}e^{\|\delta(A^f)\|_{L^\infty}} 
\cdot\left\|r_t^{\frac12}\right\|_{L^\infty}\cdot\left\|(f-f(\nu_t))\right 
\|_{L^2}\vphantom{\int} \\ 
&&\hspace{.5cm}\le\|\delta(A^f)\|_{L^\infty}e^{\frac32\|\delta(A^f)\|_{ 
L^\infty}}\cdot\left(\|f\|_{L^2}+\|f(\nu_t)\|_{L^2}\right)\vphantom{\int}
 \\ 
&&\hspace{.5cm}\le\|\delta(A^f)\|_{L^\infty}e^{\frac32\|\delta(A^f)\|_{ 
L^\infty}}\cdot\|f\|_{L^2}\left(1+e^{\frac12\|\delta(A^f)\|_{L^\infty}} 
\right)\vphantom{\int}\, ,\quad t\in [0,1], 
\end{eqnarray*} 
which implies the last convergence holds for all $\vp\in L^2(E,\bmu)$. 
Using again the identity $T^\ast_t\1\cdot\chi_G(\nu_{-t})=T^\ast_t\1$, 
cf. (\ref{4.20}), we have weak convergence 
\begin{eqnarray}\label{4.24*} 
&&\hspace{-.5cm}(f\circ\nu_{-t}-f)\cdot\frac{T^\ast_t\1-\1\cdot\chi_G 
(\nu_{-t})}{t}\nonumber \\ 
&&\hspace{.5cm}=(f\circ\nu_{-t}\cdot\chi_G(\nu_{-t})-f)\cdot\frac{ 
T^\ast_t\1-\1\cdot\chi_G(\nu_{-t})}{t}\rightharpoonup 0\quad \mbox{\rm 
in}\ L^2(E,\bmu).   
\end{eqnarray} 
Now observe that by (jjj'') and (vjj) we have $r_t\stack{t\downarrow 0} 
{\lra}\1$ in $L^\infty(E,\bmu)$. Provided the following limit exists, 
we get in addition 
\begin{eqnarray}\label{4.25*}
&&\hspace{-.5cm}\lim_{t\downarrow 0}\int_G\left(\frac{f\circ\nu_{-t}-f} 
{t}\cdot\chi_G(\nu_{-t})-\t a\right)\vp\, d\bmu\nonumber \\ 
&&\hspace{.5cm}=\lim_{t\downarrow 0}\int_G\left(\frac{f-f(\nu_t)}{t}-\t 
a(\nu_t)\right)\vp(\nu_t)\, r_t\, d\bmu=\lim_{t\downarrow 0}\int_G\left( 
\frac{f-f(\nu_t)}{t}-\t a(\nu_t)\right)\vp(\nu_t)\, d\bmu\nonumber \\ 
&&\hspace{.5cm}=\lim_{t\downarrow 0}\int_G\left(\frac{f-f(\nu_t)}{t}-\t 
a\right)\vp\, d\bmu\, ,\quad\vp\in L^2(E,\bmu),
\end{eqnarray} 
the second last equality sign because of Theorem \ref{Theorem2.3} (a) 
as well as (jjj''), and the last equality sign because of Proposition 
\ref{Proposition4.3} (a). 

As a consequence of (\ref{4.22*})-(\ref{4.25*}), 
\begin{eqnarray}\label{4.26*}
\frac{T^\ast_tf-f\cdot\chi_G(\nu_{-t})}{t}\rightharpoonup\t a-\delta 
(A^f)\cdot f\quad\mbox{\rm weakly in}\ L^2(E,\bmu) 
\end{eqnarray} 
if and only if 
\begin{eqnarray}\label{4.27*}
\frac{T_tf-f}{t}\rightharpoonup -\t a\quad\mbox{\rm weakly in}\ L^2 
(E,\bmu).  
\end{eqnarray} 
{\it Step 3 } From (jjj'') and (vjj) it follows that $r_t\stack{t 
\downarrow 0}{\lra}\1$ in $L^\infty(E,\bmu)$ and 
\begin{eqnarray*} 
\lim_{t\downarrow 0}\left\|\frac1t\left(\exp\left\{-\int_{s=0}^t\delta 
(A^f)(\nu_s)\, ds\right\}-\1\right)\cdot f+\delta (A^f)\cdot f\right 
\|_{L^2}=0\, .  
\end{eqnarray*} 
Keeping this in mind and taking into consideration that $T^\ast_t\1\cdot 
\chi_G(\nu_{-t})=T^\ast_t\1$ by (\ref{4.20}) we obtain  
\begin{eqnarray}\label{4.28*}
&&\hspace{-.5cm}\left\|\frac{T^\ast_tf-f\cdot\chi_G(\nu_{-t})}{t}\right 
\|^2_{L^2}=\int\left(\frac{T^\ast_t\1(\nu)\cdot f(\nu_{-t})-f}{t}\cdot 
\chi_G(\nu_{-t})\right)^2\, d\bmu\nonumber \\ 
&&\hspace{.5cm}=\int_G\frac1{t^2}\left(\exp\left\{-\int_{s=0}^t\delta 
(A^f)(\nu_s)\, ds\right\}\cdot f-f(\nu_t)\right)^2\, \bmu(d\nu_t) 
\nonumber \\ 
&&\hspace{.5cm}=\int_G\frac1{t^2}\left(\left(\exp\left\{-\int_{s=0}^t 
\delta (A^f)(\nu_s)\, ds\right\}-1\right)\cdot f-(f(\nu_t)-f)\right)^2 
r_t\, d\bmu\nonumber \\ 
&&\hspace{-.0cm}\stack{t\downarrow 0}{\lra}\|-\delta (A^f)\cdot f-Af 
\|^2_{L^2}\quad\mbox{\rm if and only if}\quad f\in D(A)\, .\vphantom 
{\int}
\end{eqnarray}
Let $f\in L^2(E,\bmu)$. Assuming $f\in D(A)$, relation (\ref{4.28*}) and 
the implication from (\ref{4.27*}) to (\ref{4.26*}) says $f\in D(\t A)$ 
and 
\begin{eqnarray}\label{4.29*}
\t Af=\lim_{t\downarrow 0}\frac{T^\ast_tf-f\cdot\chi_G(\nu_{-t})}{t}= 
-Af-\delta (A^f)\cdot f  
\end{eqnarray} 
in $L^2(E;\bmu)$. Next we remind of the fact that the weak infinitesimal 
generator of a strongly continuous semigroup of bounded linear operators 
coincides with its strong infinitesimal generator, see \cite{Pa83}, 
Theorem 1.3 of Chapter 2. Assuming now $f\in D(\t A)$, the implication 
from (\ref{4.26*}) to (\ref{4.27*}) says that $f\in D(A)$. 
\qed 
\begin{theorem}\label{Theorem4.6} 
Assume (j)(jj),(jjj''), (vj), and (vjj). For all $f,g\in D(A)$ it holds that 
\begin{eqnarray*} 
\langle -Af,g\rangle_{L^2}+\langle -Ag,f\rangle_{L^2}=\left\langle\delta (A^f) 
\cdot f,g\right\rangle_{L^2}+\int fg\, d\bmu_\gamma\, . 
\end{eqnarray*} 
\end{theorem} 
Proof. Recalling $T^\ast_t\1\cdot\chi_G(\nu_{-t})=T^\ast_t\1$ and $T^\ast_t\1 
(\nu_t)\cdot r_t=1$ on $G$ by (\ref{4.20}), for $f,g\in D(A)$ we may deduce 
from (\ref{4.29*}) 
\begin{eqnarray*}
&&\hspace{-.5cm}\langle -Af,g\rangle_{L^2}-\left\langle\delta (A^f)\cdot f,g 
\right\rangle_{L^2}=\lim_{t\downarrow 0}\left\langle\frac{T^\ast_tf-f\cdot\chi_G 
(\nu_{-t})}{t}\, ,\, g\right\rangle_{L^2} \\ 
&&\hspace{.5cm}=\lim_{t\downarrow 0}\int\left(\frac{T^\ast_tf-f}{t}\right)\cdot g 
\, d\bmu+\lim_{t\downarrow 0}\int\frac{1-\chi_G(\nu_{-t})}{t}\cdot fg\, d\bmu 
 \\ 
&&\hspace{.5cm}=\lim_{t\downarrow 0}\int f\cdot\frac{T_tg-g}{t}\, d\bmu
+\lim_{t\downarrow 0}\frac1t\int\chi_{G\setminus G_t}\cdot fg\, d\bmu \\ 
&&\hspace{.5cm}=\langle f,Ag\rangle_{L^2}+\lim_{t\downarrow 0}\frac1t\int_{ 
\nu\in G_{-t}\setminus G} f(\nu_t)g(\nu_t)\, \bmu(d\nu_t)\, . 
\end{eqnarray*} 
Noting that for all $B\in {\cal B}(E)$ and $t\in [0,1]$ we have $\bmu(\nu_t: 
\nu\in B)-\bmu\circ\nu_t(B)=\bmu(\nu_t:\nu\in B\setminus G)$ it follows from 
Theorem \ref{Theorem2.7} (b) (Case 2), (jjj'), (vj), and (vjj) that 
\begin{eqnarray*}
\langle -Af,g\rangle_{L^2}-\left\langle\delta (A^f)\cdot f,g\right\rangle_{L^2} 
=\langle f,Ag\rangle_{L^2}+\int fg\, d\bmu_\gamma\, . 
\end{eqnarray*} 
\qed 

%% file: Lobus_-_April_2023.bbl
\begin{thebibliography}{1} {\small
\bibitem{AGS05}{Ambrosio, L., Gigli, N., Savar\'e, G.:} {Gradient flows 
in metric spaces and in the space of probability measures}. Lectures 
in Mathematics ETH Z\"urich. Birkh\"auser Verlag, Basel Boston Berlin (2005) 

\bibitem{AK06}{Albiac, F., Kalton, N. J.:} {Topics in Banach space 
theory}. Graduate texts in Mathematics {\bf 233}. Springer, New York 
(2006) 


\bibitem{Be85}{Bell, D.:} A quasi-invariance theorem for measures 
on Banach spaces. {Trans. Amer. Math. Soc.} {\bf 290}(2), 851-855 
(1985) 

\bibitem{Be06}{Bell, D.:} {The Malliavin calculus}. Dover Publications, 
Mineola, NY (2006) 

\bibitem{BB99}{van den Berg, M., Bolthausen, E.:} Estimates for Dirichlet 
eigenfunctions. {J. London Math. Soc.} {\bf 59}(2), 607-619 (1999)

\bibitem{Bo97}{Bogachev, V.:} Differentiable measures and the 
Malliavin calculus. {J. Math. Sci.} {\bf 87}(4), 3577-3731 (1997) 

\bibitem{Bo07}{Bogachev, V.:} {Measure theory} Vol. I and II. Springer, 
Berlin (2007) 

\bibitem{Bo10}{Bogachev, V.:} {Differentiable measures and the 
Malliavin calculus}. Mathematical Surveys and Monographs {\bf 164}. 
AMS, Providence, RI (2010) 

\bibitem{BM-W99}{Bogachev, V., Mayer-Wolf, E.:} Absolutely continuous 
flows generated by Sobolev class vector fields in finite and infinite 
dimensions. {J. Funct. Anal.} {\bf 167}(1), 1-68 (1999) 

\bibitem{Bo95}{Bourbaki, N.:} General topology. Chapters 1-4. Elements of 
Mathematics. Springer, Berlin Heidelberg New York London Paris Tokyo (1995)

\bibitem{BHM00}{Burdzy, K., Holyst, R., March. P.:} A Fleming-Viot 
particle representation of the Dirichlet Laplacian. {Comm. Math. 
Phys.} {\bf 214}(3), 679-703 (2000) 

\bibitem{CPW98}{Caprino, S., Pulvirenti, M., Wagner, W.:} Stationary 
particle systems approximating stationary solutions to the Boltzmann 
equation. {SIAM J. Math. Analysis} {\bf 29}(4), 913-934 (1998) 

\bibitem{CDGR20}{C\'erou, F., Delyon, B., Guyader, A., Rousset, M.:} 
A central limit theorem for Fleming-Viot particle systems. {Ann. Inst. 
Henri Poincar\'e  Probab. Stat.} {\bf 56}(1), 637-666 (2020)

\bibitem{CM03}{Chernov, N., Markarian, R.:} Introduction to the ergodic 
theory of chaotic billiards. Second edition. Facultad de Ingenier\'{i}a
Universidad de la Rep\'{u}blica - Uruguay (2003) \\
{\tt https://www.fing.edu.uy/imerl/grupos/ssd/publicaciones/pdfs/articulos/ 
2003/2003CheMarIntrod.pdf} 

\bibitem{Cr83}{Cruzeiro, A. A.:} \'Equations diff\'erentielles sur l'espace 
de Wiener et formules de Cameron-Martin non-lin\'eaires. {J. Funct. Anal.} 
{\bf 54}(2), 206-227 (1983)

\bibitem{CBH05}{Cybulski, O., Babin, V., Holyst, R.:} Minimization of 
the Renyi entropy production in the space-partitioning process. {Physical 
Review E} {\bf 71}, 046130--1-046130--10 (2005)

\bibitem{DS88}{Daletskii, Yu. L., Sokhadze, G. A.:} Absolute continuity 
of smooth measures. {Funct. Anal. Appl.} {\bf 22}(2), 149-150 (1988) 

\bibitem{Da73}{Davies, E. B.:} Properties of the Green's function of 
some Schr\"odinger operators. {J. London Math. Soc.} {\bf 7}, 
483-491 (1973)

\bibitem{DS84}{Davies, E. B., Simon, B.:} Ultracontractivity and the 
heat kernel for Schr\"odinger operators and Dirichlet Laplacians. {J. 
Funct. Anal.} {\bf 59}(2), 335-395 (1984) 

\bibitem{DM82}{Dellacherie, C., Meyer, P.-A.:} Probabilities and 
Potential, B: Theory of Martingales. North Holland, Amsterdam New  
York Oxford (1982) 

\bibitem{DT21}{Debussche, A., Tsutsumi, Y.:} Quasi-invariance of 
Gaussian measures transported by the cubic NLS with third-order 
dispersion on {\bf T}. {J. Funct. Anal.} {\bf 281}(3), Paper No. 
109032, 23 pp (2021) 

\bibitem{EE87}{Edmunds, D. E., Evans, W.:} D. Spectral theory and 
differential operators. Oxford Mathematical Monographs. Oxford Science 
Publications. The Clarendon Press, Oxford University Press, New York 
(1987) 

\bibitem{EN06}{Engel, K.-J., Nagel, R.:} A short course on operator 
semigroups. Universitext. Springer, New York (2006). 

\bibitem{EK86}{Ethier, S. N., Kurtz, T.:} {Markov processes, 
Characterization and convergence}. John Wiley, New York Chichester 
Brisbane Toronto Singapore (1986) 


\bibitem{FM07}{Ferrari, P. A., Maric, N.:} Quasi stationary 
distributions and Fleming-Viot processes in countable spaces. 
{Electron. J. Probab.} {\bf 12}(24) 684-702 (2007). 

\bibitem{FT19}{Forlano, J., Trenberth, W. J.:} On the transport of 
Gaussian measures under the one-dimensional fractional nonlinear 
Schr\"odinger equations. {Ann. Inst. H. Poincar\'{e} Anal. Non 
Lin\'{e}aire} {\bf 36}(7), 1987-2025 (2019). 

\bibitem{FL07}{Fonseca, I., Leoni, G.:} {Modern methods in the 
calculus of variations: $L^p$ spaces}. Springer Monographs in 
Mathematics. Springer, New York (2007) 


\bibitem{GK04}{Grigorescu, I., Kang, M.:} Hydrodynamic limit for a 
Fleming-Viot type system. {Stochastic Process. Appl.} {\bf 110}(1), 
111-143 (2004) 

\bibitem{GOTW22}{Gunaratnam, T., Oh, T., Tzvetkov, N., Weber, H.:}   
Quasi-invariant Gaussian measures for the nonlinear wave equation in three 
dimensions. {Probab. Math. Phys.} {\bf 3}(2), 343-379 (2022)

\bibitem{JM22}{Journel, L., Monmarch\'e, P.:} Convergence of a particle 
approximation for the quasi-stationary distribution of a diffusion 
process: Uniform estimates in a compact soft case. ESAIM Probab. Stat. 
{\bf 26}, 1-25 (2022) 

\bibitem{Ka95}{Kato, T.:} Perturbation Theory for Linear Operators. 
Springer, Berlin Heidelberg New York (1995) 

\bibitem{K06}{Klenke, A.:} Probability theory. A comprehensive course.  
Universitext. Springer, London (2006) 


\bibitem{Ko19}{Kolokoltsov, V.:} Differential equations on measures 
and functional spaces. Birkh\"auser Advanced Texts: Basler Lehrb\"ucher. 
Birkh\"auser/Springer, Cham (2019)

\bibitem{LPR18}{Leli\`evre, T., Pillaud-Vivien, L., Reygner, J.:} 
Central limit theorem for stationary Fleming-Viot particle systems 
in finite spaces. {ALEA Lat. Am. J. Probab. Math. Stat.} {\bf 15}(2), 
1163-1182 (2018)

\bibitem{Lo051}{L\"obus, J.-U.:} A stationary Fleming-Viot type 
Brownian particle system. {Math. Z.} {\bf 263}(3), 541-581 (2009) 

\bibitem{Lo20}{L\"obus, J.-U.:} Knudsen type group for time in 
$\mathbb{R}$ and related Boltzmann type equations. {Commun. Contemp. 
Math.} {\bf 25}(3), Paper No. 2150072, 72 pp (2023)


\bibitem{Ma02}{Matveev, O. V.:} Bases in Sobolev spaces on bounded 
domains with a Lipschitz boundary. {Math. Notes} {\bf 72}(3-4), 
373-382 (2002)

\bibitem{Ma03}{Matveev, O. V.:} On bases in Sobolev spaces. {Dokl. 
Akad. Nauk} {\bf 393}(6), 740-743 (2003) 

\bibitem{OST18}{Oh, T., Sosoe, P., Tzvetkov, N.:} An optimal regularity 
result on the quasi-invariant Gaussian measures for the cubic fourth 
order nonlinear Schr\"odinger equation., {J. \'Ec. polytech. Math.} 
{\bf 5}, 793-841 (2018). 

\bibitem{OT17}{Oh, T., Tzvetkov, N.:}  Quasi-invariant Gaussian 
measures for the cubic fourth order nonlinear Schr\"odinger equation. 
{Probab. Theory Related Fields} {\bf 169}(3-4), 1121-1168 (2017)

\bibitem{OT20}{Oh, T., Tzvetkov, N.:} Quasi-invariant Gaussian measures 
for the two-dimensional defocusing cubic nonlinear wave equation. 
{J. Eur. Math. Soc.} {\bf 22}(6), 1785-1826 (2020)

\bibitem{Pa83}{Pazy, A.:} Semigroups of linear operators and 
applications to partial differential equations. Springer, New York 
(1983) 

\bibitem{PTV20}{Planchon, F., Tzvetkov, N., Visciglia, N.:} Transport 
of Gaussian measures by the flow of the nonlinear Schr\"odinger equation. 
{Math. Ann.} {\bf 378}(1-2), 389-423 (2020) 

\bibitem{Re13}{Rehmann, U.:} Convergence of measures. {Encyclopedia 
of Mathematics}. Springer and EMS, 
{\tt http://www.encyclopediaofmath.org/index.php/Convergence{\_}of{\_}measures } 
(2013) 

\bibitem{Ry02}{Ryan, R. A.:} Introduction to Tensor Products of Banach 
Spaces, Springer, London Berlin Heidelberg (2002)

\bibitem{SvW93}{Smolyanov, O. G., von Weizs\"acker, H.:} Differentiable 
families of measures. {J. Funct. Anal.} {\bf 118}(2), 454-476 (1993) 

\bibitem{STX20}{Sosoe, P., Trenberth, W. J., Xian, T.:} Quasi-invariance of 
fractional Gaussian fields by the nonlinear wave equation with polynomial 
nonlinearity. {Differential Integral Equations} {\bf 33}(7-8), 393-430 
(2020)

\bibitem{SV09}{Song, V., Vondra$\check{\rm c}$ek, Z.:} Potential theory of 
subordinate Brownian motion. In: Graczyk, P., Stos, A. (eds.), Potential 
Analysis of Stable Processes and its Extensions, Lecture Notes in Mathematics
{\bf 1980}, 87-183. Springer, Berlin Heidelberg (2009)

\bibitem{Ta11}{Tao, T.:} An Introduction to Measure Theory. Graduate Studies 
in Mathematics {\bf 126}. AMS, Providence, RI (2011)  

\bibitem{Tz15}{Tzvetkov, N.:} Quasiinvariant Gaussian measures for 
one-dimensional Hamiltonian partial differential equations. {Forum Math. 
Sigma} {\bf 3}, Paper No. e28, 35 pp. (2015) 

\bibitem{UZ00}{\"Ust\"unel, A. S., Zakai, M.:} Transformation of Measure 
on Wiener Space. Springer Monographs in Mathematics. Springer, Berlin (2000)

\bibitem{Va92}{van den Berg, M.:} On the spectral counting function for the 
Dirichlet Laplacian. {J. Funct. Anal.} {\bf 107}(2), 352-361 (1992) 

\bibitem{Vi11}{Villemonais, D.:} Interacting particle systems and Yaglom 
limit approximation of diffusions with unbounded drift. {Electron. J. 
Probab.} {\bf 16}(61), 1663-1692 (2011) 

\bibitem{Vi14}{Villemonais, D.:} General approximation method for the 
distribution of Markov processes conditioned not to be killed. {ESAIM 
Probab. Stat.} {\bf 18}, 441-467 (2014) 

\bibitem{Ze86}{Zeidler, E.:} Nonlinear Functional Analysis and Its Applications 
I: Fixed-Point Theorems. Springer, New York (1986) 

}
\end{thebibliography}
